\documentclass[twoside, 10pt]{article}

\usepackage[hmargin = 0.75in, height =9.5in]{geometry}
\usepackage{amssymb,amsthm, amsfonts,amsmath,mathrsfs,  color, amsxtra, url, hyperref}

\usepackage{fancyhdr}
\fancypagestyle{plain}{
 \fancyhf{}

}

\renewcommand{\theequation}{\thesection.\arabic{equation}}
 \numberwithin{equation}{section}

\newtheorem {thm}{Theorem}[section]
\newtheorem {prop}{Proposition}[section]
\newtheorem {lemm}{Lemma}[section]
\newtheorem {deff}{Definition}[section]
\newtheorem {cor}{Corollary}[section]
\newtheorem {rem}{Remark}[section]

\def\ba{\begin{array}}
\def\ea{\end{array}}
\def\bea{\begin{eqnarray}}
\def\eea{\end{eqnarray}}
\def\beas{\begin{eqnarray*}}
\def\eeas{\end{eqnarray*}}
\def\bi{\begin{itemize}}
\def\ei{\end{itemize}}
\def\bc{\begin{cases}}
\def\ec{\end{cases}}

\def\bhe{\begin{highlightequation}  }
\def\ehe{\end{highlightequation}  }

\def\ba{\begin{array}}
\def\ea{\end{array}}
\def\bea{\begin{eqnarray}}
\def\eea{\end{eqnarray}}
\def\beas{\begin{eqnarray*}}
\def\eeas{\end{eqnarray*}}
\def\bi{\begin{itemize}}
\def\ei{\end{itemize}}
\def\bc{\begin{cases}}
\def\ec{\end{cases}}

\def\bhe{\begin{highlightequation}  }
\def\ehe{\end{highlightequation}  }

\def\a{\alpha}
\def\ga{\gamma}

\def\d{\delta}
\def\e{\varepsilon}

\def\z{\zeta}
\def\k{\kappa}
\def\l{\lambda}
\def\vr{\varrho}
\def\si{\sigma}

\def\th{\theta}
\def\o{\omega}
\def\vf{\varphi}
\def\vth{\vartheta}

\def\D{\Delta}
\def\Ga{\Gamma}
\def\L{\Lambda}
\def\O{\Omega}

\def\Th{\Theta}
\def\U{\Upsilon}

\def\bC{{\bf C}}

\def\bF{{\bf F}}
\def\bG{{\bf G}}

\def\bK{{\bf K}}

\def\bQ{{\bf Q}}

\def\bW{{\bf W}}

\def\bz{{\bf 0}}

\def\bw{{\bf w}}
\def\bx{{\bf x}}

\def\bu{{\bf u}}
\def\bx{{\bf x}}

\def\cA{{\cal A}}
\def\cB{{\cal B}}
\def\cC{{\cal C}}
\def\cD{{\cal D}}
\def\cE{{\cal E}}
\def\cF{{\cal F}}
\def\cG{{\cal G}}
\def\cH{{\cal H}}
\def\cI{{\cal I}}

\def\cK{{\cal K}}
\def\cL{{\cal L}}
\def\cM{{\cal M}}
\def\cN{{\cal N}}
\def\cO{{\cal O}}
\def\cP{{\cal P}}
\def\cQ{{\cal Q}}

\def\cS{{\cal S}}

\def\cU{{\cal U}}

\def\cW{{\cal W}}
\def\cX{{\cal X}}

\def\hJ{\mathbb{J}}

\def\hN{\mathbb{N}}

\def\hQ{\mathbb{Q}}
\def\hR{\mathbb{R}}

\def\hT{\mathbb{T}}
\def\hU{\mathbb{U}}

\def\hX{\mathbb{X}}
\def\hY{\mathbb{Y}}
\def\hZ{\mathbb{Z}}

\def\sB{\mathscr{B}}
\def\sC{\mathscr{C}}
\def\sD{\mathscr{D}}

\def\sH{\mathscr{H}}
\def\sI{\mathscr{I}}

\def\sM{\mathscr{M}}
\def\sN{\mathscr{N}}
\def\sO{\mathscr{O}}
\def\sP{\mathscr{P}}

\def\sS{\mathscr{S}}
\def\sT{\mathscr{T}}
\def\sU{\mathscr{U}}
\def\sV{\mathscr{V}}
\def\sW{\mathscr{W}}
\def\sX{\mathscr{X}}

\def\fA{\mathfrak{A}}
\def\fB{\mathfrak{B}}
\def\fC{\mathfrak{C}}

\def\fE{\mathfrak{E}}
\def\fF{\mathfrak{F}}

\def\fI{\mathfrak{I}}

\def\fN{\mathfrak{N}}

\def\fP{\mathfrak{P}}

\def\fS{\mathfrak{S}}
\def\fT{\mathfrak{T}}
\def\fU{\mathfrak{U}}

\def\fX{\mathfrak{X}}
\def\fY{\mathfrak{Y}}

\def\fp{\mathfrak{p}}
\def\fq{\mathfrak{q}}
\def\ff{\mathfrak{f}}
\def\fm{\mathfrak{m}}
\def\fu{\mathfrak{u}}

\def\fy{\mathfrak{y}}
\def\fz{\mathfrak{z}}
\def\fx{\mathfrak{x}}
\def\fl{\mathfrak{l}}

\def\fw{\mathfrak{w}}
\def\fra{\mathfrak{a}}
\def\fri{\mathfrak{i}}
\def\friJ{\fri_{\overset{}{\hJ}}}
\def\fk{\mathfrak{k}}
\def\fn{\mathfrak{n}}

\def\ft{\mathfrak{t}}
\def\fs{\mathfrak{s}}

\def\fc{\mathfrak{c}}

\def\fr{\mathfrak{r}}
\def\fu{\mathfrak{u}}
\def\fl{\mathfrak{l}}

\def\ti{\n \times \n}
\def\oti{\n \otimes \n}

\def\df{\n := \n}
\def\ls{\n \le \n}
\def\gs{\n \ge \n}
\def\={\n = \n}
\def\+{\n + \n}
\def\-{\n - \n}
\def\ins{\n \in \n}
\def\ld{\n \land \n}
\def\ve{\n \vee \n}
\def\sb{\n \subset \n}
\def\>{\n > \n}
\def\<{\n < \n}
\def\Cp{\n \cap \n}
\def\cp{\n \cup \n}

\def\nxi{ \rule[-0.5mm]{0.45pt}{2.6mm}\hspace{-0.4pt}\rule[-0.5mm]{1.8mm}{0.45pt}\hspace{-1.4mm}\rule[0.8mm]{1mm}{0.45pt}\hspace{-1.4mm}\rule[2mm]{1.8mm}{0.45pt}\hspace{-0.4pt}\rule[-0.5mm]{0.4pt}{2.6mm}\,}
\def\snxi{ \rule[-0.4mm]{0.4pt}{1.7mm}
\hspace{-0.35pt}\rule[-0.5mm]{1.2mm}{0.4pt}
\hspace{-0.85mm}\rule[0.37mm]{0.6mm}{0.4pt}
\hspace{-0.85mm}\rule[1.15mm]{1.2mm}{0.4pt}
\hspace{-0.35pt}\rule[-0.5mm]{0.4pt}{1.7mm}\,}
\def\hAe{\wh{A}^{\phantom{\rule[-0.8mm]{1pt}{2.6mm}} \e}}

\def\({\textnormal{(}}
\def\){\textnormal{)}}
\def\[{[\n[}
\def\]{]\n]}
\def\lan{\langle}
\def\ran{\rangle}

\def\no{\noindent}

\def\ss{\smallskip}

\def\q{\quad}
\def\qq{\qquad}
\def\n{\negthinspace}
\def\dn{\n \n}
\def\tn{\n \n \n}

\def\ol{\overline}
\def\ul{\underline}
\def\ua{\mathop{\uparrow}}
\def\da{\mathop{\downarrow}}

\def\Ra {\mathop{\Rightarrow }}

\def\wt{\widetilde}
\def\wh{\widehat}

\def\oMoo{\oM^{\raisebox{2.5pt}{\scriptsize \hb{$t_\oo,\ol{\mu}^\oo$}}}}

\def\loo{{\overset{}{\oo}}}

\def\hb{\hbox}
\def\dis{\displaystyle}
\def\cd{\cdot}
\def\cds{\cdots}

\def\fa{\,\forall \,}
\def\pa{\partial}
\def\es{\emptyset}

\def\b1{{\bf 1}}
\def\qed{\hfill $\Box$ \medskip}

\def\liminf{\mathop{\ul{\rm lim}}}
\def\limsup{\mathop{\ol{\rm lim}}}
\newcommand{\Sup}[1]{\underset{#1}{\sup}\,}

\newcommand{\lsup}[1]{ \underset{#1}{\limsup}}
\newcommand{\linf}[1]{ \underset{#1}{\liminf}}
\newcommand{\lmt}[1]{ \underset{#1}{\lim}}
\newcommand{\lmtu}[1]{ \underset{#1}{\lim} \n \ua \,}
\newcommand{\lmtd}[1]{ \underset{#1}{\lim} \n \da \,}

\newcommand{\ccap}[2]{\underset{#1}{\overset{#2}{\cap}}}
\newcommand{\ccup}[2]{\underset{#1}{\overset{#2}{\cup}}}

\newcommand{\Rho}[1]{\rho_{\overset{}{#1}}}

\def\oP{{\ol{P}}}
\def\oQ{{\ol{Q}}}
\def\oA{{\ol{A}}}
\def\oeta{{\ol{\eta}}}
\def\otau{{\ol{\tau}}}
\def\oga{{\ol{\ga}}}
\def\oz{{\ol{\z}}}

\def\ocP{{\ol{\cP}}}
\def\ocI{{\ol{\cI}}}

\def\ocA{{\ol{\cA}}}

\def\ocF{{\ol{\cF}}}
\def\obF{{\ol{\bF}}}
\def\ocG{{\ol{\cG}}}
\def\obG{{\ol{\bG}}}
\def\obK{{\ol{\bK}}}

\def\oo{{\ol{\o}}}
\def\oO{{\ol{\O}}}

\def\oJ{{\ol{J}}}
\def\oW{{\ol{W}}}
\def\oK{{\ol{K}}}

\def\oX{{\ol{X}}}
\def\oY{{\ol{Y}}}
\def\oZ{{\ol{Z}}}
\def\oU{{\ol{U}}}
\def\oV{{\ol{V}}}
\def\oR{{\ol{R}}}
\def\oT{{\ol{T}}}
\def\oD{{\ol{D}}}
\def\oTh{{\ol{\Th}}}

\def\oM{{\ol{M}}}
\def\osX{\ol{\sX}}

\def\wA{{\wt{A}}}

\def\mto{\n \mapsto \n}
\def\nto{\n \to \n}

\def\bul{\no $\bullet$ }

\def\oxi{{\ol{\xi}}}

\def\ocM{{\ol{\cM}}}

\def\ocN{{\ol{\cN}}}

\def\ofN{{\ol{\fN}}}

\def\ocW{{\ol{\cW}}}
\def\osW{{\ol{\sW}}}
\def\osU{{\ol{\sU}}}
\def\obW{{\ol{\bW}}}

\def\ooV{{\ol{\oV}}}
\def\oocP{{\ol{\ocP}}}

\def\Wtzo{{\obW^t_{\otau,\oo}}}
\def\Wtgo{{\obW^t_{\oga,\oo}}}
\def\Ktgo{{\obK^t_{\oga,\oo}}}

\def\gP{{\big[ \n \big[\,\ocP\,\big] \n \big]}}
\def\gcP{{\big\{ \n \big\{ \ocP \big\} \n \big\}}}
\def\nci{\n \circ \n}
\def\tx{(t,\bx) \ins [0,\infty) \ti \OmX}

\def\OmX{\O_{\overset{}{X}}}
\def\omX{\o_{\overset{}{X}}}

\def\has{\exists \,}
\def\aand{\q \hb{and} \q}

\def\nne{\n \ne \n}

\def\cad{c\`adl\`ag }

\def\ogaP{{\oga_{\overset{}{\oP}}}}
\def\btau{\rule[1.45mm]{2.1mm}{0.7pt}\hspace{-3.35pt} \rule[-0.15mm]{0.8pt}{1.8mm}~}
\def\bbtau{\rule[1.3mm]{1.6mm}{0.6pt}\hspace{-2.6pt} \rule[0.1mm]{0.7pt}{1.3mm}\;}

\def\usa{upper semi-analytic}

\begin{document}

 \title{\bf Stochastic Control/Stopping Problem \\  with Expectation Constraints
    }

\author{
 Erhan Bayraktar\thanks{ \noindent Department of  Mathematics,
  University of Michigan, Ann Arbor, MI 48109; email: {\tt erhan@umich.edu}.}
 \thanks{E. Bayraktar is supported in part  by the National Science Foundation under  DMS-2106556,
 and in part by the Susan M. Smith Professorship.
 Any opinions, findings, and conclusions or recommendations expressed in this material are
 those of the authors and do not necessarily reflect the views of the National Science Foundation.} $\,\,$,
 $~~$Song Yao\thanks{
 \noindent Department of  Mathematics,
  University of Pittsburgh, Pittsburgh, PA 15260; email: {\tt songyao@pitt.edu}. }
  \thanks{S. Yao is supported in part  by the National Science Foundation under DMS-1613208.
} }

\date{}

\maketitle

 \begin{abstract}

 We study  a stochastic  control/stopping problem  with a series of inequality-type and  equality-type expectation constraints in a general non-Markovian framework.
We demonstrate that the stochastic  control/stopping problem  with expectation constraints (CSEC)
is independent of a specific probability setting
and   is equivalent to the constrained stochastic  control/stopping problem in weak formulation (an optimization over joint laws of Brownian motion, state dynamics, diffusion controls  and stopping rules   on an enlarged canonical space).
Using a martingale-problem formulation of controlled SDEs in spirit of \cite{Stroock_Varadhan}, we  characterize the probability classes  in   weak formulation by countably many actions of canonical processes, and thus obtain the upper semi-analyticity of   the CSEC value function.
  Then we employ  a measurable selection argument to establish a dynamic programming principle (DPP)  in weak formulation for the CSEC value function,
 in which the conditional expected costs  act  as   additional states for constraint levels at the intermediate horizon.

 This article extends 
 \cite{Elk_Tan_2013b} to the expectation-constraint case.
We extend our previous work
  \cite{OSEC_stopping} to the more complicated setting where the diffusion is controlled. Compared to that paper
 the topological properties of   diffusion-control spaces  and the corresponding measurability are more technically involved which complicate the arguments especially for the measurable selection for the super-solution side  of DPP in the weak formulation.

 \ss \no {\bf MSC:}\; 93E20, 60G40,  49L20, 60G44


 \ss \no   {\bf Keywords:}\; Stochastic control/stopping problem  with expectation constraints,
  martingale-problem formulation, enlarged canonical space, Polish space of diffusion controls,
  Polish space of stopping times,
   dynamic programming principle, regular conditional probability distribution,
 measurable selection.

\end{abstract}

  \section{Introduction}

  In this paper,  we analyze a continuous-time  stochastic control/stopping problem with
  a series of inequality-type and  equality-type expectation constraints   in a general non-Markovian framework.

 Let the game start  from time $t \ins [0,\infty)$ with a historical   path of state $\bx|_{[0,t]}$.
 The player can choose an open-loop control $\mu \= \{\mu_s\}_{s \in [t,\infty)} $
 to make the state  of the  game  evolve according to  some controlled  SDE 
 on a probability space $(\cQ, \cF, \fp)$
 whose   drift and diffusion coefficients depend on the past trajectories   of the solution.
 Let $\cX^{t,\bx,\mu} \= \big\{\cX^{t,\bx,\mu}_s\big\}_{s \in [t,\infty)} $ denote this controlled state process.
 The player also  decides an exercise time $\tau  $ 
 to maximize     the expectation of   her accumulative reward $\int_t^\tau   \n  f  \big(r, \cX^{t,\bx,\mu}_{r \land \cd}, \mu_r  \big) \, dr $ plus her terminal reward $   \pi \big(  \tau  , \cX^{t,\bx,\mu}_{ \tau   \land \cd}  \big) $
   while she is subject to a series of constraints:  for $i \ins \hN$,
 the expectation  of  some   accumulative   cost  $\int_t^\tau \n  g_i \big( r,\cX^{t,\bx,\mu}_{r \land \cd}, \mu_r \big) dr $ 
   should not exceed certain level $y_i$ and the expectation  of  some other  accumulative   cost  $\int_t^\tau \n  h_i \big( r,\cX^{t,\bx,\mu}_{r \land \cd}, \mu_r \big) dr $ should exactly reach  certain level $z_i$.
 This   stochastic control/stopping problem with expectation constraints (SCEC for short)
 has many  applications in  economy, engineering, finance, management,  etc.

 Let $V (t,\bx,y,z)$ be the SCEC value with   $(y,z) \df \big(\{y_i\},\{z_i\}\big)$.
 We aim to 
   establish an associated  dynamic programming principle (DPP) of this value function 
without imposing any  
regularity condition on  reward/cost functions in time, state and control variables.  
A dynamic programming principle   of a stochastic optimization problem allows one  to maximize/minimize  the problem stage by stage
in a backward recursive way.  
 It entails the problem value function to be  measurable   so that one can do optimization first at an intermediate horizon.

 To obtain  the measurability of the  SCEC value function,
 we first study  the topological and measurable properties  of the path space $\hJ$ of diffusion control $\mu$:
 More precisely, we show that  $\hJ$ is a Borel space under a weak topology  and attain a representation of 
 its Borel sigma-field $\sB(\hJ)$ (Lemma \ref{lem_Nov25_03} and Lemma \ref{lem_M29_01}).
    Inspired by \cite{EHJ_1987,Elk_Tan_2013b},   we then  embed   
   diffusion control $\mu$ and stopping rule  $\tau$ together with the Brownian and state information into an enlarged canonical space $\oO$ and regard their joint distribution as  a new type of controls.
   The optimization of the total expected reward  over  constrained diffusion controls/stopping times  transforms into
 a maximal  expectation of   reward functional  over a class $ \ocP_{t,\bx}(y,z) $ of probability measures  on $\oO$  under which four canonical coordinates $(\oW,\oU,\oX,\oT)$  serve as  Brownian motion, 
  diffusion control, state process and 
 stopping rules respectively.
  We demonstrate that   such a transformation is equivalent (Theorem \ref{thm_V=oV}), namely,
  the  value $V (t,\bx,y,z)$ of SCEC in   strong formulation (i.e., on $\cQ$) is equal to the value $\oV (t,\bx,y,z)$ of SCEC in   weak formulation (i.e., over $\oO$).  
  Hence, the SCEC value 
  is   a  robust value,    independent of a specific probability model.

 For the measurability of   SCEC value functions,
 we next take advantage of the martingale-problem formulation from \cite{Stroock_Varadhan} to
 describe the   probability class   $\ocP_{t,\bx}(y,z)$ 
  as a series of probabilistic tests on   stochastic behaviors of the canonical coordinates  of $\oO$.
 With such  a countable characterization,
  we employ a Polish space of diffusion control processes (Lemma \ref{lem_082020_15}) 
  and a Polish space of stopping times 
  constructed in \cite{OSEC_stopping}
 to deduce that   the set-valued mapping   $(t,\bx,y,z) \mto \ocP_{t,\bx}(y,z)$ has Borel-measurable graph
 and the SCEC value function  $V \= \oV$   is thus  upper semi-analytic in $(t,\bx,y,z)$, (Theorem \ref{thm_V_usa}).

 Our main achievement    is to derive a DPP for $\oV$ in weak formulation, 
 which takes conditional expectations of the remaining costs as additional states for constraint levels at the intermediate horizon (Theorem \ref{thm_DPP1}).
 For the subsolution side of this DPP, we use the regular conditional probability distribution
  to show that the probability classes $  \ocP_{t,\bx}(y,z) $, $ \fa (t,\bx,y,z)   $ are stable under conditioning.

 For the supersolution side of the DPP, we exploit a measurable selection theorem in the analytic-set theory
 to paste  a class of locally $\e-$optimal probability measures.
 We make a delicate analysis to demonstrate that the second  canonical coordinate $\oU$ serves as a constrained diffusion control  under the pasted probability measure,
 and we apply the martingale-problem formulation again to indicate that  the canonical coordinates $(\oW,\oX)$ are still Brownian motion and the state process under the pasted probability measure. Similar to the arguments in \cite{OSEC_stopping}, the fourth  canonical coordinate $\oT$  is a constrained stopping time   under the pasted probability measure.
  To wit, the probability classes $  \ocP_{t,\bx}(y,z) $'s are also stable under   concatenation. 

 \ss \no {\bf Relevant Literature.}

  Kennedy \cite{Kennedy_1982} employed  a {\it Lagrange multiplier} method
  to reformulate a discrete-time optimal stopping problem with first-moment constraint as a minimax problem
  and   showed that the optimal value of the dual problem
  is   equal to that of the primal problem.
    The Lagrangian technique was later adopted in many economic/financial applications of
    optimal stopping problems with expectation constraints,  see e.g. \cite{Pontier_Szpirglas_1984,LSMS_1995,Horiguchi_2001b,Balzer_Jansen_2002,Urusov_2005,Makasu_2009,Peskir_2012,Pedersen_Peskir_2017,Tanaka_2019}.
   Pfeiffer et al. \cite{PTZ_2020} recently took  a Lagrange relaxation approach
  to obtain a duality result  for general stochastic control problems with expectation constraints.

 In their study of a   continuous-time stochastic optimization problem of controlled Markov processes,
 El Karoui, Huu Nguyen and Jeanblanc-Picqu\'e \cite{EHJ_1987}
  viewed   joint laws of    state and control processes
 as {\it control rules} on the product space of canonical state  space and control space.
 They utilized a measurable selection theorem in the analytic-set theory  to establish a DPP
  without assuming any regularity on the reward functional.
    Nutz et al.  \cite{HN_2012,Neufeld_Nutz_2013} came up with a similar idea
  to analyze a superheging problem  under volatility uncertainty.
  They modeled the ``uncertainty" by path-dependent classes of controlled-diffusion laws  
   and explored the analytic measurability    of these classes. Using the measurable selection techniques,
    the authors obtained  DPP result in a form of   time-consistency 
     of a sub-linear expectation and they thus established  a duality formula for the robust superhedging of measurable claims.
  The approach of \cite{HN_2012,Neufeld_Nutz_2013} was later developed by e.g. \cite{PRT_2013,PTZ_2018} to derive  DPPs of various non-Markovian   control problems.
  \if{0}

  see  \cite{PRT_2013} for a dual formulation of  robust semi-static trading and its application  to   martingale optimal transportation and  see \cite{PTZ_2018} for stochastic control of a class of nonlinear kernels
  and its relation to second-order backward stochastic differential equations.

 \fi
    Yu et al. \cite{CYZ_2020} 
   took a similar measurable selection argument to analyze
  the DPP of a stochastic control problem  with certain expectation constraint 
  in which they dynamically relaxed  the expectation constraint by a family of auxiliary supermartingales. 

     El Karoui and Tan \cite{Elk_Tan_2013a,Elk_Tan_2013b} used the measurable selection argument to attain the DPP for a general stochastic control/stopping problem by embedding diffusion  controls and stopping times   into an enlarged canonical space  in the spirit of \cite{EHJ_1987}. However, the probability class they considered in weak formulation is not suitable for stochastic control/stopping with expectation constraints, see our discussion in Remark \ref{rem_032322}.
   Instead,   we additional   require in (D4) of Definition \ref{def_ocP} that under each $\oP$ of $\ocP_{t,\bx}(y,z)$
   the time canonical coordinate $\oT$ acts as some   stopping time
   (it turns out that such a restriction does not affect the unconstrained stochastic control/stopping problem in weak formulation).
     By  constructing a Polish space of diffusion control processes   
  and utilizing a Polish space of stopping times 
  from \cite{OSEC_stopping}, we manage to derive the Borel measurability of   graph $ \gP$ and thus obtain the measurability of the SCEC value functions.
    Because of   condition (D4) and   expectation constraints, it is more technically involved to
    verify the stability of our  probability classes $\ocP_{t,\bx}(y,z)$ under conditioning and concatenation
    and thus establish a DPP for the SCEC value function $\oV$.

  As to the optimal stopping problems with expectation constraints,  Ankirchner et al. \cite{AKK_2015} and Miller  \cite{Miller_C_2017a} took   different approaches   by transforming the constrained   optimal stopping problems for diffusion processes to stochastic optimization problems with martingale controls.
  The former   characterizes the value function 
  in terms of a Hamilton-Jacobi-Bellman   equation and obtains a verification theorem, while the latter embeds the optimal stopping problem with first-moment constraint   into  a time-inconsistent (unconstrained) 
  stopping problem. However, the authors only postulate   dynamic programming principles  for  their corresponding  problems.
  In contrast, our previous work \cite{OSEC_stopping} exploited a measurable selection method to rigorously establish
    a dynamic programming principle for the  optimal stopping problem   with expectation constraints.

  An interesting related topic to our research is optimal stopping   with constraint on the distribution of   stopping times.
 Bayraktar and Miller \cite{Bayraktar_Miller_2016}
 studied the problem of optimally stopping a Brownian motion with the restriction that
   the distribution of the stopping time must equal a given measure with finitely many atoms,
 and obtained a dynamic programming result which relates each of the sequential optimal control problems.
  K\"allblad \cite{Kallblad_2017}  used measure-valued martingales   to transform
  the distribution-constrained optimal stopping problem to a stochastic
 control problem and derived a DPP by measurable selection arguments.
  From the perspective of optimal transport,  Beiglb\"ock et al. \cite{BEES_2016} gave a geometric interpretation  of optimal stopping times  of a Brownian motion with distribution constraint.

   Moreover, for stochastic control problems  with state constraints,   stochastic target problems with controlled losses and   related geometric DPP,  
 see \cite{BEI_2009,BET_2009,Bouchard_Nutz_2012,Soner_Touzi_2002_GDP,Soner_Touzi_2002_SDV,
 Soner_Touzi_2009,Bouchard_Vu_2010,Bouchard_Dang_2013,BMN_2014,BDK_2017}.

 The rest of the paper is organized as follows:
 Section \ref{sec_genprob} introduces the stochastic control/stopping problem with expectation constraints
 in a generic probabilistic setting.
   Section \ref{sec_weak_form}  shows that  the   stochastic control/stopping problem with expectation constraints
 can be equivalently embedded into an enlarged canonical space: i.e.,
 the SCEC 
 in   strong formulation has the same value as the SCEC in   weak formulation.
 In Section \ref{sec_Mart_prob}, we  use the martingale-problem formulation to make a countable characterization of the probability class in weak formulation. With such a characterization,
   we employ   a Polish space of diffusion control processes and a Polish space of stopping times to demonstrate that  the SCEC  value function
 is upper semi-analytic. Then in Section \ref{sec_DPP},  we utilize  a measurable selection argument to establish a dynamic programming principle in weak formulation for the SCEC value function. 
 We defer the proofs of our results to  Section \ref{sec_proof} and put some technical lemmata in the appendix.

  We close this section by a description of our notation and a review of the martingale-problem formulation of controlled SDEs.

\subsection{Notation and Preliminaries}

\label{subsec:preliminary}

  Throughout this paper, let us denote  $a^+   \df a \ve 0$ and $a^-   \df (  -  a) \ve 0$ for any   $ a \ins \hR$.
  We set $\hQ_+ \df \hQ \cap [0,\infty)$, $ \hQ^{2,<}_+ \df \\ \big\{ (s,r) \ins \hQ_+ \ti \hQ_+ \n :   s \< r \big\} $
   and set   $\Re \df   (-\infty,\infty]^\hN $ as the product of countably many copies of $(-\infty,\infty]$.
 On $\hT \df [0,\infty]$  we define a metric   $\Rho{+}(t_1,t_2) \df \big|\arctan(t_1) \- \arctan(t_2)\big|$, $\fa t_1,t_2 \ins \hT$
 and consider the induced topology by $\Rho{+}$.

 For a generic topological space $\big( \hX,\fT(\hX) \big)$,
    denote its  Borel sigma-field by $\sB(\hX)$. We let $\fP(\hX)$ be the set of  all probability measures on $\big(\hX,\sB(\hX)\big)$
  and equip $\fP(\hX)$ with the topology of weak convergence $\fT_\sharp\big(\fP(\hX)\big)$.
   By e.g. Corollary 7.25.1 of \cite{Bertsekas_Shreve_1978}, $\big(\fP(\hX), \fT_\sharp(\fP(\hX)) \big)$
   is a {\it Borel space} (i.e., homeomorphic to a Borel subset of a complete separable metric space).

  Let $ n \ins \hN$.  For any $x \ins \hR^n$ and $\d \ins (0,\infty)$,   let $O_\d(x)$ denote
 the open ball   centered at $x$ with radius $\d$  and let  $\ol{O}_\d(x)$ be its closure.
  For any $x,\wt{x} \ins \hR^n$
  we denote the usual inner  product by   $x   \cd   \wt{x} \df \sum^n_{i=1} x_i \wt{x}_i$, and for any  $ n \ti n -$real matrices $A,\wt{A}  $
  we denote the Frobenius inner  product by $A \n : \n \wt{A} := trace\big(A \wt{A}^T\big)$, where $\wt{A}^T$ is the transpose of $\wt{A}$.
 Let $ \big\{\cE^n_i\big\}_{i \in \hN}$ be a countable subbase of  the Euclidean topology   $\fT (\hR^n) $
 on $\hR^n$.
 Then $ \sO (\hR^n) \df  \Big\{ \ccap{i=1}{n}  \cE^n_{k_i} \n : \{ k_i  \}^n_{i=1} \sb \hN  \Big\}
 \cp \{\es,\hR^n\}$
 forms a countable base of $\fT (\hR^n) $ and thus $\sB(\hR^n) \= \si\big(\sO (\hR^n)\big)$.
 We also set $\wh{\sO} (\hR^n) \df \ccup{k \in \hN}{} \big( \hQ_+ \ti \sO (\hR^n)   \big)^k $.
  For any $\vf \ins C^2(\hR^n)$,
  let $ D \vf   $ be its gradient, $ D^2 \vf $ be its Hessian matrix and denote $D^0\vf \df \vf$.
 \if{0}
 For any $f \ins C^1(\hR^n)$ and $x_o \ins \hR^n$, let $Df(x_o) \ins \hR^n$ be the gradient of $f$ at $x_o$
 and let $D^2 f(x_o) \ins \hR_{n \ti n}$ be the Hessian matrix of $f$ at $x_o$, i.e.,
 the $i-$the component of $Df(x_o)$ is $D_i f(x_0) \df \frac{\pa f}{\pa x_i} (x_o)$ and
 the element on $i-$th row and $j-$th column of $D^2 f(x_o)$ is  $D_{ij} f(x_0) \df \frac{\pa^2 f}{\pa x_i \pa x_j} (x_o)$
 for $i,y \= 1, \cds, n$.
 \fi
 For $ i  \= 1, \cds, n $,    define $\vf_i(x) \df x_i$,  $\fa x \= (x_1,\cds \n ,x_n) \ins \hR^n $.
 We let  $\fC(\hR^n) $ collect these coordinate  functions   and their products, i.e.,
 $\fC(\hR^n) \df \{\vf_i\}^n_{i=1} \cp \{\vf_i \vf_j\}^n_{i,j=1} $.

  Let  $(\O,\cF,P )$ be a generic probability  space. For   subsets $A_1,A_2$ of $   \O$, we  denote $A_1 \D A_2 \df ( A_1 \Cp A^c_2 ) \cp ( A_2 \Cp A^c_1 ) $.
  For a random variable $\xi$ on $\O$ with values in a measurable space $(\cQ,\cG)$,
  we say   $\xi$ is $\cF/\cG-$measurable if its induced sigma-field $\xi^{-1}(\cG) \df \{\xi^{-1}(\cA) \n : \fa \cA \ins \cG \} $ is
  included in $ \cF $.
 For  a  sub-sigma-field $\fF  $ of $\cF$,
 define $\sN_P(\fF) \df \big\{ \cN \sb \O \n : \cN \sb A $ for some $ A \ins \fF \hb{ with } P (A) \=0  \big\}$,
 which collects all $P-$null sets with respect to   $\fF  $.
  For two sub-sigma-fields $\fF_1, \fF_2  $ of $\cF$,
 we denote   $\fF_1 \ve \fF_2 \df \si(\fF_1 \cp \fF_2)$.
 Let $t \ins [0,\infty)$.
 For a  filtration $\bF \= \{\cF_s\}_{s \in [t,\infty)}$ of $ \cF  $, we decree $\cF_{t-} \df \cF_t $ and define $ \cF_{s-} \df \si \Big(\underset{r \in [t,s)}{\cup} \cF_r \Big) $, $\fa s \ins (t,\infty)$; we also set $\cF_\infty \df \si\Big(\underset{s \in [t,\infty)}{\cup}\cF_s\Big)$
 and refer to   filtration $\bF^P \n \= \big\{\cF^P_s \n \df \si \big( \cF_s \cp \sN_P ( \cF_\infty ) \big) \big\}_{s \in [t,\infty)} $ as the $P-$augmentation of $ \bF $.
 For  a process $X\= \{X_s\}_{s \in [t,\infty)}$  on $ \O  $ with values in a topological space,
   its raw filtration is $\bF^X \=  \big\{ \cF^X_s \df \si (X_r ; r \ins [t,s])  \big\}_{s \in [t,\infty)}$.
 We denote the $P-$augmentation of $\bF^X $ by $\bF^{X,P} \n \= \big\{\cF^{X,P}_s \n \df \si \big( \cF^X_s \cp \sN_P ( \cF^X_\infty ) \big) \big\}_{s \in [t,\infty)} $ and let $\sP^{X }$ be  the $\bF^X-$predictable sigma$-$field  of $[t,\infty) \ti \O $.
 We call $X$ a continuous process if its paths are all continuous.
 When the time variable $s$ of $X$ has complicated form, we may write $X(s,\o)$ as $X_s(\o)$   for readability.
 By default, a Brownian motion $\{B_s\}_{s \in [t,\infty)}$ on $(\O,\cF,P )$  is with respect to its raw filtration
 $\bF^B$ unless stated otherwise.

 Fix    $d, l \ins \hN$.  Let $ \O_0 \=  \big\{ \o \ins  C   ( [0,\infty)  ; \hR^d ) \n : \o(0) \= 0 \big\}   $
 be the   space of   all $\hR^d-$valued continuous paths 
 starting  from $\bz$, which is a Polish space under the topology of locally uniform convergence.
  Let $ P_0  $ be the Wiener measure on $\big(\O_0,  \sB(\O_0) \big)$, under which
   the canonical process $ W \=\{W_s\}_{s \in [0,\infty)} $ of $\O_0$ is a   $d-$dimensional standard Brownian motion.
     For any $t \ins [0,\infty)$,   $W^t_s \df W_s \- W_t$, $  s \ins [t,\infty)$ is also a   Brownian motion on $\big( \O_0,\sB(\O_0),P_0 \big)$.
 Let $\OmX \= C  ( [0,\infty)  ; \hR^l  )$ be the  space of all  $\hR^l-$valued   continuous paths 
 endowed  with the topology of locally uniform convergence.
 The  function $\fl_1 (t,\o_0) \df \o_0(t \ld \cd)  $ is continuous in $ (t,\o_0) \ins  [0,\infty) \ti \O_0  $
 while the  function $\fl_2 (t,\omX) \df \omX(t \ld \cd)  $ is continuous in $ (t,\omX) \ins  [0,\infty) \ti \OmX  $.

  Let  $\hU$ be a  Polish space  with a compatible metric $ \Rho{\hU} $ and let $u_0 \ins \hU$.
    \if{0}

  Let  $\hU$ be a  Polish space  with  $u_0 \ins \hU$ and let $ \Rho{\hU} $ be  a    metric on $\hU$ that is compatible   with the topology of $\hU$.

   There exists a homeomorphism $  \fI$ from $\hU$ to a    complete metric space
 $\big(\hZ, \Rho{\hZ}\big)$ with a countable dense subset $\{z_i\}_{i \in \hN}$.
 We can induce a metric $\Rho{\hU}$ on $\hU$ by: $\Rho{\hU} (u_1,u_2) \= \Rho{\hZ} \big( \fI(u_1), \fI(u_2)\big)$, $\fa u_1,u_2 \ins \hU$. Then $\hU$ is also a separable  complete   metric space under $\Rho{\hU}$,
 namely, $\big(\hU,\Rho{\hU}\big)$ is a Polish metric space.

    \fi
  As a Polish space, $\hU$ is   homeomorphic to a Borel subset $ \fE $ of $[0,1]$,
  we denote this  homeomorphism by $\sI \n : \hU \mto \fE$.
  Let $C_b\big([0,\infty) \ti \hU\big) $ \big(resp. $\wh{C}_b\big([0,\infty) \ti \hU\big) $\big) collect  all real-valued bounded continuous (resp. bounded uniformly continuous) functions   on $[0,\infty) \ti \hU$.
  For any $(\d, \fm, \phi )  \ins (0,\infty) \ti  \fP \big([0,\infty) \ti \hU\big) \ti C_b\big([0,\infty) \ti \hU\big)$,  set
$ O_\d (\fm,\phi ) \df \big\{ \fm' \ins \fP \big([0,\infty) \ti \hU\big) \n : \big| \int_0^\infty
 \n   \int_\hU    \phi (t, u) \big(  \fm'(dt,du)  \-   \fm(dt,du)  \big)   \big|  \< \d  \big\} $.

 \begin{lemm} \label{lem_082020_11}

 There exist  $\{\fm_k\}_{k \in \hN } \sb \fP \big([0,\infty) \ti \hU\big) $ and
 $ \{\phi_j\}_{j \in \hN} \sb \wh{C}_b\big([0,\infty) \ti \hU\big)  $ such that
  $ \big\{ O_{\frac{1}{n}} (\fm_k,\phi_j ) \n : n,k,j \ins \hN \big\} $
 forms  a countable subbase of 
 $\fT_\sharp\big(\fP \big([0,\infty) \ti \hU\big)\big)$.
 \end{lemm}

  Let  $\hJ \df L^0\big([0,\infty);\hU\big)$ denote  the equivalence classes of $\hU$-valued Borel measurable functions on $[0,\infty)$
  in the sense that $\fu_1, \fu_2 \ins \hJ $ are equivalent if  $\fu_1(t) \= \fu_2(t) $ for a.e. $t \ins (0,\infty)$.
  It will serve as the path space of $\hU-$valued diffusion controls.
  We can embed $\hJ$ into $\fP \big([0,\infty) \ti \hU\big)$ via a  mapping
  $\friJ \n :  \hJ \n \ni \n \fu \mto e^{-t}\d_{\fu(t)} (du) dt \ins \fP \big([0,\infty) \ti \hU\big) $.
  Let $\fT_\sharp(\hJ)$ be the    topology induced by $\fT_\sharp\big(\fP \big([0,\infty) \ti \hU\big)\big)$  via $\friJ$.  
  According to  Lemma \ref{lem_082020_11},
    $\fT_\sharp(\hJ)$  is  generated by a countable subbase
  \bea \label{J04_01}
  \, \fri^{-1}_\hJ \big( O_{\frac{1}{n}} (\fm_k,\phi_j ) \big)  \=  \Big\{ \fu  \ins \hJ \n :
  \Big| \int_0^\infty e^{-t}   \phi_j (t,\fu(t)) dt \- \int_0^\infty \n \int_\hU
  \phi_j (t,u) \fm_k (dt,du)    \Big|  \< \frac{1}{n}    \Big\}, \q \fa n,k,j \ins \hN . \q 
  \eea

\begin{lemm} \label{lem_Nov25_03}
$ \big(\hJ,\fT_\sharp(\hJ)\big) $ is a Borel space.
\end{lemm}

 Let $L^0 \big((0,\infty) \ti \hU;\hR  \big) $ collect all real$-$valued  Borel-measurable functions   on $(0,\infty) \ti \hU$.
 For any $\vf \ins L^0 \big((0,\infty) \ti \hU ;\hR \big)$,   define
  $ I_\vf (\fu)  \df \int_0^\infty \vf (s,  \fu (s) ) ds  \= \int_0^\infty \vf^+ (s, \fu (s) ) ds
   \- \int_0^\infty \vf^- (s, \fu (s) ) ds  $,  $\fa \fu \ins \hJ$.
   The Borel sigma-field $\sB(\hJ)$ of  $\fT_\sharp(\hJ)$ can be generated by these random variables $I_\vf$  on $\hJ$.

\begin{lemm} \label{lem_M29_01}

1\) We have  $  \sB(\hJ)    \=  \si \big(I_\vf  ;  \vf \ins L^0 \big((0,\infty) \ti \hU ; \hR \big) \big)   $.

\no 2\) Let $\psi \n: (0,\infty) \ti \OmX \ti \hR^{d+l} \ti \hU \mto [-\infty,\infty]$ be a Borel-measurable function. Then
the mapping $ \Psi (t,\fs,\o_0,\omX,\fu) \df \int_t^{t+\fs} \psi  \big(r,\fl_2(r,\omX),\o_0(r), \omX(r),\fu(r)\big) dr
\= \int_t^{t+\fs} \psi^+ \big(r,\fl_2(r,\omX),\o_0(r), \omX(r),\fu(r)\big) dr
\- \int_t^{t+\fs} \psi^- \big(r,\fl_2(r,\omX),\o_0(r),   \omX(r), \\ \fu(r)\big) dr $, $(t,\fs,\o_0,\omX,\fu) \ins [0,\infty) \ti [0,\infty) \ti \O_0 \ti \OmX \ti \hJ$ is  $\sB[0,\infty) \oti \sB[0,\infty) \oti \sB(\O_0) \oti \sB(\OmX) \oti \sB(\hJ)-$measurable.

\end{lemm}

 \ss     Let   $ b  \n :   (0,\infty) \ti \OmX \ti \hU \mto  \hR^l  $ and $ \si  \n :   (0,\infty) \ti \OmX \ti \hU \mto   \hR^{l \times d} $ be two Borel-measurable functions such that   for any  $t \ins (0,\infty)$
    \bea
 &&  \big|b(t,\omX,u)\-b(t,\omX',u) \big| \+ \big|\si(t,\omX,u)\-\si(t,\omX',u) \big| \ls \k(t)   \big\|\omX\-\omX' \big\|_t ,
 \q \fa   \omX, \omX' \ins \OmX, ~ \fa  u \ins \hU , \q \qq   \label{coeff_cond1} \\
 && \hspace{2.5cm} \hb{and} \q \int_0^t \Sup{u \in \hU} \big( |b (r,\bz,u)|^2 \+ |\si (r,\bz,u)|^2 \big) dr \< \infty ,   \label{coeff_cond2}
 \eea
 \if{0}

 Let   $ \ff  \n :   (0,\infty)   \ti \hU \mto  [0,\infty]  $   be a   Borel-measurable function.
 By Proposition 7.47 of \cite{Bertsekas_Shreve_1978},   $  \Sup{u \in \hU}   \ff (s, u) $ is upper semianalytic  and is thus universally measurable in $s \ins (0,\infty)$.

 \fi
  where $\k \n : (0,\infty) \mto (0,\infty)$ is some non-decreasing   function and $\big\|\omX\-\omX' \big\|_t \df \Sup{s \in [0,t]} \big|\omX(s)\-\omX'  (s) \big| $.
  Under   conditions \eqref{coeff_cond1} and \eqref{coeff_cond2},   controlled  SDEs with coefficients $(b,\si)$   are well-posed
   (see e.g.  Theorem V.7 of \cite{SIDE_Protter}):

 \begin{prop} \label{prop_122021}
  Let $ (\O,\cF,  P  )$ be  a  probability space.
  Given $t \ins [0,\infty)$, let $\{B^t_s\}_{s \in  [t,\infty)}$ be a $d-$dimensional Brownian motion   on $ (\O,\cF,  P  )$
  and let $\mu \= \{\mu_s\}_{s \in [t,\infty)}$ be a $\hU-$valued, $\bF^{B^t,P}-$progressively measurable process.
  For any $ \bx \ins \OmX $,   the SDE with the open-loop control $\mu$
   \bea \label{121621_11}
  X_s  \= \bx(t) + \int_t^s b (r,  X_{r \land \cd},\mu_r) dr \+ \int_t^s  \si (r,  X_{r \land \cd},\mu_r) d B^t_r , \; ~ \; \fa s \in [t,\infty)
  \; \hb{ with   initial condition }  X\big|_{[0,t]} \= \bx|_{[0,t]}
 \eea
  admits a unique strong solution $ X^{t,\bx,\mu} \= \{ X^{t,\bx,\mu}_s\}_{s \in  [0,\infty)}$
  \big(i.e., $X^{t,\bx,\mu}$ is an   $\{\cF^{B^t,P}_{s \vee t}\}_{s \in  [0,\infty) }-$adapted continuous process satisfying \eqref{121621_11}
  and $P\big\{ X^{t,\bx,\mu}_s \= \wt{X}^{t,\bx,\mu}_s , \fa s \ins [0,\infty) \big\} \= 1$
  if $ \big\{ \wt{X}^{t,\bx,\mu}_s \big\}_{s \in  [0,\infty)}$ is another  $\{\cF^{B^t,P}_{s \vee t}\}_{s \in  [0,\infty) }-$adapted continuous process satisfying \eqref{121621_11}\big).

 \end{prop}

 Let $\sH_o $ collect all $(-\infty,\infty]-$valued  Borel-measurable functions $\psi$ on $(0,\infty) \ti \OmX \ti \hU $ such that
  for any $(t,\bx) \ins [0,\infty) \ti \OmX$
  and any $\hU-$valued  $ \bF^{W^t} -$predictable process  $ \mu^o \=   \{\mu^o_s\}_{s \in [t,\infty)} $, 
 one has  $ E_{P_0} \big[ \int_t^\infty \n  \psi^- (r, \,   X^{t,\bx,\mu^o}_{r \land \cd},   \mu^o_r )  dr\big]   \< \infty$,
 where $    \{   X^{t,\bx,\mu^o}_s\}_{s \in  [0,\infty)}$
 is  the unique strong solution of \eqref{121621_11} on $  (\O,\cF,  P   ) \= \big( \O_0,\sB(\O_0),   P_0   \big) $ with $  (B^t, \mu ) \= \big( W^t  , \mu^o \big) $.

 Moreover, we   take the conventions $ \inf \es \df \infty$, $ \sup \es \df - \infty$   and $(+\infty)\+(-\infty) 
 \= -\infty$.
  In particular, on a measure space  $(\O,\cF,\fm) $,
   one can  define the integral  $\int_\O  \xi \, d \fm
   \df \int_\O  \xi^+ \, d \fm  \- \int_\O  \xi^- \, d \fm  $ for  any  $[-\infty,\infty]-$valued $\cF-$measurable random variable $\xi$  on  $\O$.


\subsection{Review of Martingale-Problem Formulation of Controlled SDEs}

 In this subsection, we consider  a general measurable space $(\O, \cF)$.
 Let $  \{B_s\}_{s \in [0,\infty)}$ be an $\hR^d-$valued continuous process on $\O$ with $B_0 \= \bz$
 and let $X  \= \{X_s\}_{s \in [0,\infty)}$ be an $\hR^l-$valued continuous  process on $\O$  such that
 $(B_s,X_s)$ is $\cF-$measurable for each $s \ins [0,\infty)$.

 Given $(t,\bx) \ins [0,\infty) \ti \OmX $,  let $P$ be a probability measure on $(\O,\cF)$ such that
  $P \big\{X_s\=\bx(s), \fa s \ins [0,t] \big\} \= 1$.
  We set  $B^t_s \df B_s \- B_t$, $\fa s \ins [t,\infty)$
  and let   $\mu \= \{\mu_s\}_{s \in [t,\infty)}$ be a $\hU-$valued, $\bF^{B^t,P}-$progressively measurable process.
 Define  filtration $\bF^t  \= \{\cF^t_s\}_{s \in [t,\infty)}$ by $\cF^t_s \df
 \cF^{B^t}_s \ve \cF^X_s \=   \si(B^t_r; r \ins [t,s])   \ve \si (X_r; r \ins [0,s] )   $ 
 and filtration $\bF^{t,P}  \= \{\cF^{t,P}_s\}_{s \in [t,\infty)}$ by $\cF^{t,P}_s \df
 \cF^{B^t,P}_s \ve \cF^X_s \= \si \big( \cF^{B^t}_s \cp \sN_P \big( \cF^{B^t}_\infty \big)\big) \ve \cF^X_s  $.  
 For any $\vf \ins C^2(\hR^{d+l})$,  we  define
 \beas    
   \q   M^{t,\mu}_s(\vf)   \df   \vf \big(B^t_s  , X_s \big)
    \- \n \int_t^s     \ol{b}   ( r, X_{r \land \cd}, \mu_r) \n \cd \n D \vf \big( B^t_r  , X_r \big) dr
    \-   \frac12 \n \int_t^s     \ol{\si} \, \ol{\si}^T   ( r,  X_{r \land \cd}, \mu_r) \n : \n D^2 \vf  ( B^t_r, X_r  )   dr  ,
    \q  \fa s \ins [t,\infty),
 \eeas
   where  $ \dis \ol{b}  (r,\omX,u) \df \binom{0}{  b(r,\omX,u)} \n \ins \hR^{d+l} $, $ \dis  \ol{\si}  (r,\omX,u) \df \binom{ I_{d \times d}}{   \si(r,\omX,u)} \n \ins \hR^{ (d+l)  \times d  } $,   $ \fa  (r,\omX,u) \ins   (0,\infty) \ti  \OmX \ti \hU $.
      Clearly, $\big\{  M^{t,\mu}_s(\vf) \big\}_{s \in [t,\infty)}$ is an $\bF^{t,P}-$adapted continuous process
      and it is even   $\bF^t-$adapted 
      if the control process $ \mu $ is only $\bF^{B^t}-$progressively measurable.
   For any $n \ins \hN   $ and $\fra \ins \hR^{d+l}$,     set $\tau^t_n (\fra) \df  \inf\big\{s \ins [t,\infty) \n : |(B^t_s,X_s)\-\fra|   \gs n  \big\} \ld (t\+n) $, which  is an $\bF^t-$stopping time. In particular, we   denote $\tau^t_n(\bz) $ by $\tau^t_n$.


  In virtue of \cite{Stroock_Varadhan},  we have the following  martingale-problem formulation of controlled  SDEs  with coefficients $(b,\si)$ on $\O$.

\begin{prop}  \label{prop_MPF1}
 Under the probabilistic setup of this subsection,
 the process $ \big\{ M^{t,\mu}_{s \land \tau^t_n(\fra) } (\vf) \big\}_{s \in [t,\infty)} $ is   bounded   
  for any $(\vf,n,\fra) \ins  C^2(\hR^{d+l}) \ti \hN \ti \hR^{d+l} $ and
  the following statements are equivalent  on  $(\O, \cF, P)$:

\no \(i\) The process $B^t$ is a   Brownian motion 
    and $P\{ X_s \=   X^{t,\bx,\mu}_s,   \fa s \ins [0,\infty)\}=1$, where  $ \big\{ X^{t,\bx,\mu}_s \big\}_{s \in [0,\infty)}$
  is the unique $ \big\{ \cF^{B^t,P}_{s \vee t} \big\}_{s \in [0,\infty)}   -$adapted continuous process solving   SDE \eqref{121621_11}.

\no \(ii\)   $ \big\{ M^{t,\mu}_{s \land \tau^t_n (\fra) } (\vf) \big\}_{s \in [t,\infty)} $  is a bounded  $\bF^{t,P}-$martingale for any $(\vf,n,\fra) \ins  C^2(\hR^{d+l}) \ti \hN \ti \hR^{d+l}$.

\no \(iii\)   $ \big\{ M^{t,\mu}_{s \land \tau^t_n } (\vf) \big\}_{s \in [t,\infty)} $  is a bounded $\bF^{t,P}-$martingale for any $(\vf,n) \ins \fC(\hR^{d+l}) \ti \hN $.

Moreover, if the control process $ \mu $ is   $\bF^{B^t}-$progressively measurable,   the $\bF^{t,P}-$martingales  mentioned  in   \(ii\) and \(iii\) should be $\bF^t-$martingales.

\end{prop}

 We   have the following consequence of Proposition \ref{prop_MPF1}.

 \begin{prop} \label{prop_122021b}

  Let $(\O, \cF, P)$ be   a    probability space. Given  $ t \ins [0,\infty)$,
 let $  \{B_s\}_{s \in [0,\infty)}$ be an $\hR^d-$valued continuous process on $\O$ with $B_0 \= \bz$ such that
 the process $B^t$
 is a  Brownian motion on $(\O, \cF, P)$,
 and let $\mu \= \{\mu_s\}_{s \in [t,\infty)}$ be a $\hU-$valued, $\bF^{B^t,P}-$progressively measurable  process.

 Let  $(t,\bw) \ins [0,\infty)   \ti \O_0  $  and define
 $  B^{t,\bw}_s  (\o) \df   \bw(s \ld t)  \+  B^t_{s \vee t} (\o)     $, $ \fa (s,\o) \ins [0,\infty) \ti \O $.
 There exists   a   $\hU-$valued, $ \bF^{W^t} -$predictable process $ \mu^o \=   \{ \mu^o_s \}_{s \in [t,\infty)} $
 on $\O_0$ and an $\cN_\mu  \ins \sN_P(\cF^{B^t}_\infty)$ such that for any $\o \ins \cN^c_\mu$,
  $  \mu_s (\o) \=  \mu^o_s \big(B^{t,\bw}(\o)\big)  $ for a.e. $s \ins (t,\infty)   $.
  It also holds for any $\bx \ins \OmX$ and $\psi \ins \sH_o$ that
   $P \big\{ X^{t,\bx,\mu}_s \=   X^{t,\bx,\mu^o}_s ( B^{t,\bw} )     , \, \fa s \ins [0,\infty) \big\} \= 1$
  and   $ E_P  \big[ \int_t^\infty \n  \psi^- (r, X^{t,\bx,\mu}_{r \land \cd}, \mu_r )   dr\big] \= E_{P_0} \big[ \int_t^\infty \n \psi^-(r,X^{t,\bx,\mu^o}_{r \land \cd},\mu^o_r ) dr \big]   \< \infty$.

 \end{prop}

\section{Stochastic Control/Stopping Problem with Expectation Constraints}

\label{sec_genprob}

  Let $(\cQ, \cF, \fp)$ be   a    probability space
  equipped with a $d-$dimensional standard  Brownian motion $  \{\cB_s\}_{s \in [0,\infty)}$.

  Let $t \ins [0,\infty)$. We set
  $\cB^t_s  \= \cB_s \- \cB_t $, $\fa s \ins [t,\infty)$,
  which  is also a   Brownian motion on $(\cQ, \cF, \fp)$.
  Let $\cU_t$ collect all $\hU-$valued, $\bF^{\cB^t,\fp}-$progressively measurable processes $\mu \= \{\mu_s\}_{s \in [t,\infty)}$
  and let $\cS_t$ denote the set of all $[t,\infty]-$valued  $\bF^{\cB^t,\fp}-$stopping times.
  For any $(\bx,\mu) \in   \OmX \ti \cU_t$,
  Proposition \ref{prop_122021} shows that the SDE with the open-loop control $\mu$
 \bea \label{FSDE1}
 \cX_s \= \bx(t) + \int_t^s b \big( r, \cX_{r \land \cd}, \mu_r \big) dr
 \+ \int_t^s  \si \big( r, \cX_{r \land \cd}, \mu_r \big) d \cB_r , ~ \fa  s \in [t,\infty)
 \; \hb{ with   initial condition }  \cX\big|_{[0,t]} \= \bx|_{[0,t]}
 \eea
 admits   a unique strong  solution $ \cX^{t,\bx,\mu} \= \big\{\cX^{t,\bx,\mu}_s\big\}_{s \in [0,\infty)}$ on $ \big(\cQ, \cF, \bF^{\cB^t,\fp}, \fp\big) $
 \big(i.e., $\cX^{t,\bx,\mu}$ is the unique $\big\{\cF^{\cB^t,\fp}_{s \vee t}\big\}_{s \in [0,\infty)}-$ adapted continuous process solving SDE  \eqref{FSDE1}\big).

 Let $f \ins \sH_o $, $ \{ g_i, h_i \}_{i \in \hN} \sb \sH_o$
 and let $\pi \n : [0,\infty) \ti \OmX \mto (-\infty,\infty]$ be a Borel-measurable function bounded from below
 by some $c_\pi \ins (-\infty,0)  $.

 Given a historical path $\bx |_{[0,t]}$, the state of the game then  evolves along process $\big\{\cX^{t,\bx,\mu}\big\}_{s \in [t,\infty)}$
 if  the player chooses a  control process $\mu \ins \cU_t$.
 The player  also determines an exercise time $\tau \ins \cS_t$   to cease the game,
 at which she will  receive
 an accumulative reward $\int_t^\tau   \n f \big( r, \cX^{t,\bx,\mu}_{r \land \cd},\mu_r   \big) \, dr $ plus
 a terminal reward $   \pi \big(  \tau  , \cX^{t,\bx,\mu}_{ \tau  \land \cd}  \big) $
  (both random rewards can take negative values).
 The player intends to maximize   the expectation of her    total wealth,
  but her choice of $(\mu,\tau)$  is subject to a series of expectation constraints
 \bea \label{111920_11}
  E_\fp \Big[ \int_t^\tau  g_i ( r,\cX^{t,\bx,\mu}_{r \land \cd}, \mu_r ) dr  \Big] \ls y_i,   \q
  E_\fp \Big[ \int_t^\tau  h_i ( r,\cX^{t,\bx,\mu}_{r \land \cd}, \mu_r ) dr  \Big] \= z_i,
  \q \fa i \ins \hN
 \eea
  for some  $(y,z) \= \big(\{y_i\}_{i \in \hN}, \{z_i\}_{i \in \hN}\big)  \ins \Re \ti \Re$.
One can regard each $\int_t^\tau  g_i ( r,\cX^{t,\bx,\mu}_{r \land \cd}, \mu_r ) dr$ or $\int_t^\tau h_i ( r,\cX^{t,\bx,\mu}_{r \land \cd}, \mu_r ) dr$ as   certain accumulative cost.
 So the value of this stochastic control/stopping problem with expectation constraints (SCEC for short) is
\bea  \label{081820_11}
 V (t,\bx,y,z)   \df    \Sup{(\mu,\tau) \in \cC_{t,\bx}(y,z) } E_\fp \Big[ \int_t^\tau f \big( r, \cX^{t,\bx,\mu}_{r \land \cd},\mu_r  \big) dr
 \+ \b1_{\{\tau < \infty\}} \pi \big(  \tau  , \cX^{t,\bx,\mu}_{ \tau  \land \cd}  \big)  \Big] ,
\eea
where $ \cC_{t,\bx}(y,z) \df \big\{ (\mu,\tau) \ins \cU_t \ti \cS_t \n : E_\fp \big[ \int_t^\tau  g_i ( r,\cX^{t,\bx,\mu}_{r \land \cd},\mu_r ) dr  \big] \ls y_i,   \, E_\fp \big[ \int_t^\tau  h_i ( r,   \cX^{t,\bx,\mu}_{r \land \cd},\mu_r ) dr  \big] \= z_i, \, \fa i \ins \hN \big\} $.

\if{0}

\begin{rem} \label{rem_112220}
 Let $(t,\bx) \ins [0,\infty) \ti \OmX$.

  \no 1\) \(finitely many constraints\) For   $i \ins \hN$,  the constraint $E_\fp \big[ \int_t^\tau  g_i ( r,\cX^{t,\bx,\mu}_{r \land \cd},\mu_r ) dr  \big] \ls y_i$ holds for any $(\mu,\tau) \ins \cU_t \ti \cS_t$ 
  if  $ y_i \= \infty $, and the constraint $E_\fp \big[ \int_t^\tau  h_i ( r,\cX^{t,\bx,\mu}_{r \land \cd},\mu_r ) dr  \big] \= z_i$
  holds for any $(\mu,\tau) \ins \cU_t \ti \cS_t$ 
  if $\big(h_i(\cd,\cd,\cd),z_i\big) \= (0,0)$.

  \no 1a\) If we take  $\big(y_i,h_i(\cd,\cd,\cd),z_i\big) \= (\infty,0,0)$, $\fa i \ins \hN$, there is no expectation constraint at all.

  \no 1b\) If one takes $y_i \= \infty$, $\fa i \gs 2$ and
   $\big(h_i(\cd,\cd,\cd),z_i\big) \= (0,0)$, $\fa i \ins \hN$,   \eqref{111920_11} reduces  to
  a single constraint $ E_\fp \big[ \int_t^\tau  g_1 ( r,\cX^{t,\bx,\mu}_{r \land \cd}, \\ \mu_r ) dr  \big]   \ls y_1 $.
  In addition, if   $y_1 \gs 0$,
  then $(\mu,t) \ins  \cC_{t,\bx}(y,\bz) $ for any $\mu \ins \cU_t $.

   \no 1c\) If one takes $y_i \= \infty$, $\fa i \ins \hN $ and
   $\big(h_i(\cd,\cd,\cd),z_i\big) \= (0,0)$, $\fa i \gs 2$,   \eqref{111920_11} degenerates to  
  $ E_\fp \big[ \int_t^\tau  h_1 ( r,    \cX^{t,\bx,\mu}_{r \land \cd},\mu_r ) dr  \big]    \= z_1 $.

    \no 1d\) If we take  $\big(y_i,h_i(\cd,\cd,\cd),z_i\big) \= (\infty,0,0)$, $\fa i \gs 2$, \eqref{111920_11} becomes a couple of constraints
  $ E_\fp \big[ \int_t^\tau  g_1 ( r,\cX^{t,\bx,\mu}_{r \land \cd},\mu_r ) dr  \big] \ls y_1 $ and $ E_\fp \big[ \int_t^\tau  h_1 ( r,\cX^{t,\bx,\mu}_{r \land \cd},\mu_r ) dr  \big] \= z_1$.

  \no 1e\) If we take  $g_2 \= -g_1$, $ y_2 \gs - y_1$; $y_i \= \infty$, $\fa i \gs 3$ and
   $\big(h_i(\cd,\cd,\cd),z_i\big) \= (0,0)$, $\fa i \ins \hN$,   \eqref{111920_11} becomes a range constraint  $ - y_2 \ls E_\fp \big[ \int_t^\tau  g_1 ( r,   \cX^{t,\bx,\mu}_{r \land \cd},\mu_r ) dr  \big]   \ls y_1 $.

\no 2\) \(moment constraints\) Let $i \ins \hN$, $ a  \ins (0,\infty)$ and $q \ins [1,\infty) $.
If   $g_i(s,\bx,u)  \= a q s^{q-1}   $, $\fa (s,\bx,u) \ins (0,\infty) \ti \OmX  \ti \hU $
\big(resp.  $h_i(s,\bx,u)  \= a q s^{q-1}   $, $\fa (s,\bx,u) \ins (0,\infty) \ti \OmX \ti \hU $\big),
then the expectation constraint  $E_\fp \big[ \int_t^\tau  g_i ( r,\cX^{t,\bx,\mu}_{r \land \cd},\mu_r ) dr  \big] \ls y_i$
\big(resp.  $E_\fp \big[ \int_t^\tau  h_i ( r,\cX^{t,\bx,\mu}_{r \land \cd},\mu_r ) dr  \big] \= z_i$\big)
specifies  as a moment constraint   $E_\fp \big[ a  ( \tau^q \- t^q )    \big] \ls y_i$
\big(resp.  $E_\fp \big[ a  (\tau^q \- t^q )   \big] \= z_i$\big).

\if{0}
If   $g_i(s,\bx)  \= a q s^{q-1}   $, $\fa (s,\bx) \ins (0,\infty) \ti \OmX  $,
then the expectation constraint  $E_\fp \big[ \int_t^\tau  g_i ( r,\cX^{t,\bx,\mu}_{r \land \cd} ) dr  \big] \ls y_i$
specify as a moment constraint   $E_\fp \big[ a  ( \tau^q \- t^q )    \big] \ls y_i$;
If    $h_i(s,\bx)  \= a q s^{q-1}   $, $\fa (s,\bx) \ins (0,\infty) \ti \OmX  $,
then the expectation constraint   $E_\fp \big[ \int_t^\tau  h_i ( r,\cX^{t,\bx,\mu}_{r \land \cd} ) dr  \big] \= z_i$
specify as a moment constraint    $E_\fp \big[ a  (\tau^q \- t^q )   \big] \= z_i$.
\fi

\end{rem}

\fi

 To study the measurability of   value function $V$  and  derive a  dynamic programming principle for $V$ without imposing any   continuity condition on 
  functions $f$, $\pi$, $g_i$'s and $h_i$'s in time and state variables, we follow \cite{EHJ_1987}'s approach
  to embed the controls, the stopping rules as well as the Brownian/state information into an enlarged canonical space
  via  a mapping $   \o \mto \big(\cB_\cd(\o),\mu_\cd(\o),\cX^{t,\bx,\mu}_\cd(\o),\tau(\o)\big)  $
  and consider their joint law 
  as a new type of controls.

\section{Weak Formulation}
\label{sec_weak_form}

  In this section, we study the stochastic control/stopping problem with expectation constraints in a weak formulation or over an  enlarged canonical space
 \beas
 \oO \df \O_0 \ti   \hJ  \ti \OmX  \ti \hT  .
 \eeas
 As $ \big(\hJ,\fT_\sharp(\hJ)\big) $ is a Borel space by Lemma  \ref{lem_Nov25_03},
 $\oO$ is also a Borel space  under the product topology.
  \if{0}

  We see from \eqref{J04_01} that  the  induced  topology $\fT_\sharp(\hJ)$
  has a countable subbase and thus has a countable base \big(i.e., $ \big(\hJ,\fT_\sharp(\hJ)\big) $ is a second-countable space\big).  As $\O_0$ and $\OmX$ are two complete separable metric spaces,
  Lemma \ref{lem_prod_Bsf} renders that
  \bea \label{051321_11}
  \sB\big(\oO\big) \= \sB (\O_0 \ti \hJ \ti \OmX \ti \hT) \=   \sB (\O_0) \oti \sB(\hJ) \oti \sB(\OmX) \oti \sB(\hT)  .
  \eea

  \fi
   Let $\fP(\oO)$ be  the space  of all probability measures      on $\big(\oO, \sB(\oO) \big)$
 equipped with the topology of weak convergence, which is also a Borel space (see e.g. Corollary 7.25.1 of \cite{Bertsekas_Shreve_1978}).
 For any $\oP \ins \fP(\oO)$,   set $\sB_\oP(\oO) \df \si \big( \sB(\oO) \cp \sN_\oP  (\sB(\oO) ) \big)$.
 We define the canonical coordinates on $\oO$ by
\beas
\big(\oW_{\n s} (\oo),\oU_{\n s} (\oo),\oX_{\n s} (\oo)\big) \df \big(\o_0(s) ,   \fu(s) ,   \omX(s) \big) ,  \q   s \ins [0,\infty)
 \q \hb{and} \q \oT(\oo) \df   \ft , \q \fa \oo \= \big(\o_0,\fu,\omX,\ft\big) \ins \oO ,
\eeas
 in which one can regard  $\oW$ as a canonical coordinate for Brownian motion,
 $\oU$ as a canonical coordinate for the control process,
 $\oX$ as a canonical coordinate for the state process, and $\oT$ as a canonical coordinate for stopping rules.
   Given $t \ins [0,\infty)$, we  define
   \beas
   \oW^t_{\n s} (\oo) \df \oW_{\n s} (\oo) \- \oW_{\n t} (\oo) \aand
    \ol{\U}^t_s (\oo) \df 
    \int_t^s e^{-r} \sI \big(  \oU_r(\oo)\big) dr \ins [0,1)  , \q \fa (s,\oo) \ins [t,\infty) \ti \oO  .
  \eeas

  \if{0}

 By Lemma \ref{lem_WUX_basic},
 the random variable $\oW^t   (\oo) \=   W^t(\oW(\oo)) \ins \O_0 $, $ \fa \oo \ins \oO $  is $\sB(\oO)   / \sB(\O_0)  -$measurable.

Applying Lemma \ref{lem_M29_01} (2) with  $\psi(r,\fx,\nxi,u) \df e^{-r} \sI(u) $, $\fa (r,\fx,\nxi,u) \ins (0,\infty) \ti \OmX \ti \hR^{d+l} \ti \hU $  shows that the mapping   $   \Psi   (t,\fs,\fu) \df  \Psi (t,\fs,\bz,\bz,\fu)
   \=  \int_t^{t+\fs} \psi \big(r, \fl_2(r,\bz), 0, 0 ,   \fu(r)\big) dr
   \=  \int_t^{t+\fs} e^{-r} \sI \big(\fu(r)\big) dr    $,
  $  (t,\fs,  \fu) \ins [0,\infty) \ti [0,\infty) \ti   \hJ$ is $\sB[0,\infty) \oti \sB[0,\infty) \oti   \sB(\hJ)-$measurable.
As the random variable $\oU$ on $\oO$ is $\sB(\oO)/\sB(\hJ)-$measurable,
$ \ol{\U}^t_{t+\fs}(\oo) \=  \Psi   \big(t,\fs,\oU(\oo)\big) $, $(t,\fs,\oo) \ins [0,\infty) \ti [0,\infty) \ti \oO$
is $\sB[0,\infty) \ti \sB[0,\infty) \ti \sB(\oO)-$measurable.

  \fi

 The weak formulation of the SCEC 
 relies on the following probability classes of $\fP\big(\oO\big)$.

\begin{deff} \label{def_ocP}
 For any $\tx$,
 let $\ocP_{t,\bx}$ be the collection of  all probability measures $ \oP \ins \fP\big(\oO\big) $ satisfying:

 \no \(D1\) There exists   a   $\hU-$valued, $ \bF^{W^t} -$predictable process
 $ \wh{\mu} \= \{\wh{\mu}_s\}_{s \in [t,\infty)} $   on $\O_0$ such that
 $ \oP \big\{  \oU_s  \=  \ol{\mu}_s   \hb{ for a.e. }  s \ins (t,\infty) \big\} \= 1$,
 where $\ol{\mu}_s \df  \wh{\mu}_s  (\oW   )$, $\fa s \ins [t,\infty)$.

\no \(D2\) The process $\oW^t$ is a  $d-$dimensional 
Brownian motion on $\big(\oO , \sB(\oO) , \oP\big)$.

\no \(D3\)  $  \oP\big\{ \oX_s \= \osX^{t,\bx,\ol{\mu}}_s, ~ \fa s \ins [0,\infty) \big\} \= 1$,
  where  $ \big\{\osX^{t,\bx,\ol{\mu}}_s\big\}_{s \in [0,\infty)}$
  is an $\big\{\cF^{\oW^t,\oP}_{s \vee t}\big\}_{s \in [0,\infty)}-$adapted continuous process
  that   uniquely solves   the following SDE with the open-loop control $\ol{\mu}$ on $\big(\oO , \sB\big(\oO\big) , \oP\big) \n : $
 \bea \label{Ju01_01}
 \osX_s = \bx(t) + \int_t^s b \big( r, \osX_{r \land \cd}, \ol{\mu}_r  \big)dr \+ \int_t^s \si \big( r, \osX_{r \land \cd}, \ol{\mu}_r  \big) d \oW_r, ~ \fa  s \ins [t,\infty) \hb{  with initial condition $\osX \big|_{[0,t]} \=   \bx \big|_{[0,t]} $. }
 \eea

\no \(D4\)  There exists   a $[t,\infty]-$valued  $ \bF^{W^t,P_0} -$stopping time $\wh{\tau}$   on $\O_0$ such that
 $ \oP \big\{  \oT  \=    \wh{\tau} (\oW   )   \big\} \= 1$.

 \end{deff}

   \if{0}

 \begin{rem} \label{rem_082020_11}
  \(1a\)  Let $ \wh{\mu} \= \{\wh{\mu}_s\}_{s \in [t,\infty)} $ be   a   $\hU-$valued, $ \bF^{W^t,P_0} -$predictable process
    on $\O_0$ such that $ \oP \big\{  \oU_s  \=  \ol{\mu}_s   \hb{ for a.e. }  s \ins (t,\infty) \big\} \= 1$,
 where $\ol{\mu}_s \df  \wh{\mu}_s  (\oW   )$, $\fa s \ins [t,\infty)$. Under \(D2\), applying Lemma \ref{lem_predict}   with $(\O,\cF,P,B) \= \big(\oO,\sB\big(\oO\big),\oP,\oW \big)$
 shows  that
 $\ol{\mu}_s (\oo) \df   \wh{\mu}_s \big(\oW (\oo)\big)   $, $\fa (s,\oo) \ins [t,\infty) \ti \oO$ is a $\hU-$valued,  $\bF^{\oW^t,\oP}-$predictable process on $\oO$.
  By the continuity of   $\sI$,
   $  \big\{ \sI  ( \ol{\mu}_s  ) \big\}_{s \in [t,\infty)} $ is an $[0,1]-$valued, $\bF^{\oW^t,\oP}-$predictable process
 and thus  $ \ol{\Phi}^t_s(\oo) \df \int_t^s  e^{-r}  \sI  \big( \ol{\mu}_r (\oo) \big) dr   \ins [0,1]
 $, $ \fa (s,\oo) \ins [0,\infty) \ti \oO $    is an   $ \bF^{\oW^t,\oP}  -$adapted continuous process.
    Then one can deduce from the injective of $\sI$ that
 \bea
 \q && \hspace{-1.4cm} \big\{\oo \ins \oO \n : \oU_s(\oo) \= \ol{\mu}_s (\oo) \hb{ for a.e. }  s \ins (t,\infty) \big\}
  \= \big\{\oo \ins \oO \n :   \sI (\oU_s(\oo)) \=   \sI \big(\ol{\mu}_s (\oo)\big) \hb{ for a.e. } s \ins (t,\infty) \big\} \nonumber \\
 &&   \= \bigg\{\oo \ins \oO \n : \ol{\U}^t_s(\oo) \= \int_t^s e^{-r} \sI (\oU_r(\oo)) dr  \= \int_t^s e^{-r} \sI \big(\ol{\mu}_r (\oo)\big)  dr \= \ol{\Phi}^t_s(\oo) , \, \fa s \ins (t,\infty) \bigg\} \nonumber \\
 &&  \= \underset{s \in \hQ \in (t,\infty)}{\cap} \big\{\oo \ins \oO \n :  \ol{\U}^t_s(\oo) \= \ol{\Phi}^t_s(\oo)  \big\}
  \ins \cF^{\oW^t,\ol{\U}^t,\oP}_\infty .   \label{041421_11}
  \eea

\no   \(1b\) Let $ \wh{\mu} \= \{\wh{\mu}_s\}_{s \in [t,\infty)} $ be   a   $\hU-$valued, $ \bF^{W^t} -$predictable process
    on $\O_0$ such that $ \oP \big\{  \oU_s  \=  \ol{\mu}_s   \hb{ for a.e. }  s \ins (t,\infty) \big\} \= 1$,
 where $\ol{\mu}_s \df  \wh{\mu}_s  (\oW   )$, $\fa s \ins [t,\infty)$.
    Lemma \ref{lem_predict0} implies that
 $\ol{\mu}_s (\oo) \df   \wh{\mu}_s \big(\oW (\oo)\big)   $, $\fa (s,\oo) \ins [t,\infty) \ti \oO$ is a $\hU-$valued,  $\bF^{\oW^t}-$predictable process on $\oO$.
  By the continuity of mapping $\sI$,
   $  \big\{ \sI  ( \ol{\mu}_s  ) \big\}_{s \in [t,\infty)} $ is an $[0,1]-$valued, $\bF^{\oW^t }-$predictable process
 and thus  $ \ol{\Phi}^t_s(\oo) \df \int_t^s  e^{-r}  \sI  \big( \ol{\mu}_r (\oo) \big) dr   \ins [0,1]
 $, $ \fa (s,\oo) \ins [0,\infty) \ti \oO $    is an   $ \bF^{\oW^t }  -$adapted continuous process.
    Then one can deduce from the injective of $\sI$ that
 \beas
 \q && \hspace{-1.4cm} \big\{\oo \ins \oO \n : \oU_s(\oo) \= \ol{\mu}_s (\oo) \hb{ for a.e. }  s \ins (t,\infty) \big\}
  \= \big\{\oo \ins \oO \n :   \sI (\oU_s(\oo)) \=   \sI \big(\ol{\mu}_s (\oo)\big) \hb{ for a.e. } s \ins (t,\infty) \big\} \nonumber \\
 &&   \= \bigg\{\oo \ins \oO \n : \ol{\U}^t_s(\oo) \= \int_t^s e^{-r} \sI (\oU_r(\oo)) dr  \= \int_t^s e^{-r} \sI \big(\ol{\mu}_r (\oo)\big)  dr \= \ol{\Phi}^t_s(\oo) , \, \fa s \ins (t,\infty) \bigg\} \nonumber \\
 &&  \= \underset{s \in \hQ \in (t,\infty)}{\cap} \big\{\oo \ins \oO \n :  \ol{\U}^t_s(\oo) \= \ol{\Phi}^t_s(\oo)  \big\}
  \ins \cF^{\oW^t,\ol{\U}^t }_\infty .
  \eeas

\no \(2\)
 In \(D1\) of Definition \ref{def_ocP} ,
   we can assume without loss of generality  that  $\wh{\mu}_\cd(\o_0) \ins \hJ$   for any $\o_0 \ins \O_0$:
   As an $ \bF^{W^t,P_0} -$predictable process, $ \big\{\wh{\mu}(s,\o_0)\big\}_{(s,\o_0) \in [t,\infty) \times \O_0}$
   is $\sB[t,\infty) \oti \cF^{W^t,P_0}_\infty-$measurable. By Fubini Theorem,
   \bea \label{050421_25}
   \hb{$ \wh{\cN} \df \big\{\o_0 \ins \O_0 \n : \wh{\mu}_\cd  (\o_0) \n \notin \n  \hJ \big\} $ is a $\cF^{W^t,P_0}_\infty-$measurable
   set with zero $P_0-$measure  or $ \wh{\cN} \ins \sN_{P_0}\big(\cF^{W^t}_\infty\big)$.}
   \eea
    Then $ \breve{\mu}_s \df \wh{\mu}_s    \b1_{   \wh{\cN}^c  }  \+ u_0 \b1_{   \wh{\cN}  } \ins \hU  $,
   $  s  \ins  [t,\infty)$    is an  $ \bF^{W^t,P_0} -$predictable process with all paths in $\hJ$.
 Under \(D2\), since an application of Lemma \ref{lem_122921_11} (ii) with $t_0 \= t$,  $(\O_1, \cF_1, P_1,B^1)   \= \big(\oO,\sB\big(\oO\big),\oP,\oW \big) $, $(\O_2, \cF_2, P_2,B^2) \= \big(\O_0, \sB(\O_0) , P_0 , W\big) $ and $\Phi \= \oW$  yields that
    $  \oW^{-1} (\wh{\cN}) \ins \sN_\oP\big(\cF^{\oW^t}_\infty\big)   $, we see from \eqref{041421_11} that
 \beas
 && \hspace{-1.5cm} \big\{\oo \ins \oO \n : \oU_s(\oo) \= \breve{\mu}_s \big(\oW (\oo)\big) \hb{ for a.e. }  s \ins (t,\infty) \big\} \Cp \oW^{-1} \big( \wh{\cN}^c \big) \= \big\{\oo \ins \oO \n : \oU_s(\oo) \= \wh{\mu}_s \big(\oW (\oo)\big) \hb{ for a.e. }  s \ins (t,\infty) \big\} \Cp \oW^{-1} \big( \wh{\cN}^c \big)
 \eeas
  is a $\cF^{\oW^t,\ol{\U}^t,\oP}_\infty-$measurable set whose $\oP-$measure   equals   $1$.
   It follows that $\big\{\oo \ins \oO \n : \oU_s(\oo) \= \breve{\mu}_s \big(\oW (\oo)\big) \hb{ for a.e. }  s \ins (t,\infty) \big\}$
   is also a $\cF^{\oW^t,\ol{\U}^t,\oP}_\infty-$measurable set whose $\oP-$measure   equals   $1$.

  \no \(3\)   Since $\big\{ \oX_s \= \osX^{t,\bx,\ol{\mu}}_s \big\} \ins \cF^\oX_s \ve \cF^{\oW^t,\oP}_{s \vee t}   $ for any
  $ s \ins [0,\infty)$, one has
  $  \big\{ \oX_s   \=   \osX^{t,\bx,\ol{\mu}}_s,   \fa s \ins [0,\infty) \big\}
  \=   \ccap{s \in  \hQ \cap [0,\infty)}{} \big\{    \oX_s \= \osX^{t,\bx,\ol{\mu}}_s   \big\}
  \ins  \cF^\oX_\infty \ve \cF^{\oW^t,\oP}_\infty   $.

 \end{rem}

    \fi

 Let $ t  \ins [0,\infty)   $.
  For any $s \in [t,\infty)$,   define
 $\ocF^t_{\n s} \df \cF^{\oW^t}_s \ve \cF^\oX_s \= \si \big(\oW^t_{\n r}; r \ins [t,s]\big) \ve \si \big(\oX_r;r \ins [0,s]\big)$, which is countably generated by
 $ \big\{\oX^{-1}_r (\cO) \n: r \ins \hQ \Cp [0,t], \cO \ins \sO(\hR^l) \big\} \cp \big\{ (\oW^t_r,\oX_r)^{^{-1}} (\cO' ) \n: r \ins \hQ \Cp (t,s], \cO'  \ins \sO(\hR^{d+l}) \big\}   $.
  We denote the filtration  $  \big\{\ocF^t_{\n s}   \big\}_{s \in [t,\infty)}$ by $\obF^t$.
  Let $ \wh{\mu} \= \{\wh{\mu}_s\}_{s \in [t,\infty)} $  be a
   $\hU-$valued, $ \bF^{W^t} -$predictable process  on $\O_0$.
  Then  $\ol{\mu}_s (\oo) \df   \wh{\mu}_s \big(\oW (\oo)\big)   $, $\fa (s,\oo) \ins [t,\infty) \ti \oO$ is a $\hU-$valued,  $\bF^{\oW^t }-$predictable process on $\oO$.
 For any  $ (\vf,n,\fra )  \ins   C^2 (\hR^{d+l}) \ti \hN \ti \hR^{d+l} $,
 \beas  
 \q \oM^{t,\ol{\mu}}_{\n s}(\vf)   \df  \vf \big(\oW^t_{\n s}   , \oX_s \big)
    \- \n \int_t^s  \n  \ol{b}  \big( r, \oX_{r \land \cd},\ol{\mu}_r \big) \n \cd \n D \vf \big( \oW^t_{\n r}  , \oX_r \big) dr
    \-   \frac12 \n \int_t^s  \n  \ol{\si} \, \ol{\si}^T  \big( r,  \oX_{r \land \cd},\ol{\mu}_r  \big) \n : \n D^2 \vf  ( \oW^t_{\n r}  , \oX_r  )   dr  ,
    \q \fa s \ins [t,\infty)
 \eeas
 is an $\obF^t-$adapted continuous process  and   $\otau^t_n (\fra) \df  \inf\big\{s \ins [t,\infty) \n : \big|(\oW^t_{\n s} ,   \oX_s) \- \fra \big|    \gs n   \big\} \ld (t\+n) $ is an $\obF^t-$stopping time.
 We will simply denote $ \otau^t_n (\bz) $ by $\otau^t_n  $.

 Let us also define a shifted canonical  process on $\oO$ by $\osW^t_{\n \fs} (\oo) \df \oW_{t+\fs} (\oo) \- \oW_t (\oo) \= \oW^t_{t+\fs} (\oo)$, $\fa (\fs,\oo) \ins [0,\infty) \ti \oO$.  \big(Note: the subscript $\fs \ins [0,\infty)$ of $\osW^t$  is the relative time after $t$ while the subscript $s \ins [t,\infty)$ of $\oW^t $ is the real time.\big)
 Given $\fs \ins [0,\infty]$, one has $\cF^{\osW^t}_\fs \= \si\big(\osW^t_\fr;   \fr \ins [0,\fs] \Cp \hR\big)
\=  \si\big(\oW^t_{t+\fr};   \fr \ins [0,\fs] \Cp \hR\big)  \=   \si\big(\oW^t_r;   r \ins [t,t\+\fs] \Cp \hR\big) \= \cF^{\oW^t}_{t+\fs} $.
In particular, $ \cF^{\osW^t}_\infty \= \cF^{\oW^t}_\infty$. It then holds for any $\oP \ins \fP\big(\oO\big)$ that
$ \cF^{\osW^t,\oP}_\fs \= \si \big( \cF^{\osW^t}_\fs \cp \sN_\oP(\cF^{\osW^t}_\infty )\big)
\= \si \big( \cF^{\oW^t}_{t+\fs} \cp \sN_\oP(\cF^{\oW^t}_\infty )\big) \= \cF^{\oW^t,\oP}_{t+\fs} $.

 According to the martingale-problem formulation of controlled  SDEs (Proposition \ref{prop_MPF1}), we have an alternative description  of
 the probability class  $\ocP_{t,\bx}$:

 \begin{rem} \label{rem_ocP}

 Let $(t,\bx) \ins [0,\infty) \ti \OmX$. In definition \ref{def_ocP} of $\ocP_{t,\bx}$,

\no \(i\) \(D1\)  is equivalent to

 \no  \(D1\,$'$\) There exists   a   $\hU-$valued, $ \bF^W -$predictable process
 $ \ddot{\mu} \= \{\ddot{\mu}_\fs\}_{\fs \in [0,\infty)} $   on $\O_0$ such that
 $ \oP \big\{  \osU^t_\fs   \=  \ddot{\mu}_\fs (\osW^t)   \hb{ for a.e. }  \fs \ins (0,\infty) \big\} \= 1$,
 where $\osU^t_\fs \df \oU_{t+\fs}$, $\fa \fs \ins [0,\infty)$.

\no \(ii\) Under \(D1\)$=$\(D1\,$'$\), \(D2\)+\(D3\) is equivalent to

 \no  \(D2\,$'$\)  $\oP\{\oX_s \= \bx(s), \fa s \ins [0,t]\} \= 1$ and
  $\big\{ \oM^{t,\ol{\mu}}_{s \land \otau^t_n } (\vf) \big\}_{s \in  [t,\infty) } $  is a bounded $ \big( \obF^t , \oP \big)   -$martingale,   $\fa (\vf,n) \ins \fC(\hR^{d+l}) \ti \hN $.

 \no \(iii\)   \(D4\) is equivalent to

    \no   \(D4\,$'$\)  There exists   a $[0,\infty]-$valued   $ \bF^{W,P_0} -$stopping time $\ddot{\tau}$   on $\O_0$ such that
 $ \oP \big\{  \oT  \=  t \+ \ddot{\tau} \big(\osW^t  \big)   \big\} \= 1$.


   \end{rem}

 \begin{rem} \label{rem_ocP2}

 Let $(t,\bx) \ins [0,\infty) \ti \OmX$
 and let $\oP \ins \fP(\oO)$ satisfy  \(D1\)+\(D2\)+\(D3\) of $\ocP_{t,\bx}$.

\no \(1\) For any $\psi \ins \sH_o$, Proposition \ref{prop_122021b} shows that  $ E_\oP  \big[ \int_t^\infty \n  \psi^- (r, \oX_{r \land \cd}, \oU_r )   dr\big] \= E_\oP  \big[ \int_t^\infty \n  \psi^- \big(r, \osX^{t,\bx,\ol{\mu}}_{r \land \cd}, \ol{\mu}_r \big)   dr\big]   \< \infty$.

\no \(2\) Let $ (\vf,n,\fra )  \ins   C^2 (\hR^{d+l}) \ti \hN \ti \hR^{d+l} $. As $\big\{ \oM^{t,\ol{\mu}}_{s \land \otau^t_n (\fra)} (\vf) \big\}_{s \in  [t,\infty) } $  is a bounded $ \big( \obF^t , \oP \big)   -$martingale,
the optional sampling theorem implies that
  for any   two $[t,\infty]-$valued $\obF^t-$stopping times $ \oz_1, \oz_2 $ with $\oz_1 \ls \oz_2$, $\oP-$a.s.,
 \bea    \label{020522_17}
  E_\oP  \Big[  \big( \, \oM^{t,\ol{\mu}}_{\oz_2 \land \otau^t_n(\fra)} (\vf ) \-  \oM^{t,\ol{\mu}}_{\oz_1 \land \otau^t_n(\fra)} (\vf ) \big) \b1_\oA  \Big]
  \=  E_\oP  \Big[  E_\oP  \Big[ \oM^{t,\ol{\mu}}_{\oz_2 \land \otau^t_n(\fra)} (\vf ) \-  \oM^{t,\ol{\mu}}_{\oz_1 \land \otau^t_n(\fra)} (\vf ) \Big| \ocF^t_{ \oz_1 }\Big] \b1_\oA  \Big]
  \= 0 , ~ \;  \fa \oA \ins \ocF^t_{ \oz_1} . \q
  \eea

  \end{rem}

  Let $(t,\bx) \ins [0,\infty) \ti \OmX$, $(y,z) \= \big(\{y_i\}_{i \in \hN}, \{z_i\}_{i \in \hN}\big) \ins \Re \ti \Re$
  and set   $\oR  (t) \df \int_{\oT \land t}^\oT \n  f   (r, \oX_{r \land \cd},\oU_r  ) dr  \+ \b1_{\{\oT < \infty\}} \pi   \big(\oT, \oX_{\oT \land \cd}\big) $.
   Given a historical state path   $\bx|_{[0,t]}$,  the   value  of the stochastic control/stopping problem with    expectation constraints
 \bea \label{111920_14}
 E_\oP \Big[ \int_t^\oT   g_i ( r,\oX_{r \land \cd} , \oU_r ) dr  \Big] \ls y_i,   \q  E_\oP \Big[ \int_t^\oT   h_i ( r,\oX_{r \land \cd} , \oU_r ) dr  \Big] \= z_i, \q \fa i \ins \hN
 \eea
  in   weak formulation  is
  \beas
  \oV (t,\bx,y,z) \df \underset{\oP \in \ocP_{t,\bx}(y,z)}{\sup} E_\oP \big[ \, \oR  (t) \big]
  \= \underset{\oP \in \ocP_{t,\bx}(y,z)}{\sup} E_\oP \Big[ \int_t^\oT \n  f   (r, \oX_{r \land \cd},\oU_r  ) dr  \+ \b1_{\{\oT < \infty\}} \pi   \big(\oT, \oX_{\oT \land \cd}\big) \Big] ,
 \eeas
  where  $ \ocP_{t,\bx}(y,z) \df \Big\{\oP \ins \ocP_{t,\bx} \n :
 E_\oP \big[ \int_t^\oT   g_i ( r,\oX_{r \land \cd} , \oU_r ) dr  \big] \ls y_i,   \, E_\oP \big[ \int_t^\oT   h_i ( r,\oX_{r \land \cd} , \oU_r ) dr  \big] \= z_i, \, \fa i \ins \hN  \Big\} $.    We will simply call $ \oV (t,\bx,y,z) $ the {\it weak} value of the  stochastic control/stopping problem with expectation constraints.
   In case $ \ocP_{t,\bx}(y,z) \=  \es$, $\oV (t,\bx,y,z)   \= -\infty $ by the convention   $ \sup \es \df -\infty$.

   \if{0}

Let $\ff \n: (0,\infty) \ti \OmX \ti \hU \mto [0,\infty]$ be a Borel-measurable function.
 Taking $\psi(r,\fx,\nxi,u) \df \ff (r,\fx,u) $, $\fa (r,\fx,\nxi,u) \ins (0,\infty) \ti \OmX \ti \hR^{d+l} \ti \hU $
 in Lemma \ref{lem_M29_01} (2)
  shows that the mapping   $   \Psi_\ff  (t,\fs,\omX,\fu) \df  \Psi (t,\fs,\bz,\omX,\fu)
   \=  \int_t^{t+\fs} \psi \big(r, \fl_2(r,\omX), 0,\omX(r) , \\ \fu(r)\big) dr
  \=    \int_t^{t+\fs} \ff \big(r, \fl_2(r,\omX),\fu(r)\big) dr  $,
  $  (t,\fs, \omX,\fu) \ins [0,\infty) \ti [0,\infty) \ti \OmX \ti \hJ$ is $\sB[0,\infty) \oti \sB[0,\infty) \oti \sB(\OmX) \oti  \sB(\hJ)-$measurable.

 Since the random variables $(\oX,\oU)$ on $\oO$ are $  \sB(\OmX) \oti \sB(\hJ)-$measurable, it follows that
 the mapping
 \bea
 \ol{\Psi}_\ff (t,\oo) & \tn \df & \tn  
   \int_t^{\oT(\oo) \vee t } \ff \big(r, \oX_{r \land \cd}(\oo),\oU_r(\oo)\big) dr
 \= \lmt{n \to \infty}  \int_t^{(\oT(\oo) \land n) \vee  t}  \ff \big(r, \oX_{r \land \cd}(\oo),\oU_r(\oo)\big) dr  \nonumber  \\
 & \tn \=  & \tn \lmt{n \to \infty}  \int_t^{t+(\oT(\oo) \land n -t )^+  } \ff \big(r, \oX_{r \land \cd}(\oo),\oU_r(\oo)\big) dr
  \= \lmt{n \to \infty} \Psi_\ff \big(t , (\oT(\oo) \ld n \-t)^+ ,  \oX (\oo) ,\oU (\oo)\big)  , \q \fa (t,\oo)  \ins [0,\infty) \ti \oO
  \label{081222_11_a}
 \eea
   is $\sB[0,\infty) \oti \sB(\oO)-$measurable.

   \fi

 We can consider another weak value function of the SCEC:
 Let  $(\bw,\bu) \ins   \O_0  \ti \hJ   $ and  define
   $ \ocP_{t,\bw,\bu,\bx} \df \big\{\oP \ins \ocP_{t,\bx} \n : \oP  \big\{  \oW_s  \=   \bw(s)   , \,
    \fa s \ins [0,t] ; \, \oU_s \= \bu(s) \hb{ for a.e. } s \ins (0,t)  \big\} \= 1 \big\}$
     as the subclass of $\ocP_{t,\bx}$ given the historical Brownian path $\bw \big|_{[0,t]}$
     and the historical control trajectory  $ \bu \big|_{[0,t]}$.
 \if{0}

 Set  $\ol{\U}_s (\oo) \df \ol{\U}^0_s (\oo)  $, $\fa (s,\oo) \ins [0,\infty) \ti \oO$ and
 $\cI(s,\fu) \= \int_0^s e^{-r} \sI(\fu(r))dr$, $\fa (s,\fu) \ins [0,\infty) \ti \hJ $. We have
 \beas
 && \hspace{-1cm} \big\{\oo \ins \oO \n :  \oW_s (\oo)  \=   \bw(s)   , \,
    \fa s \ins [0,t] ; \; \oU_s (\oo)\= \bu(s) \hb{ for a.e. } s \ins (0,t)  \big\}
     \=  \big\{\oo \ins \oO \n :  \oW_s (\oo) \=   \bw(s)   , \, \ol{\U}_s (\oo) \= \cI(s,\bu) ,\,
    \fa s \ins [0,t]   \big\}  \\
 && \= \ccap{q \in \hQ \cap [0,t]}{} \big\{\oo \ins \oO \n :  \oW_q (\oo) \=   \bw(q)   , \ol{\U}_q (\oo) \= \cI(q,\bu)
       \big\}
       \=   \ccap{q \in \hQ \cap [0,t]}{} \Big( \oW^{-1}_q \big\{\bw(q)\big\} \Cp \ol{\U}^{-1}_q \big\{\cI(q,\bu)\big\}  \Big)
       \ins \cF^{\oW,\ol{\U}}_t \sb \sB(\oO) .
    \eeas
    Similarly,
 \bea
 \oO_{t,\bw ,\bu, \bx} & \tn   \df  & \tn  \big\{\oo \ins \oO \n :   (\oW_s , \oX_s) (\oo)   \=   \big(\bw(s),\bx(s)\big)   , \,
    \fa s \ins [0,t] ; \; \oU_s (\oo)\= \bu(s) \hb{ for a.e. } s \ins (0,t)  \big\} \nonumber \\
  & \tn    \= & \tn   \big\{\oo \ins \oO \n :  \oW_s (\oo) \=   \bw(s)   , \, \oX_s (\oo) \=   \bx(s)   , \, \ol{\U}_s (\oo) \= \cI(s,\bu) ,\,
    \fa s \ins [0,t]   \big\}   \nonumber \\
  & \tn  \= & \tn  \ccap{q \in \hQ \cap [0,t]}{} \big\{\oo \ins \oO \n :  \oW_q (\oo) \=   \bw(q)   ,\,  \oX_q (\oo) \=   \bx(q)   ,\, \ol{\U}_q (\oo) \= \cI(q,\bu)       \big\}  \nonumber \\
   & \tn     \=  & \tn   \ccap{q \in \hQ \cap [0,t]}{} \Big( \oW^{-1}_q \big\{\bw(q)\big\} \Cp \oX^{-1}_q \big\{\bx(q)\big\} \Cp \ol{\U}^{-1}_q \big\{\cI(q,\bu)\big\}  \Big)
       \ins \cF^{\oW,\ol{\U},\oX}_t \sb \sB(\oO) . \label{070121_14}
    \eea
    It is clear that
    \bea \label{070121_19}
    \oP \big( \oO_{t,\bw ,\bu, \bx} \big) \= 1 , \q \fa (t,\bw,\bu,\bx) \ins [0,\infty) \ti \O_0  \ti \hJ \ti \OmX , ~ \fa \oP \ins \ocP_{t,\bw,\bu,\bx} .
    \eea

 \fi
 The weak value of the stochastic control/stopping problem with expectation constraints \eqref{111920_14} given
 $(\bw,\bu,\bx)\big|_{[0,t]}$ is
  $ 
  \oV (t,\bw,\bu,\bx,y,z) \df \underset{\oP \in \ocP_{t,\bw,\bu,\bx}(y,z)}{\sup} E_\oP \big[ \, \oR  (t) \big] $,
  where $ \ocP_{t,\bw,\bu,\bx}(y,z) \df \big\{\oP \ins \ocP_{t,\bx}(y,z) \n : \oP  \big\{  \oW_s  \=   \bw(s)   , \,
    \fa s \ins [0,t] ; \, \oU_s \= \bu(s) \hb{ for a.e. } s \ins (0,t)  \big\} \= 1 \big\}$.

  One of our main results   in the next theorem demonstrates that
the value $V(t,\bx,y,z)$ in \eqref{081820_11} coincides with  the weak value $\oV(t,\bx,y,z)$,
and is even equal to   $ \oV (t,\bw,\bu,\bx,y,z) $.

 \begin{thm} \label{thm_V=oV}
 Let $(t, \bw, \bu, \bx, y,z ) \ins [0,\infty) \ti \O_0 \ti   \hJ  \ti \OmX \ti \Re \ti \Re$.
 Then $   V (t,\bx,y,z ) \= \oV (t,\bx,y,z ) = \oV (t, \bw,\bu,\bx, y,z ) $, and
 $\cC_{t,\bx}(y,z) \nne \es   \Leftrightarrow \ocP_{t,\bx}(y,z) \nne \es \Leftrightarrow \ocP_{t,\bw,\bu,\bx}(y,z) \nne \es$.

 \end{thm}

Theorem \ref{thm_V=oV} indicates that the value of the SCEC is independent of a specific probabilistic setup
and is also indifferent to the Brownian/control  history. This result even allows us to deal with the robust case:

 \begin{rem}

 Let $  \big\{(\cQ_\a, \cF_\a, \fp_\a)\big\}_{\a \in \fA} $ be a family of probability spaces,
 where $\fA$ is a countable or uncountable  index set \big(e.g. one can consider   a non-dominated class $\{\fp_\a\}_{\a \in \fA}$ of probability measures on a measurable space $(\cQ, \cF)$\big).

 Let $\a \ins \fA$ and let $\cB^\a \= \big\{\cB^\a_s\big\}_{s \in [0,\infty)}$ be a $d-$dimensional standard  Brownian motion on $(\cQ_\a, \cF_\a, \fp_\a)$. Given $t \ins [0,\infty)$, set $\cB^{\a,t}_s \df \cB^\a_s \- \cB^\a_t$, $s \ins [t,\infty) $,
 let $\cU^\a_t$ collect all $\hU-$valued, $\bF^{\cB^{\a,t},\fp_\a}-$progressively measurable processes $\mu^\a \= \{\mu^\a_s\}_{s \in [t,\infty)}$
  and let $\cS^\a_t$ denote the set of all $[t,\infty]-$valued  $\bF^{\cB^{\a,t},\fp_\a}-$stopping times.
    For any $ (\bx, \mu^\a)  \ins   \OmX \ti \cU^\a_t$,
 let $\cX^{t,\bx,\mu^\a} \= \big\{ \cX^{t,\bx,\mu^\a}_s  \big\}_{s \in [0,\infty)}$ be
 the unique $\big\{\cF^{\cB^{\a,t},\fp_\a}_{s \vee t}\big\}_{s \in [0,\infty)}-$adapted continuous  process
 solving the SDE with the open-loop control $\mu^\a$
  \beas
 \cX_s \= \bx(t) + \int_t^s b ( r, \cX_{r \land \cd}, \mu^\a_r) dr \+ \int_t^s  \si ( r, \cX_{r \land \cd}, \mu^\a_r) d \cB^\a_r , ~ \fa s \in [t,\infty) \hb{ with initial condition $ \cX  \big|_{[0,t]} \= \bx \big|_{[0,t]} $ }
 \eeas
 on $\big(\cQ_\a, \cF_\a, \bF^{\cB^{\a,t},\fp_\a},  \fp_\a \big)$.

  Then we know from Theorem \ref{thm_V=oV} that  for any $(t,\bx) \ins [0,\infty)    \ti \OmX  $ and
    $(y,z) \= \big(\{y_i\}_{i \in \hN}, \{z_i\}_{i \in \hN}\big)  \ins \Re \ti \Re$
 \beas
 \;\; \oV (t,\bx,y,z) \= \Sup{\a \in \fA} \, \Sup{(\mu^\a,\tau_\a) \in \cC^\a_{t,\bx}(y,z) } E_{\fp_\a} \Big[ \int_t^{ \tau_\a} f \big( r, \cX^{t,\bx,\mu^\a}_{r \land \cd} , \mu^\a_r   \big) dr
 \+ \b1_{\{\tau_\a < \infty\}} \pi \big(  \tau_\a, \cX^{t,\bx,\mu^\a}_{ \tau_\a  \land \cd}  \big)  \Big] ,
 \eeas
  where $ \cC^\a_{t,\bx}(y,z) \df \big\{(\mu^\a,\tau_\a) \ins \cU^\a_t \ti \cS^\a_t \n : E_{\fp_\a} \big[ \int_t^{ \tau_\a} g_i (r,\cX^{ t,\bx,\mu^\a}_{r \land \cd} , \mu^\a_r  ) dr  \big] \ls y_i, \, E_{\fp_\a} \big[ \int_t^{ \tau_\a} h_i (r,\cX^{ t,\bx,\mu^\a}_{r \land \cd} , \mu^\a_r  ) dr  \big] \= z_i, ~ \fa i \ins \hN   \big\} $.
  To wit, the weak value  $ \oV (t,\bx,y,z)$ is also equal to the robust value of the SCEC
  under model uncertainty.

\end{rem}

  The equivalence between   strong and weak formulation of an  (unconstrained) stochastic control/stopping problem was discussed in
\cite{Elk_Tan_2013b}. However, their argument may not be applicable to
   stochastic control/stopping problem with expectation constraints:

\begin{rem} \label{rem_032322}

 When $\big(y_i,h_i(\cd,\cd),z_i\big) \= (\infty,0,0)$, $\fa i \ins \hN$, 
the unconstrained version of Theorem \ref{thm_V=oV} states that for any $(t,\bx ) \ins [0,\infty) \ti  \OmX   $,
$ V (t,\bx ) \df  \Sup{(\mu,\tau) \in \cU_t \times \cS_t } E_\fp \big[ \int_t^\tau \n f \big( r, \cX^{t,\bx,\mu}_{r \land \cd},\mu_r  \big) dr
 \+ \b1_{\{\tau < \infty\}} \pi \big( \tau  , \cX^{t,\bx,\mu}_{ \tau  \land \cd}  \big)  \big] $ is equal to $ \oV (t,\bx  ) \df \underset{\oP \in \ocP_{t,\bx} }{\sup} E_\oP \big[\, \oR(t)\big]$.
On the other hand,  \cite{Elk_Tan_2013b} showed  that for any $(t,\bx ) \ins [0,\infty) \ti  \OmX   $,
$ V (t,\bx ) $   equals
  $  \ooV (t,\bx) \df \underset{\oP \in \oocP_{t,\bx}}{\sup} E_\oP \big[\, \oR(t)\big] $,
  where $\oocP_{t,\bx}$ collects all $ \oP \ins \fP(\oO)$ satisfying \(D1\), \(D2\), \(D3\) and  ``\,$\oP \big\{ \oT \gs t \big\} \= 1$"
  \(We summarize \cite{Elk_Tan_2013b}'s result in our terms 
  for an easy comparison with our work\).
As $  \ocP_{t,\bx} \sb \oocP_{t,\bx}$,
the equality $   V (t,\bx ) \= \underset{\oP \in \ocP_{t,\bx} }{\sup} E_\oP \big[\, \oR(t)\big]
\= \underset{\oP \in \oocP_{t,\bx}}{\sup} E_\oP \big[\, \oR(t)\big] $ indicates   that the probability classes $\ocP_{t,\bx}$'s are more accurate than 
 $\oocP_{t,\bx}$'s  to describe the \(unconstrained\)  stochastic control/stopping problem in weak formulation.

 The condition \(D4\) on $\ocP_{t,\bx}$ is necessary for  the expectation-constraint case.
 Without it, the weak value $\ooV(t,\bx,y,z) \\ \df \underset{\oP \in \oocP_{t,\bx}(y,z) }{\sup} E_\oP \big[\, \oR(t)\big]$
 \Big(with $ \oocP_{t,\bx}(y,z)   \df \Big\{\oP \ins \oocP_{t,\bx} \n : E_\oP \big[ \int_t^\oT   \n  g_i ( r,\oX_{r \land \cd},\oU_r ) dr  \big] \ls y_i,   \,
   E_\oP \big[ \int_t^\oT   \n  h_i ( r,\oX_{r \land \cd},\oU_r ) dr  \big]   \= z_i, \, \fa i \ins \hN  \Big\} $\Big)
 may not be equal to $V (t,\bx,y,z )  $ for the following reason:

 In Proposition 4.3 of \cite{Elk_Tan_2013b}, the key  to show $\ooV (t,\bx) \ls  V (t,\bx) $ or $E_\oP \big[\, \oR(t)\big] \ls  V (t,\bx) $ for a given $\oP \ins \oocP_{t,\bx} $,
 relies on   transforming  the hitting times of process $\big\{E_\oP   \big[   \b1_{\{\oT \in [t,s]\}} \big| \cF^{\oW^t,\oP}_\infty \big]\big\}_{s \in  [t,\infty)}$
 to a member of $\cS_t$. More precisely,
 the so-called {\it Property \(K\)} assures
   an $\bF^{W^t,P_0}-$adapted \cad   process $ \wh{\vth}_\cd $  such that
$\wh{\vth}_s  ( \oW  ) \=  E_\oP   \big[   \b1_{\{\oT \in [t,s]\}} \big| \cF^{\oW^t}_s \big] \=  E_\oP   \big[   \b1_{\{\oT \in [t,s]\}} \big| \cF^{\oW^t,\oP}_\infty \big]  $, $\oP-$a.s. for any $s \ins [t,\infty)$.
 It follows  that
 $   E_\oP \big[ \b1_{\{\oT \in [t,s]\}}  \n  \b1_{\{\osX^{t,\bx,\ol{\mu}}   \in A\}}   \big| \cF^{\oW^t,\oP}_\infty \big]
   \= \b1_{\{\osX^{t,\bx,\ol{\mu}}_\cd \in A\}}  \wh{\vth}_s  ( \oW ) \\
    \=   \int_t^s   \b1_{\{\osX^{t,\bx,\ol{\mu}}  \in A\}}   \wh{\vth}  (dr ,\oW ) $, $\oP-$a.s.  for any $(s,A) \ins [t,\infty)  \ti  \sB(\OmX)$,
    where $\ol{\mu} \= \big\{\ol{\mu}_s \= \wh{\mu}_s(\oW)\big\}_{s \in [t,\infty)}$ is the $\hU-$valued, $\bF^{\oW^t}-$predictable process in \(D1\)
    and $ \osX^{t,\bx,\ol{\mu}} \= \big\{\osX^{t,\bx,\ol{\mu}}_s\big\}_{s \in [0,\infty)} $ is the unique solution of SDE \eqref{Ju01_01}.
 Let $\Phi $ be a nonnegative Borel-measurable  function  on $ [0,\infty) \ti \hJ \ti \OmX$.
     Then   
     a standard approximation argument and the ``change-of-variable" formula 
     yield that
    $ E_\oP   \big[  \Phi (\oT  , \oU, \oX )   \big| \cF^{\oW^t,\oP}_\infty \big]
  \= \int_t^\infty \n \Phi  (  r, \oU, \oX )    \wh{\vth}  (dr ,\oW )
  \= \int_0^1   \Phi (   \vr(\oW ,\l) , \wh{\mu} (\oW),  \oX )   d \l $, $\oP-$a.s.,
  where $\vr (\o_0,\l) \df \inf\big\{s \ins [t,\infty) \n : \wh{\vth}_s (  \o_0) \>  \l \big\}  $, $ \fa  (\o_0,\l) \ins \O_0 \ti (0,1) $.
  Set $\mu_s \df \wh{\mu}_s(\cB)$, $\fa s \ins [t,\infty)$.
   Since the  joint $\oP-$distribution of $ \big(\oW , \wh{\mu}_\cd(\oW) \big)$ is equal to the joint $\fp-$distribution of $\big(\cB,\wh{\mu}_\cd(\cB)\big)$, we can deduce that
    the  joint $\oP-$distribution of $ \big(\oW ,\ol{\mu}_\cd, \osX^{t,\bx,\ol{\mu}} \big)$ is equal to the joint $\fp-$distribution of $(\cB,\mu_\cd,\cX^{t,\bx,\mu})$ and thus
   \bea \label{031922_11}
    E_\oP   \big[  \Phi (\oT  , \oU, \oX )    \big]
   \= \int_0^1  E_\oP  \big[ \Phi \big(   \vr(\oW ,\l) , \wh{\mu} (\oW), \osX^{t,\bx,\ol{\mu}} \big) \big]  d \l
   \=   \int_0^1  E_\fp  \big[ \Phi (   \vr(\cB ,\l) , \wh{\mu} (\cB),  \cX^{t,\bx,\mu} ) \big]  d \l  .
   \eea
  As   $\tau_\l \df \vr(\cB ,\l) \ins \cS_t $ for each $\l \ins (0,1)$, taking $\Phi$ to  be the total reward  function
  implies  that
  \bea \label{031922_14}
  E_\oP \big[\, \oR(t)\big] \= \int_0^1 E_\fp \Big[ \int_t^{\tau_\l}  \n  f \big( r, \cX^{t,\bx,\mu}_{r \land \cd},\mu_r  \big) dr
 \+ \b1_{\{\tau_\l < \infty\}} \pi \big( \tau_\l  , \cX^{t,\bx,\mu}_{ \tau_\l  \land \cd}  \big)  \Big]   d \l
 \ls \int_0^1 V(t,\bx) d \l \= V(t,\bx) .
 \eea

 However, this argument is not applicable to the expectation-constraint case: Given a $\oP \ins \oocP_{t,\bx} (y,z)$,
 since   $ (\mu,\tau_\l) $ may not belong to $\cC_{t,\bx} (y,z) $ for a.e. $\l \ins (0,1)$,   one can not get $E_\oP \big[\, \oR(t)\big]  \ls V(t,\bx,y,z)$ like \eqref{031922_14}.
  Actually, for each $\l \ins (0,1)$, $ (\mu, \tau_\l) $ is only of $\cC_{t,\bx} (y_\l,z_\l) $ with
 $(y_\l,z_\l) \= \big(\{y^i_\l\}_{i \in \hN} , \{z^i_\l\}_{i \in \hN} \big)$ and
 $(y^i_\l,z^i_\l) \df \Big(  E_\fp \big[ \int_t^{\tau_\l} \n  g_i ( r,\cX^{t,\bx,\mu}_{r \land \cd},\mu_r ) dr  \big] ,   E_\fp \big[  \\  \int_t^{\tau_\l} \n  h_i ( r,    \cX^{t,\bx,\mu}_{r \land \cd},\mu_r ) dr  \big]  \Big) $.
 For $i \ins \hN$, choosing   accumulative cost functions for $\Phi$ in \eqref{031922_11}
 renders that
 \beas
  \int_0^1   E_\fp \Big[ \int_t^{\tau_\l} \n  g_i ( r,\cX^{t,\bx,\mu}_{r \land \cd},    \mu_r ) dr  \Big]   d \l
  \= E_\oP \bigg[ \int_t^\oT   g_i ( r,\oX_{r \land \cd},\oU_r ) dr  \bigg]   \ls y_i
  \eeas
  and similarly
   $ \int_0^1   E_\fp \big[ \int_t^{\tau_\l}    h_i ( r,\cX^{t,\bx,\mu}_{r \land \cd},\mu_r ) dr  \big]  d \l
    \= z_i $, so $ V \big(t,\bx,\{ \int_0^1 y_\l d \l \}_{i \in \hN} ,    \{ \int_0^1 z_\l d \l \}_{i \in \hN}\big)  \ls V(t,\bx,y,z)$.
  Then  the attempt to show   $ E_\oP \big[\, \oR(t)\big] \ls V(t,\bx,y,z) $ reduces to deriving  a Jensen-type inequality:
 \beas
  \int_0^1 V(t,\bx,y_\l,z_\l)  d \l \ls V \Big(t,\bx, \Big\{ \int_0^1 y_\l d \l \Big\} , \Big\{ \int_0^1 z_\l d \l \Big\}\Big)   .
  \eeas
 But this does not hold  since the value function $V $ is not concave in level $z$ of equality-type expectation constraints.

\end{rem}

\section{The Measurability of SCEC Values}
\label{sec_Mart_prob}

   In this section, using the martingale-problem formulation of controlled SDEs, 
   we   characterize the probability class $ \ocP_{t,\bx} $ by countably many stochastic behaviors of the
   canonical coordinators $(\oW,\oU,\oX,\oT)$ of $\oO$.
   This will enable us to analyze   the measurability  of   value functions of the  stochastic control/stopping problem
   with expectation constraints.

 Let $ \fU $ be the equivalence classes of  all   $\hU-$valued, $ \bF^{W,P_0} -$predictable processes
 on $\O_0$ in the sense that $\mu^1, \mu^2 \ins \fU $ are equivalent if
 $  \big\{(s,\o_0) \ins [0,\infty) \ti \O_0 \n : \mu^1_s(\o_0) \nne \mu^2_s(\o_0) \big\} $
 is a  $ds \ti dP_0  -$null set.
    Given $ \mu  \ins \fU$,  Fubini Theorem shows that
    $  \cN  \df \big\{\o_0 \ins \O_0 \n :  \mu_\cd  (\o_0) \n \notin \n  \hJ \big\} $ is a $\cF^{W,P_0}_\infty-$measurable
   set with zero $P_0-$measure or $  \cN  \ins \sN_{P_0}\big(\cF^W_\infty\big)$.
   By modifying $\mu$ on $\cN$, we obtain an $ \bF^{W,P_0} -$predictable process with all paths in $\hJ$:
   $     \{ \mu_s    \b1_{   \cN^c  }  \+ u_0 \b1_{   \cN  }  \}_{  s  \in  [0,\infty) }  $,
   which is   in the same   equivalence class 
   as $ \mu $.
   To wit, one can assume without loss of generality that for all $\mu \in \fU$,
   $ \mu_\cd(\o_0) \ins \hJ $   for any $\o_0 \ins \O_0$.
  We equip $\fU$  with the topology of local convergence in measure, i.e., the metric
  \beas
   \Rho{\fU} (\mu^1,\mu^2) \df    E_{P_0} \Big[ \int_0^\infty \n  e^{-  s} \big( 1 \ld  \Rho{\hU} (\mu^1_s,\mu^2_s)  \big)   ds \Big]  , \q \fa \mu^1, \, \mu^2 \ins \fU \, .
  \eeas

    Also, let $ \fS $ be the equivalence classes of  all   $[0,\infty]-$valued,
 $ \bF^{W,P_0} -$stopping times on $\O_0$
 in the sense that $\tau_1, \tau_2 \ins \fS $ are equivalent if
 $ P_0 \{   \tau_1  \= \tau_2  \} \= 1 $.
 We  endow $\fS$   with the   metric
 \beas
 \Rho{\fS} (\tau_1,\tau_2) \df E_{P_0} \big[ \Rho{+}  ( \tau_1 , \tau_2 )  \big]  ,    \qq  \fa \tau_1, \, \tau_2 \ins \fS .
 \eeas

 \begin{lemm} \label{lem_082020_15}
  $\big(\fU, \Rho{\fU}\big)$ and $\big(\fS, \Rho{\fS}\big)$ are two complete separable metric  spaces, i.e., Polish spaces.

\end{lemm}

 For any $(\mu,\tau) \ins \fU \ti \fS$,   define their joint distribution with $W$ under $P_0$ by
 $\Ga  (  \mu,\tau   ) \df P_0 \nci  ( W,  \mu,\tau    )^{-1} \ins \fP \big( \O_0 \ti \hJ \ti \hT \big)$.

\begin{lemm} \label{lem_082020_17}

The mapping $ \Ga \n : \fU \ti \fS \mto \fP \big( \O_0 \ti \hJ \ti \hT \big) $ is a continuous injection
    from $\fU \ti \fS$  into $ \fP \big(\O_0 \ti \hJ \ti \hT \big) $.

\end{lemm}

For any $(t,\oP) \ins [0,\infty) \ti \fP\big(\oO\big)$,  define a probability measure on $\big(\O_0 \ti \hJ \ti \hT,\sB ( \O_0 \ti \hJ \ti \hT) \big)$ by $\oQ_{t,\oP} (D) \df \oP \big\{(\osW^t, \osU^t, \oT\-t) \ins D \big\}$, $\fa D \ins \sB ( \O_0 \ti \hJ \ti \hT)$.

   \if{0}

 Similar to the metric $\Rho{+}$ on $\hT \= [0,\infty]$, we define a metric $\ddot{\rho}$ on $ \ddot{\hR} \df [-\infty,\infty]$ by
 $\ddot{\rho}(a_1,a_2) \df \big|\arctan(a_1) \- \arctan(a_2)\big|$, $\fa a_1,a_2 \ins \ddot{\hR}$.

\begin{lemm} \label{lem040221}
The mapping $(t,\oo) \mto \oT(\oo)\-t$ is continuous from $[0,\infty) \ti \oO$ to $ \ddot{\hR} $.
\end{lemm}

\no {\bf Proof:} Let $ t  \ins [0,\infty)$, $\oo  \= (\o_0,\fu,\omX,\ft) \ins  \oO$ and   $\e \ins (0,1)$. We discuss by   two cases:

\no (i) When $\ft \< \infty$, we set $\dis \d \= \d(\e,\ft) \df \frac{\e}{2 (1\+(\ft\+1)^2)} \ld \big(\arctan(\ft\+1) \- \arctan(\ft)\big)   $.
Let $ t' \ins (t \- \d,t\+ \d)\Cp [0,\infty)$ and $\oo'  \= (\o'_0,\fu',\omX',\ft') \ins  \oO$ with
 $ \rho^2_{\overset{}{\oO}}   (\oo,\oo') \=   \rho^2_{\overset{}{\O_0}} (\o_0,\o'_0)  \+  \rho^2_{\overset{}{\hJ}} (\fu,\fu') \+ \rho^2_{\overset{}{\OmX}} (\omX,\omX') \+ \rho^2_{\overset{}{+}}(\ft,\ft')  \< \d^2 $.
 As $ \arctan(\ft\+1) \- \arctan(\ft) \gs \d \> \Rho{+}(\ft,\ft') \= \big|\arctan(\ft) \- \arctan(\ft')\big|
 \gs \arctan(\ft') \- \arctan(\ft) $, one has that $\ft' \< \ft\+1$.
  Since
 \bea \label{040521_11}
 \frac{b \- a}{1+b^2} \ls \arctan(b) \- \arctan(a) \= \int_a^b \frac{1}{1+x^2}dx \ls \frac{b \- a}{1+a^2} \ls b\-a  , \q \hb{ for } 0 \ls a \ls b \< \infty ,
\eea
we can also deduce that
 $\dis  \frac{\e}{2 (1\+(\ft\+1)^2)} \gs \d \>   \big|\arctan(\ft) \- \arctan(\ft')\big|
 \gs \frac{ |\ft'\-\ft|}{1+(\ft'\ve \ft)^2} \> \frac{ |\ft'\-\ft|}{1+(  \ft\+ 1)^2}$ and thus $  |\ft'\-\ft| \< \e/2$.
 As   $\d \< \e/2$,  \eqref{040521_11} also implies that
$\ddot{\rho} \big( \oT(\oo')\-t' , \oT(\oo)\-t\big)
\= \big|   \arctan \big(\ft'\-t'\big) \- \arctan \big(\ft\-t\big) \big|
   \ls    \big| \ft'\-t' \- \ft \+ t \big|  \ls | t'\- t| \+ |\ft'\-\ft|
\< \d \+ \e/2 \< \e $.

\no (ii) When $\ft \= \infty$, we set $\dis \wt{\d} \= \wt{\d}(\e,t ) \df \frac{\e}{2} \ld \Big( \frac{\pi}{2} \- \arctan \big( t\+1\+ \sqrt{  2(t\+1) \e^{-1}\-1 } \big) \Big)$.  Let $ t' \ins \big(t \- \wt{\d},t\+ \wt{\d} \big)\Cp [0,\infty)$ and $\oo'  \= (\o'_0,\fu',\omX',\ft') \ins  \oO$ with
 $ \rho^2_{\overset{}{\oO}}   (\oo,\oo') \=   \rho^2_{\overset{}{\O_0}} (\o_0,\o'_0)  \+  \rho^2_{\overset{}{\hJ}} (\fu,\fu') \+ \rho^2_{\overset{}{\OmX}} (\omX,\omX') \+ \rho^2_{\overset{}{+}}(\ft,\ft')  \< \wt{\d}^2 $.
 Since
$\frac{\pi}{2} \- \arctan \big( t\+1\+ \sqrt{  2(t\+1) \e^{-1}\-1 } \big) \gs \wt{\d} \> \Rho{+}(\ft,\ft') \= \big|\arctan(\infty) \- \arctan(\ft')\big| \= \frac{\pi}{2} \- \arctan(\ft') $ and since $t' \< t\+ \wt{\d} \ls t \+ \frac{\e}{2} \< t \+ 1$,
we see that
$\ft' \> t\+1\+ \sqrt{  2(t\+1) \e^{-1}\-1 } \> t' \+ \sqrt{  2(t\+1) \e^{-1}\-1 }$.
 By \eqref{040521_11} again,
 $  \arctan  (\ft' ) \-  \arctan \big(\ft'\-t'\big)
\ls \frac{t'}{ 1+ (\ft'\-t')^2 } \< \frac{t\+1}{1+ ( 2(t\+1) \e^{-1}\-1) } \= \frac{\e}{2} $.
 Adding it to the inequality  $\frac{\pi}{2} \- \arctan(\ft') \= \Rho{+}(\ft,\ft') \< \wt{\d} \ls \e/2 $ yields that
  $ \ddot{\rho} \big( \oT(\oo')\-t' , \oT(\oo)\-t\big)
 \= \big|   \arctan \big(\ft'\-t'\big) \- \arctan \big(\ft\-t\big) \big|
 \= \big|   \arctan \big(\ft'\-t'\big) \- \arctan  (\infty ) \big|
 \= \frac{\pi}{2} \- \arctan \big(\ft'\-t'\big) \< \e $.    \qed

\begin{rem} \label{rem_061222}
Given $t \ins [0,\infty)$, 
$ \big( \osW^t,\osU^t, \oT \- t\big)  \n : \oO \mto \O_0 \ti \hJ \ti \ddot{\hR}$
is $\sB(\oO) \big/ \sB(\O_0) \oti \sB(\hJ) \oti \sB\big(\ddot{\hR}\big) -$measurable.
  So $(\osW^t, \osU^t, \oT\-t)^{-1} ( D ) \ins \sB(\oO)$ for any $  D \ins \sB ( \O_0 \ti \hJ \ti \hT)$.
  The probability measure $  \oP \nci \big( \osW^t, \osU^t,\oT \- t\big)^{-1} $ is well-defined.
\end{rem}

\no {\bf Proof of Remark \ref{rem_061222}:}
 Define $\sW^t_\fs (\o_0) \df W_{t+\fs} (\o_0) \- W_t (\o_0) \= \o_0(t\+\fs) \- \o_0(t) $, $(\fs,\o_0) \ins [0,\infty) \ti \O_0$
 and set $T \df \lceil t \rceil$.

  Let $\e \ins (0,1)$.
  \if{0}
  Using the inequality
 \bea \label{040321_11}
  a \ld (b\+c) \ls a \ld b \+ a \ld c, \q \fa a,b,c \ins (0,\infty),
  \eea
  \fi
 For any $\o_0 ,\o'_0 \ins  \O_0$ with $\dis \Rho{\O_0}  (\o_0,\o'_0) \< \d   \df \frac{\e}{(1 \+ 2^T)} $\,,
 \beas
 \q  && \hspace{-1cm} \Rho{\O_0} \big(\sW^t (\o_0), \sW^t (\o'_0)\big)
  \= \sum_{n \in \hN} \Big( 2^{-n} \ld \Sup{\fs \in [0,n]}  \big|\sW^t_\fs (\o_0)\- \sW^t_\fs (\o'_0) \big| \Big)
   \= \sum_{n \in \hN} \Big( 2^{-n} \ld \Sup{\fs \in [0,n]}   \big|  \o_0(t\+\fs) \- \o_0(t) \-    \o'_0 (t\+\fs) \+  \o'_0(t) \big| \Big) \\
   &&  \ls  \sum_{n \in \hN}     \Big( 2^{-n} \ld  \Sup{\fs \in [0,n]}|  \o_0(t\+\fs) \- \o'_0(t\+\fs) | \Big) \+  \sum_{n \in \hN}  \Big( 2^{-n} \ld \Sup{\fs \in [0,n]} |    \o_0(t) \-  \o'_0(t) |\Big) \\
   &&  \=  \sum_{n \in \hN}     \Big( 2^{-n} \ld  \Sup{s \in [t,t+n]}|  \o_0(s) \- \o'_0(s) | \Big) \+  \Rho{\O_0}  (\o_0,\o'_0)
      \ls  \sum_{n \in \hN}   2^T  \Big( 2^{-n-T} \ld  \Sup{s \in [0, n+T]}|  \o_0(s) \- \o'_0(s) | \Big) \+  \Rho{\O_0}  (\o_0,\o'_0) \\
  && \< \big(1 \+ 2^T \big) \Rho{\O_0}  (\o_0,\o'_0)  \< \e ,
 \eeas
 which shows that $\sW^t$ is a continuous mapping from $\O_0$ to $\O_0$. So $\osW^t \= \sW^t(\oW)$ is $\sB(\oO) / \sB(\O_0)-$measurable.

 Let   $\vf \ins L^0 \big((0,\infty) \ti \hU;\hR\big)$. Clearly, $\vf_t(s,u) \df \b1_{\{s>t\}} \vf  (s\-t,u) $, $(s,u) \ins (0,\infty) \ti \hU$
 is also of $L^0 \big((0,\infty) \ti \hU;\hR\big)$.
  For any $  \cE \ins \sB(\hR)$, since the mapping $\oU$ is  $\sB(\oO) / \sB(\hJ) -$measurable,
  \beas
 && \hspace{-1.2cm} \big(\osU^t\big)^{-1} \big((I_\vf)^{-1}(\cE)\big)
   \=   \big\{\oo \ins \oO \n : \osU^t_\cd (\oo) \ins (I_\vf)^{-1}(\cE) \big\}
 \= \big\{\oo \ins \oO \n : I_\vf \big(\osU^t_\cd (\oo)\big) \ins  \cE  \big\}
 \= \big\{\oo \ins \oO \n : \int_0^\infty \vf (\fs, \osU^t_\fs (\oo)) d\fs \ins  \cE  \big\} \\
  && \hspace{-0.5cm} \=   \big\{\oo \ins \oO \n : \int_t^\infty \vf (s\-t, \oU_s (\oo)) ds \ins  \cE  \big\}
 \= \big\{\oo \ins \oO \n : \int_0^\infty \vf_t (s, \oU_s (\oo)) ds \ins  \cE  \big\}
 \= \big\{\oo \ins \oO \n : I_{\vf_t} \big(\oU  (\oo)\big)   \ins  \cE  \big\}
 \= \oU^{-1} \big( (I_{\vf_t})^{-1}(\cE) \big) \ins \sB(\oO) ,
 \eeas
   which together with Lemma \ref{lem_M29_01} (1) shows that the sigma-field 
 $    \big\{ A \sb \hJ \n : \big(\osU^t\big)^{-1} (A) \ins \sB(\oO) \big\}$
 includes all generating sets of $ \sB(\hJ)$ and thus contains $ \sB(\hJ)$.
 So $\big(\osU^t\big)^{-1} (A) \ins \sB(\oO)$ for any $A \ins \sB(\hJ)$. Namely,
 $ \osU^t $ is $\sB(\oO) / \sB(\hJ)-$measurable.

 Moreover, we see from Lemma \ref{lem040221}   that  the function $ \oo  \mto \oT(\oo)\-t$ is continuous from $  \oO$ to $ \ddot{\hR} $.
 Hence, the mapping $ \big( \osW^t,\osU^t, \oT \- t\big)  \n : \oO \mto \O_0 \ti \hJ \ti \ddot{\hR}$
is $\sB(\oO) \big/ \sB(\O_0) \oti \sB(\hJ) \oti \sB\big(\ddot{\hR}\big) -$measurable. \qed

   \fi

\begin{lemm} \label{lem_082020_19}
The mapping $   [0,\infty) \ti \fP\big(\oO\big) \ni (t,\oP) \mto  \oQ_{t,\oP}
\ins \fP \big( \O_0 \ti \hJ \ti \hT \big) $  is continuous.
\end{lemm}

  For $t \ins [0,\infty)$ and  $\vf \ins C^2 (\hR^{d+l})  $,    define process
 \beas
   \oM^t_s(\vf)   \df   \vf \big(\oW^t_s,\oX_s \big)
 \- \n \int_t^s    \ol{b}  \big( r,\oX_{r \land \cd} ,\oU_r \big) \n \cd \n D \vf \big( \oW^t_r,\oX_r \big) dr
    \-   \frac12 \int_t^s   \ol{\si} \, \ol{\si}^T  \big( r, \oX_{r \land \cd} , \oU_r \big) \n : \n D^2 \vf \big( \oW^t_r,\oX_r \big)   dr  , \q \fa s \ins [t,\infty) .
 \eeas

We can use Remark \ref{rem_ocP} and Lemma \ref{lem_082020_17} to decompose the probability class $\ocP_{t,\bx}$ as the intersection of countably
many action sets of processes $(\oW,\oU,\oX,\oT)$:

\begin{prop} \label{prop_Ptx_char}
For any $(t,\bx) \ins [0,\infty) \ti \OmX$,
the probability class $\ocP_{t,\bx}$ is the intersection of the following three subsets of  $\fP\big(\oO\big)$:

\no i\)    $ \ocP^1_{t,\bx} \df \big\{ \oP \ins \fP\big(\oO\big) \n : \oP  \{   \oX_s \= \bx(s), \fa s \ins [0,t]   \} \= 1
     \big\}$.

\no  ii\)    $ \ocP^2_t \df \big\{ \oP \ins \fP\big(\oO\big) \n :
 \oQ_{t,\oP} \ins \Ga (\fU \ti \fS)   \big\}$.

\no  iii\)
 $ \ocP^3_t \df \Big\{ \oP \ins \fP\big(\oO\big) \n :
  E_\oP \Big[  \Big( \oM^t_{\otau^t_n \land (t+\fr)} (\vf )  - \oM^t_{\otau^t_n \land (t+\fs)} (\vf ) \Big) \underset{i=1}{\overset{k}{\prod}}   \b1_{    \{(\oW^t_{ t+s_i    },\oX_{ t+s_i  }) \in \cO_i\}    }
      \Big]   \= 0  ,  ~ \fa (\vf,n) \ins \fC(\hR^{d+l}) \ti \hN , \,\fa (\fs,\fr) \\ \ins \hQ^{2,<}_+  , \,\fa
  \{(s_i,   \cO_i )\}^k_{i=1} \sb \big(\hQ \Cp [0,\fs]\big) \ti \sO (\hR^{d+l})   \Big\}$.

\end{prop}

 Based on the countable decomposition of the probability class $\ocP_{t,\bx}$ by Proposition \ref{prop_Ptx_char},
 the next proposition shows that   the graph of   probability classes $ \{\ocP_{t,\bx}\}_{(t,\bx) \in [0,\infty) \times \OmX} $
 is a Borel subset of   $ [0,\infty) \ti \OmX \ti \fP\big(\oO\big) $,
 which is crucial for the measurability of the value functions $V \= \oV $.

\begin{prop} \label{prop_graph_ocP}
 The graph  $\big\lan\n\big\lan  \ocP  \big\ran\n\big\ran   \df \big\{\big(t,\bx,\oP \big)  \ins [0,\infty) \ti \OmX \ti \fP\big(\oO\big) \n :  \oP \ins \ocP_{t,\bx} \big\}$ is a Borel subset of $[0,\infty) \ti \OmX \ti \fP\big(\oO\big)$.
\end{prop}

  Set $ D_\ocP \df \big\{(t,\bx,y,z) \ins [0,\infty) \ti \OmX \ti \Re \ti \Re \n :   \ocP_{t,\bx}(y,z) \nne \es \big\} $
 and   $ \cD_\ocP \df \big\{(t,\bw,\bu,\bx,y,z) \ins [0,\infty) \ti \O_0 \ti \hJ \ti \OmX \ti \Re \ti \Re \n :   \ocP_{t,\bw,\bu,\bx}(y,z) \nne \es \big\} $.

  \if{0}
    According to Remark \ref{rem_112220} (1b) and Theorem \ref{thm_V=oV}, if   $ h_i  \= 0$ for all $  i \ins \hN$,
     $\hb{Proj}  \big(D_\ocP\big) \df \{(t,\bx) \ins [0,\infty) \ti \OmX \n : (t,\bx,y,z) \ins D_\ocP \hb{ for some } (y,z) \ins \Re \ti \Re  \} \= [0,\infty) \ti \OmX$ and $\hb{Proj}  \big(\cD_\ocP\big) \df \{(t,\bw,\bu,\bx) \ins [0,\infty) \ti \O_0 \ti \hJ \ti \OmX \n : (t,\bw,\bu,\bx,y,z) \ins \cD_\ocP \hb{ for some } (y,z) \ins \Re \ti \Re \} \= [0,\infty) \ti \O_0 \ti \OmX$.
  \fi

 \begin{cor} \label{cor_graph_ocP}

   The graph $\gP \df \big\{ \big(t,\bx,y,z, \oP\big) \ins D_\ocP \ti \fP\big(\oO\big)  \n :
 \oP \ins \ocP_{t,\bx}(y,z)   \big\}$   is a Borel subset of $ D_\ocP \ti \fP\big(\oO\big)$
  and the graph $\gcP \df \big\{ \big(t,\bw,\bu,\bx,y,z, \oP\big) \ins \cD_\ocP \ti \fP\big(\oO\big)  \n :
 \oP \ins \ocP_{t,\bw,\bu,\bx}(y,z)   \big\}$   is a Borel subset of $\cD_\ocP \ti \fP\big(\oO\big)$.


 \end{cor}

By Corollary \ref{cor_graph_ocP}, the value function $\oV$ is   \usa ~ and is thus universally measurable.

\begin{thm} \label{thm_V_usa}

  The value function $\oV (t,\bx,y,z)$ is \usa ~ on $   D_\ocP$
  and the value function $\oV (t,\bw,\bu,\bx,y,z) $ is \usa ~ on $   \cD_\ocP $.

\end{thm}

\section{Dynamic Programming Principle for $\oV$}
\label{sec_DPP}

 In this section, we    explore a   dynamic programming principles (DPP)   for the value function  $\oV$ in   weak formulation,
 which takes the conditional expected integrals of constraint functions as additional states.

  Given $t \ins [0,\infty)$, let  $\oga$ be  a  $[t,\infty)-$valued  $\bF^{\oW^t}-$stopping time and let $\oP \ins \fP\big(\oO\big)$.
  According to Lemma 1.3.3 and  Theorem 1.1.8 of \cite{Stroock_Varadhan},
    $ \cF^{\oW^t}_\oga  $ is  countably generated and  there   thus exists a family $\big\{ \oP^t_{\oga,\oo} \big\}_{\oo \in \oO}$
 of probability measures in $\fP\big(\oO\big)$, called  the {\it regular conditional probability distribution} (r.c.p.d.) of $\oP$ with respect to $\cF^{\oW^t}_\oga$,  such that

\no (R1) for any $\oA \ins \sB(\oO)$, the mapping $\oo \mto \oP^t_{\oga,\oo} \big(\oA\big)$ is $\cF^{\oW^t}_\oga-$measurable;

 \no  (R2) for   any $ (-\infty,\infty]-$valued, $\sB_\oP (\oO) -$measurable random variable $\oxi$ that is bounded from below under $\oP$,
 it holds for any $\oo \ins \oO$ except on a $ \ocN_\oxi   \ins \sN_\oP\big(\cF^{\oW^t}_\oga\big)   $ that
 $ \oxi $ is $ \sB_{\oP^t_{\oga,\oo}}(\oO)-$measurable and
 $E_{\oP^t_{\oga,\oo}} \big[ \, \oxi \, \big] \= E_\oP \big[ \, \oxi \, \big| \cF^{\oW^t}_\oga \big] (\oo)$;

\no  (R3) for some   $ \ocN_0   \ins \sN_\oP \big(\cF^{\oW^t}_\oga\big) $,
   $ \oP^t_{\oga,\oo} \big(\oA\big) \= \b1_{\{\oo \in \oA\}} $,  $\fa \big(\oo,\oA\big) \ins \ocN^c_0   \ti \cF^{\oW^t}_\oga $.

\if{0}

 \ss  We have the following extension of (R2).

\begin{cor} \label{cor_070420_11}
  Given $t \ins [0,\infty)$, let $\ocF$ be a sub-sigma field of $\sB(\oO)$, let  $\oga$ be  a  $[0,\infty)-$valued, $\bF^{\oW^t}-$stopping time and let $\oP \ins \fP\big(\oO\big)$.

  \no  1\)  For any    $  \ofN   \ins  \sN_\oP\big(\ocF\big)      $,
  it holds  for all $\oo \ins \oO$ except on a $ \ocN_\ofN   \ins \sN_\oP\big(\cF^{\oW^t}_\oga\big)   $ that $\ofN  \ins  \sN_{\oP^t_{\oga,\oo}} \big(\ocF\big) $.

\no  2\)     For any $ [0,\infty]-$valued, $\si \big( \ocF \cp \sN_\oP\big(\ocF\big) \big) -$measurable random variable $\oxi$ on $\oO$,
 it holds  for all $\oo \ins \oO$ except on a $ \ocN_\oxi   \ins \sN_\oP\big(\cF^{\oW^t}_\oga\big)   $ that
$ \oxi$ is $ \si \big( \ocF \cp \sN_{\oP^t_{\oga,\oo}} \big(\ocF\big) \big)  -$measurable and $ E_{\oP^t_{\oga,\oo}} \big[ \, \oxi \, \big] \= E_\oP \big[ \, \oxi \, \big| \cF^{\oW^t}_\oga \big] (\oo) $.

\end{cor}

\fi

    Let $\oo \ins \oO  $ and set $   \Wtgo   \df \big\{\oo' \ins \oO \n : \oW^t_r (\oo') \= \oW^t_r (\oo), ~ \fa r \ins [t,\oga(\oo)] \big\} $.
    We know from  Galmarino's test  that
   \bea \label{090520_11}
   \oga(\oo') \= \oga(\oo) , \q \fa \oo' \ins \Wtgo ,
  \eea
  and $\Wtgo$ is thus $ \cF^{\oW^t}_\oga -$measurable. 
  Since $\oo \ins \Wtgo$ for any $\oo \ins \oO$,  (R3) shows that
  \bea \label{Jan11_03}
   \oP^t_{\oga,\oo} \big(\Wtgo  \big) \= \b1_{\big\{\oo \in \Wtgo\big\}} \= 1  , \q \fa  \oo \ins \ocN^c_0.
  \eea

 For any $i \ins \hN$, define
  $   \oY^i_{\n \oP} (\oga)   \df   E_\oP \Big[   \int_{\oT \land  \oga }^\oT    g_i (r,\oX_{r \land \cd},\oU_r ) dr \Big| \cF^{\oW^t}_\oga \Big]
  $ and
    $  \oZ^i_\oP (\oga)  \df   E_\oP \Big[   \int_{\oT \land  \oga }^\oT    h_i (r,\oX_{r \land \cd},\oU_r ) dr \Big| \cF^{\oW^t}_\oga \Big] $.
 \if{0}

      \beas
    \big(\oY^i_\oP (\oga)\big)^\pm \df \lsup{n \to \infty} E_\oP \bigg[ n \ld \int_{\oT \land  \oga }^\oT    g^\pm_i (r,\oX_{r \land \cd},\oU_r ) dr \Big| \cF^{\oW^t}_\oga \bigg] , \q
     \big(\oZ^i_\oP (\oga)\big)^\pm \df \lsup{n \to \infty} E_\oP \bigg[ n \ld \int_{\oT \land  \oga }^\oT    h^\pm_i (r,\oX_{r \land \cd},\oU_r ) dr \Big| \cF^{\oW^t}_\oga \bigg] .
   \eeas

 \fi
   So  $ \big(\oY_{\n \oP} (\oga), \oZ_\oP (\oga) \big)  \df \Big( \big\{\oY^i_{\n \oP} (\oga)\big\}_{i \in \hN}, \big\{\oZ^i_\oP (\oga)\big\}_{i \in \hN} \Big)  $   is   an $ \Re \ti \Re -$valued $ \cF^{\oW^t}_\oga -$measurable  random variable. 

 In terms of the r.c.p.d. $\big\{ \oP^t_{\oga,\oo} \big\}_{\oo \in \oO}$,
 the probability class $ \big\{ \ocP_{t,\bx}(y,z) \n : (t,\bx,y,z) \ins D_\ocP \big\}$
 is stable under conditioning as follows.
 It will play an important role in deriving the sub-solution side  of the DPP for $\oV$.

 \begin{prop} \label{prop_flow}

 Given $ (t,\bx ) \ins [0,\infty) \ti \OmX   $,
 let  $\oga $ be a  $[t,\infty)-$valued $\bF^{\oW^t}-$stopping time and let $\oP \ins \ocP_{t,\bx}$.
 There exists a $\oP-$null set $\ocN    $ such that
  \bea   \label{062821_11}
  \oP^t_{\oga,\oo} \ins \ocP_{ \oga(\oo),\oX_{\oga \land \cd}  (\oo)} \Big(  \big(\oY_{\n \oP}  (\oga )\big)    (\oo), \big(\oZ_\oP  (\oga )\big)    (\oo) \Big), \q \fa \oo \ins  \big\{\oT \gs   \oga  \big\}  \Cp  \ocN^c  .
  \eea

 \end{prop}

 Now, we are ready to present a dynamic programming principle in weak formulation for
the value function $\oV$, 
in which  $ \big(\oY_{\n \oP} (\oga),\oZ_\oP (\oga)\big) $  act  as   additional states for constraint levels at the intermediate horizon $\oga$.

\begin{thm} \label{thm_DPP1}

 Given $ (t,\bx,y,z) \ins D_\ocP  $,
 let  $\big\{\ogaP   \big\}_{\oP \in \ocP_{t,\bx}(y,z)}$ be a   family of $[t,\infty)-$valued  $\bF^{\oW^t}  - $stopping times. Then
 \beas
 &&\hspace{-1.5cm}
 \oV(t,\bx,y,z)   \=   \Sup{\oP \in \ocP_{t,\bx}(y,z)} \n E_\oP \bigg[
 \b1_{\{\oT <  \ogaP \}} \Big( \n \int_t^\oT \n f(r,\oX_{r \land \cd},\oU_r) dr  \+ \pi \big(\oT, \oX_{\oT \land \cd}\big) \Big) \nonumber \\
 &&\hspace{2.3cm} + \b1_{\{\oT \ge  \ogaP \}} \bigg( \n \int_t^\ogaP \n f(r,\oX_{r \land \cd},\oU_r) dr  \+ \oV \Big(  \ogaP ,\oX_{ \ogaP  \land \cd} ,    \oY_{\n \oP}  \big( \ogaP \big) , \oZ_\oP  \big( \ogaP \big)   \Big) \bigg) \bigg] .     
  \eeas

 \end{thm}

\section{Proofs}   \label{sec_proof}

 \if{0}

   \begin{lemm} \label{Jul28_01}
  Let $(X,\fT)$ be a topological space with a subbase $\fB$   of the topology $\fT$
  and let $ \wt{\fB} $ be a sub-collection of $\fT$. If each set in $\fB$ is a union of sets in $ \wt{\fB} $,
  then $ \wt{\fB} $ is another subbase of $\fT$.
  \end{lemm}

  \ss \no {\bf Proof:}  Let $\cO \ins \fT$, so  $\cO \= \underset{i \in \cI}{\cup} \,
  \underset{j=1}{\overset{N_i}{\cap}} A^i_j$ for some $ \ccup{i \in  \cI}{}\big\{A^i_j \n :    j \= 1,\cds, N_i \big\} \sb  \fB$.
  For $i \ins \cI$ and $j \= 1,\cds, N_i$,  $ A^i_j \= \underset{k \in \cK^i_j}{\cup} \wA^{i,j}_k $
  for some $ \big\{\wA^{i,j}_k \n : k \ins \cK^i_j \big\} \sb \wt{\fB}$. Then
  \beas
  \hspace{3.8cm}
  \cO \=  \underset{i \in \cI}{\cup} \Big( \underset{j=1}{\overset{N_i}{\cap}} \, \underset{k \in \cK^i_j}{\cup} \wA^{i,j}_k \Big)
  \= \underset{i \in \cI}{\cup} \,\underset{(k_1,\cds,k_{N_i}) \in \cK^i_1 \times \cds \cK^i_{N_i} }{\cup}
  \wA^{i,1}_{k_1} \Cp \cds \Cp \wA^{i,N_i}_{k_{N_i}} . \hspace{3.8cm} \hb{\qed}
  \eeas

 \fi

     \ss \no {\bf Proof of Lemma \ref{lem_082020_11}:}
       Proposition 7.19 of  \cite{Bertsekas_Shreve_1978} shows that
      the topology of weak convergence $\fT_\sharp\big(\fP \big([0,\infty) \ti \hU\big)\big)$ on $\fP \big([0,\infty) \ti \hU\big)$
     can be generated by a subbase $\L  \df \big\{ O_\d (\fm, \phi_j ) \n : \d \ins (0,\infty)  , \fm \ins \fP \big([0,\infty) \ti \hU\big), j \ins \hN  \big\}  $,    where $ \{ \phi_j\}_{j \in \hN}$ is a countable dense subset of $\wh{C}_b\big([0,\infty) \ti \hU\big)$.
    As the Borel space $\big(\fP  ([0,\infty) \ti \hU ), \fT_\sharp\big(\fP  ([0,\infty) \ti \hU )\big)\big) $ is separable, it has a   countable dense subset $\{\fm_k\}_{k \in \hN}$.

  To show that $\wt{\L}  \df \big\{ O_{\frac{1}{n}} (\fm_k, \phi_j ) \n : n, k, j \ins \hN  \big\}  $
  is another subbase of $\fT_\sharp\big(\fP \big([0,\infty) \ti \hU\big)\big)$, it suffices to verify that
  any member of $\L$ is a union of some members in $\wt{\L}$:
    Let $ (\d,\fm,j) \ins (0,\infty) \ti \fP \big([0,\infty) \ti \hU\big) \ti \hN $ and let $\fm' \ins O_\d (\fm, \phi_j )$.
    There exists $n   \ins \hN $ such that $ \frac{2}{n}     \< \d \- \big| \int_0^\infty
 \n   \int_\hU    \phi_j (t, u) \big(  \fm'(dt,du) \-   \fm(dt,du)  \big)    \big| $,
 and one can find $\fm_k \ins O_{\frac{1}{n}} (\fm';\phi_j) $ for some  $k   \ins \hN$.
 For any $\fm'' \ins O_{\frac{1}{n}} (\fm_k,\phi_j)$,  we can deduce  that
   $ \big| \int_0^\infty
 \n   \int_\hU    \phi_j (t, u) \big(  \fm''(dt,du) \-   \fm(dt,du)  \big)\big|
   \ls   \big| \int_0^\infty
 \n   \int_\hU    \phi_j (t, u) \big(  \fm''(dt,du) \-   \fm_k (dt,du)  \big)\big|
\+  \big| \int_0^\infty
 \n   \int_\hU    \phi_j (t, u) \big(  \fm_k(dt,du) \-   \fm'(dt,du)  \big)\big|
 \+  \big| \int_0^\infty
 \n   \int_\hU    \phi_j (t, u) \big(  \fm'(dt,du) \-   \fm(dt,du)  \big)\big|
  \ls \frac{2}{n} \+ \big| \int_0^\infty
 \n   \int_\hU    \phi_j (t, u) \big(  \fm'(dt, \\ du) \-   \fm(dt,du)  \big)\big| \< \d $,
    which implies that  $\fm' \ins O_{\frac{1}{n}} (\fm_k,\phi_j)   \sb O_\d (\fm,\phi_j )$.
   \qed

\no {\bf Proof of Lemma   \ref{lem_Nov25_03}: 1)}   We first show that $\hJ$ is a complete separable space under the metric
\beas 
\Rho{\hJ} (\fu,\fu') \df  \int_0^\infty \n  e^{-s} \big( 1 \ld  \Rho{\hU} (\fu(s),\fu'(s))  \big)   ds, \q \fa \fu,\fu' \ins \hJ .
\eeas

    \if{0}

    Given $\fu,\fu' \ins \hJ$,   $ \Rho{\hJ} (\fu,\fu') \= 0 $
   iff $ \Rho{\hU} (\fu(s),\fu'(s)) \= 0 $ for a.e. $s \ins (0,\infty)$
   iff $  \fu(s) \= \fu'(s)   $ for a.e. $s \ins (0,\infty)$ or $\fu \= \fu' \ins \hJ$.
  Also, let $ \fu'' \ins \hJ$.  Since $ 1 \ld (a\+b) \ls 1 \ld a \+ 1 \ld b$ for $a,b \ins [0,\infty)$,  
  one has
 \beas
 1 \ld  \Rho{\hU} (\fu(s),\fu''(s))
 \ls 1 \ld \big(  \Rho{\hU} (\fu(s),\fu'(s)) \+ \Rho{\hU} (\fu'(s),\fu''(s)) \big)
 \ls 1 \ld \big(  \Rho{\hU} (\fu(s),\fu'(s))   \big)
 \+ 1 \ld \big(   \Rho{\hU} (\fu'(s),\fu''(s)) \big)
 , \q \fa s \ins (0,\infty) .
 \eeas
 It follows that
 \beas
 \Rho{\hJ} (\fu,\fu'') & \tn \=  & \tn     \int_0^\infty \n  e^{-  s} \big( 1 \ld  \Rho{\hU} (\fu(s),\fu''(s))  \big)   ds
 \ls    \int_0^\infty \n  e^{-  s} \big( 1 \ld  \Rho{\hU} (\fu(s),\fu'(s))  \big)   ds
 \+    \int_0^\infty \n  e^{-  s} \big( 1 \ld  \Rho{\hU} (\fu'(s),\fu''(s))  \big)   ds   \\
  & \tn \= & \tn  \Rho{\hJ} (\fu,\fu') \+ \Rho{\hJ} (\fu',\fu'') .
 \eeas
 So $ \Rho{\hJ} $ is a metric on $\hJ$.

    \fi

     \ss \no {\bf 1a)} Let $\{\fu_n\}_{n \in \hN}$ be a Cauchy sequence in $\big(\hJ, \Rho{\hJ}\big)$.
  We can find a subsequence $\{n_k\}_{k \in \hN} $ of $\hN$ such that
  $  \Rho{\hJ} \big(\fu_{n_k}, \fu_{n_{k+1}} \big) \< 2^{-k} $, $\fa k \ins \hN$.
   Fix $k  \ins \hN$. Since
   \bea \label{060722_11}
    1 \ld (a\+b) \ls 1 \ld a \+ 1 \ld b  ,   \q  \fa a,b \ins [0,\infty)  ,
    \eea 
      it holds for any $\ell \ins \hN$ that
  $ 1 \ld \Rho{\hU} (\fu_{n_k}(s),\fu_{n_{k+\ell}}(s))
   \ls \sum^\ell_{i=1}   1 \ld  \Rho{\hU} (\fu_{n_{k+i-1}}(s),\fu_{n_{k+i}}(s))
  $,   $ s \ins (0,\infty) $,
  taking $ \underset{\ell \in \hN}{\sup} $   yields that
  \bea  \label{050921_23}
  1 \ld \Big( \underset{\ell \in \hN}{\sup} \,   \Rho{\hU} \big(\fu_{n_k}(s),\fu_{n_{k+\ell}}(s)\big) \Big)
  \= \underset{\ell \in \hN}{\sup} \Big( 1 \ld \Rho{\hU} \big(\fu_{n_k}(s),\fu_{n_{k+\ell}}(s)\big) \Big)
    \ls \sum_{i \in \hN}   1 \ld  \Rho{\hU} \big(\fu_{n_{k+i-1}}(s),\fu_{n_{k+i}}(s)\big)   , \q s \ins (0,\infty) .
  \eea 
  The monotone convergence theorem then implies that
  $  \int_0^\infty e^{-s} \Big( 1 \ld \Big( \underset{\ell \in \hN}{\sup} \,   \Rho{\hU} (\fu_{n_k}(s),\fu_{n_{k+\ell}}(s)) \Big) \Big)  ds
   \ls  
   \sum_{i \in \hN} \Rho{\hJ} \big(\fu_{n_{k+i-1}}, \\ \fu_{n_{k+i}}  \big)
  \ls \sum_{i \in \hN} 2^{-k-i+1} \= 2^{-k+1} $.
 So $\big\{  n_k  \big\}_{k \in \hN}$ has a subsequence $\big\{  n_{k_m}  \big\}_{m \in \hN}$ such that
 $ \lmt{m \to \infty} \, e^{-s} \Big( 1 \land \Big( \underset{\ell \in \hN}{\sup} \,   \Rho{\hU} (\fu_{n_{k_m}}(s), \\ \fu_{n_{k_m+\ell}}(s)) \Big) \Big)  \= 0 $
 or $ \lmt{m \to \infty}    \Big( \underset{\ell \in \hN}{\sup} \,   \Rho{\hU}   \big(\fu_{n_{k_m}}(s),   \fu_{n_{k_m+\ell}}(s)\big) \Big)   \= 0 $
   for all $s \ins (0,\infty) $ except on a   $ds  -$null set  
   $\cN$ of $(0,\infty)$.
   Given $s \ins (0,\infty)\backslash \cN $, one has
   $   \lmt{m \to \infty}    \bigg( \underset{j \in \hN}{\sup} \,   \Rho{\hU} \big(    \fu_{n_{k_m}}(s) ,
      \fu_{n_{k_{j+m}}}(s)   \big) \bigg) \= 0$,
   i.e., $ \big\{  \fu_{n_{k_m}}(s)  \big\}_{m \in \hN} $ is a Cauchy sequence
   in $\hU$. Let $ \fu_o (s) $ be the limit of $ \big\{  \fu_{n_{k_m}}(s)  \big\}_{m \in \hN} $ in $\big(\hU, \Rho{\hU}\big)$.

   Define    $ \ul{\mu}  (s) \df \linf{m \to \infty} \sI \big( \fu_{n_{k_m}} (s) \big) \ins [0,1] $, $s \ins [0,\infty) $,
   which is   a  Borel measurable function  on $[0,\infty)$. Then
 $ \fu_* (s)\df  \sI^{-1} \big(\ul{\mu}(s) \big)
 \b1_{ \{\ul{\mu}(s) \in \fE\}} \+  u_0  \b1_{ \{\ul{\mu}(s) \notin \fE\}}
   $, $ s \ins [0,\infty)$
  is   a $\hU-$valued Borel measurable function  on $[0,\infty)$, i.e., $\fu_* \ins \hJ$.
 For $s \ins (0,\infty)\backslash \cN$, 
 the continuity of   mapping $ \sI   $ shows
 $ \ul{\mu} (s) \= \lmt{m \to \infty} \sI \big(\fu_{n_{k_m}}(s)\big) \= \sI   \big( \fu_o  (s) \big) \ins \fE $, so   $ \fu_* (s) \= \sI^{-1} \big(\ul{\mu}(s) \big) \= \fu_o  (s)  $. 
The dominated convergence theorem implies   $ \lmt{m \to \infty} \Rho{\hJ} \big( \fu_{n_{k_m}} , \fu_*\big)
  \=  \lmt{m \to \infty}  \int_0^\infty \n \b1_{\{s \in \cN^c\}}    e^{-  s}   \big( 1 \ld  \Rho{\hU}  (  \fu_{n_{k_m}}(s) ,  \fu_o  (s)   )  \big)   ds \= 0  $.

   Let $\e \ins (0,1)$. There exists a $N \ins  \hN$ such that $ \Rho{\hJ} \big( \fu_{\fn } , \fu_{\fn'}\big) \< \e/2 $ for any $\fn ,\fn' \gs N$. We can also find a $m \ins \hN$ such that $n_{k_m} \gs N $ and that
 $ \Rho{\hJ} \big( \fu_{n_{k_m}} , \fu_*\big) \< \e/2 $.
 It then  holds  for any $n \ins \hN$ with $n \gs   N  $ that
 \bea \label{050921_25}
 \Rho{\hJ} \big( \fu_n , \fu_*\big)
 \ls \Rho{\hJ} \big( \fu_n , \fu_{n_{k_m}} \big) \+
 \Rho{\hJ} \big( \fu_{n_{k_m}} , \fu_*\big) \< \e  .
 \eea
  Hence,    $\fu_*$ is the limit of $\{\fu_n\}_{n \in \hN}$ in $\big(\hJ, \Rho{\hJ}\big)$.

 \ss \no {\bf 1b)} In this step, we   demonstrate   that $\hJ$ is separable under $ \Rho{\hJ} $.

 Let $\{u_i\}_{i \in \hN}$ be a countable dense subset of $\big(\hU, \Rho{\hU}\big)$ and
 let $\{O_i\}_{i \in \hN}$ be a countable base of the Euclidean topology on $[0,\infty)$.
      Given $n \in \hN$,  let us enumerate the $2^n$ elements of $   \big\{ \underset{i \in  I}{\cup} O_i \n :  I \sb \{1,\cds,n\} \big\}$   by $ \big\{ \breve{O}^n_1, \cds , \breve{O}^n_{2^n} \big\}$ and consider
     the following countable collections of  $\sB[0,\infty)/\sB[-1,1]-$measurable functions on $[0,\infty)$:
     \beas
     \sC_n \df \Big\{    \sI   ( u_j) \b1_{\breve{O}^n_i} \- \b1_{(\breve{O}^n_i)^c}   \n : i \in \{1,\cds,2^n\},~ j  \ins \hN   \Big\} , \q
     \wt{\sC}_n \df \Big\{ \underset{i = 1,\cds, k}{ \max }   \fx_i    \n :  k \ins \hN, ~ \{ \fx_1,   \cds, \fx_k \} \sb \sC_n \Big\} .
     \eeas
     Clearly, each   $\fx \ins \wt{\sC}_n$ takes   values  in a finite subset of $ \big\{\sI ( u_i)\big\}_{i \in \hN}   \cp \{-1\} $.
     By additionally assigning $\sI^{-1} (-1) \df u_0 $, one has
     $ \wt{\bC} \df \big\{ \sI^{-1} (\fx) \n : \fx \ins \wt{\sC}_n \hb{ for some }  n \ins \hN  \big\}
     \sb \big\{\hb{$\{u_i\}^\infty_{i = 0}-$valued $\sB[0,\infty)-$measurable  functions}\big\}
     \n =:\n \wh{\bC} $.

   We claim that   $ \wh{\bC}  $ is a dense subset of $ \hJ $ under $\Rho{\hJ}$.
 To see this,  for any $i,n \ins \hN$ we set $o^n_i \df \big\{u \ins \hU: \Rho{\hU} (u,u_i) \< 2^{-n} \big\}$
 as  the open ball centered at $u_i$ with radius $2^{-n}$.
 We also set   $\wt{o}^n_1 \df o^n_1 $ and $\wt{o}^n_i \df o^n_i \big\backslash \Big( \underset{j < i }{\cup} o^n_j \Big) $ for $i \gs 2$.
  Given $ \fu  \ins \hJ  $ and $ n  \ins  \hN$,   define a member of  $\wh{\bC}$
  by    $ \fu_n (s) \df   \sum_{i \in \hN} \b1_{\{  s   \in  \cE^n_i \}}  u_i $,
  $  s \ins [0,\infty)  $, where
    $\cE^n_i \df \big\{s \ins [0,\infty) \n :   \fu  (s) \ins \wt{o}^n_i \big\} \ins \sB[0,\infty)$.
  As   $  \Rho{\hU} \big(\fu_n (s), \fu (s) \big)
    \=   \sum_{i \in \hN} \b1_{\{  s   \in \cE^n_i \}}  \Rho{\hU} \big( u_i , \fu (s)   \big)
   \< 2^{-n} $ for any  $      s \ins [0,\infty) $,
 one has $ \Rho{\hJ}   \big(\fu_n  , \fu   \big) \=  \int_0^\infty \n  e^{-s} \big( 1 \ld  \Rho{\hU} (\fu_n(s),\fu(s))  \big)   ds
 \ls 
 2^{-n} $. So  $ \wh{\bC}  $ is a dense subset of $ \hJ $ under $\Rho{\hJ}$.

    We then show that the countable collection $\wt{\bC}$ is dense in $\wh{\bC}$ 
 and is thus dense in $\hJ$ under $\Rho{\hJ}$:
  Let $\wh{\fu} \ins \wh{\bC} $ and  $\e \ins (0,1)$.
  We set    $\wh{\cE}_i \= \big\{s \ins [0,\infty) \n : \wh{\fu}(s) \= u_i  \big\} \ins \sB[0,\infty)$, $\fa i \ins \hN$.
  Since $ \l(\cE) \df   \int_{  s  \in \cE } \n  e^{-  s}    ds $
  is a probability measure  on  $\big([0,\infty), \sB[0,\infty)\big)$, there exists  $N \ins \hN$ such that
  $  \l \Big( \underset{i > N}{\cup} \wh{\cE}_i \Big) < \e/3 $.

  Given $i \=  1,\cds \n , N$,     Proposition 7.17 of \cite{Bertsekas_Shreve_1978} shows that
  $\l( \cO_i \backslash \wh{\cE}_i) \< \frac{\e}{3N} $ for  some open subset $ \cO_i$ of $[0,\infty)$ containing $\wh{\cE}_i$.
  Since $  \cO_i \= \underset{n \in \hN}{\cup} O_{\ell^{\,i}_n} $ for some  subsequence $\big\{\ell^{\,i}_n\big\}_{n \in \hN}$
  of $\hN$, one can find $ n_i \ins \hN$ such that
  $\l \big( \cO_i \backslash \breve{\cO}_i \big) <  \frac{\e}{3N} $,
  where $\breve{\cO}_i \df \underset{n  = 1}{\overset{ n_i}{\cup}} O_{\ell^{\,i}_n} \ins \sB[0,\infty) $.
 As $   \breve{\cO}_i    \ins \big\{ \breve{O}^\fn_j   \big\}^{2^\fn}_{j=1}$ for $ \fn \df  \underset{i  = 1,\cds,N}{\max} \ell^{\,i}_{n_i}  $,
  we see that $ \fx_i \df \sI   ( u_i ) \b1_{\breve{\cO}_i} \-    \b1_{\breve{\cO}^c_i}  $ belongs to  $\sC_\fn$.
 Define     $\fu \df \sI^{-1} \Big( \underset{i = 1,\cds, N}{\max} \,  \fx_i \Big) \ins \wt{\bC} $.

 For  $i \= 1,\cds  \n , N $, if   $  \cA_i \df \big( \wh{\cE}_i \Cp \breve{\cO}_i   \big) \big\backslash \Big( \underset{j \le N;j \ne i}{\cup}\breve{\cO}_j \Big) \ins \sB[0,\infty)$ is not empty,   it holds for any $s \ins \cA_i$ that
    $  \fx_j (s) \= \b1_{\{j = i\}} \sI   ( u_i) \- \b1_{\{j \ne i\}}$   for $j \ins \{1,\cds  \n ,N\}$ and thus
  $ \fu  (s) \= \sI^{-1} \big( \fx_i (s)  \big) \=   u_i  \= \wh{\fu} (s)$.
 Also, if $   \cA_0 \df \Big( \underset{i \in \hN}{\cup} \wh{\cE}_i \Big)^c \bigcap \Big( \underset{j = 1}{\overset{N}{\cup}} \breve{\cO}_j \Big)^c $ is not empty, it holds for any  $s \ins \cA_0$
 that  $ \fu(s) \= \sI^{-1}(-1)   \= u_0 \= \wh{\fu}(s) $.   Then
  $ \Rho{\hJ} \big(  \fu  , \wh{\fu} \big)
  \=    \int_{s \in \cA}   e^{-  s}   \big( 1 \ld \Rho{\hU} \big( \fu(s),\wh{\fu}(s) \big) \big)    ds $
  for  $\cA \df \Big( \underset{i=0}{\overset{N}{\cup}} \cA_i \Big)^c \= \Big(\underset{i=1}{\overset{N}{\cup}} (\wh{\cE}_i \backslash \cA_i)\Big) \bigcup  \Big( \underset{i > N}{\cup} \wh{\cE}_i\Big) \bigcup \bigg( \Big( \underset{i \in \hN}{\cup} \wh{\cE}_i \Big)^c \bigcap \Big( \underset{j = 1}{\overset{N}{\cup}} \breve{\cO}_j \Big)\bigg) \sb [0,\infty) $.
  Since
  \beas
    \wh{\cE}_i \backslash \cA_i \= \big(\wh{\cE}_i \Cp \breve{\cO}^c_i \big) \hb{$\bigcup$}
  \Big(   \underset{j \le N; j \ne i}{\cup} \big(\wh{\cE}_i  \Cp \breve{\cO}_i  \Cp \breve{\cO}_j \big)  \Big)
  \sb \big(\cO_i \backslash   \breve{\cO}_i \big) \hb{$\bigcup$} \Big(   \underset{j \le N; j \ne i}{\cup}   \breve{\cO}_j \backslash \wh{\cE}_j    \Big),
  \q  i \= 1,\cds, N
  \eeas
  and since $\Big( \underset{i \in \hN}{\cup} \wh{\cE}_i \Big)^c \bigcap \Big( \underset{j = 1}{\overset{N}{\cup}} \breve{\cO}_j \Big)
  \= \underset{j = 1}{\overset{N}{\cup}} \Big( \breve{\cO}_j \Cp \Big( \underset{i \in \hN}{\cup} \wh{\cE}_i \Big)^c \Big)
  \sb \underset{j = 1}{\overset{N}{\cup}} \big( \breve{\cO}_j \Cp   \wh{\cE}^c_j   \big) \sb \underset{j = 1}{\overset{N}{\cup}} \big( \cO_j \backslash  \wh{\cE}_j   \big) $,
  we can deduce that
  $ \cA 
  \sb   \Big(   \underset{i  = 1}{\overset{N}{\cup}}  \cO_i \backslash   \breve{\cO}_i   \Big) \bigcup \Big(   \underset{i  = 1}{\overset{N}{\cup}}   \cO_i \backslash \wh{\cE}_i    \Big) \bigcup  \Big( \underset{i > N}{\cup} \wh{\cE}_i\Big) $.
  It follows that    $  \Rho{\hJ} \big(  \fu  , \wh{\fu} \big)   \ls \l ( \cA  )  \< \e $.
  Namely, the countable collection $ \wt{\bC} $ is   dense in $\wh{\bC}$  under $ \Rho{\hJ} $.

 Therefore,   $\big(\hJ, \Rho{\hJ}\big)$ is a complete separable metric space.

 \if{0}

  Let $\fT(\hJ)$ be the   topology induced by $\Rho{\hJ}$,   which is stronger than $\fT_\sharp(\hJ)$.

 It suffices to show that each member of the countable subbase of $\fT_\sharp(\hJ)$ in \eqref{J04_01} is included in $\fT(\hJ)$:
 Let $n,k,j  \ins \hN$ and $\fu \ins \fri^{-1}_\hJ \big( O_{\frac{1}{n}} (\fm_k , \phi_j ) \big) $.
 As $ \phi_j \ins \wh{C}_b\big([0,\infty) \ti \hU\big) $, there exists $c_j \ins (0,\infty)$ such that
 $\big|\phi_j(t,u)\big| \ls c_j$ and $ \big|\phi_j(t,u)\- \phi_j(t',u') \big| \ls c_j \big(|t\-t'|\+ \Rho{\hU}(u,u') \big) $
 for any $(t,u),(t',u') \ins [0,\infty) \ti \hU $. In particular, one has
 \bea \label{050822_21}
  \big|\phi_j(t,u)\- \phi_j(t,u') \big| \ls 2c_j \big( 1 \ld \Rho{\hU}(u,u') \big) , \q \fa t  \ins [0,\infty) , ~ \fa u,u' \ins \hU .
  \eea

  Set $\e \df \frac{1}{n} \- \big| \int_0^\infty e^{-t}   \phi_j (t,\fu(t)) dt \- \int_0^\infty \n \int_\hU
  \phi_j (t,u) \fm_k (dt,du)     \big| \> 0$.
  Given $\fu' \ins \hJ$ with $\Rho{\hJ} (\fu',\fu) \< \frac{\e}{2 c_j   } $,   we can   deduce from \eqref{050822_21} that
 \beas
 && \hspace{-2cm}  \bigg| \int_0^\infty e^{-t}   \phi_j (t,\fu'(t)) dt \- \int_0^\infty \n \int_\hU
  \phi_j (t,u) \fm_k (dt,du)     \bigg|  \\
 && \ls   \int_0^\infty e^{-t}  \big| \phi_j (t,\fu'(t))   \-    \phi_j (t,\fu(t)) \big| dt
   \+ \bigg| \int_0^\infty e^{-t}   \phi_j (t,\fu(t)) dt \- \int_0^\infty \n \int_\hU
  \phi_j (t,u) \fm_k (dt,du)     \bigg|  \\
 && \ls 2 c_j \int_0^\infty e^{-t}  \big( 1 \ld \Rho{\hU}(\fu'(t),\fu(t)) \big) dt  \+   \frac{1}{n}  \- \e
    \< \frac{1}{n} ,
 \eeas
i.e., $ \fu' \ins \fri^{-1}_\hJ \big( O_{\frac{1}{n}} (\fm_k , \phi_j ) \big)$. So  $\fri^{-1}_\hJ \big( O_{\frac{1}{n}} (\fm_k , \phi_j ) \big) \ins \fT(\hJ)$.

 \fi

  \ss \no {\bf 2)} We next  show that the mapping $\friJ \n : \big(\hJ, \Rho{\hJ}\big) \mto \Big(\fP \big([0,\infty) \ti \hU\big), \fT_\sharp\big(\fP \big([0,\infty) \ti \hU\big)\big)\Big)$  is a continuous injection.

 \ss \no {\bf 2a)} Let $\fu ,\fu' \ins \hJ$ such that $\friJ(\fu) \= \friJ(\fu')$. Since   $ \int_0^t e^{-s}   \sI \big(\fu(s)\big)   ds \=
\int_0^\infty \int_\hU \b1_{\{s \le t\}} \sI (u) \friJ(\fu) (ds,du)
\= \int_0^\infty \int_\hU \b1_{\{s \le t\}} \sI (u) \friJ(\fu') \\ (ds,du) \= \int_0^t e^{-s}   \sI \big(\fu'(s)\big)   ds $ for any $t \ins (0,\infty)$,
 we see that $ e^{-s}   \sI \big(\fu(s)\big) \= e^{-s}   \sI \big(\fu'(s)\big) $ and thus
$ \fu(s) \= \fu'(s) $ for a.e. $s \ins (0,\infty)$, namely,  $\fu \= \fu'$ in $\hJ$.  So $ \friJ $ is an injection.

 \ss \no {\bf 2b)} 
 Let $\{\fu_n \}_{n \in \hN}  $ be a sequence of $\hJ$ converging to a $\fu \ins \hJ$ under $\Rho{\hJ}$
 and  we show  that $ \friJ(\fu_n)$ converges to $ \friJ(\fu)$ under $\fT_\sharp\big(\fP \big([0,\infty) \ti \hU\big)\big)$.
  By e.g. Lemma 7.6 of \cite{Bertsekas_Shreve_1978},
 this is equivalent to verify that
 \bea \label{050921_21}
  \lmt{n \to \infty} \int_0^\infty \n \int_\hU  \phi (s,u) \friJ(\fu_n)(ds,du)
  \=  \int_0^\infty \n \int_\hU  \phi (s,u) \friJ(\fu)(ds,du)
 \eea
  for any bounded continuous function $\phi \n : [0,\infty) \ti \hU \mto \hR$.

    Let $\phi$ be such a continuous function on $ [0,\infty) \ti \hU $.
  For the limit \eqref{050921_21}, it suffices to show that   any subsequence $\big\{\fu_{n_k}\big\}_{k \in \hN}$ of
  $\{\fu_n \}_{n \in \hN}  $ has in turn a subsequence  $\big\{\fu_{n'_k}\big\}_{k \in \hN}$ 
  satisfying \eqref{050921_21}:
  As $   \lmt{k \to \infty}   \int_0^\infty \n  e^{-  s} \big( 1 \ld  \Rho{\hU} (\fu_{n_k}(s),\fu(s))  \big)   ds
   \=  \lmt{k \to \infty} \Rho{\hJ} (\fu_{n_k},\fu ) \= 0  $,
  there exists a subsequence $\big\{\fu_{n'_k}\big\}_{k \in \hN}$ of   $\big\{\fu_{n_k}\big\}_{k \in \hN}$ such that
   $ \lmt{k \to \infty} \Rho{\hU} \big(\fu_{n'_k}  (s),\fu   (s) \big) \= 0 $ for a.e.
   $s \ins (0,\infty)   $.
  Then   the continuity of $ \phi $ and the dominated convergence theorem imply that
   $ \lmt{k \to \infty} \int_0^\infty \n \int_\hU  \phi (s,u) \friJ \big(\fu_{n'_k}\big)(ds,du)
   \=   \lmt{k \to \infty} \int_0^\infty     e^{-s} \phi \big(s, \fu_{n'_k}(s)  \big)      ds
   \= \int_0^\infty     e^{-s} \phi \big(s, \fu(s)  \big)      ds
   \= \int_0^\infty \n \int_\hU  \phi (s,u) \friJ(\fu)(ds,du) $.

 So   $\friJ $ is a continuous injection from the complete separable metric space $\big(\hJ, \Rho{\hJ}\big)$   to the Borel space $ \Big(\fP \big([0,\infty) \ti \hU\big), \fT_\sharp\big(\fP \big([0,\infty) \ti \hU\big)\big)\Big)$.
 (In particular,   the   topology induced by $\Rho{\hJ}$  is stronger than $\fT_\sharp(\hJ)$.)
 Then  the image $\friJ(\hJ)$ is   a Lusin subset and thus a Borel subset
 of $ \Big( \fP \big([0,\infty) \ti \hU\big), \fT_\sharp \big( \fP \big([0,\infty) \ti \hU\big) \big)\Big) $,
 see e.g. Theorem A.6 of \cite{Takesaki_1979}.

  As the embedding mapping $\friJ$ is clearly  a  homeomorphism between
  $\big(\hJ, \fT_\sharp(\hJ)\big)$ and $\friJ (\hJ)$ with the relative topology to $ \fT_\sharp\big(\fP \big([0,\infty) \ti \hU\big)\big) $,
  we obtain that $\big(\hJ, \fT_\sharp(\hJ)\big)$ is   a Borel space.        \qed

 \if{0}

  As $\hJ$ is a Polish space and thus a Lusin space,   there exist a Polish space $\hX$
 and a continuous bijection $   \phi_{\overset{}{\hX}} $ from $\hX$   onto $\hJ$  such that
 every point $\fx$ of $\hX$ has a fundamental system of neighborhood of open and closed sets.
 Then the image of the continuous composite injection $ \friJ \circ \phi_{\overset{}{\hX}}$ is a Lusin subset of $ \Big( \fP \big([0,\infty) \ti \hU\big), \fT_\sharp \big( \fP \big([0,\infty) \ti \hU\big) \big)\Big) $.
 In light of  Theorem A.6 of \cite{Takesaki_1979},
 $\friJ \circ \phi_{\overset{}{\hX}} (\hX) \= \friJ(\hJ)$ is a Borel subset of $ \Big( \fP \big([0,\infty) \ti \hU\big), \fT_\sharp \big( \fP \big([0,\infty) \ti \hU\big) \big)\Big) $.

  Lusin space is defined in the book "General Topology" (Part. II, page 205) or in the notes "Polish Spaces and Standard Borel Spaces".  It can   alternatively be characterized as a   topological space that is homeomorphic to a Borel subset of a compact metric space
  (see https://en.wikipedia.org/wiki/Polish_space#Lusin_spaces)
  (see also line 3 in the proof of Theorem A.13 in "Polish Spaces and Standard Borel Spaces", pg 382 therein).

  A metrizable topological space $\hY$ is called a Lusin space if there exist a Polish space $\hX$
  and a bijective continuous mapping from $\hX$   onto $\hY$  such that
  every point $\fx$ of $\hX$ has a fundamental system of neighborhood of open and
  closed sets (i.e., a family $\{\cE_i\}_{i \in \cI}$ of open and closed sets such that
  each $\cE_i$ is a neighborhood of $\fx$ and that every neighborhood of $\fx$ contains some  $\cE_i$).

  A subset of a topological space is called a Lusin set if it is a Lusin space as  a topological space (with the relative topology). Any countable intersection of Lusin sets in a topological  space is a Lusin set.

 https://proofwiki.org/wiki/Metrizable_Space_is_Hausdorff

 https://math.stackexchange.com/questions/2880098/continuous-bijection-is-homeomorphism

 http://mathonline.wikidot.com/homeomorphisms-between-compact-and-hausdorff-spaces

 a continuous bijection that is not Homeomorphism

 https://joemathjoe.wordpress.com/2018/05/19/homeomorphism-is-not-just-continuous-bijection/

  \fi

\no {\bf Proof of Lemma \ref{lem_M29_01}: 1)}
 For any $t \ins (0,\infty)$ and $\phi \ins C_b \big([0,\infty) \ti \hU\big)$,
 set $\cI^\phi_t(\fu) \df \int_0^t \phi (s,  \fu (s) ) ds$, $\fa \fu \ins \hJ$.
   We first show that $\fF \df \si \big(\cI^\phi_t  ; t \ins (0,\infty), \phi \ins C_b ([0,\infty) \ti \hU) \big) \= \sB(\hJ)$.

 \no {\bf 1a)}   Let $t \ins (0,\infty)$, $\phi \ins C_b\big([0,\infty) \ti \hU\big)$ and $a \ins \hR$.
  To see that $A \df \{\fu \ins \hJ \n : \cI^\phi_t (\fu) \< a \}$ belongs to  $\fT_\sharp(\hJ)$,
  we let $\fu \ins A$ and set $\e \df \frac12 \big(a\-\cI^\phi_t (\fu)\big)$.
  There exists  a positive bounded continuous function $\beta $ on $[0,\infty)$ such that
 $\beta(s) \= e^s$, $\fa s \ins [0,t]$ and $\int_t^\infty e^{-s} \beta(s) ds \ls  \frac12 \e \big(\|\phi\|_\infty \ve 1\big)^{-1} $,
 where $\|\phi\|_\infty \df \Sup{(s,u) \in [0,\infty) \times \hU} |\phi(s,u)| $.
 Set  $   \cO \df \fri^{-1}_\hJ \big( O_\e (\friJ(\fu),\beta\phi) \big)$, which clearly contains $\fu$.
 Let $\fu' \ins \cO $. We can deduce that
 \beas
   \e    \>     \Big| \int_0^\infty \n e^{-s} \beta(s) \big(  \phi (s, \fu'(s) )   \-  \phi (s, \fu(s) ) \big)  ds\Big|
   \gs   \Big| \int_0^t   \big(  \phi (s, \fu'(s) )   \-  \phi (s, \fu(s) ) \big)   ds \Big| \-
  \Big| \int_t^\infty e^{-s} \beta(s) \big(  \phi (s, \fu'(s) )   \-  \phi (s, \fu(s) ) \big)  ds   \Big| .
 \eeas
 Since
 $ \big| \int_t^\infty \n e^{-s} \beta(s) \big(  \phi (s, \fu'(s) )   \-  \phi (s, \fu(s) ) \big)  ds   \big|
 \ls 2 \|\phi\|_\infty \int_t^\infty e^{-s} \beta(s) ds \ls \e $,
 one has  $ \big| \cI^\phi_t (\fu') \- \cI^\phi_t (\fu) \big|
 \< 2 \e $, which implies that $ \cI^\phi_t (\fu') \< a$ or $\fu' \ins A$.
 So $\fu \ins \cO \sb A$.  As the induced topology $\fT_\sharp(\hJ)$   is   generated by the subbase
 $ \Big\{ \fri^{-1}_\hJ \big( O_\d (\fm,\phi ) \big)   $,
 $ \fa  (\d, \fm, \phi )  \ins (0,\infty) \ti \fP \big([0,\infty) \ti \hU\big) \ti C_b \big([0,\infty) \ti \hU\big) \Big\}$,
 we see that $A \ins \fT_\sharp(\hJ) \sb \sB(\hJ) $.
 It follows that $\{\fu \ins \hJ \n : \cI^\phi_t (\fu) \ins \cE \} \ins \sB(\hJ)$ for any $\cE \ins \sB(\hR)$
 and  thus   $ \fF \= \si \big(\cI^\phi_t  ; t \ins (0,\infty), \phi \ins C_b([0,\infty) \ti \hU)\big)   \sb \sB(\hJ) $.

 \no {\bf 1b)} Let $\big\{\fri^{-1}_\hJ \big(O_{\frac{1}{n}} (\fm_k,\phi_j )\big)\big\}_{n,k,j \in \hN}$
be the   countable subbase of $\fT_\sharp(\hJ)$ in \eqref{J04_01}.
 Given $n,k,j  \ins  \hN$, let $\{I^j_t\}_{t \in (0,\infty)}$ denote $\{\cI^\phi_t\}_{t \in (0,\infty)}$ with $\phi (s,u) \= e^{-s} \phi_j (s,u)$,
 $\fa (s,u) \ins [0,\infty) \ti \hU$  and set $\k_{k,j} \df \int_0^\infty   \int_\hU   \phi_j (t,u) \fm_k (dt,du)      \ins \hR$.
 Since the function $ I^j_*(\fu) \df \linf{t \to \infty} I^j_t (\fu) $, $\fu \ins \hJ$ is $\fF-$measurable and
 since $ \int_0^\infty e^{-s}   \phi_j (s,\fu(s))   ds  \= \lmt{t \to \infty} \int_0^t  e^{-s}   \phi_j (s,\fu(s))   ds
 \= \lmt{t \to \infty} I^j_t (\fu) \= I^j_*(\fu) $ holds for any $  \fu \ins \hJ$,
 \eqref{J04_01} shows that
$ \fri^{-1}_\hJ \big(O_{\frac{1}{n}} (\fm_k,\phi_j )\big)
  \= \big\{\fu \ins \hJ \n : \big| I^j_*(\fu) \- \k_{k,j} \big| \< 1/n \big\} \ins \fF $.
 Then $ \fT_\sharp(\hJ) \sb \fF $ and thus $\sB(\hJ) \= \fF $. 

 \no {\bf 2)} We demonstrate the second statement in several steps. Then   $  \sB(\hJ)    \=  \si \big(I_\vf  ;  \vf \ins L^0 \big((0,\infty) \ti \hU ; \hR \big) \big)   $ easily follows.

  Let $\{\omX^i\}_{i \in \hN}$ be a countable dense subset of $\OmX$
  and let $\big\{\,\nxi_j \= (w_j,x_j)\big\}_{j \in \hN}$ be a countable dense subset of $\hR^{d+l}$.
 Given $n \ins \hN$, we set $O^n_i \df \big\{\omX \ins \OmX \n : \Rho{\OmX}(\omX,\omX^i) \< 1/n \big\} \ins \sB(\OmX)$
   and  $\wt{O}^n_i \df O^n_i \big\backslash \Big( \ccup{i'<i}{} O^n_{i'} \Big) \ins \sB(\OmX)$ for any   $i \ins   \hN$.
 We also denote $\cE^n_j \df O_{\frac{1}{n}} (\nxi_j) \Big\backslash \Big( \ccup{j'<j}{} O_{\frac{1}{n}} (\nxi_{j'}) \Big) $
  for any $j \ins \hN$.

\no {\bf 2a)} Let $\phi \n: [0,\infty) \ti \OmX \ti \hR^{d+l} \ti \hU \mto [0,c_\phi] $  be a  continuous function
  for some $c_\phi \ins (1,\infty)$.

 Let $T \ins (0,\infty)$ and $m \ins \hN$. We set  $t^m_k \df k 2^{-m} T$, $\fa k \ins  \{0, 1, \cds, 2^m\}$.
 Then, we let $k \ins \{0, 1, \cds, 2^m\-1\}$.

  For any  $i,j \ins \hN$, as the function $\phi_{i,j} (s,u) \df   \phi (s,\omX^i , \nxi_j ,u) $, $   (s,u) \ins [0,\infty) \ti \hU $   is of  $C_b \big([0,\infty) \ti \hU\big)$,
    we know from Part (1)  that
    $I^{m,k}_{i,j} (\fu) \df \cI^{\phi_{ij}}_{t^m_{k+1}}(\fu) \- \cI^{\phi_{ij}}_{t^m_k}(\fu)
  \= \int_{t^m_k}^{t^m_{k+1}}  \phi \big(r,\omX^i , \nxi_j ,\fu(r)\big) dr 
  $, $\fa \fu \ins \hJ$ is $\sB(\hJ)-$measurable.
   Since the function $   \fl_2 (t^m_k,\omX) $ is continuous in $\omX \ins \OmX$
 and since the function $ \cW^m_k(\o_0,\omX)\df \big(\o_0(t^m_k), \omX(t^m_k)\big)  $ is continuous in $(\o_0,\omX) \ins \O_0 \ti \OmX$,
 we can deduce from the continuity of $\phi$ and the bounded convergence theorem   that the mapping
   \beas
 \hspace{-1.2cm} \Phi^m_k (\o_0,\omX,\fu) & \tn   \df & \tn   \int_{t^m_k}^{t^m_{k+1}} \n \phi \big(r, \fl_2(t^m_k,\omX),\o_0(t^m_k), \omX(t^m_k) ,\fu(r)\big) dr  \nonumber \\
  & \tn \=  & \tn  \lmt{n \to \infty} \n  \int_{t^m_k}^{t^m_{k+1}}  \n \sum_{i,j \in \hN} \b1_{ \{\fl_2 (t^m_k,\omX) \in \wt{O}^n_i \} } \b1_{\{\cW^m_k(\o_0,\omX) \in  \cE^n_j \}} \phi \big(r,\omX^i , \nxi_j, \fu(r)\big) dr \nonumber \\
   & \tn   \=  & \tn   \lmt{n \to \infty} \sum_{i,j \in \hN} \b1_{ \{\fl_2 (t^m_k,\omX) \in \wt{O}^n_i \} } \b1_{\{\cW^m_k(\o_0,\omX) \in  \cE^n_j \}}   I^{m,k}_{i,j} (\fu) , \q  \fa (\o_0,\omX,\fu) \ins \O_0 \ti \OmX \ti \hJ \q
  \eeas
   is $\sB(\O_0) \oti \sB(\OmX) \oti \sB(\hJ)-$measurable.
 Then the continuity of $\fl_2(r,\omX)$ in $r \ins [0,\infty)$ 
  and the  bounded convergence theorem imply  that the mapping
   \beas
  \hspace{-0.7cm}  \Phi  (\o_0,\omX,\fu) & \tn \df & \tn  \int_0^T \phi \big(r,\fl_2(r,\omX),\o_0(r), \omX(r),\fu(r)\big) dr
   \= \lmt{m \to \infty} \int_0^T \sum^{2^m-1}_{k = 0 } \b1_{\{r \in [t^m_k,t^m_{k+1}) \}} \phi \big(r,\fl_2(t^m_k,\omX),\o_0(t^m_k), \omX(t^m_k),\fu(r)\big) dr \\
     & \tn  \= & \tn 
      \lmt{m \to \infty} \sum^{2^m-1}_{k = 0 } \Phi^m_k (\o_0,\omX,\fu)   , \q  \fa (\o_0,\omX,\fu) \ins \O_0 \ti \OmX \ti \hJ
   \eeas
     is also $\sB(\O_0) \oti \sB(\OmX) \oti \sB(\hJ)-$measurable.

 \no {\bf 2b)} Set $\hX_1 \df [0,\infty)$,  $\hX_2  \df \OmX$, $\hX_3  \df \hR^{d+l}$   and $\hX_4  \df \hU$. They are all metric spaces.
 Let $i \= 1,2,3,4$.  We denote the corresponding metric of $\hX_i$ by $\Rho{\hX_i}$ and Let $ C_i $ be a closed subset of $\hX_i$
 Given $n \ins \hN$, we define an open subset  $O^i_n$ of $\hX_i$ by  $O^i_n \df \Big\{x \ins \hX_i \n :  dist_{\overset{}{\hX_i}}(x,C_i) \df \underset{x' \in C_i}{\inf} \Rho{\hX_i}(x,x') \< 1/n \Big\}$.  In light of Urysohn's Lemma,
there exist continuous function   $\phi^i_n \n : \hX_i \to [0,1]$  such that
$ \phi^i_n (x)   \= 1 $, $\fa x \ins C_i$ and $ \phi^i_n (x)   \= 0 $, $\fa x \ins (O^i_n)^c  $.
 So  $ \wh{\phi}_n (s,\omX,\nxi,u) \df \phi^1_n (s) \phi^2_n (\omX) \phi^3_n (\,\nxi) \phi^4_n(u) \ins [0,1] $, $(s,\omX,\nxi,u) \ins [0,\infty) \ti \OmX \ti \hR^{d+l} \ti \hU$
is a continuous function    satisfying   
 $ \lmt{n \to \infty} \wh{\phi}_n (s,\omX,\nxi,u)
 \= \b1_{\{(s,\omX,\,\snxi,u) \in C_1 \times C_2 \times C_3 \times C_4\}} $, $\fa (s,\omX,\nxi,u) \ins [0,\infty) \ti \OmX \ti \hR^{d+l} \ti \hU$.

 Let $T \ins (0,\infty)$.     Since the mappings $\int_0^T \wh{\phi}_n \big(r,\fl_2(r,\omX),\o_0(r), \omX(r),\fu(r)\big) dr  $,
 $(\o_0,\omX,\fu) \ins \O_0 \ti \OmX \ti \hJ$ are $\sB(\O_0) \oti \sB(\OmX) \oti \sB(\hJ)-$measurable for all $n \ins \hN$ by Part (2a),
    the bounded convergence theorem shows that  the function
  \beas
   \int_0^T \b1_{\{(r,\fl_2(r,\omX),\o_0(r), \omX(r),\fu(r)) \in C_1 \times C_2 \times C_3 \times C_4\}}  dr
   \= \lmt{n \to \infty}  \int_0^T  \wh{\phi}_n \big(r,\fl_2(r,\omX),\o_0(r), \omX(r),\fu(r)\big)  dr  ,
   \eeas
  $(\o_0,\omX,\fu) \ins \O_0 \ti \OmX \ti \hJ$  is also $\sB(\O_0) \oti \sB(\OmX) \oti \sB(\hJ)-$measurable.
  Let $\cH$ collect all real-valued  Borel-measurable functions $\psi$  on $[0,\infty) \ti \OmX \ti \hR^{d+l} \ti \hU$
    such that the mapping $\O_0 \ti \OmX \ti \hJ \ni (\o_0,\omX,\fu)   \mto \int_0^T \psi \big(r,\fl_2(r,\omX),\o_0(r), \omX(r),\fu(r)\big) dr  $  is $\sB(\O_0) \oti  \sB(\OmX) \oti \sB(\hJ)-$measurable.
 Clearly, $\cH$ is closed under linear combination.
 If $\{\psi_n\}_{n \in \hN} \sb \cH $ is a   sequence of non-negative functions that increases to a bounded function $\psi$
 on $[0,\infty) \ti \OmX \ti \hR^{d+l} \ti \hU$,
 the bounded convergence theorem implies that   $ \psi $ is also of $\cH$.
 Since $\{C_1 \ti C_2 \ti C_3 \ti C_4 \n: \hb{$ C_1 \sb [0,\infty)$, $ C_2 \sb \OmX$, $ C_3 \sb \hR^{d+l}$  and $C_4 \sb \hU$ are closed}\}$
 is a Pi-system that 
 generates $ \sB[0,\infty) \oti \sB(\OmX) \oti \sB(\hR^{d+l}) \oti \sB(\hU) $, we know from the monotone class theorem that
 $\cH$ includes all bounded Borel-measurable functions on $[0,\infty) \ti \OmX \ti \hR^{d+l} \ti \hU$.

 \no {\bf 2c)} Let $\psi \n: (0,\infty) \ti \OmX \ti \hR^{d+l} \ti \hU \mto [-\infty,\infty]$ be a Borel-measurable function
 and set $s^m_k \df k 2^{-m}  $ for any $m \ins \hN$ and $  k \ins \hN \cp \{0\}$.
 Given $n \ins \hN$, the bounded convergence theorem renders that
 the mapping  $  \int_t^{t+\fs} n \ld \psi^\pm \big(r,   \fl_2(r,\omX),\o_0(r), \\ \omX(r),   \fu(r)\big)  dr  \=   \lmt{m  \to \infty} \sum^\infty_{k,\ell = 0} \b1_{\{t \in [s^m_k,s^m_{k+1}) \}} \b1_{\{\fs \in [s^m_\ell,s^m_{\ell+1}) \}} \int_{s^m_k}^{s^m_{k+\ell}} n \ld \psi^\pm \big(r,\fl_2(r,\omX),\o_0(r), \omX(r),\fu(r)\big)  dr
 $, $  (t,\fs,\o_0,  \omX, \\ \fu) \ins [0,\infty) \ti [0,\infty) \ti \O_0 \ti \OmX \ti \hJ $  is $\sB[0,\infty) \oti \sB[0,\infty) \oti \sB(\O_0) \oti \sB(\OmX) \oti \sB(\hJ)-$measurable.
 Then it follows from  the monotone convergence theorem   that
  the mapping $ \Psi^\pm (t,\fs,\o_0,\omX,\fu) \df \int_t^{t+\fs} \psi^\pm  \big(r,\fl_2(r,\omX),\o_0(r), \omX(r),\fu(r)\big) dr \=   \lmt{n  \to \infty}  \int_t^{t+\fs} n \ld \psi^\pm  \big(r,\fl_2(r,\omX),\o_0(r), \omX(r),\fu(r)\big)  dr  $
 is Borel-measurable in $(t,\fs,\o_0,\omX,\fu) \ins [0,\infty) \ti [0,\infty) \ti \O_0 \ti \OmX \ti \hJ$,
 proving the second statement.

 Clearly, $\fF \sb  \si \big(I_\vf  ;  \vf \ins L^0 \big([0,\infty) \ti \hU ; \hR \big) \big)
 \= \si \big(I_\vf  ;  \vf \ins L^0 \big((0,\infty) \ti \hU ; \hR \big) \big)$,
 where $L^0 \big([0,\infty) \ti \hU;\hR  \big) $ collect all real$-$valued  Borel-measurable functions   on $[0,\infty) \ti \hU$.
 Let $\vf$ be a non-negative function in $   L^0 \big((0,\infty) \ti \hU ; \hR \big)$.
 Taking $\psi(r,\fx,\nxi,u) \df \vf (r,u) $, $\fa (r,\fx,\nxi,u) \ins (0,\infty) \ti \OmX \ti \hR^{d+l} \ti \hU $.
  shows that the mapping   $ I_\vf (\fu) \= \lmt{T \to \infty} \int_0^T \vf \big(r,\fu(r)\big) dr
  \= \lmt{T \to \infty} \int_0^T \psi \big(r,\bz,0,\fu(r)\big) dr
  \= \lmt{T \to \infty}  \Psi (0,T,\bz,\bz,\fu)$, $\fa \fu \ins \hJ$ is $\sB(\hJ)-$measurable.
  Hence, we have  $\sB(\hJ) \= \fF   \=  \si \big(I_\vf  ;  \vf \ins L^0 \big((0,\infty) \ti \hU ; \hR \big) \big)  $.
   \qed


\no {\bf Proof of Proposition \ref{prop_MPF1}: 1)} Set $\cN \df \big\{\o \ins \O \n: X_s (\o) \nne \bx(s) \hb{ for some } s \ins [0,t] \big\} \ins \sN_P(\cF^X_t) $ and let $(\vf,n,\fra) \ins  C^2(\hR^{d+l}) \ti \hN \ti \hR^{d+l}$.
 We set  $ c^\vf_n(\fra) \df \underset{|( w, x )| \le  n + \fra }{\sup} \big( \sum^2_{i=0} | D^i \vf ( w, x ) |   \big)  \+   \big|\vf (0,\bx(t) )\big|    \< \infty$ and $c^n_{t,\bx}(\fra) \df \big[ d / 2 \+  \k( t\+n)  ( \|\bx\|_t  \+  n \+ \fra    )  \+     \k^2( t\+n)  ( \|\bx\|_t  \+  n \+ \fra     )^2      \big] n \+     \int_t^{t+n}  \Sup{u \in \hU} \big(  |b( r,\bz,u)| \+ |\si( r,\bz,u)|^2 \big)  dr \< \infty $.
   Given $\o \ins \cN^c$, since $\big\|X_{r \land \cd} (\o)\big\|_r \= \Sup{r' \in [0,r]} \big|X_{r'} (\o) \big|
    \ls \|\bx\|_t  \ve ( n \+ \fra)  $, $\fa r \ins \big[t,  (\tau^t_n (\fra)) (\o)\big]$,
   we can deduce from   \eqref{coeff_cond1}, \eqref{coeff_cond2}  and Cauchy-Schwarz inequality
   \bea
   && \hspace{-1.2cm} \Sup{s \in [t,(\tau^t_n (\fra)) (\o)]} \big| \big(M^{t,\mu}_s(\vf)\big)(\o) \big|  \nonumber \\
    & & \hspace{-0.7cm}  \ls   \Sup{s \in [t,(\tau^t_n (\fra)) (\o)]} \big| \vf  \big( B^t_s (\o) , X_s (\o) \big) \big| \+ c^\vf_n(\fra) \n \int_t^{(\tau^t_n (\fra)) (\o)} \n \Big( \big|   b  (r, X_{r \land \cd}(\o),\mu_r(\o) ) \big| \+  \frac12   \big(d\+\big|    \si   (r, X_{r \land \cd}(\o),\mu_r(\o) ) \big|^2\big) \Big)  dr \nonumber  \\
   & & \hspace{-0.7cm}  \ls c^\vf_n(\fra)   \+ c^\vf_n(\fra)  \n \int_t^{(\tau^t_n (\fra)) (\o)} \Big(     \k( r)   \big\|X_{r \land \cd}(\o)\big\|_r
    \+   |b( r,\bz,\mu_r(\o) )|    \+        d/2   \+    \k^2( r)   \big\|X_{r \land \cd}(\o)\big\|^2_r
    \+   |\si( r,\bz,\mu_r(\o) )|^2 \Big) dr \nonumber  \\
   & & \hspace{-0.7cm}    \ls c^\vf_n(\fra)  \big( 1 \+   c^n_{t,\bx}(\fra) \big)     . \q   \label{122721_11}
    \eea
     So the process $ \big\{ M^{t,\mu}_{s \land \tau^t_n (\fra) } (\vf)   \big\}_{s \in [t,\infty)} $ is   bounded.    

 \no {\bf 2)} We next show that  (i) implies  (ii):  Suppose that (i) holds and let   $ (\vf,n,\fra)  \ins  C^2(\hR^{d+l}) \ti \hN \ti \hR^{d+l}  $. Setting $\th \df (t,\bx,\mu)$,
 we simply denote   $\Xi^\th_s \df (B^t_s ,  X^\th_s) $, $ \fa s \ins [t,\infty)$
 and set $\tau^\th_n (\fra)  \df  \inf\big\{s \ins [t,\infty) \n :  |\Xi^\th_s \- \fra|    \gs n   \big\} \ld (t\+n) $,
 which is an $\bF^{B^t,P} -$stopping time.
 Applying  It\^o's formula yields   that   $P-$a.s.
\beas
 M^\th_s(\vf)  & \tn \df & \tn   \vf  ( \Xi^\th_s  )
 \- \n \int_t^s \ol{b}  \big( r, X^\th_{r \land \cd}, \mu_r \big)  \n \cd \n D \vf  ( \Xi^\th_r  ) dr
    \-   \frac12 \int_t^s   \ol{\si} \, \ol{\si}^T  \big( r, X^\th_{r \land \cd}, \mu_r \big) \n : \n D^2 \vf ( \Xi^\th_r  )   dr \nonumber \\
 & \tn \= & \tn  \vf \big(0,\bx(t)\big)   \+ \int_t^s    D  \vf ( \Xi^\th_r  )  \n \cd \n
 \ol{\si}  \big(r, X^\th_{r \land \cd}, \mu_r \big) d B_r , \q s \ins [t,\infty) . 
\eeas

 For any $\o \ins \O$, an analogy to \eqref{122721_11} shows that
    $\Sup{s \in [t,(\tau^\th_n(\fra))(\o)]} \big| \big(M^\th_s(\vf)\big)(\o) \big| \ls  c^\vf_n(\fra)  \big( 1 \+   c^n_{t,\bx}(\fra) \big)  $ and   $  \int_t^{(\tau^\th_n(\fra))(\o)}   \big|  D  \vf (   \Xi^\th_r (\o)  )   \cd
 \ol{\si}  (r, X^\th_{r \land \cd} (\o),\mu_r(\o) ) \big|^2 d r
  \ls    (c^\vf_n(\fra))^2   \big[ d   \+ 2    \k^2( t\+n)  ( \|\bx\|_t \+ n \+ \fra    )^2      \big] n \+ 2   (c^\vf_n(\fra))^2 \n   \int_t^{t+n}  \Sup{u \in \hU} |\si( r,\bz,u)|^2  dr    \< \infty $.
 So
 \bea \label{122921_14}
  \hb{$\big\{ M^\th_{s \land \tau^\th_n (\fra)   }(\vf) \big\}_{s \in [t,\infty)}$ is   a    bounded  $  \bF^{B^t,P}  -$martingale.}
  \eea
 Set $ \cN_\th \df \{\o \ins \O \n: X_s (\o) \nne   X^\th_s (\o)  \hb{ for some } s \ins [0,\infty)\} \ins \sN_P(\cF^{t,P}_\infty) $.
 For any $(s,\o) \ins [0,\infty) \ti \cN^c_\th $,
 \bea \label{122021_17}
 X^\th_s (\o) \= X_s (\o)   ,\q \big( M^\th_{s \vee t}(\vf) \big) (\o) \=   \big( M^{t,\mu}_{s \vee t}(\vf) \big) (\o)
 \q \hb{and thus} \q (\tau^\th_n (\fra)) (\o) \= (\tau^t_n (\fra)) (\o) .
 \eea

  Fix $t_1,t_2 \ins [t,\infty)$ with $t_1 \< t_2$.
  Let $ \big\{(s_i,\cE_i)\big\}^m_{i=1} \sb [t,t_1] \ti \sB(\hR^d)$ and $ \big\{(r_j,A_j)\big\}^k_{j=1} \sb [0,t_1] \ti \sB(\hR^l)$.
 We can derive from \eqref{122921_14} and \eqref{122021_17}   that
 $
    E_P \Big[ \b1_{\cN^c_\th} \big( M^{t,\mu}_{t_2 \land \tau^t_n(\fra)}(\vf) \-  M^{t,\mu}_{t_1 \land \tau^t_n(\fra)}(\vf) \big)  \prod^m_{i=1}   \b1_{(B^t_{s_i})^{-1}(\cE_i)} \prod^k_{j=1} \b1_{ X_{r_j}^{-1}(A_j)}   \Big]
     \=  E_P \Big[ \b1_{\cN^c_\th} \big( M^\th_{t_2 \land \tau^\th_n(\fra)}(\vf) \\  \-  M^\th_{t_1 \land \tau^\th_n(\fra)}(\vf) \big)  \prod^m_{i=1} \b1_{(B^t_{s_i})^{-1}(\cE_i)} \prod^k_{j=1} \b1_{ (X^\th_{r_j})^{-1}(A_j)}  \Big]
 \= 0 $.
   So   the Lambda-system $ \L \df \big\{ A \ins \cF^{t,P}_\infty \n : E_P \big[   \big( M^{t,\mu}_{t_2 \land \tau^t_n(\fra) } (\vf) \- M^{t,\mu}_{t_1 \land \tau^t_n(\fra) } (\vf)  \big)   \b1_A \big] \= 0  \big\}$   contains the  Pi-system
 $ \Big\{  \Big( \underset{i=1}{\overset{m}{\cap}}  (B^t_{s_i})^{-1}(\cE_i) \Big) \Cp \Big( \underset{j=1}{\overset{k}{\cap}} X_{r_j}^{-1}(A_j) \Big) \n :  \big\{(s_i,\cE_i)\big\}^m_{i=1} \sb [t,t_1] \ti \sB(\hR^d) , \, \big\{(r_j,A_j)\big\}^k_{j=1} \sb [0,t_1] \ti \sB(\hR^l)   \Big\} $,
  which generates $\cF^t_{t_1}  $.    Dynkin's Pi-Lambda Theorem (see e.g Theorem 3.2 of \cite{Billingsley_PM}) renders   $ \cF^t_{t_1} \sb \L $.
  As $   \cF^t_{t_1} \cp \sN_P \big( \cF^t_\infty \big)$ is another Pi-system included in $\L$. Applying Dynkin's Pi-Lambda Theorem again yields that $ \cF^{t,P}_{t_1} \sb   \si\big(\cF^t_{t_1} \cp \sN_P ( \cF^t_\infty ) \big) \sb \L $, i.e.,
  $E_P \big[   \big( M^{t,\mu}_{t_2 \land \tau^t_n(\fra) } (\vf) \- M^{t,\mu}_{t_1 \land \tau^t_n(\fra) } (\vf)   \big)   \b1_A \big] \= 0 $, $\fa A \ins \cF^{t,P}_{t_1}$.
  Hence, $\big\{ M^{t,\mu}_{s \land \tau^t_n(\fra) } (\vf) \big\}_{s \in [t,\infty)}$ is a bounded $ \bF^{t,P}-$martingale,
  which is further a $ \bF^t-$martingale if    process $ \mu $ is   $\bF^{B^t}-$progressively measurable.

 \no {\bf 3)} As $\fC(\hR^{d+l}) \sb C^2(\hR^{d+l}) $, (ii) \n $\Ra$ \n (iii)   is straightforward. It remains to show that (iii) gives rise to (i).

 If the  process $ \mu $ is   $\bF^{B^t}-$progressively measurable and if  $ \big\{ M^{t,\mu}_{s \land \tau^t_n } (\vf) \big\}_{s \in [t,\infty)} $  is a bounded $\bF^t-$martingale for any $(\vf,n) \ins \fC(\hR^{d+l}) \ti \hN $, then Dynkin's Pi-Lambda Theorem implies that
 $ \big\{ M^{t,\mu}_{s \land \tau^t_n } (\vf) \big\}_{s \in [t,\infty)} $  is further an $\bF^{t,P}-$martingale for any $(\vf,n) \ins \fC(\hR^{d+l}) \ti \hN $.

 \no {\bf 3a)} Define $  \cF^{t,P}_{s+} \df \ccap{\e>0}{} \cF^{t,P}_{s+\e} $, $\fa s \ins   [t,\infty)  $
   and set $\bG^{t,P} \= \big\{ \cG^{t,P}_s \df \cF^{t,P}_{s+}\big\}_{s \in   [t,\infty)}$.

  Let $i,j \ins \{ 1,\cds \n ,d \}$. We  set $  \phi_i(w,x) \df w_i $ and $\phi_{ij}(w,x) \df w_i w_j $
 for any $w \= (w_1,\cds \n ,w_d) \ins \hR^d$ and $  x \ins \hR^l$. Clearly,  $ \phi_i , \phi_{ij} \ins \fC(\hR^{d+l})$.
 One can calculate that   $ M^{t,\mu}_s(\phi_i) \= B^{t,i}_s   $,   $ M^{t,\mu}_s(\phi_{ij}) \=  B^{t,i}_s B^{t,j}_s \- \d_{ij} (s\-t) $, $\fa s \ins [t,\infty)$,
 where  $B^t_s \= \big( B^{t,1}_s ,\cds \n , B^{t,d}_s  \big)$ and $\d_{ij}$ is the $(i,j)-$element of the identity matrix $I_{d \times d}$.

 Let $n \ins \hN $.  By (iii),   $ \big\{ M^{t,\mu}_{s \land \tau^t_n } (\phi_i) \big\}_{s \in [t,\infty)} $ and $ \big\{ M^{t,\mu}_{s \land \tau^t_n } (\phi_{ij}) \big\}_{s \in [t,\infty)} $ are bounded $   \bF^{t,P}-$martingales.  
   The {\it optional sampling} theorem (e.g. Theorem 1.3.22 of \cite{Kara_Shr_BMSC}) implies that
  they are further  $   \bG^{t,P}-$martingales.
  Since $\lmtu{n \to \infty} \tau^t_n \= \infty$,
  we see that $\big\{M^{t,\mu}_s(\phi_i) \= B^{t,i}_s   \big\}_{s \in [t,\infty)}  $
  and $\big\{M^{t,\mu}_s(\phi_{ij}) \= B^{t,i}_s B^{t,j}_s \- \d_{ij} (s\-t) \big\}_{s \in [t,\infty)}$   are $   \bG^{t,P}-$local martingales.
    L\'evy's characterization theorem then yields that $ B^t   $ is a    Brownian motion with respect to  filtration $   \bG^{t,P} $
   and is thus a Brownian motion with respect to   filtration $ \bF^{B^t} $. 

 \no {\bf 3b)}
 We simply denote   $\Xi_s \df (B^t_s ,  X_s) $, $\beta_s \df \ol{b}  \big(s, X_{s \land \cd},\mu_s \big)$ and
 $ \a_s \df \ol{\si} \, \ol{\si}^T  \big( s,  X_{s \land \cd},\mu_s \big) $, $ \fa  s \ins [t,\infty)$.
 Let  $i,j  \ins \{ 1,\cds \n ,d\+l\}$. We set $  \psi_i(\nxi) \df \nxi_i  $
 and $\psi_{ij}(\nxi) \df \nxi_i  \nxi_j  $
 for any  $  \nxi \= \big( \nxi_1   ,\cds \n , \nxi_{d + l} \big)  \ins \hR^{d+l}$.
 Similar to $ M^{t,\mu}_\cd (\phi_i) $ and $ M^{t,\mu}_\cd   (\phi_{ij}) $,
  the processes  $  M^{t,\mu}_s(\psi_i) \=    \Xi^{(i)}_s   \-  \int_t^s \beta^{(i)}_r  dr$  and
$ M^{t,\mu}_s(\psi_{ij}) \= \Xi^{(i)}_s  \Xi^{(j)}_s
 \- \int_t^s \beta^{(i)}_r \Xi^{(j)}_r dr
 \- \int_t^s  \beta^{(j)}_r   \Xi^{(i)}_r dr
 \- \int_t^s  ( \a_r  )_{ij} dr $, $ s \ins [t,\infty)$
 are  $ \bG^{t,P}-$local martingales.
 Using the {\it integration by parts} formula, we obtain   that  $P-$a.s.
 \beas
 && \hspace{-1.5cm}
 \Xi^{(i)}_s  \Xi^{(j)}_s  \-  M^{t,\mu}_s(\psi_i) M^{t,\mu}_s(\psi_j)
   \=      M^{t,\mu}_s(\psi_i) \int_t^s  \beta^{(j)}_r    dr \+M^{t,\mu}_s(\psi_j) \int_t^s \beta^{(i)}_r  dr \+ \int_t^s \beta^{(i)}_r  dr \cd  \int_t^s  \beta^{(j)}_r    dr \\
  &&  \hspace{-0.5cm} \=    \int_t^s M^{t,\mu}_r(\psi_i)  \beta^{(j)}_r    dr
  \+  \int_t^s   \Big(  \int_t^r \beta^{(j)}_{r'}  dr' \Big) d M^{t,\mu}_r(\psi_i)
  \+ \int_t^s M^{t,\mu}_r(\psi_j) \beta^{(i)}_r  dr
  \+  \int_t^s   \Big(  \int_t^r \beta^{(i)}_{r'}  dr' \Big) d M^{t,\mu}_r(\psi_j) \\
  && \hspace{-0.5cm} \q
   +  \int_t^s \Big(  \int_t^r \beta^{(i)}_{r'}  dr' \Big)   \beta^{(j)}_r    dr
  \+ \int_t^s  \Big(  \int_t^r \beta^{(j)}_{r'}  dr' \Big) \beta^{(i)}_r  dr \\
  && \hspace{-0.5cm} \=       \int_t^s \Big[ \Xi^{(i)}_r  \beta^{(j)}_r
  \+ \Xi^{(j)}_r \beta^{(i)}_r \Big] dr
  \+  \int_t^s  \n \Big(  \int_t^r \beta^{(i)}_{r'}  dr' \Big) d M^{t,\mu}_r(\psi_i)
  \+  \int_t^s  \n \Big(  \int_t^r \beta^{(j)}_{r'}  dr' \Big) d M^{t,\mu}_r(\psi_j)  , \q  s \ins [t,\infty) .
 \eeas
   \if{0}
 Note
 \beas
 \int_t^s \beta^{(i)}_r  dr \int_t^s  \beta^{(j)}_r
 \= \int_t^s \Big(  \int_0^r \beta^{(i)}_{r'}  dr' \Big)   \beta^{(j)}_r    dr
  \+ \int_t^s  \Big(  \int_0^r \beta^{(j)}_{r'}  dr' \Big) \beta^{(i)}_r  dr ,
   \q s \ins [t,\infty) .
 \eeas
   \fi
 So $ M^{t,\mu}_s(\psi_i) M^{t,\mu}_s(\psi_j) \- \int_t^s  ( \a_r  )_{ij} dr
 \= M^{t,\mu}_s(\psi_{ij}) \-  \int_t^s   \big(  \int_t^r \beta^{(i)}_{r'}  dr' \big) d M^{t,\mu}_r(\psi_i)
\- \int_t^s   \big(  \int_t^r \beta^{(j)}_{r'}  dr' \big) d M^{t,\mu}_r(\psi_j) $,
  $s \ins [t,\infty)$  is also an $ \bG^{t,P}-$local martingale, which implies that
 the quadratic variation of  the $ \bG^{t,P}-$local martingale
  $  M^{t,\mu}_s  \df \big( M^{t,\mu}_s(\psi_1), \cds \n ,    M^{t,\mu}_s(\psi_{d+l})\big)
  \= \Xi_s  \-  \int_t^s \beta_r   dr $, $ s \ins [t,\infty)$   is
  $ \big\lan  M^{t,\mu} ,  M^{t,\mu}    \big\ran_s  \= \int_t^s   \a_r   dr $, $ s \ins [t,\infty) $.

  Let $ n \ins \hN $, $a \ins \hR^l$ and set $ \dis \cH^a_s \df   \binom{-   \si^T \big( s, X_{s \land \cd},\mu \big) a }{a}$, $ \fa s \ins (t,\infty)$.
   The stochastic exponential   of the  $ \bG^{t,P}-$ martingale
  $\big\{\int_t^{\tau^t_n \land s} \cH^a_r \n \cd \n d M^{t,\mu}_r\big\}_{s \in [t,\infty)}$ is
  \beas
 && \hspace{-1.2cm} \exp\Big\{\int_t^{\tau^t_n \land s} \cH^a_r \n \cd \n d M^{t,\mu}_r \- \frac12 \int_t^{\tau^t_n \land s}   (\cH^a_r)^T  \a_r  \cH^a_r  dr \Big\}
  \= \exp\Big\{\int_t^{\tau^t_n \land s} \cH^a_r \n \cd \n d \Xi_r  \-  \int_t^{\tau^t_n \land s} \cH^a_r \n \cd \n  \beta_r   dr  \Big\} \\
 && \= \exp\bigg\{ a \n \cd \n \Big( \int_t^{\tau^t_n \land s}        d X_r  \- \int_t^{\tau^t_n \land s}    \si \big( r, X_{r \land \cd},\mu_r \big)   d B_r  \-  \int_t^{\tau^t_n \land s}      b \big( r, X_{r \land \cd},\mu_r \big)   dr \Big)   \bigg\} , \q s \ins [t,\infty) .
  \eeas
  \if{0}
  We can deduce from \eqref{coeff_cond1} that  $\big\{\int_0^{\tau^t_n \land s} \cH^a_r \n \cd \n d M^{t,\mu}_r\big\}_{s \in [t,\infty)}$ is a $  \bF^{t,P}-$BMO martingale
  and thus its stochastic exponential is a  $   \bF^{t,P}-$uniformly integrable martingale.
  \fi
  Letting $a $ vary over $\hR^l$ yields that $P-$a.s.,
  $ X_{ \tau^t_n \land s} \=   \bx(t)  \+  \int_t^{\tau^t_n \land s}     b \big( r, X_{r \land \cd},\mu_r \big)   dr \+ \int_t^{\tau^t_n \land s}    \si \big( r, X_{r \land \cd},\mu_r \big)   d B_r   $,  $ \fa s \ins [t,\infty) $.
   Sending $n \nto \infty$ then renders  that      $P-$a.s.,
  $ X_s \=   \bx(t)  \+  \int_t^s    b \big( r, X_{r \land \cd},\mu_r \big)   dr \+ \int_t^s    \si \big( r, X_{r \land \cd},\mu_r \big)   d B_r   $,  $ \fa s \ins [t,\infty) $.
  Viewing  SDE \eqref{121621_11} with the control process $\{\mu_s\}_{s \in [t,\infty)}$ on   $\big(\O,\cF, \bG^{t,P},P\big)$,
  we know from Proposition \ref{prop_122021} that  there is a unique   $ \big\{\cG^{t,P}_{s \vee t}   \big\}_{s \in [0,\infty)}  -$adapted continuous process satisfying   \eqref{121621_11}. Hence, $P\{X_s \= X^{t,\bx,\mu}_s , \;  \fa s \ins [0,\infty)\} \= 1$.  \qed


\no {\bf Proof of Proposition \ref{prop_122021b}: 1a)}    Given a sigma-field $\cG$ of $\O_0$,   define
    $ 
  \sS_{t,\bw}(\cG) \df  \big\{  A \sb   \O   \n:
     \has \cA \ins  \cG   $ and $\cN \ins \sN_P(\cF^{B^t}_\infty)  $ s.t.
      $ \b1_{\{ B^{t,\bw} (\o)  \in \cA \} }  \=  \b1_{ \{ \o  \in A \} }, \, \fa \o \ins \cN^c \big\} $.
    As   $ \b1_{\{ B^{t,\bw} (\o)  \in  \O \}} \= 1 \= \b1_{\{  \o  \in  \O\}} $  for any  $  \o   \ins   \O $,
 it is clear that    $\O \ins \sS_{t,\bw}\big(\cG\big)$.
  When $A \ins \sS_{t,\bw}\big(\cG\big)$,   there exist   $\cA \ins  \cG $ and $\cN \ins \sN_P(\cF^{B^t}_\infty)$ such that
    $ \b1_{\{ B^{t,\bw} (\o)  \in \cA \} }   \=   \b1_{ \{ \o  \in A \} } $,   $ \fa  \o \ins \cN^c $.
    Then $\cA^c \ins  \cG $  satisfies that
   $ \b1_{\{ B^{t,\bw} (\o)  \in \cA^c \} }   \= 1\- \b1_{\{ B^{t,\bw} (\o)  \in \cA \} }   \=
   1\- \b1_{ \{ \o  \in A \} } \=  \b1_{ \{ \o  \in A^c \} } $, $ \fa \o \ins \cN^c $,
    so $A^c $ also belongs to $   \sS_{t,\bw}\big(\cG\big) $.
   If $\{A_k\}_{k \in \hN} \sb \sS_{t,\bw}\big(\cG\big)$, for any $k \ins \hN$
    there exist   $\cA_k \ins  \cG $ and $\cN_k \ins \sN_P(\cF^{B^t}_\infty)$ such that
    $ \b1_{\{ B^{t,\bw} (\o)  \in \cA_k \} }   \=   \b1_{ \{ \o  \in A_k \} } $, $ \fa \o \ins \cN^c_k $. Then $ \underset{k \in \hN}{\cap} \cA_k \ins \cG $ satisfies that
    $ \b1_{\big\{ B^{t,\bw} (\o)  \in \underset{k \in \hN}{\cap} \cA_k \big\} }
     \= \underset{k \in \hN}{\prod} \b1_{\{ B^{t,\bw} (\o)  \in \cA_k \} }
     \= \underset{k \in \hN}{\prod}  \b1_{ \{ \o  \in A_k \} }
     \= \b1_{\big\{ \o  \in \underset{k \in \hN}{\cap}  A_k \big\} } $, $ \fa \o \ins \underset{k \in \hN}{\cap} \cN^c_k $,
   which shows that $ \underset{k \in \hN}{\cap}  A_k \ins \sS_{t,\bw}\big(\cG\big) $ with $(\cA,\cN) \= \Big(\underset{k \in \hN}{\cap} \cA_k ,  \ccup{k \in \hN}{} \cN_k\Big) $.
   Hence $\sS_{t,\bw}\big(\cG\big)$ is a sigma$-$field of $\O$.

  Let $s \ins [t,\infty)$.   For any $r \ins [t,s]$ and $\cE \ins \sB(\hR^d) $, since
  $  \o \ins    (B^t_r)^{-1} (\cE) $  iff $  W^t_r (B^{t,\bw}  (\o)) \= W_r(B^{t,\bw} (\o))\-W_t(B^{t,\bw} (\o)) 
  \= B^t_r(\o) \ins   \cE $
  iff $ B^{t,\bw} (\o) \ins (W^t_r)^{-1} (\cE)   $,
   one has
 $ \b1_{\{ B^{t,\bw} (\o) \in   (W^t_r)^{-1} (\cE) \}} \= \b1_{\{ \o \in   (B^t_r)^{-1} (\cE) \}}    $, $\fa \o \ins  \O $, which means
  $(B^t_r)^{-1} (\cE) \ins \sS_{t,\bw}(\cF^{W^t}_s)  $ with $\cA \= (W^t_r)^{-1} (\cE)
 \ins \cF^{W^t}_s$ and $\cN \= \es$.  So  $ \cF^{B^t }_s \sb \sS_{t,\bw}(\cF^{W^t}_s)$\,.
 When $s \> t $, as $  \cF^{B^t }_r \sb \sS_{t,\bw}(\cF^{W^t}_r) \sb \sS_{t,\bw}(\cF^{W^t}_{s-})$ for any $r \ins [t,s)$,
 one further has
 $ 
  \cF^{B^t }_{s-} \= \si \Big(\underset{r \in [t,s)}{\cup} \cF^{B^t }_r \Big) \sb \sS_{t,\bw}(\cF^{W^t}_{s-}) $.

  \if{0}

     It is clear that for any $(s,\o) \ins [t,\infty) \ti \O$ that
  $ 
  \b1_{\{ (s,B^{t,\bw}(\o))  \in [t,\infty) \times \O_0 \} } \= 1 \=   \b1_{ \{ (s,\o)  \in [t,\infty) \times \O \} } $,
    so 
    $[t,\infty) \ti \O$ belongs to $\wh{\sS}_{t,\bw}$.
    When $D \ins \wh{\sS}_{t,\bw}$,   there exist  an $\bF^{W^t}-$predictable set $ \cD   $
   and $\cN \ins \sN_P(\cF^{B^t}_\infty) $ such that
    $   \b1_{\{ (s,B^{t,\bw}(\o))  \in \cD \} }   \=  \b1_{ \{ (s,\o)  \in D \} }  $
    for any $(s,\o) \ins [t,\infty) \ti \cN^c$.
    Then the $\bF^{W^t}-$predictable set $ \cD^c   $  satisfies that
   \beas
    \b1_{\{ (s,B^{t,\bw}(\o))  \in \cD^c \} }   \= 1\- \b1_{\{ (s,B^{t,\bw}(\o)) \in \cD \} }   \=
   1\- \b1_{ \{ (s,\o)  \in D \} } \=  \b1_{ \{ (s,\o)  \in D^c \} } , \q \fa (s,\o) \ins [t,\infty) \ti \cN^c ,
   \eeas
   which shows  $D^c \ins \wh{\sS}_{t,\bw}$.

   If $\{D_k\}_{k \in \hN} \sb \wh{\sS}_{t,\bw}$, for any $k \ins \hN$
    there exist   an $\bF^{W^t}-$predictable set $ \cD_k   $ and $\cN_k \ins \sN_P(\cF^{B^t}_\infty) $  such that
    $   \b1_{\{ (s,B^{t,\bw}(\o))  \in \cD_k \} }   \=  \b1_{ \{ (s,\o)  \in D_k \} }  $
    for any $(s,\o) \ins [t,\infty) \ti \cN^c_k$.
    Then the $\bF^{W^t}-$predictable set $ \underset{k \in \hN}{\cap} \cD_k   $ satisfies that
    \beas
     \b1_{\big\{ (s,B^{t,\bw}(\o))  \in \underset{k \in \hN}{\cap} \cD_k \big\} }
     \= \underset{k \in \hN}{\prod} \b1_{\{ (s,B^{t,\bw}(\o))  \in \cD_k \} }
     \= \underset{k \in \hN}{\prod}  \b1_{ \{ (s,\o) \in D_k \} }
     \= \b1_{\big\{ (s,\o) \in \underset{k \in \hN}{\cap}  D_k \big\} } ,  \q \fa (s,\o) \ins
     [t,\infty) \ti \Big(\underset{k \in \hN}{\cup} \cN_k \Big)^c ,
     \eeas
  which shows that $ \underset{k \in \hN}{\cap} D_k \ins \wh{\sS}_{t,\bw} $ with $\cN \=  \underset{k \in \hN}{\cup} \cN_k   $.
  Hence $\wh{\sS}_{t,\bw}$ is a sigma$-$field of $[t,\infty) \ti \O$.

  \fi

 \no {\bf 1b)} Similar to $\sS_{t,\bw}\big(\cG\big)$ in Part (1a),
   $ \wh{\sS}_{t,\bw}   \df  \big\{   D \sb [t,\infty) \ti \O   \n:
  \has \cD \ins \sP^{W^t}   $ and $\cN \ins \sN_P(\cF^{B^t}_\infty) $ s.t. $    \b1_{ \{(s,B^{t,\bw}(\o)) \in  \cD   \}}  \=    \b1_{\{(s, \o ) \in   D  \}} , \fa  (s,\o)  \ins [t,\infty) \ti \cN^c   \big\} $ is a sigma$-$field of $[t,\infty) \ti \O$.

Set $ \L_t  \df \big\{\{t\} \ti A_o \n : A_o \ins \cF^{B^t}_t \big\} \cp \big\{ (s,\infty) \ti A \n :    s \ins [t,\infty) \Cp   \hQ, \,
   A  \ins \cF^{B^t}_{s-}  \big\} $,
 which generates the $\bF^{B^t}-$predictable sigma$-$field $\sP^{B^t}$. 
  For any $A_o \ins \cF^{B^t}_t \sb \sS_{t,\bw}(\cF^{W^t}_t)$,   there exist    $\cA_o \ins   \cF^{W^t}_t $ and $\cN  \ins \sN_P(\cF^{B^t}_\infty)$ such that
    $ \b1_{\{ B^{t,\bw}(\o)  \in \cA_o \} }   \=   \b1_{ \{ \o  \in A_o \} } $,   $\fa \o \ins \cN^c $. We can deduce that
    $ \b1_{ \{(s,B^{t,\bw}(\o)) \in  \{t\} \ti \cA_o  \}}
       \=   \b1_{\{(s, \o ) \in  \{t\} \ti A_o \}}   $,
      $ \fa (s,\o) \ins [t,\infty) \ti \cN^c $.
  So $\{t\} \ti A_o$   is  of $\wh{\sS}_{t,\bw}$ with $ \cD \= \{t\} \ti \cA_o$.
  For any $s \ins [t,\infty) \cap \hQ$
  and $ A  \ins \cF^{B^t}_{s-} \sb \sS_{t,\bw}(\cF^{W^t}_{s-}) $, one can find some    $\cA \ins   \cF^{W^t}_{s-} $
  and  $\cN \ins \sN_P(\cF^{B^t}_\infty)$ such that
     $\b1_{\{ B^{t,\bw}(\o)  \in \cA \} }   \=   \b1_{ \{ \o  \in A \} }  $,    $ \fa    \o \ins \cN^c $.
   It then holds for any $  (r,  \o) \ins [t,\infty) \ti \cN^c $ that
    $\b1_{\{ (r,B^{t,\bw}(\o))  \in (s,\infty) \times \cA \} }
    \=   \b1_{ \{ (r,\o)  \in (s,\infty) \times A \} }  $,  which implies   $ (s,\infty) \ti  A \ins \wh{\sS}_{t,\bw} $ with $ \cD \= (s,\infty) \ti  \cA $.
   Hence,  
    $\wh{\sS}_{t,\bw}$ contains $\L_t$
    and thus includes  $\sP^{B^t}$.
    \if{0}

    \no {\bf 1c)}  Let $\{\triangle_i\}^\infty_{i=0} $ be an arbitrary countable set   and
     let $  \{\nu_s\}_{s \in [t,\infty)}$ be an $ \bF^{B^t} - $predictable process on $\O$ taking
      values in  $\{\triangle_i\}^\infty_{i=1}$. In this step, we shall construct
      a  $\{\triangle_i\}^\infty_{i=0}-$valued, $\bF^{W^t}-$predictable process $  \big\{\wh{\nu}_s\big\}_{s \in [t,\infty)}$
    on $\O_0$    and an $\wh{\cN} \ins \sN_P(\cF^{B^t}_\infty)$  such that
    $     \wh{\nu}_s (B^{t,\bw}(\o)) \= \nu_s (\o) $ for any $  (s,\o)  \ins [t,\infty) \ti \wh{\cN}^c $.

  Let $i \ins \hN$ and  set  an $\bF^{B^t}-$predictable set $D_i \df \{(s,\o) \ins [t,\infty) \ti \O \n : \nu_s(\o) \= \triangle_i \}$.
  As $\sP^{B^t} \sb \wh{\sS}_{t,\bw} $ by Part (1b),
    there exist  an $\bF^{W^t}-$predictable set $ \cD_i $ and an $\cN_i \ins \sN_P(\cF^{B^t}_\infty) $ such that
  $   \b1_{ \{(s,B^{t,\bw}(\o)) \in  \cD_i   \}} \= \b1_{\{(s, \o ) \in   D_i  \}} $,   $ \fa (s,\o) \ins [t,\infty) \ti \cN^c_i$.
 Set $\wh{\cN} \df \underset{i \in \hN}{ \cup }  \cN_i \ins \sN_P(\cF^{B^t}_\infty) $ and define an    $\bF^{W^t}-$predictable process $\wh{\nu}$ on $\O_0$ by
 $ 
 \wh{\nu}_s(\o_0) \df  \b1_{\big\{ (s,\o_0) \in \underset{i \in \hN}{ \cap } \cD^c_i \big\}} \triangle_0
 \+   \sum_{i \in \hN}\b1_{\{(s,\o_0) \in \cD'_i\}} \triangle_i $, $ (s,\o_0) \ins [t,\infty) \ti \O_0 $,
 where   $\cD'_1 \= \cD_1 $ and $\cD'_i \df \cD_i \big\backslash \Big( \underset{j < i}{\cup} \cD_j\Big)$ for $i \gs 2$.

   Let  $(s,\o) \ins  [t,\infty) \ti \wh{\cN}^c$ and let $i \gs 2 $.
   Since $D_i \Cp D_j \= \es$ for $j \= 1 \cds \n , i \-1$,
    one can deduce that $  0      \ls \b1_{ \big\{(s,B^{t,\bw}(\o)) \in  \cD_i \cap (\underset{j < i}{\cup} \cD_j )  \big\} }
    \ls  \sum^{i-1}_{j=1} \b1_{ \{(s,B^{t,\bw}(\o)) \in  \cD_i \cap   \cD_j   \}}
    \=  \sum^{i-1}_{j=1} \b1_{\{(s,\o ) \in D_i \cap D_j\}} \= 0 $, which implies
       $      \b1_{ \{(s,B^{t,\bw}(\o)) \in  \cD'_i\}} \=  \b1_{ \{(s,B^{t,\bw}(\o)) \in  \cD_i\}}   \= \b1_{\{(s,\o ) \in D_i  \}}  $. It follows that
 $   \b1_{\big\{ (s,B^{t,\bw}(\o)) \in \underset{i \in \hN}{ \cap } \cD^c_i \big\}}
 \=  1 \- \b1_{\big\{ (s,B^{t,\bw}(\o)) \in  \underset{i \in \hN}{ \cup } \cD_i   \big\}}
 \=  1 \- \b1_{\big\{ (s,B^{t,\bw}(\o)) \in  \underset{i \in \hN}{ \cup } \cD'_i   \big\}}
 \= 1 \- \sum_{i \in \hN}\b1_{\big\{ (s,B^{t,\bw}(\o)) \in   \cD'_i \big\}} \= 1 \- \sum_{i \in \hN} \b1_{\{(s,\o ) \in D_i  \}} \= 0 $    and thus
  $ \wh{\nu}_s(B^{t,\bw}(\o)) \= \sum_{i \in \hN}  \b1_{ \{(s,B^{t,\bw}(\o)) \in  \cD'_i\}}  \triangle_i
 \= \sum_{i \in \hN} \b1_{\{(s,\o ) \in D_i\}} \triangle_i \= \nu_s(\o) $.  

\ss \no {\bf 1d)} Let $  \{a_s\}_{s \in [t,\infty)}$ be a $[0,\infty)-$valued, $ \bF^{B^t} - $predictable process  on $\O$.
 Given $n \ins \hN$, we  define an $ \bF^{B^t} - $predictable process 
 by  $ a^n_s(\o) \df \sum_{i \in \hN} \b1_{\{ (i-1) 2^{-n} \le a_s(\o)   < i 2^{-n}\}} i 2^{-n} $, $ (s,\o) \ins [t,\infty) \ti \O $.
   Part (1c) 
   assures  a $\{i2^{-n}\}^\infty_{i=0}-$valued, $\bF^{W^t}-$predictable process $  \big\{\wh{a}^n_s\big\}_{s \in [t,\infty)}$   on $\O_0$
    and an $\cN_n \ins \sN_P(\cF^{B^t}_\infty)$ such that
 $ \wh{a}^n_s(B^{t,\bw}(\o)) \= a^n_s(\o)  $, $ \fa (s,\o) \ins [t,\infty) \ti \cN^c_n $.
  Define  $  \a_s (\o_0) \df \linf{n \to \infty}  \wh{a}^n_s (\o_0)  $, $   (s,\o_0) \ins [t,\infty) \ti \O_0$,
   which is a $[0,\infty)-$valued, $\bF^{W^t}-$predictable process. 

    Set $ \cN_a \df \Big(\underset{n \in \hN}{ \cup }  \cN_n \Big)   \ins \sN_P(\cF^{B^t}_\infty) $.
 It holds for any $(s,\o) \ins [t,\infty) \ti \cN^c_a $ that      $ a_s(\o)  \= \lmt{n \to \infty} a^n_s( \o )
\= \lmt{n \to \infty}  \wh{a}^n_s(B^{t,\bw}  (\o)) \= \a_s (B^{t,\bw}  (\o)) $.

    \fi

  \no {\bf 1c)}
Now, let $  \{\mu_s\}_{s \in [t,\infty)}$ be a general $\hU-$valued, $ \bF^{B^t,P} - $progressively measurable process  on $\O$.
Similar  to   Lemma 2.4 of \cite{STZ_2011a}, 
one can construct   a   $[0,1]-$valued  $  \bF^{B^t}  -$predictable  process  $  \big\{ \nu_s \big\}_{s \in [t,\infty)}$ on $\O$
   such that $ \nu_s  (\o)   \= \sI \big(  \mu_s (\o )\big) $ for $ds \ti dP-$a.s. $(s,\o) \ins [t,\infty) \ti \O$.
     Fubini Theorem yields that    for all $\o \ins \O$ except on a $\cN^1_\nu \ins \sN_P \big(\cF^{B^t}_\infty\big)$,
   $\nu_s  (\o)   \= \sI \big(  \mu_s (\o ) \big)  $ for a.e. $s \ins [t,\infty)$.
  By Part (1b) and a standard approximation scheme, we can find a $[0,1]-$valued $\bF^{W^t}-$predictable  process
   $  \big\{ \nu^o_s \big\}_{s \in [t,\infty)}$ on $\O_0$
   and an $ \cN^2_\nu \ins \sN_P(\cF^{B^t}_\infty)$ such that
   $\nu_s(\o) \= \nu^o_s(B^{t,\bw}(\o))$ for any $(s,\o) \ins [t,\infty) \ti (\cN^2_\nu)^c $.

   Set $\cN_\mu \= \cN^1_\nu \cp \cN^2_\nu \ins \sN_P(\cF^{B^t}_\infty) $
   and define $  \mu^o_s(\o_0 ) \df \sI^{-1} \big(\nu^o_s (\o_0)\big) \b1_{\{\nu^o_s (\o_0) \in \fE\}}  \+ u_0 \b1_{\{\nu^o_s (\o_0) \notin \fE\}}
  $, $ (s,\o_0) \ins [t,\infty) \ti \O_0 $, which is  a $\hU-$valued $ \bF^{W^t}-$predictable process. 
  Given $\o \ins \cN^c_\mu$, it holds for a.e. $s \ins (t,\infty)$ that
  $\nu^o_s(B^{t,\bw}(\o)) \= \nu_s(\o) \= \sI \big(  \mu_s (\o ) \big)  $ and thus   $ \mu^o_s (B^{t,\bw}(\o)) \= \mu_s(\o) $.

 \ss \no {\bf 2)}  Let $\bx \ins \OmX$ and let $(\vf,n) \ins \fC(\hR^{d+l}) \ti \hN $.  On $\O_0$,
  $   
    M^{t,\mu^o}_s(\vf)   \df   \vf \big(W^t_s  , \, X^{t,\bx,\mu^o}_s \big)
    \- \n \int_t^s  \n  \ol{b}  \big( r, \, X^{t,\bx,\mu^o}_{r \land \cd},\mu^o_r \big)   \cd   D \vf \big( W^t_r  , \, X^{t,\bx,\mu^o}_r \big) dr
    \\ -   \frac12 \n \int_t^s  \n  \ol{\si} \, \ol{\si}^T  \big( r, \, X^{t,\bx,\mu^o}_{r \land \cd},\mu^o_r  \big) \n : \n D^2 \vf  ( W^t_r, \, X^{t,\bx,\mu^o}_r  )   dr  $,   $ \fa s \ins [t,\infty) $
  is an $\bF^{W^t,P_0}-$adapted continuous process
 and  $\tau^{t,\mu^o}_n   \df  \inf\big\{s \ins [t,\infty) \n : \big|(W^t_s  , \, X^{t,\bx,\mu^o}_s ) \big|   \gs n  \big\} \ld (t\+n) $  is an $\bF^{W^t,P_0}-$stopping time.
 As   $ X^{t,\bx,\mu^o}  $ is
  the unique strong solution of SDE \eqref{121621_11} on $   \big( \O_0,\sB(\O_0),   P_0   \big) $ with $  (B^t, \mu ) \= \big( W^t  , \mu^o \big) $,
 taking $(\O,\cF,P,B,X,\mu) \= \big(\O_0,\sB(\O_0),P_0,W,X^{t,\bx,\mu^o},\mu^o \big) $ in Part \(iii\) of Proposition \ref{prop_MPF1} shows that     $ \Big\{ M^{t,\mu^o}_{s \land \tau^{t,\mu^o}_n } (\vf) \Big\}_{s \in [t,\infty)} $  is a bounded $\bF^{W^t,P_0}-$martingale.

 Since $ W^t_s(B^{t,\bw}(\o)) 
  \= B^t_s(\o) $, $\fa (s,\o) \ins [t,\infty) \ti \O$,
   applying Lemma \ref{lem_122921_11} with $t_0 \= t$,  $(\O_1, \cF_1, P_1,B^1)   \= \big(\O, \cF, P , B \big) $, $(\O_2, \cF_2, P_2,B^2) \= \big(\O_0, \sB(\O_0) , P_0 , W\big) $ and $\Phi \= B^{t,\bw} $ implies that
     $ \breve{X}_s \df X^{t,\bx,\mu^o}_s  ( B^{t,\bw} )$, $  \breve{M}_s(\vf) \df  \big( M^{t,\mu^o}_s(\vf) \big) (B^{t,\bw} )$, $s \ins [t,\infty)$
     are $\bF^{B^t,P}-$adapted continuous processes and  $\breve{\tau}_n \df \tau^{t,\mu^o}_n  (B^{t,\bw})   $  is an $\bF^{B^t,P }-$stopping time.

      Let $ t \ls s \< r \< \infty$ and $ \big\{(s_i,\cE_i)\big\}^k_{i=1} \sb [t,s] \ti \sB(\hR^d)  $.
     We can also deduce from  Lemma \ref{lem_122921_11}   that
    \beas
      0 & \tn \= & \tn   E_{P_0} \Big[ \Big(  M^{t,\mu^o}_{  r \land \tau^{t,\mu^o}_n   } (\vf) \- M^{t,\mu^o}_{ s \land  \tau^{t,\mu^o}_n   } (\vf) \Big) \b1_{\ccap{i=1}{k} (W^t_{s_i})^{-1} (\cE_i)} \Big] \\
      & \tn  \= & \tn  \int_{\o_0 \in \O_0}  \Big( \big( {M^{t,\mu^o}} (\vf) \big) \big(r \ld \tau^{t,\mu^o}_n (\o_0),\o_0 \big) \- \big( {M^{t,\mu^o}} (\vf) \big) \big(s \ld \tau^{t,\mu^o}_n (\o_0),\o_0 \big) \Big)   \b1_{\ccap{i=1}{k} \{  W^t_{s_i} (\o_0) \in \cE_i \}} (P \nci (B^{t,\bw})^{-1}) (d \o_0) \\
      & \tn  \= & \tn  \int_{\o  \in \O }  \Big( \big( {M^{t,\mu^o}} (\vf) \big) \big(r \ld \tau^{t,\mu^o}_n (B^{t,\bw}(\o)),B^{t,\bw}(\o) \big) \- \big( {M^{t,\mu^o}} (\vf) \big) \big(s \ld \tau^{t,\mu^o}_n (B^{t,\bw}(\o)),B^{t,\bw}(\o) \big) \Big)  \b1_{\ccap{i=1}{k} \{  W^t_{s_i} (B^{t,\bw}(\o)) \in \cE_i \}}  P    (d \o ) \\
       & \tn  \= & \tn  \int_{\o  \in \O }  \Big( \big( {M^{t,\mu^o}} (\vf) \big) \big(r \ld  \breve{\tau}_n  (\o) ,B^{t,\bw}(\o) \big) \- \big( {M^{t,\mu^o}} (\vf) \big) \big(s \ld \breve{\tau}_n  (\o),B^{t,\bw}(\o) \big) \Big)  \b1_{\ccap{i=1}{k} \{  B^t_{s_i}  (\o)  \in \cE_i \}}  P    (d \o ) \\
      & \tn  \= & \tn  E_P \Big[ \Big( \breve{M}_{  r \land \breve{\tau}_n   } (\vf) \- \breve{M}_{ s \land \breve{\tau}_n   } (\vf) \Big) \b1_{\ccap{i=1}{k} (B^t_{s_i})^{-1} (\cE_i)} \Big] .
      \eeas
      So the Lambda-system $\breve{\L}_{s,r} \df \big\{ A \ins \cF^{B^t,P}_\infty \n : E_P \big[ \big( \breve{M}_{  r \land \breve{\tau}_n   } (\vf) \- \breve{M}_{ s \land \breve{\tau}_n   } (\vf) \big) \b1_A \big] \= 0 \big\}$ includes
      the Pi-system $\Big\{ \ccap{i=1}{k} (B^t_{s_i})^{-1} (\cE_i) \n : \big\{(s_i,\cE_i)\big\}^k_{i=1} \sb [t,s] \ti \sB(\hR^d) \Big\} \cp \sN_P(\cF^{B^t}_\infty)$.
      An application of Dynkin's Pi-Lambda Theorem (see e.g Theorem 3.2 of \cite{Billingsley_PM}) shows that
      $ \cF^{B^t,P}_s \sb \breve{\L}_{s,r} $, i.e., $  E_P \big[ \big( \breve{M}_{  r \land \breve{\tau}_n   } (\vf) \- \breve{M}_{ s \land \breve{\tau}_n   } (\vf) \big) \b1_A \big] \= 0  $ for any $A \ins \cF^{B^t,P}_s$.
      To wit, the $\bF^{B^t,P}-$adapted continuous process
  $\breve{M}_s(\vf) \=    \vf \big( B^t_s  , \breve{X}_s \big)
    \- \n \int_t^s  \n  \ol{b}  \big( r, \breve{X}_{r \land \cd},\mu_r \big) \n \cd \n D \vf \big( B^t_r  , \breve{X}_r \big) dr
    \-   \frac12 \n \int_t^s  \n  \ol{\si} \, \ol{\si}^T  \big( r, \breve{X}_{r \land \cd},\mu_r  \big) \n : \n D^2 \vf  ( B^t_r, \breve{X}_r  )   dr$, $\fa s \ins [t,\infty)$ stopped by the $\bF^{B^t,P}-$stopping time
    $ \breve{\tau}_n   \=  \inf\big\{s \ins [t,\infty) \n : \big|(B^t_s  , \breve{X}_s ) \big|   \gs n  \big\} \ld (t\+n)$
    is a bounded $\bF^{B^t,P}-$martingale.  Using Part (i) of Proposition \ref{prop_MPF1} yields that $P\{   \breve{X}_s \=   X^{t,\bx,\mu}_s,   \fa s \ins [0,\infty)\}=1$.

   Moreover, for any $ \psi \ins \sH_o $, one has
   $  E_P  \big[ \int_t^\infty \n \psi^-(r,  X^{t,\bx,\mu}_{r \land \cd},\mu_r  ) dr \big]
  \= E_P  \big[ \int_t^\infty \n \psi^- \big(r, X^{t,\bx,\mu^o}_{r \land \cd} (B^{t,\bw}),\mu^o_r(B^{t,\bw}) \big) dr \big]
   \= E_{P_0} \big[ \int_t^\infty \n \psi^-(r,X^{t,\bx,\mu^o}_{r \land \cd},\mu^o_r ) dr \big]  \< \infty   $.    \qed

   \if{0}

 \no {\bf Proof of Remark \ref{rem_ocP}: 1)} By Lemma \ref{lem_061222}, $\big\{\wh{\mu}_s\big\}_{s \in [t,\infty)}$ is a $\hU-$valued, $\bF^{W^t}-$predictable process on $\O_0$
if and only if there exists a $\hU-$valued, $\bF^{W}-$predictable process $\big\{\ddot{\mu}_\fs\big\}_{\fs \in [0,\infty)}$ on $\O_0$ such that
$ \ddot{\mu}_\fs \big(\sW^t(\o_0)\big) \= \wh{\mu}_{t+\fs} (\o_0)    $ for any $  (\fs,\o_0) \ins [0,\infty) \ti \O_0$.
So it holds for any $  (\fs,\oo) \ins [0,\infty) \ti \oO $ that
$ \ddot{\mu}_\fs \big(\osW^t (\oo) \big) \= \ddot{\mu}_\fs \big(\sW^t\big(\oW(\oo)\big)\big) \= \wh{\mu}_{t+\fs} \big(\oW(\oo)\big)    $.

   Also, the equivalence between \(D2\)+\(D3\) and  \(D2\,$'$\) under (D1) is straightforward  by the martingale-problem formulation of controlled  SDEs on $\oO$.

 \no {\bf  2)} Next, we   verify the equivalence between (D4) and (D4') under (D2).
 Define two processes on $\O_0$ by:
\beas
    \beta_s (\o_0) \df W_{s \vee t} (\o_0) \- W_t (\o_0) \aand
  \ddot{\beta}_s(\o_0) \df W_{(s-t)^+}(\o_0)   ,  \q \fa (s,\o_0) \ins [0,\infty) \ti \O_0 .
\eeas

\no {\bf  2a)} We first show   (D4) \n $\Ra$ \n (D4') under (D2):   Let $\wh{\tau}$ be   a $[t,\infty]-$valued $\bF^{W^t,P_0} -$stopping time   on $\O_0$. We claim that
\beas
  \wh{\tau}\= \wh{\tau}(\beta)  , \q \hb{$P_0- $a.s.}
  \eeas
  Since it holds for any $(r,\cE) \ins [t,\infty) \ti \sB(\hR^d)$ that    $\beta^{-1} \big(\{W^t_r \ins \cE\}\big) \= \big\{W^t_r(\beta) \ins \cE \big\} \=  \{ W_r(\beta) \- W_t(\beta) \ins \cE \big\}
 \=  \{  \beta_r \-  \beta_t \ins \cE \big\}
\= \{W_r  \- W_t  \ins \cE\} \= \{W^t_r   \ins \cE\} $,
  the sigma-field $\{A_0 \sb \O_0 : \beta^{-1}(A_0) \= A_0 \} $ contains all generating sets of $\cF^{W^t}_\infty$ and thus includes
$\cF^{W^t}_\infty$.
 For any $\cN_0 \ins \sN_{P_0} \big(\cF^{W^t}_\infty\big)$, there exists $A_0 \ins \cF^{W^t}_\infty$ such that $\cN_0  \sb A_0  $ and  $ P_0(A_0) \= 0$.
As $ \beta^{-1}(\cN_0) \sb \beta^{-1}(A_0) \= A_0 $, we see that $ \beta^{-1}(\cN_0) \ins \sN_{P_0}\big(\cF^{W^t}_\infty\big)$. Hence,
 \bea \label{010722_11}
 \beta^{-1}(A_0) \= A_0 , \q \fa A_0 \ins \cF^{W^t}_\infty \aand
 \beta^{-1}(\cN_0) \sb \sN_{P_0} \big(\cF^{W^t}_\infty\big) , \q \fa  \cN_0 \ins \sN_{P_0} \big(\cF^{W^t}_\infty \big) .
 \eea

 Let $n \ins \hN$ and set $s^n_i \df t \+ i2^{-n}$, $\fa i \ins \hN \cp \{0\}   $.
  We denote $A^n_i  \df  \{ s^n_{i-1}   \ls  \wh{\tau}  \< s^n_i  \}  \ins \cF^{W^t,P_0}_{ s^n_i }  $, $\fa i \ins \hN   $
  and $A^n_\infty  \df  \{    \wh{\tau}  \= \infty  \}  \ins \cF^{W^t,P_0}_\infty  $.
 Define $\wh{\tau}_n \df \sum_{i \in \hN} s^n_i \b1_{A^n_i}   \+ \infty \b1_{A^n_\infty}$.
 Clearly,
    $\wh{\tau}_n(\beta) \= \sum_{i \in \hN} s^n_i \b1_{\beta^{-1}(A^n_i)}  \+ \infty \b1_{\beta^{-1}(A^n_\infty)} $
    is equal to $ \wh{\tau}_n $  on $\cA_n \df \ccup{i \in \hN \cup \{\infty\}}{} \big(  A^n_i  \Cp \beta^{-1}(A^n_i)   \big) $.
  For any $i \ins \hN \cp \{\infty\}$, since there exists  some    $\wA^n_i   \ins  \cF^{W^t}_{ s^n_i }$
   such that   $ \cN^n_i \df  A^n_i   \D  \wA^n_i  \ins  \sN_{P_0}\big(\cF^{W^t}_\infty\big)  $ (see e.g.   Problem 2.7.3 of \cite{Kara_Shr_BMSC}),
one can deduce from \eqref{010722_11} that
 $ (A^n_i)^c \Cp \beta^{-1}(A^n_i)   \sb 
\big( (A^n_i)^c \Cp \beta^{-1} \big(\wA^n_i  \big)   \big) \cp \big( (A^n_i)^c \Cp \beta^{-1}(  \cN^n_i)   \big)
\= \big(  (A^n_i)^c \Cp \wA^n_i    \big) \cp \big( (A^n_i)^c \Cp \beta^{-1}(  \cN^n_i)   \big)
\sb \cN^n_i \cp \beta^{-1}(  \cN^n_i) \ins \sN_{P_0} \big(\cF^{W^t}_\infty \big)  $.
  So $\cA^c_n \= \ccup{i \in \hN \cup \{\infty\}}{} \big( ( A^n_i )^c \Cp \beta^{-1}(A^n_i)   \big) \ins \sN_{P_0} \big(\cF^{W^t}_\infty\big)$.

 Set $\cN_* \df \ccup{n \in  \hN}{} \cA^c_n   \ins \sN_{P_0} \big(\cF^{W^t}_\infty\big)$.
 Given $\o_0 \ins \cN^c_* \= \ccap{n \in  \hN}{} \cA_n $, as $\wh{\tau}   \= \lmtd{n \to \infty} \wh{\tau}_n   $, we have
 $ \wh{\tau} \big(\beta(\o_0)\big) \= \lmtd{n \to \infty} \wh{\tau}_n  \big(\beta(\o_0)\big) \= \lmtd{n \to \infty} \wh{\tau}_n  ( \o_0 )
 \= \wh{\tau} (\o_0) $, proving the claim.

 Since $\oW^t$ is a Brownian motion under $\oP$ by (D2), applying Lemma \ref{lem_122921_11} with $t_0 \= t$, $(\O_1, \cF_1, P_1,B^1)   \= \big(\oO ,  \sB(\oO ),  \oP, \oW \big) $, $(\O_2, \cF_2, P_2,B^2) \= \big(\O_0,  \sB(\O_0),  P_0,  W \big) $  and $\Phi \= \oW$ yields that
 $  \ocN_* \df \oW^{-1} (\cN_*) \ins \sN_\oP \big(\cF^{\oW^t}_\infty\big) $. So
 \bea \label{010722_14}
 \wh{\tau} \big(\beta(\oW(\oo))\big)  \= \wh{\tau} \big(\oW(\oo)\big) , \q \fa \oo \ins \ocN^c_* \=  \oW^{-1} (\cN^c_*) .
 \eea

  As a shifted canonical process on $\O_0$,   $\sW^t_\fs(\o_0) \df W_{t+\fs}(\o_0) \- W_t (\o_0) \= W^t_{t+\fs}(\o_0) $, $  (\fs,\o_0) \ins [0,\infty) \ti \O_0$  is also a Brownian motion under $P_0$.  
   Since it holds for any $(\fs,\o_0) \ins [0,\infty) \ti \O_0$ that
   $ \sW^t_\fs \big( \ddot{\beta} (\o_0)\big) \= W_{t+\fs} \big( \ddot{\beta} (\o_0)\big) \- W_t \big( \ddot{\beta} (\o_0)\big)
    \= \ddot{\beta}_{t+\fs} (\o_0)  \-   \ddot{\beta}_t (\o_0) 
    \=  W_\fs(\o_0) $,
 using Lemma \ref{lem_122921_11} again with $t_0 \= 0$, $(\O_1, \cF_1, P_1,B^1)   \= \big(\O_0,  \sB(\O_0),  P_0,W\big) $, $(\O_2, \cF_2, P_2,B^2) \= \big(\O_0,  \sB(\O_0),  P_0, \sW^t\big) $ and $\Phi \= \ddot{\beta}$   shows that
  \bea  \label{010622_11}
  \ddot{\beta}^{-1} \big(\cF^{\sW^t,P_0}_\fs\big) \sb \cF^{W,P_0}_\fs  , \q \fa \fs \ins [0,\infty] .
  \eea
 As  $\cF^{\sW^t}_\fs \= \si\big(\sW^t_\fr; \fr \ins [0,\fs] \Cp \hR \big) 
 \= \si\big( W^t_r; r \ins [t,t\+\fs]  \Cp \hR \big) \= \cF^{W^t}_{t+\fs} $ for any $\fs \ins [0,\infty]$, we also obtain that
 $  
  \cF^{\sW^t,P_0}_\fs \= \si\big( \cF^{\sW^t}_\fs \cp \sN_{P_0} (\cF^{\sW^t}_\infty) \big)
 \= \si\big( \cF^{W^t}_{t+\fs} \cp \sN_{P_0} \big(\cF^{W^t}_\infty\big) \big) \= \cF^{W^t,P_0}_{t+\fs} $, $ \fa \fs \ins [0,\infty) $.

 Define $\ddot{\tau}(\o_0) \= \wh{\tau} \big( \ddot{\beta} (\o_0)\big) \- t$, $\fa \o_0 \ins \O_0$.
 For any  $\fs \ins [0,\infty)$, since $ \big\{\wh{\tau} \ins [t,t\+\fs]\big\} \ins \cF^{W^t,P_0}_{t+\fs} \= \cF^{\sW^t,P_0}_\fs $,
   \eqref{010622_11} renders that $\big\{ \ddot{\tau} \ins [0,\fs]\big\} \= \big\{\wh{\tau} \big(\ddot{\beta}) \ins [t,t\+\fs] \big\} \= \ddot{\beta}^{-1}\big(\{\wh{\tau} \ins [t,t\+\fs]\}\big)  \ins \cF^{W,P_0}_\fs$.
 So $ \ddot{\tau} $ is a  $[0,\infty]-$valued $\bF^{W,P_0}-$stopping time on $\O_0$.
 For any $ s \ins [0,\infty) $,
 $ \ddot{\beta}_s (\osW^t)   \= \osW^t_{(s-t)^+} \= \oW_{t + (s-t)^+} \- \oW_t  \= \oW_{s \vee t} \- \oW_t \=    \beta_s (\oW)   $.
Then \eqref{010722_14} implies that
  $ \ddot{\tau} \big(\osW^t(\oo)\big) \= \wh{\tau} \big( \ddot{\beta} \big(\osW^t(\oo)\big)\big) \- t
   \= \wh{\tau} \big(\beta (\oW(\oo))\big) \- t \= \wh{\tau}  \big( \oW (\oo) \big) \- t $, $ \fa \oo \ins \ocN^c_* $.
 If   $   \oP \big\{ \oT \=   \wh{\tau}(\oW ) \big\} \= 1$, we also have
 $  \oP \big\{ \oT \= t \+ \ddot{\tau}  \big(\osW^t \big)   \big\} \= 1 $. Hence,  (D4) gives rise to   (D4') under (D2).

\no {\bf  2b)}  We then  show  (D4') $\Ra$  (D4): Let $ \ddot{\tau} $ be a  $[0,\infty]-$valued $\bF^{W,P_0}-$stopping time on $\O_0$.

 Since   ${\ddot{\beta}}^t_s \df \ddot{\beta}_s \- \ddot{\beta}_t \= W_{s-t}$, $s \ins [t,\infty)$ is   a Brownian motion on $\big(\O_0,  \sB(\O_0),  P_0\big)$
 and since $  {\ddot{\beta}}^t_s\big(\sW^t(\o_0)\big) \= \sW^t_{s-t}(\o_0) \= W_s(\o_0) \- W_t(\o_0) \= W^t_s (\o_0)$ for any $(s,\o_0) \ins [t,\infty) \ti \O_0 $,
 applying Lemma \ref{lem_122921_11} with $t_0 \= t$,  $(\O_1, \cF_1, P_1,B^1)   \= \big(\O_0,  \sB(\O_0),  P_0,W\big) $, $(\O_2, \cF_2, P_2,B^2) \= \big(\O_0,  \sB(\O_0),  P_0, \ddot{\beta}\big) $ and $\Phi \= \sW^t$   shows that
  \bea \label{010722_17}
   (\sW^t)^{-1} ( \cF^{{\ddot{\beta}}^t,P_0}_s) \sb \cF^{W^t,P_0}_s , \q \fa s \ins [t,\infty) .
   \eea
       Since $\cF^{{\ddot{\beta}}^t}_s \= \si\big({\ddot{\beta}}^t_r; r \ins [t,s] \Cp \hR \big)   \= \si\big( W_{r-t}; r \ins [t,s]  \Cp \hR \big)
 \= \si\big( W_{r'}; r' \ins [0,s\-t]  \Cp \hR  \big) \= \cF^W_{s-t} $ for any $s \ins [t,\infty]$,  we see that
 $ \cF^{{\ddot{\beta}}^t,P_0}_s \= \si\big( \cF^{{\ddot{\beta}}^t}_s \cp \sN_{P_0} (\cF^{{\ddot{\beta}}^t}_\infty) \big)
 \= \si\big( \cF^W_{s-t} \cp \sN_{P_0} (\cF^W_\infty) \big) \= \cF^{W,P_0}_{s-t} $, $ \fa s \ins [t,\infty) $.

 Define $\wh{\tau} (\o_0) \df  t \+ \ddot{\tau} \big(\sW^t (\o_0)\big)   $, $\fa \o_0 \ins \O_0$. For any $s \ins [t,\infty)$,
 as $\big\{ \ddot{\tau}  \ins [0,s\-t]\big\} \ins \cF^{W,P_0}_{s-t} \= \cF^{{\ddot{\beta}}^t,P_0}_s $,  \eqref{010722_17} yields that
 $ \big\{\wh{\tau} \ins [t,s]\big\} \= \{ \ddot{\tau} (\sW^t) \ins [0,s\-t]\}
\= (\sW^t)^{-1} \big(\{ \ddot{\tau}  \ins [0,s\-t]\}\big) \ins 
\cF^{W^t,P_0}_s  $.
So $\wh{\tau}  $ is a  $[t,\infty]-$valued $\bF^{W^t,P_0}-$stopping time on $\O_0$.
Clearly, $\wh{\tau} \big(\oW(\oo)\big)   \= t \+ \ddot{\tau} \big(\sW^t (\oW(\oo))\big)
\=  t \+ \ddot{\tau} \big(\osW^t (\oo)\big) $ for any $\oo \ins \oO$. It means that (D4') directly implies (D4). \qed

   \fi

 \no {\bf Proof of Theorem \ref{thm_V=oV}:}
  Fix  $(t,\bw,\bu,\bx) \ins [0,\infty) \ti \O_0 \ti \hJ  \ti \OmX  $ and $(y,z) \= \big(\{y_i\}_{i \in \hN}, \{z_i\}_{i \in \hN}\big)
 \ins \Re \ti \Re$.

  \no {\bf 1)} We first show   that $V (t,\bx,y,z )  \ls \oV (t, \bw,\bu,\bx, y,z ) $:
 If $\cC_{t,\bx}(y,z) \= \es$, then   $V (t,\bx,y,z) \= -\infty \ls \oV (t, \bw,\bu,\bx, y,z )$. 

 So we assume   $\cC_{t,\bx}(y,z)  \nne \es$ and let $(\mu,\tau) \ins \cC_{t,\bx}(y,z) $.

 \no {\bf 1a)}  Define a  process  $ \cB^{t,\bw}_s  (\o) \df   \bw(s \ld t)  \+ \cB^t_{s \vee t} (\o)     $,
 $ \fa (s,\o) \ins [0,\infty) \ti \cQ  $.
 According to Proposition \ref{prop_122021b}, there exists   a   $\hU-$valued, $ \bF^{W^t} -$predictable process $ \wh{\mu} \=   \{ \wh{\mu}_s \}_{s \in [t,\infty)} $
 on $\O_0$ and an $\cN_\mu  \ins \sN_\fp(\cF^{\cB^t}_\infty)$ such that for any $\o \ins \cN^c_\mu$,
  $  \mu_s (\o) \= \wh{\mu}_s \big(\cB^{t,\bw}(\o)\big)  $ for a.e. $s \ins (t,\infty)   $.
  By Fubini Theorem,
   $ 
      \cN  \df \big\{\o  \ins \cQ \n :  \mu_s  (\o ) \hb{ is not Borel-measurable in }s \in [t,\infty) \big\} $ is a $\cF^{\cB^t,\fp}_\infty-$measurable
   set with zero $\fp-$measure  or $  \cN  \ins \sN_\fp \big(\cF^{\cB^t}_\infty\big)$.
  Then   process
 $ \wt{\mu}_s (\o) \df \big( \b1_{\{s \in [0,t)\}} \bu(s)\+ \b1_{\{s \in [t,\infty)\}} \mu_s (\o)\big)   \b1_{ \{\o \in  \cN^c \} }  \+ u_0 \b1_{ \{\o \in  \cN \} } $, $ \fa (s,\o) \ins [0,\infty) \ti \cQ $ satisfies that  $\wt{\mu}_\cd(\o) \ins \hJ$ for any $\o \ins \cQ$.

 Let us simply set $\th \= (t,\bx,\mu)$. Define  a mapping $\Psi  \n : \cQ \mto \oO$ by
  $ \Psi (\o) \df \big( \cB^{t,\bw}  (\o), \wt{\mu}_\cd (\o), \cX^\th  (\o) ,  \tau(\o) \big) \ins \oO $,  $ \fa \o \ins \cQ $.
 Since    $\oW^t_s(\Psi (\o)) \= \oW_s(\Psi (\o)) \-  \oW_t(\Psi (\o)) \=  \cB^{t,\bw}_s(\o) \- \cB^{t,\bw}_t(\o) \= \cB^t_s(\o)     $,
 $\fa (s,\o) \ins [t,\infty) \ti \cQ$
 and since $ 
 \big\{ \cX^\th_s\big\}_{s \in [0,\infty)} $ is an $ \big\{ \cF^{\cB^t,\fp}_{s \vee t} \big\}_{s \in [0,\infty)}   -$adapted continuous process,
 we can deduce  that    the mapping  $ \Psi $ 
 is $\cF^{\cB^t,\fp}_s \big/ \ocF^t_s-$measurable for any $s \ins [t,\infty)$.

   Let $\vf \ins L^0 \big((0,\infty) \ti \hU ; \hR \big)$.
   The $\bF^{\cB^t,\fp}-$progressive measurability of process  $\big\{\wt{\mu}_s\big\}_{s \in [t,\infty)}$ implies that
    process $\big\{\vf (s,  \wt{\mu}_s )\big\}_{s \in [t,\infty)}$ is also  $\bF^{\cB^t,\fp}-$progressively measurable
  and   the   random variable
 $ I_\vf (\wt{\mu}_\cd(\o)) 
 \= \int_0^t \vf  (s, \bu(s) ) ds
  \+ \int_t^\infty \vf  (s, \wt{\mu}_s (\o) ) ds $, $\o \ins  \cQ $ is thus $\cF^{\cB^t,\fp}_\infty-$measurable.
     Lemma \ref{lem_M29_01} (1) then renders that the mapping $\wt{\mu}_\cd \n : \cQ \mto \hJ $ is   $\cF^{\cB^t,\fp}_\infty \big/ \sB(\hJ)-$measurable,
    which together with the $\cF^{\cB^t,\fp}_\infty-$measurability of $  \cB^{t,\bw}   $, $ \cX^\th  $, $ \tau $ shows that
  the mapping  $ \Psi $ is also  $\cF^{\cB^t,\fp}_\infty \big/ \sB(\oO)-$measurable.
    \if{0}


 Let $W^X \= \{W^X_s\}_{s \in [0,\infty)}$ be the canonical process of $\OmX$.
 For any $s \ins [0,\infty)$, $\cE_1 \ins \sB(\hR^d)$ and $\cE_2 \ins \sB(\hR^l)$, we can deduce that
  $ \big\{\o \ins \cQ \n :  \cB^{t,\bw} (\o) \ins W^{-1}_s (\cE_1) \big\}
  \= \big\{\o \ins \cQ \n : W_s \big(\cB^{t,\bw} (\o)\big) \ins \cE_1 \big\}
  \= \big\{\o \ins \cQ \n : \cB^t_{s \vee t} (\o) \ins \cE^{t,\bw}_1 \big\} \ins \cF^{\cB^t}_{s \vee t} \sb \cF^{\cB^t}_\infty$
   with $\cE^{t,\bw}_1 \df \{a \-  \bw(s \ld t)   \n : a \ins \cE_1\} \ins  \sB(\hR^d)$,
  and that  $ \big\{\o \ins \cQ \n :  \cX^\th (\o) \ins (W^X_s)^{-1} (\cE_2) \big\}
  \= \big\{\o \ins \cQ \n : W^X_s \big(\cX^\th (\o)\big) \ins \cE_2 \big\}
  \= \big\{\o \ins \cQ \n : \cX^\th_s (\o) \ins \cE^{t,\bw}_2 \big\} \ins \cF^{\cB^t,\fp}_{s \vee t}  \sb \cF^{\cB^t,\fp}_\infty$.
  So the sigma-field $\big\{A \sb \cQ \n : (\cB^{t,\bw})^{-1}(A) \ins \cF^{\cB^t }_\infty  \big\}$ of $\cQ$ contains all generating sets of
  $ \cF^W_\infty \= \sB(\O_0) $
  and   the sigma-field $\big\{A \sb \cQ \n : (\cX^\th)^{-1}(A) \ins \cF^{\cB^t,\fp}_\infty  \big\}$ of $\cQ$ contains all generating sets of
  $ \cF^{W^X}_\infty \= \sB(\OmX) $. It follows that
  the mapping $\cB^{t,\bw} \n : \cQ \mto \O_0 $ is $\cF^{\cB^t}_\infty \big/ \sB(\O_0)-$measurable
  and  the mapping $\cX^\th \n : \cQ \mto \OmX $ is $\cF^{\cB^t,\fp}_\infty \big/ \sB(\OmX)-$measurable.

  Let $\vf \ins L^0 \big((0,\infty) \ti \hU ; \hR \big)$.
  The  measurability of function $\vf $ and the $\bF^{\cB^t,\fp}-$progressive measurability of process
  $\big\{\wt{\mu}_s\big\}_{s \in [t,\infty)}$ imply that
  the process $\big\{\vf (s,  \wt{\mu}_s )\big\}_{s \in [t,\infty)}$ is   $\bF^{\cB^t,\fp}-$progressively measurable
  and   the   random variable $  \int_t^\infty \vf (s, \wt{\mu}_s (\o)) ds $, $\o \ins \cQ$   is thus $\cF^{\cB^t,\fp}_\infty-$measurable.
     It follows that the random variable
 $ I_\vf (\wt{\mu}_\cd(\o)) 
 \= \int_0^t \vf  (s, \bu(s) ) ds
  \+ \int_t^\infty \vf  (s, \wt{\mu}_s (\o) ) ds $, $\o \ins  \cQ $ is $\cF^{\cB^t,\fp}_\infty-$measurable.
  Namely,     $\big\{\o \ins \cQ \n : \wt{\mu}_\cd (\o) \ins I^{-1}_\vf(\cE) \big\} \= \big\{\o \ins \cQ \n : I_\vf(\wt{\mu}_\cd (\o)) \ins  \cE  \big\}   \ins \cF^{\cB^t,\fp}_\infty $ for any $\cE \ins \sB[-\infty,\infty]$.
  Then by Lemma \ref{lem_M29_01} (1),    the sigma-field $\big\{A \sb \cQ \n :  \wt{\mu}_\cd^{-1}(A) \ins \cF^{\cB^t,\fp }_\infty  \big\}$ of $\cQ$ contains all generating sets of $\sB(\hJ)$ and the mapping $\wt{\mu}_\cd \n : \cQ \mto \hJ $ is thus $\cF^{\cB^t,\fp}_\infty \big/ \sB(\hJ)-$measurable.

     Since the  $\bF^{\cB^t,\fp}-$stopping time  $\tau$ satisfies that   $    \big\{\o  \ins \cQ \n : \tau(\o) \ins [t,s] \big\} \ins \cF^{\cB^t,\fp}_s \sb \cF^{\cB^t,\fp}_\infty$ for any $s \ins [t,\infty]$,
      the sigma field of $[t,\infty]$,
   $\L_\tau \df \big\{\cE \sb [t,\infty] \n :   \tau^{-1}(\cE) \ins \cF^{\cB^t,\fp}_\infty \big\}$ contains all closed intervals $[t,s]$, $s \ins [t,\infty]$, which generates $\sB[t,\infty]$. It follows that
   $\sB[t,\infty] \sb \L_\tau $ or $  \tau^{-1}(\cE) \ins \cF^{\cB^t,\fp}_\infty $, $ \fa \cE \ins \sB[t,\infty] $.
   So $\tau \n: \cQ \mto [t,\infty]$ is $\cF^{\cB^t,\fp}_\infty   -$measurable,
   the mapping $\Psi \n : \cQ \mto \oO$ is eventually $\cF^{\cB^t,\fp}_\infty \big/ \sB(\oO)-$measurable.

  Let $s \ins [t,\infty)$. For any $(r_1,\cE_1) \ins [t,s] \ti \sB(\hR^d)$ and $(r_2,\cE_2) \ins [0,s] \ti \sB(\hR^l)$,
  one has  $ \big\{\o \ins \cQ \n :  \Psi (\o) \ins (\oW^t_{r_1})^{-1} (\cE_1) \big\}
  \= \big\{\o \ins \cQ \n :  \oW^t_{r_1}\big(\Psi (\o)\big) \ins  \cE_1  \big\}
  \= \big\{\o \ins \cQ \n :  \cB^t_{r_1} (\o)  \ins  \cE_1  \big\} \ins \cF^{\cB^t}_s $
  and  $ \big\{\o \ins \cQ \n :  \Psi (\o) \ins  \oX_{r_2}^{-1} (\cE_2) \big\}
  \= \big\{\o \ins \cQ \n :  \oX_{r_2}\big(\Psi (\o)\big) \ins  \cE_2  \big\}
  \= \big\{\o \ins \cQ \n :  \cX^\th_{r_2} (\o)  \ins  \cE_2  \big\} \ins \cF^{\cB^t,\fp}_s $.
   So the sigma-field $\big\{A \sb \cQ \n : \Psi^{-1}(A) \ins \cF^{\cB^t,\fp}_s  \big\}$ of $\cQ$ contains all generating sets of
   $\ocF^t_s$ and thus includes $\ocF^t_s$.
   Then  the mapping $\Psi  $ is also $\cF^{\cB^t,\fp}_s \big/ \ocF^t_s -$measurable.

     \fi

 Let $\oP_\Psi \ins \fP\big(\oO\big)  $ be the probability measure induced by $\Psi$, i.e.,
 $  \oP_\Psi \big(\oA\big) \df \fp \big( \Psi^{-1}\big(\oA\big) \big) $, $ \fa \oA \ins \sB(\oO)  $.

 \no {\bf 1b)} Set $\ol{\mu}_s \df \wh{\mu}_s(\oW)$, $\fa s \ins [t,\infty)$, which is a $\hU-$valued, $ \bF^{\oW^t} -$predictable process on $\oO$.
 Since $ \cN^c \Cp \cN^c_\mu   \sb \cN^c \Cp \big\{ \o \ins \cQ \n : \mu_s(\o) \= \wh{\mu}_s(\cB^{t,\bw}(\o)) \hb{ for a.e. } s \ins (t,\infty) \big\} \= \cN^c \Cp \big\{ \o \ins \cQ \n : \wt{\mu}_s(\o) \= \wh{\mu}_s(\cB^{t,\bw}(\o)) \hb{ for a.e. }   s \ins (t,\infty) \big\}$,
 we can deduce  that
  $ \oP_\Psi \big\{\oU_s \= \ol{\mu}_s \hb{ for a.e. } s \ins (t,\infty) \big\} \= \fp \big\{\oU_s (\Psi) \= \wh{\mu}_s(\oW(\Psi)) \hb{ for a.e. } s \ins (t,\infty) \big\}
 \= \fp \big\{\o \ins \cQ \n : \wt{\mu}_s(\o) \= \wh{\mu}_s(\cB^{t,\bw}(\o)) \hb{ for a.e. } s \ins (t,\infty) \big\}
 \= 1 $. Namely, $\oP_\Psi $ satisfies  (D1) in the  definition of $\ocP_{t,\bx}$.

  Fix  $ (\vf,n)  \ins   \fC(\hR^{d+l}) \ti \hN $.    We define an $\bF^{\cB^t,\fp}-$adapted continuous process
  $   
    \cM^\th_s(\vf)   \df   \vf \big(\cB^t_s  , \cX^\th_s \big)
    \- \n \int_t^s  \n  \ol{b}  \big( r, \cX^\th_{r \land \cd},\mu_r \big) \n \cd \n D \vf \big( \cB^t_r  , \cX^\th_r \big) dr
    \-   \frac12 \n \int_t^s  \n  \ol{\si} \, \ol{\si}^T  \big( r,  \cX^\th_{r \land \cd} ,\mu_r \big) \n : \n D^2 \vf  ( \cB^t_r, \cX^\th_r  )   dr  $, $  \fa s \ins [t,\infty)$
   and define an  $\bF^{\cB^t,\fp}-$stopping time $\btau^\th_n  \df  \inf\big\{s \ins [t,\infty) \n : \big| (\cB^t_s ,\cX^\th_s) \big|    \gs n   \big\} \ld (t\+n) $.
   Applying Proposition \ref{prop_MPF1} with $(\O,\cF,P,B,X,\mu) \= (\cQ,\cF,\fp,\cB,\cX^\th,\mu) $ yields that
   $ \big\{\cM^\th_{s \land \bbtau^\th_n}(\vf)\big\}_{s \in [t,\infty)} $
   is a bounded $ \big( \bF^{\cB^t,\fp} , \fp \big) -$martingale.  
 \if{0}

 Actually, $ \big\{\cM^\th_{s \land \bbtau^\th_n}(\vf)\big\}_{s \in [t,\infty)} $
   is a  martingale with respect to the filtration $\bF^{\cB^t,\cX^\th} \=
   \Big\{ \cF^{\cB^t,\cX^\th}_s \df \si \big(\cB^t_r;r \ins [t,s]\big) \ve \si\big( \cX^\th_r;r \ins [0,s]\big)\Big\}_{s \in [t,\infty)}$ and is thus a martingale with respect to the filtration $\bF^{\cB^t,\cX^\th,\fp} \=
   \Big\{ \cF^{\cB^t,\cX^\th,\fp}_s \df \si \big(\cF^{\cB^t,\cX^\th}_s \cp \sN_\fp(\cF^{\cB^t,\cX^\th}_\infty)\big)\Big\}_{s \in [t,\infty)}$.  Since $\cF^{\cB^t,\fp}_s \=  \si \big(\cF^{\cB^t}_s   \cp \sN_\fp(\cF^{\cB^t}_\infty)\big) \sb \cF^{\cB^t,\cX^\th,\fp}_s $ for any $s \ins [t,\infty)$, we see that $ \big\{\cM^\th_{s \land \bbtau^\th_n}(\vf)\big\}_{s \in [t,\infty)} $
   is also an $ \bF^{\cB^t,\fp} -$martingale.

 \fi

  Since $\oP_\Psi \big\{   \oX_s   \= \bx(s) , \fa s \ins [0,t]  \big\} \= \fp \big\{  \oX_s (\Psi ) \= \bx(s), \fa s \ins [0,t] \big\}
   \= \fp \big\{  \cX^\th_s    \= \bx(s), \fa s \ins [0,t] \big\} \=1$, using Proposition \ref{prop_MPF1}
   with $(\O,\cF,P,B,X,\mu) \= \big(\oO,\sB(\oO),\oP_\Psi,\oW,\oX,\ol{\mu} \big) $
    shows that  $\big\{\oM^{t,\ol{\mu}}_{s \land \otau^t_n } (\vf) \big\}_{s \in  [t,\infty) }  $ is a bounded $\obF^t-$adapted continuous process under $ \oP_\Psi $.
    Given $\o \ins \cN^c_\mu$,
   since 
     $ \ol{\mu}_s \big(\Psi (\o)\big) \= \wh{\mu}_s\big(\oW(\Psi (\o))\big) \= \wh{\mu}_s\big( \cB^{t,\bw}  (\o) \big)   \= \mu_s (\o) $
   for a.e. $s \ins (t,\infty)   $,
   we see that
  $ \big(\oM^{t,\ol{\mu}}_s(\vf)\big)\big(\Psi(\o)\big) \=  \big(\cM^\th_s(\vf) \big) (\o)$, $ \fa s \ins [t,\infty) $ and
  $\otau^t_n \big(\Psi(\o)\big) \= \btau^\th_n (\o)$.  Then
  \bea
   \big(\oM^{t,\ol{\mu}}_{s \land \otau^t_n }(\vf) \big) \big(\Psi(\o)\big)
  & \tn \= & \tn  \big(\oM^{t,\ol{\mu}} (\vf) \big) \big( s \ld \otau^t_n (\Psi(\o)),  \Psi(\o)  \big)
  \= \big(\oM^{t,\ol{\mu}} (\vf) \big) \big( s \ld \btau^\th_n  (\o) ,  \Psi(\o)  \big) \nonumber \\
  & \tn  \= & \tn  \big(\cM^\th (\vf) \big) \big( s \ld \btau^\th_n  (\o) ,  \o   \big)
  \= \big(\cM^\th_{s \land \bbtau^\th_n} (\vf) \big)   (\o) , \q   \fa (s,\o) \ins [t,\infty) \ti \cN^c_\mu . \label{122621_17}
  \eea

  Let $t_1,t_2 \ins [t,\infty)$ with $t_1 \< t_2$ and let $\oA \ins \ocF^t_{t_1} $.
    As  $\Psi^{-1}(\oA) \ins \cF^{\cB^t,\fp}_{t_1}$, the $ \big(\bF^{\cB^t,\fp}, \fp\big) -$martingality of
    $ \big\{\cM^\th_{s \land \bbtau^\th_n}(\vf)\big\}_{s \in [t,\infty)} $ and \eqref{122621_17}    imply that
 $  E_{\oP_\Psi} \Big[ \big(\oM^{t,\ol{\mu}}_{t_2 \land \otau^t_n }(\vf) \-  \oM^{t,\ol{\mu}}_{t_1 \land \otau^t_n }(\vf)  \big) \b1_\oA    \Big]
   \=     E_\fp \Big[ \Big( \big(\oM^{t,\ol{\mu}}_{t_2 \land \otau^t_n }(\vf) \big) (\Psi ) \-  \big(\oM^{t,\ol{\mu}}_{t_1 \land \otau^t_n }(\vf) \big) (\Psi ) \Big) \b1_{\Psi^{-1}(\oA)}  \Big]
     \=    E_\fp \Big[ \big(  \cM^\th_{t_2 \land \bbtau^\th_n }(\vf)      \-  \cM^\th_{t_1 \land \bbtau^\th_n }(\vf)   \big)  \b1_{\Psi^{-1}(\oA)} \Big]   \= 0 $.
  So  $ \big\{ \oM^{t,\ol{\mu}}_{s \land \otau^t_n } (\vf) \big\}_{s \in [t,\infty)} $  is a bounded $ \big(\obF^t,\oP_\Psi\big)   -$martingale.
 By Remark \ref{rem_ocP} (ii), $\oP_\Psi$ satisfies (D2)+(D3) in the  definition of $\ocP_{t,\bx}$.

 \no {\bf 1c)} Since   $   W^t_s( \cB^{t,\bw} (\o))   \n \= \n   \cB^{t,\bw}_s(\o)   \-   \cB^{t,\bw}_t(\o)    \n \= \n    \cB^t_s(\o)    $ for any $(s,\o) \ins [t,\infty) \ti \O $,
 taking $(\O,\cF,P,B,\Phi) \= (\cQ,\cF,\fp,\cB,\cB^{t,\bw}) $ in Lemma \ref{lem_M31_01} (2)   shows that
 $\fp \big\{ \tau \= \wh{\tau}(\cB^{t,\bw}) \big\} \= 1$ for some  $[t,\infty]-$valued $ \bF^{W^t,P_0}-$stopping time $\wh{\tau}$ on $\O_0$,
 it follows that $ \oP_\Psi \big\{\oT \= \wh{\tau}(\oW) \big\} \= \fp \big\{\oT(\Psi) \= \wh{\tau}(\oW(\Psi))\big\}
 \= \fp \big\{\tau \= \wh{\tau}(\cB^{t,\bw})\big\} \= 1 $. 
 Since  $ \oW_{\n s} \big(\Psi (\o)\big) \= \cB^{t,\bw}_s(\o) \= \bw(s) $, $ \fa (s,\o) \ins [0,t] \ti \cQ $
 and since $\oU_s \big(\Psi(\o)\big) \= \wt{\mu}_s(\o) \= \bu(s)$, $ \fa (s,\o) \ins [0,t) \ti \cN^c $,
 it is clear that
 $ \oP_\Psi \big\{\oW_s \= \bw(s) , \fa s \ins [0,t] ; \oU_s \= \bu(s) \hb{ for a.e. } s \ins (0,t) \big\} \= \fp \big\{ \oW_s (\Psi) \= \bw(s) , \fa s \ins [0,t] ; \oU_s (\Psi) \= \bu(s) \hb{ for a.e. } s \ins (0,t) \big\}
 \= 1 $.    Thus  $\oP_\Psi \ins \ocP_{t,\bw,\bu,\bx}$.
  For any $i \ins \hN$,
  \bea
 E_{\oP_\Psi}   \bigg[ \int_t^\oT g_i \big(r,\oX_{r \land \cd},\oU_r \big) dr    \bigg]
  & \tn  \=  & \tn    E_\fp \bigg[    \int_t^{\oT (\Psi   )  } g_i \big(r,\oX_{r \land \cd} (\Psi   ), \oU_r(\Psi) \big) dr  \bigg]
  \=  E_\fp \Big[ \int_t^\tau g_i ( r, \cX^\th_{r \land \cd}, \wt{\mu}_r ) dr  \Big]  \nonumber  \\
  & \tn  \=  & \tn   E_\fp \Big[ \int_t^\tau g_i ( r, \cX^\th_{r \land \cd},  \mu_r ) dr  \Big] \ls y_i    \label{Feb12_01}
 \eea
 and   similarly
 $ E_{\oP_\Psi}   \big[ \int_t^\oT h_i \big(r,\oX_{r \land \cd},\oU_r \big) dr    \big]
   \=  E_\fp \big[ \int_t^\tau h_i ( r, \cX^\th_{r \land \cd},\mu_r ) dr  \big] \= z_i $,
 which means that  $ \oP_\Psi \ins \ocP_{t,\bw,\bu,\bx}(y,z)   $.
 Then an analogy to  \eqref{Feb12_01}   renders that
  $  E_\fp \big[ \int_t^\tau f ( r, \cX^\th_{r \land \cd},\mu_r ) dr \+ \b1_{\{\tau < \infty\}} \pi\big(\tau, \cX^\th_{\tau \land \cd} \big) \big]
 \=  E_{\oP_\Psi}  \Big[ \int_t^\oT f \big(r,\oX_{r \land \cd},\oU_r \big) dr \+ \b1_{\{\oT < \infty\}} \pi \big( \oT, \oX_{\oT \land \cd} \big) \Big] \ls \oV (t,\bw,\bu,\bx,y,z) $.
 Taking supremum over $ (\mu,\tau) \ins \cC_{t,\bx}(y,z)   $  yields   $V(t,\bx,y,z) \ls \oV (t,\bw,\bu,\bx,y,z)  $.

 \ss \no {\bf 2)} As $\ocP_{t,\bw,\bu,\bx}(y,z) \sb \ocP_{t,\bx}(y,z)$, we automatically have  $\oV (t, \bw,\bu,\bx, y,z ) \ls \oV (t,\bx,y,z )  $.
  It remains to demonstrate  that $\oV (t,\bx,y,z ) \ls V (t,\bx,y,z )   $.
  If $\ocP_{t,\bx}(y,z) \= \es$, then $\oV (t,\bx,y,z ) \= -\infty \ls V (t,\bx,y,z )   $.

   Assume   $\ocP_{t,\bx}(y,z) \nne \es$ and let $\oP \ins \ocP_{t,\bx}(y,z) $.
   By (D1) of Definition \ref{def_ocP}, there exists   a   $\hU-$valued, $ \bF^{W^t} -$predictable process
 $ \wh{\nu} \= \{\wh{\nu}_s\}_{s \in [t,\infty)} $   on $\O_0$ such that
 $ \oP \big\{  \oU_s    \= \ol{\nu}_s  \hb{ for a.e. }  s \ins (t,\infty) \big\} \= 1$,
 where $\ol{\nu}_s (\oo)  \df  \wh{\nu}_s   \big( \oW (\oo) \big)$, $\fa  (s,\oo)  \ins [t,\infty) \ti \oO $ is a $\hU-$valued, $ \bF^{\oW^t} -$predictable process  on $\oO$.

 Set $\vth \= (t,\bx,\ol{\nu})$.
   Given $ (\vf,n)  \ins  \fC(\hR^{d+l}) \ti \hN   $,
$ \ocM^\vth_s(\vf)    \df    \vf  ( \oW^t_{\n s} ,  \osX^\vth_{\n s} )
 \- \n \int_t^s \ol{b}  \big( r, \osX^\vth_{\n r \land \cd}, \ol{\nu}_r \big)  \n \cd \n D \vf  ( \oW^t_{\n r} ,  \osX^\vth_{\n r}  ) dr
    \-   \frac12 \int_t^s   \ol{\si} \, \ol{\si}^T  \big( r, \osX^\vth_{\n r \land \cd}, \ol{\nu}_r \big) \n : \n D^2 \vf ( \oW^t_{\n r} ,  \osX^\vth_{\n r} )   dr  $, $s \ins [t,\infty)$   is  an $\bF^{\oW^t,\oP}-$adapted continuous process
 and   $\otau^\vth_n  \df  \inf\big\{s \ins [t,\infty) \n :  \big|(\oW^t_s,\osX^\vth_s ) \big|   \gs n  \big\} \ld (t\+n) $
   is an $\bF^{\oW^t,\oP} -$stopping time. Since $\oW^t$ is a Brownian motion under $\oP$ by (D2) of Definition \ref{def_ocP},
 applying Proposition \ref{prop_MPF1}    with $(\O,\cF,P,B,X,\mu) \= \big(\oO,\sB(\oO),\oP,\oW,\osX^\vth,\ol{\nu} \big) $ shows that
 \bea \label{123021_11}
 \big\{\ocM^\vth_{s \land \otau^\vth_n} (\vf) \big\} \hb{ is a bounded $\bF^{\oW^t,\oP}-$martingale.}
 \eea

 Let $(\fu_o,\fx_o,\ft_o)$ be an arbitrary triplet  in $\hJ \ti \OmX \ti [t,\infty]$ and define  a mapping $\Psi_o \n : \cQ \mto \oO$ by
 $ \Psi_o (\o) \df \big( \cB (\o), \fu_o, \fx_o, \ft_o \big)  \ins \oO $, $ \fa \o \ins \cQ $.
 (Actually, we are indifferent to the second, third and fourth components  of $\Psi_o(\o)$.)
 Since   $   \oW^t_s  ( \Psi_o (\o) )    \=    \oW_s  (\Psi_o (\o)) \- \oW_t  (\Psi_o (\o)) \=   \cB^t_s (\o)    $ for any $(s,\o) \ins [t,\infty) \ti \cQ$,
    applying Lemma \ref{lem_122921_11} with $t_0 \= t$,  $(\O_1, \cF_1, P_1,B^1)   \= \big(\cQ, \cF, \fp , \cB \big) $, $(\O_2, \cF_2, P_2,B^2) \= \big(\oO , \sB(\oO ) , \oP  , \oW\big) $ and $\Phi \= \Psi_o$ yields that
   \bea
   \Psi^{-1}_o \big( \cF^{\oW^t }_s \big) \= \cF^{\cB^t }_s , ~
    \Psi^{-1}_o \big( \cF^{\oW^t,\oP}_s \big) \sb  \cF^{\cB^t,\fp}_s  , ~ \fa s \ins [t,\infty] \aand
     \big(\fp \nci \Psi_o^{-1}\big) (\oA) \= \oP(\oA)  , ~ \; \fa \oA \ins \cF^{\oW^t, \oP}_\infty   .  \label{Oct01_07b}
 \eea
 Then $ \sX^\vth_s (\o) \df  \osX^\vth_{\n s} ( \Psi_o(\oo))$, $s \ins [0,\infty)$ defines  an $\big\{\cF^{\cB^t,\fp}_{s \vee t}\big\}_{s \in [0,\infty)}-$adapted continuous process.

  Set $ \nu_s (\o) \df  \wh{\nu}_s  \big(\cB (\o)  \big)$, $\fa (s,\o) \ins [t,\infty) \ti \cQ$, which is a  $\hU-$valued, $ \bF^{\cB^t} -$predictable process  on $\cQ$.   Let $ (\vf,n)  \ins  \fC(\hR^{d+l}) \ti \hN   $. We define an $\bF^{\cB^t,\fp}-$adapted continuous process
$ \sM^\vth_s(\vf)    \df    \vf  ( \cB^t_s ,  \sX^\vth_{\n s} )
 \- \n \int_t^s \ol{b}  \big( r, \sX^\vth_{\n r \land \cd},\nu_r \big)  \n \cd \n D \vf  ( \cB^t_r ,  \sX^\vth_{\n r}  ) dr
    \-   \frac12 \int_t^s   \ol{\si} \, \ol{\si}^T  \big( r, \sX^\vth_{\n r \land \cd},\nu_r \big) \n : \n D^2 \vf ( \cB^t_r ,  \sX^\vth_{\n r} )   dr  $, $ \fa s \ins [t,\infty)$  and define an  $\bF^{\cB^t,\fp}-$stopping time $\z^\vth_n  \df  \inf\big\{s \ins [t,\infty) \n : \big|(\cB^t_s , \sX^\vth_s) \big|   \gs n   \big\} \ld (t\+n) $.
 Since $ \ol{\nu}_s (\Psi_o) \= \wh{\nu}_s   \big(\oW (\Psi_o)   \big) \= \wh{\nu}_s(\cB) \= \nu_s  $, $\fa  s  \ins [t,\infty)  $,
   applying Proposition \ref{prop_MPF1}    with $(\O,\cF,P,B,X,\mu) \= \big(\cQ,\cF,\fp,\cB,\sX^\vth,\nu \big) $ and using an analogy to \eqref{122621_17} yield that
   $ \big\{ \sM^\vth_{s \land \z^\vth_n} (\vf) \big\}_{s \in [t,\infty)}  $ is a bounded $\bF^{\cB^t,\fp}-$adapted continuous process under $ \fp  $ satisfying
   \if{0}
 It holds  for any $\o \ins \cQ$ that
  \beas
  \big(\ocM^\vth_s(\vf)\big) \big(\Psi_o(\o)\big) \= \big( \sM^\vth_s(\vf) \big) (\o) , \q \fa s \ins [t,\infty) \aand
  \otau^\vth_n \big(\Psi_o(\o)\big) \= \z^\vth_n(\o) .
  \eeas
   \fi
  \bea \label{123021_14}
 \big( \ocM^\vth_{s \land \otau^\vth_n} (\vf) \big) \big(\Psi_o(\o)\big) \= \big( \sM^\vth_{s \land \z^\vth_n} (\vf) \big) (\o) ,
 \q \fa (s,\o) \ins [t,\infty) \ti \cQ.
 \eea

 Let $t_1,t_2 \ins [t,\infty)$ with $t_1 \< t_2$ and let $A \ins \cF^{\cB^t }_{t_1}$.
 Since $\Psi_o^{-1}\big(\oA\big) \= A   $ for some  $\oA \ins \cF^{\oW^t }_{t_1}$ by \eqref{Oct01_07b},
 we can derive  from \eqref{123021_11}, \eqref{Oct01_07b} and \eqref{123021_14} that
 $ 0   \=    E_\oP  \Big[ \big(\ocM^\vth_{t_2 \land \otau^\vth_n }(\vf) \-  \ocM^\vth_{t_1 \land \otau^\vth_n }(\vf)  \big) \b1_\oA    \Big]
   \= E_\fp \Big[ \Big( \big(\ocM^\vth_{t_2 \land \otau^\vth_n }(\vf) \big) (\Psi_o) \- \big( \ocM^\vth_{t_1 \land \otau^\vth_n }(\vf)  \big)(\Psi_o) \Big) \b1_{\Psi_o^{-1}(\oA)}   \Big] \\
    \=    E_\fp  \Big[ \big(\sM^\vth_{t_2 \land \z^\vth_n }(\vf) \-  \sM^\vth_{t_1 \land \z^\vth_n }(\vf)  \big) \b1_A    \Big] $,
 which implies that $ \big\{\sM^\vth_{s \land \z^\vth_n }(\vf)\big\}_{s \in [t,\infty)}$ is a bounded $\bF^{\cB^t,\fp}-$martingale.
 Then an application of  Proposition \ref{prop_MPF1} with $(\O,\cF,P,B,X,\mu) \= \big(\cQ,\cF,\fp,\cB,\sX^\vth,\nu\big) $ shows that
    \if{0}

 Given $s \ins [t,\infty)$, we set   $\cF^\vth_s \df \cF^{\cB^t}_s \ve \cF^{\sX^\vth}_s \= \si\big(\cB^t_r; r \ins [t,s]\big) \ve \si (\sX^\vth_r; r \ins [0,s] )$.  Clearly, $\cF^\vth_s \sb \cF^{\cB^t,\fp}_s$.
 Actually, 
 $ \big\{ \sM^\vth_{s \land \z^\vth_n} (\vf) \big\}_{s \in [t,\infty)}  $ is a bounded $\{\cF^\vth_s\}_{s \in [t,\infty)}-$adapted continuous process under $ \fp  $ and is thus a bounded $\bF^{\cB^t,\fp}-$adapted continuous process under $ \fp  $.

  We see from this equality that the Lambda-system $\L \df \Big\{A \ins \cF^{\cB^t, \fp }_\infty  \sb \cQ \n :  E_\fp  \big[ \big(\sM^\vth_{t_2 \land \z^\vth_n }(\vf) \-  \sM^\vth_{t_1 \land \z^\vth_n }(\vf)  \big) \b1_A    \big] \= 0 \Big\} $ contains the Pi-system
 $ \cF^{\oW^t }_{t_1} \cp \sN_\fp \big(\cF^{\oW^t }_\infty \big)$. Then Dynkin's Pi-Lambda Theorem implies that
 $ \cF^{\oW^t }_{t_1,\fp} \= \si \Big( \cF^{\oW^t }_{t_1} \Cp \sN_\fp \big(\cF^{\oW^t }_\infty \big) \Big) \sb \L $, i.e.,
 $  E_\fp  \Big[ \big(\sM^\vth_{t_2 \land \z^\vth_n }(\vf) \-  \sM^\vth_{t_1 \land \z^\vth_n }(\vf)  \big) \b1_A    \Big] \= 0 $ for any $A \ins \cF^{\oW^t }_{t_1,\fp}$.
 Hence, $ \big\{\sM^\vth_{s \land \z^\vth_n }(\vf)\big\}_{s \in [t,\infty)}$ is a bounded $\bF^{\cB^t,\fp}-$martingale.
 In particular, $ \big\{\sM^\vth_{s \land \z^\vth_n }(\vf)\big\}_{s \in [t,\infty)}$ is a bounded $\{\cF^\vth_s\}_{s \in [t,\infty)}-$martingale.

    \fi
 \bea \label{012521_23}
 \fp\{ \sX^\vth_s \= \cX^{t,\bx,\nu}_s, ~ \fa s \ins [0,\infty) \} \= 1 .
 \eea

    By (D4) of Definition \ref{def_ocP},
 there exists   a $[t,\infty]-$valued $\bF^{W^t,P_0}-$stopping time $\wh{\ga}$  on $\O_0$
 such that $   \oP \big\{  \oT  \=    \wh{\ga} (\oW  )   \big\} \= 1 $.
 Lemma \ref{lem_M31_01} (1) renders that  $ \ga \df \wh{\ga}(\cB ) $ is an $\bF^{\cB^t,\fp}-$stopping time on $\cQ$
 while $   \wh{\ga}(\oW ) $ is an $\bF^{\oW^t,\oP}-$stopping time on $\oO$.
 For any $ i \ins \hN $,   we can deduce from (D3),  \eqref{Oct01_07b}  and  \eqref{012521_23}     that
  \bea
    y_i & \tn \gs  & \tn  E_\oP \bigg[ \int_t^\oT g_i \big( r,\oX_{r \land \cd},\oU_r \big) dr\bigg]
   \=  E_\oP \bigg[ \int_t^{ \wh{\ga}(\oW )} g_i \big( r,\osX^\vth_{\n r \land \cd},\wh{\nu}_r (\oW) \big) dr\bigg]
     \=  E_\fp \bigg[ \int_t^{ \wh{\ga}(\oW (\Psi_o))} g_i \big( r,\osX^\vth_{\n r \land \cd}(\Psi_o),\wh{\nu}_r (\oW(\Psi_o)) \big) dr\bigg] \nonumber \\
        & \tn  \=  & \tn  E_\fp \bigg[ \int_t^{ \wh{\ga}(\cB )} g_i ( r, \sX^\vth_{r \land \cd},\wh{\nu}_r(\cB)  ) dr\bigg]
         \=  E_\fp \bigg[ \int_t^{ \ga } g_i ( r, \cX^{t,\bx,\nu}_{r \land \cd},\nu_r  ) dr\bigg]  , \label{012521_15}
 \eea
 and similarly   $  E_\fp \big[ \int_t^{  \ga } h_i ( r, \cX^{t,\bx,\nu}_{r \land \cd},\nu_r  ) dr\big]
 \=  E_\oP \big[ \int_t^\oT h_i ( r,\oX_{r \land \cd},\oU_r ) dr\big] \= z_i$.
 So    $ (\nu,\ga)   \ins \cC_{t,\bx} (y,z)$. Analogous  to \eqref{012521_15},
 \beas
 \q   E_\oP   \bigg[ \int_t^\oT f \big(r,\oX_{r \land \cd},\oU_r \big) dr \+ \b1_{\{\oT < \infty\}}
 \pi \big( \oT,\oX_{\oT \land \cd} \big)  \bigg]
 \=  E_\fp \Big[ \int_t^{ \ga} f \big( r,\cX^{t,\bx,\nu}_{r \land \cd},\nu_r  \big) dr
 \+ \b1_{\{\ga < \infty\}} \pi \big(  \ga,\cX^{t,\bx,\nu}_{ \ga  \land \cd}  \big)  \Big] \ls V(t,\bx,y,z).
 \eeas
 Taking supremum over $ \oP \ins  \ocP_{t,\bx} (y,z) $   yields that   $ \oV(t,\bx,y,z) \ls V(t,\bx,y,z)$. \qed

   \no {\bf Proof of Lemma \ref{lem_082020_15}:}
 Lemma 4.1 of \cite{OSEC_stopping} shows that $\big(\fS, \Rho{\fS}\big)$ is a complete separable metric  space.
 so we only need to verify that $\big(\fU, \Rho{\fU}\big)$ is a complete separable metric  space.

 \ss \no {\bf 1)}
   Let $\mu^1,\mu^2 \ins \fU$ with $ \Rho{\fU} (\mu^1,\mu^2) \= 0 $.
   For any $ n,k \ins \hN$,
   the set $A^n_k \df    \big\{(s,\o_0) \ins [0,n] \ti \O_0 \n : \Rho{\hU} \big(\mu^1_s (\o_0),\mu^2_s (\o_0)\big)\ge 1/k \big\} \ins \sB[0,n] \oti \cF^{W,P_0}_n
  $ satisfies
  \beas
  \frac{e^{-n}}{k}  (ds \ti dP_0) (A^n_k)  \ls E_{P_0} \bigg[ \int_0^n \n \b1_{\{\Rho{\hU} (\mu^1_s,\mu^2_s)\ge 1/k\}} e^{-  s} \big( 1 \ld  \Rho{\hU} (\mu^1_s,\mu^2_s)  \big)   ds \bigg]
  \ls E_{P_0} \bigg[ \int_0^\infty \n  e^{-  s} \big( 1 \ld  \Rho{\hU} (\mu^1_s,\mu^2_s)  \big)   ds \bigg] \= 0 ,
  \eeas
  and thus $(ds \ti dP_0) (A^n_k) \= 0$.    As
  $ A \df  \{(s,\o_0) \ins [0,\infty) \ti \O_0 \n :  \mu^1_s (\o_0) \nne \mu^2_s (\o_0)\}
   \= \underset{n,k \in \hN}{\cup} A^n_k$,
   we see that $ (ds \ti dP_0) (A) \= 0 $, i.e., $\mu^1 \= \mu^2$ in $\fU$.
  Also, For $ \mu^1,\mu^2, \mu^3 \ins \fU$, 
\if{0}
 We can  discuss by three cases:

\ss \no i) When  $a\+b \ls 1$: as both $a$ and $b$ are $\ls 1$, one has   $ 1 \ld (a\+b) \= a \+ b \= 1 \ld a  \+ 1 \ld  b $;

\ss \no ii) When  $a+b \>1$ \& $a   \gs 1$, $ 1 \ld (a\+b) \= 1 \= 1 \ld a \ls 1 \ld a  \+ 1 \ld  b   $;

\ss \no iii) When $a+b \>1$ \& $b \gs 1$, $ 1 \ld (a\+b) \= 1 \= 1 \ld b \ls 1 \ld a  \+ 1 \ld  b   $;

\ss \no iv) $a+b \>1$ \& $a \ve b \< 1$,  $ 1 \ld (a\+b) \= 1  \< a\+b \= 1 \ld a \+ 1 \ld b   $.

 see also Lemma \ref{lem_Nov21_01},
\fi
   \eqref{060722_11} implies  that
 \if{0}
 \beas
 1 \ld  \Rho{\hU} (\mu^1_s,\mu^3_s)
 \ls 1 \ld \big(  \Rho{\hU} (\mu^1_s,\mu^2_s) \+ \Rho{\hU} (\mu^2_s,\mu^3_s) \big)
 \ls 1 \ld    \Rho{\hU} (\mu^1_s,\mu^2_s)    \+ 1 \ld    \Rho{\hU} (\mu^2_s,\mu^3_s)
 , \q \fa s \ins [0,\infty) .
 \eeas
 It follows that
 \fi
 $ \Rho{\fU} (\mu^1,\mu^3)   \=    E_{P_0} \Big[ \int_0^\infty \n  e^{-  s} \big( 1 \ld  \Rho{\hU} (\mu^1_s,\mu^3_s)  \big)   ds \Big]
 \ls  E_{P_0} \Big[ \int_0^\infty \n  e^{-  s} \big( 1 \ld  \Rho{\hU} (\mu^1_s,\mu^2_s)
 \+     1 \ld  \Rho{\hU} (\mu^2_s,\mu^3_s)  \big)   ds \Big]
    \=   \Rho{\fU} (\mu^1,\mu^2) \+ \Rho{\fU} (\mu^2,\mu^3)   $.
 So   $\Rho{\fU}$ is a metric on $\fU$.

  \ss \no {\bf 2)}
  Let $\{\mu^n\}_{n \in \hN}$ be a Cauchy sequence in $\big(\fU, \Rho{\fU}\big)$.
  We  can find a subsequence $\{n_k\}_{k \in \hN} $ of $\hN$ such that
  $  \Rho{\fU} \big(\mu^{n_k}, \mu^{n_{k+1}} \big) \< 2^{-k} $, $\fa k \ins \hN$.
  Fix $k  \ins \hN$.
  \if{0}
  Since it holds for any $ \ell \ins \hN$ that
  $  1 \ld \Rho{\hU} \big(\mu^{n_k}_s,\mu^{n_{k+\ell}}_s\big)
   \ls 1 \ld \Big( \sum^\ell_{i=1}   \Rho{\hU} \big(\mu^{n_{k+i-1}}_s,\mu^{n_{k+i}}_s\big) \Big)
   \ls \sum^\ell_{i=1}   1 \ld  \Rho{\hU} \big(\mu^{n_{k+i-1}}_s,\mu^{n_{k+i}}_s\big)
   \ls \sum_{i \in \hN}   1 \ld  \Rho{\hU} \big(\mu^{n_{k+i-1}}_s,\mu^{n_{k+i}}_s\big)$,
   $ s \ins [0,\infty)  $,
  taking $ \underset{\ell \in \hN}{\sup} $ on both sides yields that
  \fi
  Since an analogy to \eqref{050921_23} shows that for any $\ell \ins \hN$,
  $  1 \ld \Big( \underset{\ell \in \hN}{\sup} \,   \Rho{\hU} (\mu^{n_k}_s,\mu^{n_{k+\ell}}_s) \Big)
    \ls \sum_{i \in \hN}   1 \ld  \Rho{\hU} (\mu^{n_{k+i-1}}_s,\mu^{n_{k+i}}_s)   $, $ \fa  s \ins [0,\infty) $,
  the monotone convergence theorem   implies
  $  E_{P_0} \Big[ \int_0^\infty \n e^{-s} \Big( 1 \ld \Big( \underset{\ell \in \hN}{\sup} \,   \Rho{\hU} (\mu^{n_k}_s,\mu^{n_{k+\ell}}_s) \Big) \Big)  ds \Big]
   \ls 
      \sum_{i \in \hN} \Rho{\fU} \big(\mu^{n_{k+i-1}}, \mu^{n_{k+i}}  \big)
  \ls \sum_{i \in \hN} 2^{-k-i+1} \= 2^{-k+1} $.
  So $\big\{  n_k  \big\}_{k \in \hN}$ has a subsequence $\big\{  n_{k_m}  \big\}_{m \in \hN}$ such that
 $ \lmt{m \to \infty}    \Big( \underset{\ell \in \hN}{\sup} \,   \Rho{\hU} \\ \big(\mu^{n_{k_m}}(s,\o_0),\mu^{n_{k_m+\ell}}(s,\o_0)\big) \Big)   \= 0 $
   for all $(s,\o_0) \ins [0,\infty) \ti \O_0 $ except on a $ds \ti d P_0-$null set $\cN$.
   Given $(s,\o_0) \ins \cN^c $, one has
   $   \lmt{m \to \infty}    \bigg( \underset{j \in \hN}{\sup} \,   \Rho{\hU}   \big(\mu^{n_{k_m}}(s,\o_0) ,  \mu^{n_{k_{j+m}}}(s,\o_0) \big)   \bigg) \= 0 $
   or $ \big\{  \mu^{n_{k_m}}(s,\o_0)  \big\}_{m \in \hN} $ is a Cauchy sequence   in $\hU$.
     Let $ \mu^o  (s,\o_0) $ be the limit of $ \big\{  \mu^{n_{k_m}}(s,\o_0)  \big\}_{m \in \hN} $ in $\big(\hU, \Rho{\hU}\big)$.

  Let $\sP^{W,P_0}$  denote the predictable sigma-field on $[0,\infty) \ti \O_0$ with respect to   filtration $\bF^{W,P_0}$ \big(i.e., the sigma-field of $[0,\infty) \ti \O_0$ generated by all $\big\{\cF^{W,P_0}_{s-}\big\}_{s \in [0,\infty)}-$adapted, c\`agl\`ad processes\big).
   Define  $ \ul{\mu}  (s,\o_0) \df \linf{m \to \infty} \sI \big( \mu^{n_{k_m}} (s,\o_0) \big) \ins [0,1] $, $(s,\o_0) \ins [0,\infty) \ti \O_0 $,
   which is   a   $\sP^{W,P_0}-$measurable function  on $[0,\infty) \ti \O_0$. Then
 \beas
  \mu^* (s,\o_0)\df
 \sI^{-1} \big(\ul{\mu}(s,\o_0) \big) \b1_{ \{\ul{\mu}(s,\o_0) \in \fE\}} \+ u_0 \b1_{ \{\ul{\mu}(s,\o_0) \notin \fE\}}    \ins \hU  ,
 \q \fa (s,\o_0) \ins [0,\infty) \ti \O_0
  \eeas
  is  also   $\sP^{W,P_0}-$measurable or  the process $\{\mu^*_s\}_{s \in \hN}$ is a $\hU-$valued $\bF^{W,P_0}-$predictable
 process on $\O_0$. To wit, $\mu^* \ins \fU$.

 For $(s,\o_0) \ins \cN^c$,  the continuity of   mapping $ \sI $ shows
  $ \ul{\mu} (s,\o_0) \= \lmt{m \to \infty} \sI \big(\mu^{n_{k_m}}(s,\o_0)\big) \= \sI   \big( \mu^o  (s,\o_0) \big) \ins \fE $,
 so   $\mu^* (s,\o_0) \= \sI^{-1} \big(\ul{\mu}(s,\o_0) \big) \=   \mu^o(s,\o_0)  $.
 The dominated convergence theorem yields that $ \lmt{m \to \infty} \Rho{\fU} \big( \mu^{n_{k_m}} , \mu^*\big) \= \lmt{m \to \infty}  \int_{\o_o \in \O_0}   \int_0^\infty \n  e^{-  s} \b1_{\{(s,\o_0) \in \cN^c\}} \Big( 1 \ld  \Rho{\hU} \big(  \mu^{n_{k_m}}(s,\o_0) ,  \mu^o  (s,\o_0)   \big)  \Big)   ds P_0 (d \o_0) \= 0 $.
  Similar to \eqref{050921_25}, one can further deduce that
 $ \lmt{n \to \infty} \Rho{\fU} \big( \mu^n , \mu^*\big) \= 0 $.
   \if{0}

 Let $\e \ins (0,1)$. There exists a $N \in \hN$ such that $ \Rho{\fU} \big( \mu^{\fn_1} , \mu^{\fn_2}\big) \< \e/2 $ for any $\fn_1,\fn_2 \gs N$. We can also find a $m \ins \hN$ such that $n_{k_m} \gs N $ and that
 $ \Rho{\fU} \big( \mu^{n_{k_m}} , \mu^*\big) \< \e/2 $.
 Then it holds  for any $n \ins \hN$ with $n \gs   N  $ that
 \beas
 \Rho{\fU} \big( \mu^n , \mu^*\big)
 \ls \Rho{\fU} \big( \mu^n , \mu^{n_{k_m}} \big) \+
 \Rho{\fU} \big( \mu^{n_{k_m}} , \mu^*\big) \< \e .
 \eeas

   \fi
  Hence,   $\mu^*$ is a limit of $\{\mu^n\}_{n \in \hN}$ in $\big(\fU, \Rho{\fU}\big)$.

 \ss \no {\bf 3)}
 Let $ \cL^0 \df L^0 \big( [0,\infty) \ti \O_0 ;  \hU \big) $ be
 the equivalence classes of all $\hU-$valued,  $ \sB ([0,\infty) \ti \O_0) -$measurable
 functions $ \fX  $ on $[0,\infty) \ti \O_0$   in the sense that
 $\fX^1,\fX^2 \ins \cL^0   $ are   equivalent   if
 $\{ (s,\o_0) \ins [0,\infty) \ti \O_0 \n : \fX^1(s,\o_0) \nne \fX^2(s,\o_0)  \} $ is a $ds \ti dP_0-$null set.
 We endow $ \cL^0 $   with a metric
\beas
   \Rho{\cL^0} ( \fX^1,\fX^2 ) \df \int_{ \o_0 \in \O_0 }    \int_0^\infty \n  e^{-  s} \big( 1 \ld \Rho{\hU} \big( \fX^1(s,\o_0),\fX^2(s,\o_0) \big) \big)    ds  P_0 (d\o_0)  , \q \fa \fX^1,\fX^2 \ins \cL^0 ,
 \eeas
  and we shall  demonstrate that $\big(\cL^0,\Rho{\cL^0}\big) $ is a separable  space.

 Let $\{u_i\}_{i \in \hN}$ be  a countable dense subset of the complete metric space $\big(\hU, \Rho{\hU}\big)$.
      Since both $[0,\infty)$ and $\O_0$ are separable metric spaces (and thus second-countable spaces),
     $  [0,\infty) \ti \O_0$ is also second-countable or it  has a countable base   $\{O_i\}_{i \in \hN}$.
      Given $n \in \hN$,  let us enumerate the $2^n$ elements of $   \big\{ \underset{i \in  I}{\cup} O_i \n :  I \sb \{1,\cds,n\} \big\}$   by $ \big\{ \breve{O}^n_1, \cds , \breve{O}^n_{2^n} \big\}$ and consider
     the following countable collections of   $\sB ([0,\infty) \ti \O_0) / \sB [-1,1]-$measurable functions on $[0,\infty) \ti \O_0$:
     \beas
     \sC_n \df \Big\{    \sI   ( u_j) \b1_{\breve{O}^n_i} \- \b1_{(\breve{O}^n_i)^c}   \n : i \in \{1,\cds,2^n\},~ j  \ins \hN   \Big\} , \q
     \wt{\sC}_n \df \Big\{ \underset{i = 1,\cds, k}{ \max }   \fX^i    \n :  k \ins \hN, ~ \{ \fX^1,   \cds, \fX^k \} \sb \sC_n \Big\} .
     \eeas
     Clearly, each   $\fX \ins \wt{\sC}_n$ takes   values  in a finite subset of $    \big\{\sI (u_i) \big\}_{i\in \hN}   \cup \{-1\} $.
     By additionally assigning $\sI^{-1} (-1) \df u_0 $, one has
     $ \wt{\bC} \df \big\{ \sI^{-1} (\fX) \n : \fX \ins \wt{\sC}_n \hb{ for some }  n \ins \hN  \big\}
     \sb \big\{\hb{$\{u_i\}^\infty_{i=0}-$valued  $\sB ([0,\infty) \ti \O_0)-$measurable   functions}\big\}
     \n =:\n \wh{\bC} $.

      \ss \no {\bf 3a)} We claim that
       $ \wh{\bC}  $ is a dense subset of $ \cL^0 $ under $\Rho{\cL^0}$:
 For any $i,n \ins \hN$, we set $o^n_i$ and $\wt{o}^n_i$ as in the proof of Lemma \ref{lem_Nov25_03}.
 Given $\fX \ins \cL^0$ and $n \ins \hN$,
   define a member of $\wh{\bC}$ by  $ \fX^n (s,\o_0) \df   \sum_{i \in \hN} \b1_{\{ (s,\o_0)  \in  A^n_i \}}   u_i $,
  $ (s,\o_0)  \ins [0,\infty) \ti \O_0$, where
    $A^n_i \df \big\{ (s,\o_0)  \ins [0,\infty) \ti \O_0 \n : \fX (s,\o_0) \ins \wt{o}^n_i \big\} \ins \sB ([0,\infty) \ti \O_0)$. 
   Since
 $  \Rho{\hU} \big(\fX^n (s,\o_0), \fX (s,\o_0) \big)
   \=   \sum_{i \in \hN} \b1_{\{ (s,\o_0)  \in  A^n_i \}}  \Rho{\hU} \big( \fX (s,\o_0) , u_i \big)
   \< 2^{-n} $, $ \fa (s,\o_0) \ins [0,\infty) \ti \O_0 $,
  one can deduce that  $ \Rho{\cL^0}   \big(\fX^n  , \fX   \big) \=   \int_{\o_0 \in \O_0} \int_0^\infty \n  e^{-  s} \big( 1 \ld \Rho{\hU} \big( \fX^n(s,\o_0),\fX(s,\o_0) \big) \big)    ds P_0(d\o_0)  \ls 
  2^{-n} $.
 So $ \wh{\bC}  $ is a dense subset of $ \cL^0 $ under $\Rho{\cL^0}$.

  \ss \no {\bf 3b)} We next show that   the countable set $\wt{\bC}$ is dense in $\wh{\bC}$ 
  and is thus dense in $\cL^0$ under $ \Rho{\cL^0} $:
  Let $\wh{\fX} \ins \wh{\bC} $ and  $\e \ins (0,1)$.
  We set    $\wh{A}_i \= \big\{(s,\o_0) \ins [0,\infty) \ti \O_0 \n : \wh{\fX} (s,\o_0) \=  u_i  \big\} \ins \sB ([0,\infty) \ti \O_0)$, $\fa i \ins \hN$.
  Since $ P(A) \df \int_{\O_0} \int_0^\infty \n  e^{-  s}  \b1_{\{(s,\o_0) \in A\}}   ds P_0(d \o_0)$
  is a probability measure  on  $\big([0,\infty) \ti \O_0, \sB ([0,\infty) \ti \O_0)\big)$, there exists  $N \ins \hN$ such that
  $  P \Big( \underset{i > N}{\cup} \wh{A}_i \Big) < \e/3 $.

  Given $i \=  1,\cds, N$,     Proposition 7.17 of \cite{Bertsekas_Shreve_1978} shows that
  $P( \cO_i \backslash \wh{A}_i) \< \frac{\e}{3N} $ for  some open subset $ \cO_i$ of $[0,\infty) \ti \O_0$ containing $\wh{A}_i$.
  Since $  \cO_i \= \underset{n \in \hN}{\cup} O_{\ell^i_n} $ for some subsequence $\big\{\ell^i_n\big\}_{n \in \hN}$
  of $\hN$, we can find $ n_i \ins \hN$ such that
  $P \big( \cO_i  \backslash \breve{\cO}_i \big) <  \frac{\e}{3N} $,
  where $\breve{\cO}_i \df \underset{n  = 1}{\overset{ n_i}{\cup}} O_{\ell^i_n} \ins \sB ([0,\infty) \ti \O_0) $.
  As $   \breve{\cO}_i    \ins \big\{ \breve{O}^\fn_j   \big\}^{2^\fn}_{j=1}$ for $ \fn \df \underset{i  = 1,\cds,N}{\max} \ell^i_{n_i}$,
  we see that $ \fX^i \df \sI   ( u_i ) \b1_{\breve{\cO}_i}
 \-    \b1_{\breve{\cO}^c_i}  $ belongs to $\sC_\fn$.
 Define     $\fX \df \sI^{-1} \Big( \underset{i = 1,\cds, N}{\max} \,  \fX^i \Big) \ins \wt{\bC} $.

 For $i \= 1,\cds, N $,
 if $  \cA_i \df \big( \wh{A}_i \Cp \breve{\cO}_i   \big) \Big\backslash \Big( \underset{j \le N; j \ne i}{\cup}\breve{\cO}_j \Big) \ins \sB ([0,\infty) \ti \O_0)$ is not empty,
  it holds for any $(s,\o_0) \ins  \cA_i $ that
   $ \fX^j  (s,\o_0)   \= \b1_{\{j=i\}} \sI ( u_i) \- \b1_{\{j \ne i\}} $ for $j \ins \{1,\cds,N\}$
   and thus $ \fX  (s,\o_0) \= \sI^{-1} \big(  \fX^i (s,\o_0) \big) \=   u_i  \= \wh{\fX} (s,\o_0)$.
Also, if $ \cA_0 \df \Big( \underset{i \in \hN}{\cup} \wh{A}_i \Big)^c \bigcap \Big( \underset{j = 1}{\overset{N}{\cup}} \breve{\cO}_j \Big)^c $
is not empty, it holds for any $(s,\o_0) \ins  \cA_0$ that $ \fX (s,\o_0) \= \sI^{-1}(-1) \= u_0   \= \wh{\fX} (s,\o_0) $.
  Then
  $  \Rho{\cL^0} \big(  \fX  , \wh{\fX} \big)
  \= \int_{\o_0 \in \O_0}  \int_0^\infty \n  \b1_{ \{(s,\o_0) \in \cA    \}}  e^{-  s} \big( 1 \ld \Rho{\hU} \big( \fX(s,\o_0),\wh{\fX}(s,\o_0) \big) \big)    ds P_0(d \o_0) $
  for  $\cA \df \Big( \underset{i=0}{\overset{N}{\cup}} \cA_i \Big)^c \= \Big(\underset{i=1}{\overset{N}{\cup}} (\wh{A}_i \backslash \cA_i)\Big) \bigcup  \Big( \underset{i > N}{\cup} \wh{A}_i\Big) \bigcup \bigg( \Big( \underset{i \in \hN}{\cup} \wh{A}_i \Big)^c  \bigcap \Big( \underset{i = 1}{\overset{N}{\cup}} \breve{\cO}_i \Big)\bigg)  $.
  Since
  \beas
    \wh{A}_i \backslash \cA_i \= \big(\wh{A}_i \Cp \breve{\cO}^c_i \big) \hb{$\bigcup$}
  \Big(   \underset{j \le N; j \ne i}{\cup} \big(\wh{A}_i  \Cp \breve{\cO}_i  \Cp \breve{\cO}_j \big)  \Big)
  \sb \big( \cO_i \backslash   \breve{\cO}_i \big) \hb{$\bigcup$} \Big(   \underset{j \ne i}{\cup}   \breve{\cO}_j \backslash \wh{A}_j    \Big),
  \q  i \= 1,\cds, N
  \eeas
  and since $\Big( \underset{i \in \hN}{\cup} \wh{A}_i \Big)^c \bigcap \Big( \underset{j = 1}{\overset{N}{\cup}} \breve{\cO}_j \Big)
  \= \underset{j = 1}{\overset{N}{\cup}} \Big( \breve{\cO}_j \Cp \Big( \underset{i \in \hN}{\cup} \wh{A}_i \Big)^c \Big)
  \sb \underset{j = 1}{\overset{N}{\cup}} \big( \breve{\cO}_j \Cp   \wh{A}^c_j   \big)
  \sb \underset{j = 1}{\overset{N}{\cup}} \big(  \cO_j \backslash  \wh{A}_j   \big) $,
  we can deduce that
  $ \cA 
  \sb   \Big(   \underset{i  = 1}{\overset{N}{\cup}}  \cO_i \backslash   \breve{\cO}_i   \Big) \bigcup \Big(   \underset{i  = 1}{\overset{N}{\cup}}   \cO_i \backslash \wh{A}_i    \Big) \bigcup  \Big( \underset{i > N}{\cup} \wh{A}_i\Big)  $.
  It follows that    $  \Rho{\cL^0} \big(  \fX  , \wh{\fX} \big)   \ls P  ( \cA  )  \< \e $.
  Namely, the countable collection $ \wt{\bC} $ is   dense in $ \wh{\bC} $    under $ \Rho{\cL^0} $.

   Viewed as a subset of $\cL^0$, $\fU$ is also  separable under $ \Rho{\fU} $.
   \if{0}

  Let   $\mu^0$ be an arbitrary element of $\fU$.
  For any $\fX \ins \wt{\bC}$ and $ n \ins \hN$, if $O^\fU_{1/n}(\fX) \df \big\{\mu \ins \fU \n : \Rho{\cL^0} (\fX,\mu) \< 1/n \big\} \nne  \es$,
  we pick   a $ \mu^\fX_n $ from it; otherwise, we set $ \mu^\fX_n \df \mu^0 $.
  Let $\mu \ins \fU$ and $k \ins \hN$, since $ \wt{\bC} $ is   dense   in $\cL^0$ under $ \Rho{\cL^0} $,
  there exists $\fX \ins \wt{\bC}$ such that $\Rho{\cL^0} (\fX,\mu)
    \< \frac{1}{2k}$. As $\mu \ins O^\fU_{1/2k}(\fX)$,
     $ \mu^\fX_{2k} $ is also in  $O^\fU_{1/2k}(\fX)$.
  It  follows that $ \Rho{\fU} (\mu,\mu^\fX_{2k})
   \= \Rho{\cL^0} (\mu,\mu^\fX_{2k}) \ls \Rho{\cL^0}(\mu,\fX )
   \+ \Rho{\cL^0}(\fX ,\mu^\fX_{2k}) \< \frac{1}{2k} \+ \frac{1}{2k} \= \frac{1}{k}$.
  Hence, $\big\{\mu^\fX_n\big\}_{\fX \in \wt{\sC}, n \in \hN}$ forms a dense subset of $\fU$ under $\Rho{\fU}$
  (Note the difference between the equivalent class of $\mu$ in $\cL^0$ and the equivalent class of $\mu$ in $\fU$).

  \fi
  Therefore, $\big(\fU, \Rho{\fU}\big)$ is a complete separable space.   \qed

    \no {\bf Proof of Lemma \ref{lem_082020_17}: 1)}
     \if{0}

    Let $(\mu,\tau) \ins \fU \ti \fS$, we first show that
    the mapping $(W, \mu, \tau) \n : \O_0 \mto \O_0 \ti \hJ \ti \hT$ is
     $\cF^{W,P_0}_\infty \big/ \sB \big(\O_0 \ti \hJ \ti \hT\big)-$measurable.
     As an identical mapping from $\O_0$ to $\O_0$, $W$ is clearly $\sB(\O_0)/\sB(\O_0)-$measurable.
     So it holds for any $A_0 \ins \sB(\O_0)$ that  $W^{-1}(A_0) \ins \sB(\O_0) \=   \cF^W_\infty $.

   We   set $\mu^{-1} (A) \df \{\o_0 \ins \O_0 \n : \mu_\cd (\o_0) \ins A\}$ for any $A \sb   \hJ$.
   Recall from \eqref{J04_01} that the topology $\fT_\sharp(\hJ)$ of $\hJ$ has a countable subbase
  $\big\{ \fri^{-1}_\hJ \big(O_{\frac{1}{n}} (\fm_k,\phi_j)\big) \big\}_{n,k,j \in \hN}$.
  Given $n,k,j \ins \hN$, the continuity of function $\phi_j$ and the $ \bF^{W,P_0} -$predictability of process $\mu$
  imply that the process $ \mu^j_s(\o_0) \df \phi_j \big(s, \mu_s (\o_0)\big) $, $(s,\o_0) \ins [0,\infty) \ti \O_0$
  is also $ \bF^{W,P_0} -$predictable.
  Then the random variable $\xi_j(\o_0) \df \int_0^\infty e^{-s}   \phi_j (s, \mu_s (\o_0))   ds$, $\o_0  \ins  \O_0$ is $ \cF^{W,P_0}_\infty -$measurable and thus
 \beas
\hspace{-0.5cm} \mu^{-1} \big( \fri^{-1}_\hJ \big(O_{\frac{1}{n}} (\fm_k,\phi_j)\big) \big)
 & \tn  \=  & \tn  \bigg\{\o_0 \ins \O_0 \n :  \Big| \int_0^\infty e^{-s}   \phi_j (s, \mu_s (\o_0)) ds
 \- \int_0^\infty \n \int_\hU
  \phi_j (s,u) \fm_k (ds,du)      \Big| \< 1/n  \bigg\} \\
 & \tn \= & \tn  \big\{\o_0 \ins \O_0 \n :    \xi_j(\o_0) \ins  ( c_{k,j} \- 1/n,  c_{k,j} \+ 1/n ) \big\} \ins \cF^{W,P_0}_\infty ,
 \eeas
 where $ c_{k,j} \df  \int_0^\infty \n \int_\hU
 \phi_j (s,u) \fm_k (ds,du) \ins \hR   $ is a constant independent of $\o_0 \ins \O_0$.
 It follows that the sigma-field
 $ \L \df  \{A \sb \hJ \n : \mu^{-1} (A) \ins \cF^{W,P_0}_\infty  \} $
 of $\hJ$ contains all open sets of $\fT_\sharp(\hJ)$ and thus includes $\sB(\hJ)$.
 To wit,   $ \mu^{-1} (A) \ins \cF^{W,P_0}_\infty $ for any $A \ins \sB(\hJ)$.

 On the other hand,  since the  $\bF^{W,P_0}-$stopping time  $\tau$ satisfies that   $ \big\{\o_0 \ins \O_0 \n : \arctan(\tau(\o_0)) \ins [0,\fr] \big\} \=   \big\{\o_0 \ins \O_0 \n : \tau(\o_0) \ins [0,\tan(\fr)] \big\} \ins \cF^{W,P_0}_{\tan(\fr)} \sb \cF^{W,P_0}_\infty$ for any $\fr \ins [0,\pi/2]$,
      the sigma field of $[0,\pi/2]$,
   $\L_\tau \df \big\{\cE \sb [0,\pi/2] \n :  (\arctan(\tau))^{-1}(\cE) \ins \cF^{W,P_0}_\infty \big\}$ contains all closed intervals $[0,\fr]$, $\fr \ins [0,\pi/2]$, which generates $\sB[0,\pi/2]$. It follows that
   $\sB[0,\pi/2] \sb \L_\tau $. So $ (\arctan(\tau))^{-1}(\cE) \ins \cF^{W,P_0}_\infty $, $ \fa \cE \ins \sB[0,\pi/2] $
   or $  \tau^{-1}(\cE) \ins \cF^{W,P_0}_\infty $, $ \fa \cE \ins \sB(\hT) \= \{\tan(\cE') \n : \cE' \ins \sB[0,\pi/2]\} $.

 As $ (W,\mu,\tau)^{-1} (A_0 \ti A \ti \cE)
  \= W^{-1}(A_0) \ti \mu^{-1}(A) \ti \tau^{-1}(\cE) \ins \cF^{W,P_0}_\infty$
  for any $A_0 \ins \sB(\O_0)$, $A \ins \sB(\hJ)$ and $\cE \ins \sB(\hT)$,
 we see that
  the sigma-field
 $ \ol{\sS} \df  \big\{ \oD \sb \O_0 \ti \hJ \ti \hT \n : (W, \mu, \tau)^{-1} (\oD) \ins \cF^{W,P_0}_\infty  \big\} $
 of $\O_0 \ti \hJ \ti \hT $ contains all measurable rectangles of  the product sigma-field
 $\sB(\O_0) \oti \sB(\hJ) \oti \sB(\hT)$  and  thus includes $\sB(\O_0) \oti \sB(\hJ) \oti \sB(\hT) \= \sB (\O_0 \ti \hJ \ti \hT)$
 by  Lemma \ref{lem_prod_Bsf}.  
  Hence, $ (W, \mu, \tau)^{-1} (\oD) \ins \cF^{W,P_0}_\infty $ for any $\oD \ins \sB(\O_0 \ti \hJ \ti \hT)$.

  \fi
  We first  show that   $\Ga$ is injective:
  Set $\hQ_\pi \df \big(\hQ  \Cp [0,\pi/2) \big) \cp \{\pi/2\}$
 and let $(\mu^1,\tau_1),(\mu^2,\tau_2) \ins \fU \times  \fS$ such that $\Ga(\mu^1,\tau_1) \= \Ga(\mu^2,\tau_2)$.

 Let $q   \ins \hQ_\pi   $ and $n \ins \hN$. We set $\cE^q_n \df (q\-1/n,q\+1/n) \Cp [0,\pi/2]$.
 Also let $ k,j  \ins \hN$ and $i \= 1,2$. We know from \eqref{J04_01} and Lemma \ref{lem_061122} that  $ A^i_{n,k,j,q} \df  (\mu^i )^{-1} \big(\fri^{-1}_\hJ \big(O_{\frac{1}{n}} (\fm_k,\phi_j)\big)\big) \Cp \big\{ \arctan(\tau_i) \ins \cE^q_n \big\} $
  belongs to $ \cF^{W,P_0}_\infty$.
    Then   $A_{n,k,j,q} \df  A^1_{n,k,j,q} \Cp (A^2_{n,k,j,q})^c$  satisfies that
 \bea
  P_0 (A_{n,k,j,q} ) & \tn \= & \tn    P_0 \big\{ \o_0 \ins (A^2_{n,k,j,q})^c :  \mu^1_\cd (\o_0) \ins \fri^{-1}_\hJ \big(O_{\frac{1}{n}} (\fm_k,\phi_j)\big) , \tau_1 (\o_0)  \ins \tan(\cE^q_n) \big\} \nonumber  \\
   & \tn \= & \tn  \big( \Ga (\mu^1,\tau_1) \big) \big( (A^2_{n,k,j,q})^c \ti \fri^{-1}_\hJ \big(O_{\frac{1}{n}} (\fm_k,\phi_j)\big) \ti \tan(\cE^q_n) \big)
      \=   \big( \Ga (\mu^2,\tau_2) \big) \big( (A^2_{n,k,j,q})^c \ti \fri^{-1}_\hJ \big(O_{\frac{1}{n}} (\fm_k,\phi_j)\big) \ti \tan(\cE^q_n) \big) \nonumber \\
    & \tn  \=  & \tn    P_0 \big\{ \o_0 \ins (A^2_{n,k,j,q})^c :  \mu^2_\cd (\o_0) \ins \fri^{-1}_\hJ \big(O_{\frac{1}{n}} (\fm_k,\phi_j)\big) , \tau_2 (\o_0) \ins \tan(\cE^q_n) \big\} 
    \=  P_0 \big(   (A^2_{n,k,j,q})^c \Cp  A^2_{n,k,j,q} \big)
    \= 0   .   \label{031821_a01}
 \eea

  We claim that
  \bea \label{082020_15}
   \hb{$\cA \df \big\{\o_0 \ins \O_0 \n : \mu^1_\cd (\o_0) \nne \mu^2_\cd (\o_0) \big\} \cp \big\{\o_0 \ins \O_0 \n :  \tau_1 (\o_0) \nne \tau_2 (\o_0) \big\} $ is equal to
 $ \underset{n,k,j  \in \hN}{\cup} \, \ccup{q \in \hQ_\pi}{} A_{n,k,j,q}    $.} \q
 \eea
 Clearly, $ \underset{n,k,j  \in \hN}{\cup} \, \ccup{q \in \hQ_\pi}{} A_{n,k,j,q} \sb \cA $.
 Assume that $ \cA \Cp \Big(  \underset{n,k,j  \in \hN}{\cup} \, \ccup{q \in \hQ_\pi}{}  A_{n,k,j,q} \Big)^c $
 is not empty and has an element $\o_0$.

  Given $n,j  \ins \hN$, since the proof of Lemma \ref{lem_082020_11} selected  $\{\fm_k\}_{k \in \hN}$ as  a   countable dense subset
  of   the topological space $\big(\fP  ([0,\infty) \ti \hU ), \fT_\sharp \big(\fP   ([0,\infty) \ti \hU )\big)\big) $,
  there exist  $\fk \= \fk (n,j ) \ins \hN$ and $\fq \= \fq(n) \ins \hQ_\pi$ such that
  $\friJ \big(\mu^1_\cd (\o_0)\big) \ins
   O_{\frac{1}{n}} (\fm_\fk,\phi_j ) $
  and $\arctan\big(\tau_1(\o_0)\big) \ins \cE^\fq_n$.
  This shows $ \o_0 \ins (\mu^1 )^{-1} \big( \fri^{-1}_\hJ \big(O_{\frac{1}{n}} (\fm_\fk,\phi_j ) \big) \big) \Cp \big\{\arctan(\tau_1) \ins  \cE^\fq_n \big\}
  \= A^1_{n,\fk,j,\fq}  $.
  Since $ \o_0 \ins  A^c_{n,\fk,j,\fq}  $, we see that $ \o_0 \ins A^2_{n,\fk,j,\fq} $, i.e.,
  $ \big( \mu^2_\cd (\o_0), \arctan(\tau_2 (\o_0)) \big)  $ also belongs to $   \fri^{-1}_\hJ \big(O_{\frac{1}{n}} (\fm_\fk,\phi_j ) \big)  \ti \cE^\fq_n  $. It follows that
  $\big| \int_0^\infty e^{-s} \big[ \phi_j \big(s,\mu^1_s(\o_0)\big) \- \phi_j \big(s,\mu^2_s(\o_0)\big)   \big] ds  \big|
   \< 2/n $ and $\Rho{+} \big(\tau_1(\o_0), \tau_2(\o_0) \big) \= \big|\arctan(\tau_1(\o_0)) \- \arctan(\tau_2(\o_0)) \big| \< 2/n$.
  Letting $n \nto \infty$ yields that $ \int_0^\infty e^{-s}   \phi_j \big(s,\mu^1_s(\o_0)\big) ds \= \int_0^\infty e^{-s} \phi_j \big(s,\mu^2_s(\o_0)\big)    ds  $ and $ \tau_1(\o_0) \= \tau_2(\o_0) $.

  As $\{\phi_j\}_{j \in \hN}$ is   dense in $ \wh{C}_b\big([0,\infty) \ti \hU\big) $  by Proposition 7.20 of \cite{Bertsekas_Shreve_1978},
  the dominated convergence theorem implies that
  $ 
   \int_0^\infty e^{-s} \,  \phi \big(s,\mu^1_s(\o_0)\big) ds \= \int_0^\infty e^{-s} \, \phi \big(s,\mu^2_s(\o_0)\big)    ds $
    holds for any $\phi \ins \wh{C}_b\big([0,\infty) \ti \hU\big)$.
    By a standard approximation, 
  this equality also  holds for any bounded Borel-measurable functions $\phi$ on $ [0,\infty) \ti \hU $.
     \if{0}

    Using the convolution with smooth mollifier, 
   \eqref{051421_23} also holds for any $\phi \ins C_b\big([0,\infty) \ti \hU\big)$.
   Then   Urysohn's lemma and the dominated convergence theorem yield  that  \eqref{051421_23}   holds for any
  step $\sB\big([0,\infty) \ti \hU\big)-$measurable functions $\phi$   and thus for any bounded $\sB\big([0,\infty) \ti \hU\big)-$measurable functions $\phi$ on $ [0,\infty) \ti \hU $.
         \if{0}

    First for $ \phi \ins \wh{C}_b\big([0,\infty) \ti \hU\big) $, then for $ \phi \ins C_b\big([0,\infty) \ti \hU\big)$, next for bounded Borel-measurable functions on $ [0,\infty) \ti \hU $.

  If $(X,\fT_X)$ is a separable metrizable topological space, then   $\sB_X \= \si(\fS)$ for any subbase $\fS$ of $\fT_X$
  (see the comment  below Definition 7.6 of \cite{Bertsekas_Shreve_1978}).

  Note: $\{\phi_j\}_{j \in \hN}$ is a dense subset of  $ \wh{C}_b\big([0,\infty) \ti \hU\big) $  by Proposition 7.20 of \cite{Bertsekas_Shreve_1978}.

         \fi
    \fi
  For any $s \ins [0,\infty)$, taking $\phi (r,u) \= \b1_{\{r \in [0,s]\}}   \sI(u)$ gives that
  $\int_0^s  e^{-r} \sI \big( \mu^1_r(\o_0)\big) dr \= \int_0^s  e^{-r} \sI \big( \mu^2_r(\o_0)\big)    dr $.
  Then we obtain that $ \mu^1_s(\o_0) \= \mu^2_s(\o_0) $ for a.e. $ s \ins (0,\infty)$
  or $ \mu^1_\cd(\o_0) \= \mu^2_\cd(\o_0) $ in $\hJ$. A contradiction appears.
  So  the claim \eqref{082020_15} holds.

   By \eqref{031821_a01}, one has  $P_0 \big\{\o_0 \ins \O_0 \n : \mu^1_\cd (\o_0) \nne \mu^2_\cd (\o_0) \big\}
   \= P_0 \big\{\o_0 \ins \O_0 \n :  \tau_1 (\o_0) \nne \tau_2 (\o_0) \big\}
  \= 0$. The former   together with Fubini Theorem renders  that $ (ds \ti dP_0)  \big\{(s,\o_0) \ins [0,\infty) \ti \O_0 \n : \mu^1_s(\o_0) \nne \mu^2_s(\o_0) \big\} \= 0 $ or $\mu^1 \= \mu^2$ in $\fU$,  while  the latter directly means   $\tau_1 \= \tau_2$ in $\fS$.
   \if{0}

   Let $T \ins [1,\infty)$ and set  $ D_T \df \big\{(s,\o_0) \ins [0,T] \ti \O_0 \n : \mu^1_s(\o_0) \nne \mu^2_s(\o_0) \big\}$.
   Since $ \int_0^T \b1_{\{ \mu^1_s(\o_0)   \ne   \mu^2_s(\o_0) \}} ds \= 0$ for any $\o_0 \ins \{ \mu^1_\cd \= \mu^2_\cd \}$, we can deduce that
   \beas
   (ds \ti dP_0) (D_T) \= E \int_0^T  \b1_{\{ \mu^1_s    \ne   \mu^2_s  \}} ds
   \= \int_{\o_0 \in \{\mu^1_\cd \ne \mu^2_\cd \}} \Big( \int_0^T  \b1_{\{ \mu^1_s(\o_0)   \ne   \mu^2_s(\o_0) \}} ds \Big) P_0(d \o_0) \= 0 .
   \eeas
   Letting $T \nto \infty$ yields that
   $ (ds \ti dP_0)  \big\{(s,\o_0) \ins [0,\infty) \ti \O_0 \n : \mu^1_s(\o_0) \nne \mu^2_s(\o_0) \big\} \= 0 $
   or $\mu^1 \= \mu^2$ in $\fU$.

    \fi
   Hence, the mapping $\Ga \n : \fU \ti \fS \mto \fP \big(\O_0 \ti \hJ \ti \hT\big)$
   is injective.

\no  {\bf 2)} We next discuss the continuity of $\Ga$:
    Let $\{\mu^n \}_{n \in \hN}  $ be a sequence of $\fU$   converging to a $\mu \ins \fU$ under $\Rho{\fU}$
    and let $\{\tau_n \}_{n \in \hN}  $ be a sequence of $\fS$   converging to a $\tau \ins \fS$ under $\Rho{\fS}$
 \big(i.e., $ \lmt{n \to \infty} \Rho{\fU} (\mu^n,\mu ) \= \lmt{n \to \infty} \Rho{\fS} (\tau_n,\tau ) \= 0 $\big).
 We need to demonstrate  that $ P^n \df \Ga(\mu^n,\tau_n)$ converges to $ P \df \Ga(\mu,\tau)$ under
 the weak topology of $\fP \big(\O_0 \ti \hJ \ti \hT\big)$, i.e., 
 \bea \label{Au18_01}
  \lmt{n \to \infty} \int_{(\o_0,\fu,\ft) \in \O_0 \times \hJ \times \hT} \psi (\o_0,\fu,\ft) P^n \big(d(\o_0,\fu,\ft)\big)
  \= \int_{(\o_0,\fu,\ft) \in \O_0 \times \hJ \times \hT} \psi (\o_0,\fu,\ft) P \big(d(\o_0,\fu,\ft)\big)
 \eea
  for any bounded continuous function $\psi \n :  \O_0 \ti \hJ \ti \hT \mto \hR$.

    Let $\psi$ be   a real-valued bounded continuous function on $ \O_0 \ti \hJ \ti \hT $.
  For \eqref{Au18_01}, it suffices to show that for any subsequence $\big\{\big(\mu^{n_k},\tau_{n_k}\big)\big\}_{k \in \hN}$ of
  $\big\{(\mu^n,\tau_n) \big\}_{n \in \hN}  $, we can find a subsequence $\big\{\big(\mu^{n'_k},\tau_{n'_k}\big)\big\}_{k \in \hN}$   of   $\big\{\big(\mu^{n_k},\tau_{n_k}\big)\big\}_{k \in \hN}$ that satisfies \eqref{Au18_01}.

  Let  $\big\{\big(\mu^{n_k},\tau_{n_k}\big)\big\}_{k \in \hN}$ be an arbitrary subsequence of   $\big\{(\mu^n,\tau_n) \big\}_{n \in \hN}  $.
  Since $ 0 \=  \lmt{k \to \infty} \Rho{\fU} (\mu^{n_k},\mu ) \= \lmt{k \to \infty}  E_{P_0} \big[ \int_0^\infty \n  e^{-  s} \big( 1 \ld  \Rho{\hU} (\mu^{n_k}_s,\mu_s)  \big)   ds \big] $,
  there exists a subsequence $\Big\{\mu^{\wt{n}_k}\Big\}_{k \in \hN}$ of   $\big\{\mu^{n_k}\big\}_{k \in \hN}$ such that
  for all $\o_0 \ins \O_0 $ except on a $P_0-$null set $\cN_1$,
   $ \lmt{k \to \infty} \Rho{\hU} \Big(\mu^{\wt{n}_k}_s  ( \o_0),\mu_s   ( \o_0) \Big) \= 0 $ for a.e.
   $s \ins (0,\infty)   $. For any $\phi \ins C_b\big([0,\infty) \ti \hU\big)$,  the dominated convergence theorem implies that
  $   \lmt{k \to \infty} \int_0^\infty     \phi  (s,u  )  \friJ \big(\mu^{\wt{n}_k}_\cd  (\o_0) \big)    (ds,du)
    \= \lmt{k \to \infty} \int_0^\infty      e^{-s} \phi \Big(s, \mu^{\wt{n}_k}_s (\o_0) \Big)      ds
   \= \int_0^\infty      e^{-s} \phi \big(s, \mu_s (\o_0)  \big)      ds
   \= \int_0^\infty     \phi  (s,u  )  \friJ \big( \mu_\cd (\o_0)  \big)    (ds,du) $.
    Namely, $\Big\{ \friJ \big(\mu^{\wt{n}_k}_\cd (\o_0)\big) \Big\}_{k \in \hN} $
  converges to $\friJ \big(\mu_\cd (\o_0)\big)$ under the    weak
  topology $\fT_\sharp\big(\fP \big([0,\infty) \\ \ti   \hU\big)\big)$
  of $\fP \big([0,\infty) \ti   \hU\big)$, or equivalently, $\big\{\mu^{\wt{n}_k}_\cd (\o_0) \big\}_{k \in \hN} $
  converges to $\mu_\cd (\o_0)$ under the induced 
  topology $\fT_\sharp(\hJ)$ of $\hJ$.

  As $ \dis 0 \=  \lmt{k \to \infty} \Rho{\fS} \big(\tau_{\overset{}{\wt{n}_k}},\tau \big) \= \lmt{k \to \infty}   E_{P_0} \big[ \Rho{+} (\tau_{\overset{}{\wt{n}_k}},\tau )   \big]   $, one can extract a subsequence $\big\{  {n'_k} \big\}_{k \in \hN} $ from $ \phantom{\Big|} \big\{  {\wt{n}_k} \big\}_{k \in \hN} $ such that $ \lmt{k \to \infty}   \Rho{+} \big(\tau_{n'_k}(\o_0),\tau(\o_0) \big)   \= 0 $ for all $\o_0 \ins \O_0$  except on a $P_0-$null set $\cN_2$.
 Given $\o_0 \ins (\cN_1 \cp \cN_2)^c  $,
 since   $\Big\{\mu^{n'_k}_\cd (\o_0) \Big\}_{k \in \hN} $  also converges to $\mu_\cd (\o_0)$ under $\fT_\sharp(\hJ)$,
 \if{0}

   the induced   weak topology $\fT_\sharp(\hJ)$ of $\hJ$
   \Big(as    a subsequence of $   \Big\{  \mu^{\wt{n}_k}_\cd (\o_0) \Big\}_{k \in \hN} $\Big)
   and since  $ \lmt{k \to \infty}   \Rho{+} \big(\tau_{n'_k}(\o_0),\tau(\o_0) \big)   \= 0 $,
   the continuity of $\psi$ renders that $ \lmt{k \to \infty} \psi \Big( \o_0, \mu^{n'_k}_\cd (\o_0) , \tau_{n'_k} (\o_0) \Big)
  \= \psi \big( \o_0, \mu_\cd (\o_0), \tau (\o_0) \big) $.   Applying the bounded convergence theorem   yields  that

  \fi
   the continuity of $\psi$ and the bounded convergence theorem   yield   that
   \beas
  \hspace{2.2cm}
  && \hspace{-1.4cm} \lmt{k \to \infty} \int_{ (\o_0,\fu,\ft) \in \O_0 \times \hJ \times \hT} \psi (\o_0,\fu,\ft) P^{n'_k} \big(d(\o_0,\fu,\ft)\big)
    \=    \lmt{k \to \infty} \int_{\O_0 } \psi \Big( \o_0,\mu^{n'_k}_\cd (\o_0),\tau_{n'_k}(\o_0)\Big) P_0 (d \o_0) \\
 && \= \int_{\O_0 } \psi \big( \o_0,\mu_\cd (\o_0),\tau(\o_0)\big) P_0 (d \o_0)
     \=    \int_{ (\o_0,\fu,\ft) \in \O_0 \times \hJ \times \hT} \psi (\o_0,\fu,\ft) P \big(d(\o_0,\fu,\ft)\big) . \hspace{3cm} \hb{\qed}
 \eeas

\no {\bf Proof of Lemma \ref{lem_082020_19}:}
Let $\{ t_n \}_{n \in \hN} \sb [0,\infty) $ converge to $t \ins [0,\infty) $
 and let $\{ \oP_n \}_{n \in \hN} \sb \fP\big(\oO\big)$ converge to $\oP \ins \fP\big(\oO\big)$ under
 the weak topology of $\fP\big(\oO\big)$   \big(i.e., $\lmt{n \to \infty} \int_{ \oo \in \oO } \phi (\oo) \oP_n (d\oo)
\= \int_{ \oo \in \oO } \phi (\oo) \oP (d\oo) $ for any bounded continuous function $\phi \n : \oO \mto \hR$\big).
 To see that
$\big\{\oQ_{t_n,\oP_n}  \big\}_{n \in \hN}$ converges to $ \oQ_{t,\oP} $
under the weak topology of $\fP\big(\O_0 \ti \hJ \ti \hT\big) $,
we let $\psi \n : \O_0 \ti \hJ \ti \hT \mto \hR$ be a bounded continuous function and have to show that
$ \lmt{n \to \infty} \int_{(\o_0,\fu,\l) \in  \O_0 \times \hJ \times \hT} \psi (\o_0,\fu,\l ) \oQ_{t_n,\oP_n} \big(d (\o_0,\fu,\l)\big)
  \= \int_{(\o_0,\fu,\l) \in  \O_0 \times \hJ \times \hT} \psi (\o_0,\fu,\l ) \oQ_{t,\oP} \big(d (\o_0,\fu,\l)\big) $.

Set $\|\psi\|_\infty \df \Sup{(\o_0,\fu,\l) \in \O_0 \times \hJ \times \hT} \big| \psi (\o_0,\fu,\l) \big|$
and let $\e \ins (0,1)$. Since   the weakly convergent  sequence $\{\oP_n\}_{n \in \hN}$ 
 is relatively compact in $ \fP\big(\oO\big) $,
  Prohorov's Theorem yields that 
$ \{\oP_n\}_{n \in \hN} $ is tight, i.e., $\dis \Sup{n \in \hN}\oP_n \big(\ol{\cK}^c_\e\big) \ls \frac{\e}{ 4 \|\psi\|_\infty}  $ for some compact subset $\ol{\cK}_\e$ of $\oO$.

Since   $  \oo \mto \big( \oW(\oo),\oU(\oo) \big)   $ is a continuous mapping from $\oO$ to $ \O_0 \ti \hJ $,
Lemma \ref{lem_080422} 
implies that
\beas
 \ol{\Phi}(s,\oo) \df \big(\osW^s  (\oo),\osU^s (\oo),\oT(\oo)\-s\big)
  \= \big(\sW^s  (\oW(\oo)),\sU^s (\oU(\oo)),\oT(\oo)\-s\big) 
   , \q  \fa (s,\oo) \ins [0,\infty) \ti \oO
 \eeas
 is a continuous mapping from $[0,\infty) \ti \oO$ to $\O_0 \ti \hJ \ti \hT$.   
 There exists $\d \= \d(t,\e) \ins (0,1)$ such that
 $\big|\psi \nci \ol{\Phi} (s,\oo) \- \psi \nci \ol{\Phi} (s',\oo') \big| \\  < \n \e / 4  $ for any  $   (s,\oo) , (s',\oo') \ins [0,t\+1] \ti \ol{\cK}_\e $ with $ |s\-s'| \vee \Rho{\oO} (\oo,\oo') \< \d   $.
 And we can find   $N  \= N (t,\e) \ins \hN$ such that
 $ \big| \int_{ \oo \in \oO } \psi \nci \ol{\Phi} (t,\oo) \oP_n (d\oo)
\- \int_{ \oo \in \oO } \psi \nci \ol{\Phi} (t,\oo) \oP (d\oo) \big| \< \frac{\e}{4} $ and $ |t_n \- t| \< \d $ for any  $ n \gs N  $.

 For any $n \ins \hN$,  we can deduce   that
 \beas
 \hspace{1.3cm}
 && \hspace{-1.2cm} \bigg| \int_{(\o_0,\fu,\l) \in  \O_0 \times \hJ \times \hT} \psi (\o_0,\fu,\l ) \oQ_{t_n,\oP_n} \big(d (\o_0,\fu,\l)\big)
   \- \int_{(\o_0,\fu,\l) \in  \O_0 \times \hJ \times \hT} \psi (\o_0,\fu,\l ) \oQ_{t,\oP} \big(d (\o_0,\fu,\l)\big) \bigg| \\
 &&  \hspace{-0.7cm}  \ls \bigg| \int_{\oo \in  \oO} \Big( \psi \big( \ol{\Phi} (t_n,\oo) \big)
 \-   \psi \big( \ol{\Phi} (t,\oo) \big) \Big) \oP_n  (d\oo) \bigg|
 \+ \bigg| \int_{\oo \in  \oO} \psi \big( \ol{\Phi} (t,\oo) \big) \oP_n (d\oo)
 \- \int_{\oo \in  \oO} \psi \big( \ol{\Phi} (t,\oo) \big) \oP  (d\oo) \bigg| \\
 &&  \hspace{-0.7cm}   \<   \int_{\oo \in \ol{\cK}_\e} \Big| \psi \big( \ol{\Phi} (t_n,\oo) \big)
 \-   \psi \big( \ol{\Phi} (t,\oo) \big) \Big| \oP_n  (d\oo)
 \+   \int_{\oo \in  \ol{\cK}^c_\e} \Big| \psi \big( \ol{\Phi} (t_n,\oo) \big) \Big| \oP_n (d\oo)
 \+ \int_{\oo \in  \ol{\cK}^c_\e} \Big| \psi \big( \ol{\Phi} (t,\oo) \big) \Big| \oP_n  (d\oo) \+ \e / 4  \\
 &&  \hspace{-0.7cm}  \ls \frac{\e}{4} \oP_n (\ol{\cK}_\e) \+ 2 \|\psi\|_\infty \oP_n \big(\ol{\cK}^c_\e\big) \+ \e/4 \ls \e.
 \hspace{10.5cm} \hb{\qed}
 \eeas

\no {\bf Proof of Proposition \ref{prop_Ptx_char}:}
Fix $(t,\bx) \ins   [0,\infty) \ti \OmX$.

  \no {\bf 1)} Let $\oP \ins \ocP_{t,\bx}$, which is clearly of $\ocP^1_{t,\bx}$.
  By (D1$'$) and (D4$'$) of Remark \ref{rem_ocP},
 there exist    a   $\hU-$valued, $ \bF^W -$predictable process
 $ \ddot{\mu} \= \{\ddot{\mu}_\fs\}_{\fs \in [0,\infty)} $
 and  a $[0,\infty]-$valued   $ \bF^{W,P_0} -$stopping time $\ddot{\tau}$   on $\O_0$   such that
 $ \oP \big\{  \osU^t_\fs   \=  \ddot{\mu}_\fs (\osW^t)   \hb{ for a.e. }  \fs \ins (0,\infty) \big\}
 \= \oP \big\{  \oT  \=  t \+ \ddot{\tau} \big(\osW^t  \big)   \big\} \= 1$.
   Fubini Theorem shows that
    $  \ddot{\cN}  \df \big\{\o_0 \ins \O_0 \n :  \ddot{\mu}_\cd  (\o_0) \n \notin \n  \hJ \big\} $ is a $\cF^{W,P_0}_\infty-$measurable
   set with zero $P_0-$measure or $ \ddot{\cN} \ins \sN_{P_0}\big(\cF^W_\infty\big)$.
   Then $\breve{\mu}_\fs  \df    \ddot{\mu}_\fs    \b1_{   \ddot{\cN}^c  }  \+ u_0 \b1_{  \ddot{\cN} } $, $  \fs  \ins  [0,\infty)  $
   is an $ \bF^{W,P_0} -$predictable process with all paths in $\hJ$.

 As $\osW^t_{\n \fs} \= \oW_{t+\fs} \- \oW_t  $, $\fs \ins [0,\infty)$ is a Brownian motion under $\oP$ by (D2),
 using Lemma \ref{lem_122921_11}   with $t_0 \= 0$, $(\O_1, \cF_1, P_1,B^1)   \= \big(\oO ,  \sB(\oO ),  \oP , \osW^t\big) $, $(\O_2, \cF_2, P_2,B^2) \= \big(\O_0,  \sB(\O_0),  P_0, W\big) $ and $\Phi \= \osW^t$   yields that
  \bea \label{010922_14}
 \big(\osW^t\big)^{-1}  ( A_0 ) \ins \cF^{\osW^t,\oP}_\infty \= \cF^{\oW^t,\oP}_\infty   \aand
 \oP \nci \big(\osW^t\big)^{-1}  ( A_0 ) \= P_0  ( A_0 )   , \q  \fa A_0 \ins \cF^{W,P_0}_\infty   .
 \eea
  For any   $ A_0 \ins \sB(\O_0) \= \cF^W_\infty$, $\cA \ins \sB(\hJ)$ and   $\cE \ins \sB(\hT)$,
 since $\breve{\mu}^{-1}(\cA) \ins \cF^{W,P_0}_\infty$ by Lemma \ref{lem_061122}  and since $ \ddot{\tau}^{-1}(\cE) \ins \cF^{W,P_0}_\infty $,
 we can derive that
 \beas
  && \hspace{-1cm} \oQ_{t,\oP}   ( A_0 \ti \cA \ti \cE)
  \=   \oP \big\{     \big(\osW^t,\osU^t,\oT \- t  \big) \ins A_0 \ti \cA \ti \cE \big\}
  \= \oP \big\{    \big(\osW^t  , \ddot{\mu}  (\osW^t) ,  \ddot{\tau} (\osW^t   )   \big) \ins A_0 \ti \cA \ti \cE \big\} \\
  &&  \= \oP \big\{    \big(\osW^t  , \breve{\mu}  (\osW^t) ,  \ddot{\tau} (\osW^t   )   \big) \ins A_0 \ti \cA \ti \cE \big\}
   \=      \oP \nci (\osW^t)^{-1}  \big\{     (W, \breve{\mu},  \ddot{\tau}  ) \ins A_0 \ti \cA \ti \cE\big\}
       \=      \oP \nci (\osW^t)^{-1}  \big(      A_0 \Cp \breve{\mu}^{-1}(\cA) \Cp \ddot{\tau}^{-1}(\cE)\big) \\
  && \= P_0  \big(     A_0 \Cp \breve{\mu}^{-1}(\cA) \Cp \ddot{\tau}^{-1}(\cE) \big)
   \=  P_0 \big\{     (W, \breve{\mu}, \ddot{\tau}  ) \ins A_0 \ti \cA \ti \cE \big\}
    \= \big( \Ga (\breve{\mu}, \ddot{\tau}) \big)  (A_0 \ti \cA \ti \cE) .
 \eeas
 \if{0}

 Then the Lambda-system $\big\{ D \ins \sB(\O_0  ) \oti \sB(\hJ)  \oti \sB (\hT) 
 \n : \oP \nci \big(\osW^t,\osU^t,\oT \- t\big)^{-1}(D)
 \=   P_0 \nci (W, \breve{\mu}, \ddot{\tau})^{-1}   (D) \big\} $
 contains all measurable rectangles of $   \sB(\O_0  )   \oti \sB(\hJ)  \oti \sB (\hT) $
 and is thus equal  to $   \sB(\O_0  )   \oti \sB(\hJ)  \oti \sB (\hT) $  
 thanks to Dynkin's Pi-Lambda Theorem.
 So  $ \oP \nci \big(\osW^t,\osU^t,\oT \- t \big)^{-1}   \= P_0 \nci (W, \breve{\mu}, \ddot{\tau})^{-1}    \= \Ga (\breve{\mu}, \ddot{\tau}) \ins \Ga( \fU \ti \fS)  $, i.e., $\oP \ins \ocP^2_t$.

 \fi
  Then Dynkin's Pi-Lambda Theorem implies that
 $ \oQ_{t,\oP} \n \=    \Ga (\breve{\mu}, \ddot{\tau})    $ on $   \sB(\O_0  \ti \hJ \ti \hT) $.
  So  $\oP $   belongs to $ \ocP^2_t$.

    Let $ \wh{\mu} \= \{\wh{\mu}_s\}_{s \in [t,\infty)} $  be the $\hU-$valued, $ \bF^{W^t} -$predictable process 
 in (D1) such that the complement of $ \oO_\mu \df \big\{ \oo \ins \oO \n : \oU_s (\oo) \=  \ol{\mu}_s (\oo)   \hb{ for a.e. }  s \ins (t,\infty) \big\} $ is of $\sN_\oP \big(\sB(\oO)\big)$, where  $\ol{\mu}_s (\oo) \df  \wh{\mu}_s  \big(\oW (\oo)  \big)$, $\fa (s,\oo) \ins [t,\infty) \ti \oO $.
 Given $(\vf,n) \ins \fC(\hR^{d+l}) \ti \hN$,
     (D2$'$) of Remark \ref{rem_ocP} shows that $\big\{\oM^{t,\ol{\mu}}_{s  \land \otau^t_n  } (\vf)\big\}_{s  \in [t,\infty)}$ is a bounded $(\obF^t,\oP)-$martingale.
 For any  $(\fs,\fr) \ins \hQ^{2,<}_+  $  and $  \{(s_i,\cO_i )\}^k_{i=1} \sb \big(\hQ \cap [0,\fs]\big) \ti \sO (\hR^{d+l})  $,
 since   $\oM^t_s(\vf) \= \oM^{t,\ol{\mu}}_s(\vf)$, $\fa s \ins [t,\infty)$ on $\oO_\mu$  and since $ \big\{ (\oW^t_{ t+s_i  },\oX_{ t+s_i   }) \ins \cO_i \big\} \ins 
 \ocF^t_{t+\fs}$ for $i \= 1, \cds \n , k$,
  one has  $ E_\oP \Big[ \big(\oM^t_{  \otau^t_n \land (t+\fr) } (\vf) \- \oM^t_{  \otau^t_n \land (t+\fs) } (\vf) \big) \underset{i=1}{\overset{k}{\prod}}   \b1_{    \{(\oW^t_{ t+s_i  },\oX_{ t+s_i   }) \in \cO_i  \}    }  \Big] \= E_\oP \Big[ \big(\oM^{t,\ol{\mu}}_{  \otau^t_n \land (t+\fr) } (\vf) \- \oM^{t,\ol{\mu}}_{  \otau^t_n \land (t+\fs) } (\vf) \big) \underset{i=1}{\overset{k}{\prod}}   \b1_{    \{(\oW^t_{ t+s_i  },\oX_{ t+s_i   }) \in \cO_i  \}    }  \Big] \= 0 $.
  So $\oP $ is also of $  \ocP^3_t$,   which shows $ \ocP_{t,\bx} \sb \ocP^1_{t,\bx} \Cp \ocP^2_t \Cp \ocP^3_t   $.

\ss \no {\bf 2)} Next, let $\oP \ins \ocP^1_{t,\bx} \Cp \ocP^2_t \Cp \ocP^3_t  $.

\ss \no {\bf 2a)}
 Fix $i,j \ins \{ 1,\cds \n ,d \}$. We  set $  \phi_i(w,x) \df w_i $ and $\phi_{ij}(w,x) \df w_i w_j $
 for any $w \= (w_1,\cds \n ,w_d) \ins \hR^d$ and $  x \ins \hR^l$. Clearly,  $ \phi_i , \phi_{ij} \ins \fC(\hR^{d+l})$.
 One can calculate that   $ \oM^t_s(\phi_i) \= \oW^{t,i}_s   $,   $ \oM^t_s(\phi_{ij}) \= \oW^{t,i}_s \oW^{t,j}_s \- \d_{ij} (s\-t) $, $\fa s \ins [t,\infty)$,
 where  $ \oW^t_s \= \big( \oW^{t,1}_s ,\cds \n , \oW^{t,d}_s  \big)$.

  Let   $n \ins \hN$ and   $ (\fs,\fr) \ins \hQ^{2,<}_+  $. As $\oP \ins \ocP^3_t$,
  it holds for any   $   \big\{(s_i,\cO_i )\big\}^k_{i=1} \n \sb \big(\hQ \cap (0,\fs]\big) \ti \sO (\hR^d)   $ that
   $  E_\oP \Big[  \big( \oM^t_{\otau^t_n \land (t+ \fr)} (\phi_i ) \- \oM^t_{\otau^t_n \land (t+ \fs)} (\phi_i ) \big)     \underset{i=1}{\overset{k}{\prod}}    \b1_{    \{  \oW^t_{t+s_i}  \in  \cO_i\}    }  \Big] \= 0  $
   and  $  E_\oP \Big[  \big( \oM^t_{\otau^t_n \land (t+ \fr)} (\phi_{ij} ) \- \oM^t_{\otau^t_n \land (t+ \fs)} (\phi_{ij} ) \big)     \underset{i=1}{\overset{k}{\prod}}    \b1_{    \{  \oW^t_{t+s_i}  \in  \cO_i\}    }  \Big] \= 0  $.
 So the Lambda-system 
 $ \ol{\L}^{t,n}_{\fs,\fr}  \df \Big\{\oA \ins \sB\big(\oO\big) \n : E_\oP \Big[  \big( \oM^t_{\otau^t_n \land (t+ \fr)} ( \phi_i ) \- \oM^t_{\otau^t_n \land (t+ \fs)} (\phi_i ) \big)  \b1_{\oA} \Big]   \= 0 \hb{ and } E_\oP \Big[  \big( \oM^t_{\otau^t_n \land (t+ \fr)} (\phi_{ij} ) \- \oM^t_{\otau^t_n \land (t+ \fs)} (\phi_{ij} ) \big)  \b1_{\oA} \Big]  \\ \= 0 \Big\} $
  includes the   Pi-system 
  $   \Big\{   \Big( \underset{i=1}{\overset{k}{\cap}}   (\oW^t_{t+s_i} )^{-1}(\cO_i)   \Big)
       \n :    \big\{(s_i,\cO_i )\big\}^k_{i=1} \n   \sb \big(\hQ \Cp (0,\fs]\big) \ti \sO (\hR^d)  \Big\} $,
  which generates $\cF^{\oW^t}_{t+\fs} $.
  \if{0}

 As $\sO(\hR^d)$ is closed under intersection, $ \ol{\sC}^t_\fs \df\Big\{   \Big( \underset{i=1}{\overset{k}{\cap}}   (\oW^t_{t+s_i} )^{-1}(\cO_i)   \Big)
       \n :    \big\{(s_i,\cO_i )\big\}^k_{i=1} \n   \sb \big(\hQ \Cp (0,\fs]\big) \ti \sO (\hR^d)  \Big\} $ is a Pi-system  of $\oO$.
 It is clear that $ \si \big( \ol{\sC}^t_\fs \big) \sb \cF^{\oW^t}_{t+\fs} $.
 The   continuity of process $ \{\oW^t_r\}_{r \in [t,\infty)} $   implies that
  \bea
   \cF^{\oW^t}_{t+\fs}   \=    \si \big(  \oW^t_r   ; r \ins   (t,t\+s]\big)
   \=  \si \big(  \oW^t_{t+\fr}  ; \fr \ins  \hQ  \Cp  (0,\fs]   \big)   .   \label{Jan13_21}
  \eea

   Let $\fr \ins  \hQ  \Cp  (0,s]  $. Since
 $ \big(\oW^t_{t+\fr} \big)^{-1}(\cO)   \ins \si \big( \ol{\sC}^t_\fs \big) $, $ \fa \cO \ins \sO(\hR^d) $,
 the sigma-field 
 $ \L^t_\fs \df \big\{ \cE \sb \hR^d \n : \big(\oW^t_{t+\fr} \big)^{-1} (\cE) \ins \si \big( \ol{\sC}^t_\fs \big) \big\}$
 contains $ \sO(\hR^d) $. Then
 $ \sB(\hR^d) \= \si \big( \sO(\hR^d) \big) \sb \L^t_s $ or
 $ \big(\oW^t_{t+\fr} \big)^{-1} (\cE) \ins \si \big( \ol{\sC}^t_\fs \big) $ for any $   \cE \ins \sB(\hR^d) $.
 It follows from \eqref{Jan13_21} that
 $  \cF^{\oW^t}_{t+\fs}   \=  \si \big(  \oW^t_{t+\fr}  ; \fr \ins  \hQ  \Cp  (0,\fs]   \big)
 \sb \si \big(\ol{\sC}^t_\fs\big) \sb \cF^{\oW^t}_{t+\fs} $.

  \fi
     Dynkin's Pi-Lambda Theorem renders that $ \cF^{\oW^t}_{t+\fs} \sb \ol{\L}^{t,n}_{\fs,\fr}  $, i.e.,
  \bea \label{010922_11}
   E_\oP \Big[  \big( \oM^t_{\otau^t_n \land (t+ \fr)} (\phi_i ) \- \oM^t_{\otau^t_n \land (t+ \fs)} (\phi_i ) \big)  \b1_{\oA} \Big]
   \= 0 \hb{ and }   E_\oP \Big[  \big( \oM^t_{\otau^t_n \land (t+ \fr)} (\phi_{ij} ) \- \oM^t_{\otau^t_n \land (t+ \fs)} (\phi_{ij} ) \big)  \b1_{\oA} \Big]
   \= 0 ,  ~ \fa \oA \ins \cF^{\oW^t}_{t+\fs} . \q
  \eea

   Let $t \ls s \< r \< \infty$ and $\oA \ins \cF^{\oW^t}_s $.
     \if{0}

   Given $ m \ins \hN $,
   set $s_m \df \frac{ \lceil (s -t) 2^m \rceil}{2^m}   \ins \hQ_+ $ and $r_m \df \frac{1+\lceil (r -t) 2^m \rceil}{2^m} \ins \hQ_+ $.
   As $\oA \ins \cF^{\oW^t}_{t+s_m}$,  taking $(\fs,\fr) \= (\fs_m,\fr_m)$ in \eqref{010922_11} yields that
   $  E_\oP \Big[  \big( \oM^t_{\otau^t_n \land (t+ \fr_m)} (\vf ) \- \oM^t_{\otau^t_n \land (t+ \fs_m)} (\vf ) \big)   \b1_{\oA} \Big]   \= 0 $.
   Letting $m \nto \infty$, we can deduce from the continuity of bounded processes $\big\{\oM^t_{s \land \otau^t_n } (\phi_i) \big\}_{s \in [t,\infty)} $, $\big\{\oM^t_{s \land \otau^t_n } (\phi_{ij}) \big\}_{s \in [t,\infty)} $     and the bounded convergence theorem that
   $  E_\oP \Big[ \big( \oM^t_{\otau^t_n \land r} (\vf ) \- \oM^t_{\otau^t_n \land s} (\vf ) \big) \b1_{\oA} \Big]
   \= \lmt{m \to \infty} E_\oP \Big[  \big( \oM^t_{\otau^t_n \land (t+ \fr_m)} (\vf ) \- \oM^t_{\otau^t_n \land (t+ \fs_m)} (\vf ) \big) \b1_{\oA} \Big]   \= 0 $.

     \fi
   Taking $  (\fs,\fr) \= \Big(\frac{\lceil (s -t) 2^m \rceil}{2^m}, \frac{1+\lceil (r -t) 2^m \rceil}{2^m}\Big)$, $ m \ins \hN $  in \eqref{010922_11}
    and sending   $m \nto \infty$,  we can deduce from
      the continuity of bounded processes $\big\{\oM^t_{s \land \otau^t_n } (\phi_i) \big\}_{s \in [t,\infty)} $ and  $\big\{\oM^t_{s \land \otau^t_n } (\phi_{ij}) \big\}_{s \in [t,\infty)} $
   that   $  E_\oP \Big[ \big( \oM^t_{\otau^t_n \land r} (\phi_i ) \- \oM^t_{\otau^t_n \land s} (\phi_i ) \big) \b1_{\oA} \Big]
   \n \= 0 $ and $  E_\oP \Big[ \big( \oM^t_{\otau^t_n \land r} (\phi_{ij} ) \- \oM^t_{\otau^t_n \land s} (\phi_{ij} ) \big) \b1_{\oA} \Big] \= 0 $.
   So $\big\{\oM^t_{s \land \otau^t_n  } (\phi_i )  \big\}_{s \in [t,\infty)}$ and $\big\{\oM^t_{s \land \otau^t_n  } (\phi_{ij} )  \big\}_{s \in [t,\infty)}$ are two      $\big(\bF^{\oW^t},\oP\big)-$martingales.
  As $\lmtu{n \to \infty} \otau^t_n \= \infty$,
  we see that $\big\{\oM^t_s(\phi_i) \= \oW^{t,i}_s   \big\}_{s \in [t,\infty)}  $
  and $\big\{\oM^t_s(\phi_{ij}) \= \oW^{t,i}_s \oW^{t,j}_s \- \d_{ij} (s\-t) \big\}_{s \in [t,\infty)}$   are   $  \bF^{\oW^t}-$local martingales.
    L\'evy's characterization theorem implies that $ \oW^t   $ is a    Brownian motion on $\big(\oO,\sB(\oO),\oP\big)$.
  So $\oP$ satisfies (D2) of $\ocP_{t,\bx}$.
 We still have  \eqref{010922_14} since $\osW^t$ 
 is also a Brownian motion under $\oP$.  

\ss \no {\bf 2b)} As $\oP \ins \ocP^2_t $, there exist  a  $\hU-$valued, $ \bF^{W,P_0} -$predictable process
 $\ddot{\nu} \= \big\{ \ddot{\nu}_s \big\}_{s \in [0,\infty)} $
with all paths in $\hJ$ and a $[0,\infty]-$valued $ \bF^{W,P_0} -$stopping time $\ddot{\ga}$ on $\O_0$ such that
$ \oQ_{t,\oP} 
\= \Ga \big(\ddot{\nu},\ddot{\ga}\big) \= P_0 \nci  \big( W, \ddot{\nu},\ddot{\ga}\big)^{-1} $.
 Given $D \ins \sB(\O_0  \ti \hJ \ti \hT)$, 
 taking $A_0 \= \big( W,\ddot{\nu}, \ddot{\ga}\big)^{-1} (D) \ins \cF^{W,P_0}_\infty $ in \eqref{010922_14}  yields that
 \beas
 \oP \big\{  (\osW^t  ,\osU^t ,\oT\-t ) \ins D \big\}
   \=   \oQ_{t,\oP} (D)
 \= P_0 \nci  \big( W, \ddot{\nu},\ddot{\ga}\big)^{-1} (D)
 \= \oP \nci (\osW^t)^{-1} \big(  ( W, \ddot{\nu},\ddot{\ga} )^{-1} (D) \big)
    \=    \oP \big\{  \big(  \osW^t   , \ddot{\nu}_\cd (\osW^t  ), \ddot{\ga} (\osW^t  ) \big) \ins D \big\} ,
\eeas
 which shows   the  joint $\oP-$distribution of $ (\osW^t,\osU^t ,\oT )$ is the same as the  joint $\oP-$distribution of $ \big(\osW^t, \ddot{\nu}_\cd (\osW^t), t\+\ddot{\ga} (\osW^t  )\big)$.
 Similar to Part (2b) in the proof of \cite[Proposition 4.1]{OSEC_stopping}, we can use
 the equality $\oP \nci (\osW^t,\oT )^{-1} \n \= \oP \nci \big(\osW^t \n ,   t\+\ddot{\ga} (\osW^t  )\big)^{\n -1} $
 to derive that $\oP\{\oT \=  t \+ \ddot{\ga}(\osW^t)\} \= 1$. Namely, $\oP$ satisfies (D4$'$) or equivalently
 (D4) of $\ocP_{t,\bx}$.

  We   claim that
  \bea \label{071622_14}
 \oP \big( \oA \Cp  \big\{\osU^t \ins \cA \big\}\big) \= \oP\big( \oA \Cp \big\{ \ddot{\nu}_\cd (\osW^t) \ins \cA \big\} \big) , \q \fa \oA \ins \cF^{\oW^t,\oP}_\infty , ~ \fa \cA \ins \sB(\hJ) .
\eea
 To see this, we fix $\cA \ins \sB(\hJ)$ and define $\L \df \big\{\oA \in \sB_\oP(\oO) \n :  \oP \big( \oA \Cp  \big\{\osU^t \ins \cA \big\}\big) \= \oP\big( \oA \Cp \big\{ \ddot{\nu}_\cd (\osW^t) \ins \cA \big\} \big)  \big\}$.
 As $\oP      \big\{\osU^t \ins \cA \big\} \\  \=   \oP \{(\osW^t,\osU^t) \ins \O_0 \ti \cA\}
 \= \oP \{(\osW^t,\ddot{\nu}_\cd (\osW^t)) \ins \O_0 \ti \cA\}
\= \oP    \big\{ \ddot{\nu}_\cd (\osW^t) \ins \cA \big\}  $, we see that $\oO \ins \L$ and $\L$ is thus a  Lambda- system.

 For any $(\fs,\cE) \ins [0,\infty) \ti \sB(\hR^d)$, since $W_\fs \n : \O_0 \mto \hR^d$ is a continuous function,
 one has $ W^{-1}_\fs (\cE) \ins \sB(\O_0) $.
 Then it holds for any $\{(\fs_i,\cE_i)\}^N_{i=1} \sb  [0,\infty) \ti \sB(\hR^d)$ that
 \beas
 \oP \Big( \Big(\ccap{i=1}{N} (\osW^t_{\fs_i})^{-1} (\cE_i)\Big) \Cp \{ \osU^t \ins \cA   \}\Big)
& \tn \=& \tn  \oP \Big( \Big\{  \osW^t  \ins  \ccap{i=1}{N}   W^{-1}_{\fs_i} (\cE_i) \Big\} \Cp \{ \osU^t \ins \cA  \}\Big)
\= \oP \Big\{ \big(\osW^t,\osU^t\big) \ins  \ccap{i=1}{N}   W^{-1}_{\fs_i} (\cE_i) \ti \cA  \Big\} \\
& \tn \= & \tn   \oP \Big\{ \big(\osW^t,\ddot{\nu}_\cd(\osW^t)\big) \ins \ccap{i=1}{N}   W^{-1}_{\fs_i} (\cE_i) \ti \cA \Big\}
 \=   \oP \Big( \Big(\ccap{i=1}{N} (\osW^t_{\fs_i})^{-1} (\cE_i)\Big) \Cp \{ \ddot{\nu}_\cd(\osW^t)  \ins \cA  \}\Big) .
\eeas
So $\L$ contains the Pi-system $\Big\{ \ccap{i=1}{N} (\osW^t_{\fs_i})^{-1} (\cE_i) \n : (\fs_i,\cE_i) \ins [0,\infty) \ti \sB(\hR^d), i \= 1,\cds,N \Big\}$, which generates $\cF^{\osW^t}_\infty   \= \cF^{\oW^t}_\infty $.
 Dynkin's Pi-Lambda Theorem implies that   $ \cF^{\oW^t,\oP}_\infty \sb \L $,
\if{0}

 Dynkin's Pi-Lambda Theorem implies that  $ \cF^{\oW^t}_\infty   \sb \L $.

 For any $\fs \ins [0,\infty)$, $\cF^{\oW^t}_\fs \cp \sN_\oP\big(\cF^{\oW^t}_\infty\big)$ is another Pi-system included in $\L$,
 Dynkin's Pi-Lambda Theorem  renders that $\cF^{\oW^t,\oP}_\fs \sb \L$.
 As $\ccup{\fs \in [0,\infty)}{} \cF^{\oW^t,\oP}_\fs $ is also a Pi-system,
 applying Dynkin's Pi-Lambda Theorem again yields   $ \cF^{\oW^t,\oP}_\infty \sb \L $.

 \fi
 proving the claim \eqref{071622_14}.

Let $n,k,j \ins \hN$. By Lemma \ref{lem_061122}, $\cO_{n,k,j} \df \big\{\o_0 \ins \O_0 \n : \ddot{\nu}_\cd (\o_0) \ins \friJ^{-1} \big( O_{\frac{1}{n}}  (\fm_k,\phi_j ) \big\} \ins \cF^{W,P_0}_\infty$.  Since \eqref{010922_14} shows that
  $  \big(\osW^t\big)^{-1}  ( \cO_{n,k,j} ) \ins   \cF^{\oW^t,\oP}_\infty  $, applying
  \eqref{071622_14}  with $\oA \=\big(\osW^t\big)^{-1}  ( \cO^c_{n,k,j} ) $ and $\cA \= \friJ^{-1} \big( O_{\frac{1}{n}}  (\fm_k,\phi_j )$
  yields that
  $ \oP \Big(  (\osW^t)^{-1}  ( \cO^c_{n,k,j} ) \Cp  \big\{\osU^t \ins \friJ^{-1} \big( O_{\frac{1}{n}}  (\fm_k,\phi_j ) \big\}\Big)
  \= \oP\Big( (\osW^t)^{-1}  ( \cO^c_{n,k,j} ) \Cp \big\{ \ddot{\nu}_\cd (\osW^t) \ins \friJ^{-1} \big( O_{\frac{1}{n}}  (\fm_k,\phi_j ) \big\} \Big)
  \=  \oP\big( (\osW^t)^{-1}  ( \cO^c_{n,k,j} ) \\ \Cp (\osW^t)^{-1} (\cO_{n,k,j}) \big) \= 0 $.

 Similar to \eqref{082020_15}, we can deduce that
  $ 
     \big\{\oo \ins \oO \n : \osU^t (\oo) \nne \ddot{\nu}_\cd (\osW^t(\oo)) \big\}   \= \underset{n,k,j  \in \hN}{\cup} \Big(  (\osW^t)^{-1}  ( \cO^c_{n,k,j} ) \cap  \big\{\osU^t \ins \friJ^{-1} \big( O_{\frac{1}{n}}  (\fm_k,\phi_j ) \big\} \Big)  $.
 \if{0}

  We claim that
  $ 
     \big\{\oo \ins \oO \n : \osU^t (\oo) \nne \ddot{\nu}_\cd (\osW^t(\oo)) \big\}   $ is equal to
 $ \underset{n,k,j  \in \hN}{\cup} \Big(  (\osW^t)^{-1}  ( \cO^c_{n,k,j} ) \cap  \big\{\osU^t \ins \friJ^{-1} \big( O_{\frac{1}{n}}  (\fm_k,\phi_j ) \big\} \Big)  $.
 Clearly, $ \underset{n,k,j  \in \hN}{\cup} \Big(  (\osW^t)^{-1}  ( \cO^c_{n,k,j} ) \cap  \big\{\osU^t \ins \friJ^{-1} \big( O_{\frac{1}{n}}  (\fm_k,\phi_j ) \big\} \Big) \sb \big\{  \osU^t   \nne \ddot{\nu}_\cd (\osW^t ) \big\} $.
 Assume that $ \big\{  \osU^t   \nne \ddot{\nu}_\cd (\osW^t ) \big\} \Cp \bigg(  \underset{n,k,j  \in \hN}{\cup}  \Big(  (\osW^t)^{-1}  ( \cO^c_{n,k,j} ) \cap  \big\{\osU^t \ins \friJ^{-1} \big( O_{\frac{1}{n}}  (\fm_k,\phi_j ) \big\} \Big) \bigg)^c $
 is not empty and has an element $\oo $.

  Given $n,j  \ins \hN$, since the proof of Lemma \ref{lem_082020_11} selected  $\{\fm_k\}_{k \in \hN}$ as  a   countable dense subset
  of   the topological space $\big(\fP  ([0,\infty) \ti \hU ), \fT_\sharp \big(\fP   ([0,\infty) \ti \hU )\big)\big) $,
  there exists  $\fk \= \fk (n,j ) \ins \hN$   such that
  $\friJ \big(\osU^t   (\oo)\big) \ins
   O_{\frac{1}{n}} (\fm_\fk,\phi_j ) $.
  This implies 
  $\oo \ins  (\osW^t)^{-1}  ( \cO_{n,\fk,j} ) $
  and thus
  $\big| \int_0^\infty e^{-\fs} \big[ \phi_j \big(\fs,\osU^t_\fs(\oo)\big) \- \phi_j \big(s,\ddot{\nu}_\fs \big(\osW^t (\oo)\big)\big)   \big] d\fs  \big|
   \< 2/n $.
  Letting $n \nto \infty$ yields that $ \int_0^\infty e^{-\fs}   \phi_j \big(\fs,\osU^t_\fs(\oo)\big) d\fs \= \int_0^\infty e^{-\fs} \phi_j \big(\fs,\ddot{\nu}_\fs \big(\osW^t (\oo)\big)\big)    d\fs  $.

  As $\{\phi_j\}_{j \in \hN}$ is   dense in $ \wh{C}_b\big([0,\infty) \ti \hU\big) $  by Proposition 7.20 of \cite{Bertsekas_Shreve_1978},
  the dominated convergence theorem implies that
  $ 
  \int_0^\infty e^{-\fs}   \phi  \big(\fs,\osU^t_\fs(\oo)\big) d\fs \= \int_0^\infty e^{-\fs} \phi  \big(\fs,\ddot{\nu}_\fs (\osW^t (\oo))\big)    d\fs $
    holds for any $\phi \ins \wh{C}_b\big([0,\infty) \ti \hU\big)$.
    By a standard approximation, 
  this equality also  holds for any bounded Borel-measurable functions $\phi$ on $ [0,\infty) \ti \hU $.
      \if{0}
    Using the convolution with smooth mollifier, 
   \eqref{051421_23} also holds for any $\phi \ins C_b\big([0,\infty) \ti \hU\big)$.
   Then   Urysohn's lemma and the dominated convergence theorem yield  that  \eqref{051421_23}   holds for any
  step $\sB\big([0,\infty) \ti \hU\big)-$measurable functions $\phi$   and thus for any bounded $\sB\big([0,\infty) \ti \hU\big)-$measurable functions $\phi$ on $ [0,\infty) \ti \hU $.
          \if{0}
    First for $ \phi \ins \wh{C}_b\big([0,\infty) \ti \hU\big) $, then for $ \phi \ins C_b\big([0,\infty) \ti \hU\big)$, next for bounded Borel-measurable functions on $ [0,\infty) \ti \hU $.

  If $(X,\fT_X)$ is a separable metrizable topological space, then   $\sB_X \= \si(\fS)$ for any subbase $\fS$ of $\fT_X$
  (see the comment  below Definition 7.6 of \cite{Bertsekas_Shreve_1978}).

  Note: $\{\phi_j\}_{j \in \hN}$ is a dense subset of  $ \wh{C}_b\big([0,\infty) \ti \hU\big) $  by Proposition 7.20 of \cite{Bertsekas_Shreve_1978}.
          \fi
       \fi
  For any $\fs \ins [0,\infty)$, taking $\phi (\fr,u) \= \b1_{\{\fr \in [0,\fs]\}}   \sI(u)$ gives that
  $\int_0^\fs  e^{-\fr} \sI \big( \osU^t_\fr(\oo) \big) d\fr \= \int_0^\fs  e^{-\fr} \sI \big( \ddot{\nu}_\fr (\osW^t (\oo)) \big)    d\fr $.
   Then we obtain that $ \osU^t_\fs(\oo) \= \ddot{\nu}_\fs \big(\osW^t (\oo)\big) $ for a.e. $ \fs \ins (0,\infty)$
  or $ \osU^t (\oo) \= \ddot{\nu}_\cd \big(\osW^t (\oo)\big) $ in $\hJ$. A contradiction appears.
  So  the claim holds.  

 \fi
   It follows that
     $\oP \big\{\oo \ins \oO \n : \osU^t (\oo) \nne \ddot{\nu}_\cd (\osW^t(\oo)) \big\} \= 0$
       and thus
       \bea \label{080622_27}
       \oP \big\{\oo \ins \oO \n : \osU^t_\fs (\oo) \= \ddot{\nu}_\fs (\osW^t(\oo)) \hb{ for a.e. } \fs \ins (0,\infty)\big\}
       \= \oP \big\{\oo \ins \oO \n : \osU^t (\oo) \= \ddot{\nu}_\cd (\osW^t(\oo)) \big\} \= 1 .
       \eea

 Like  Lemma 2.4 of \cite{STZ_2011a}, 
 we can construct   a   $[0,1]-$valued  $  \bF^{W}  -$predictable  process  $  \big\{ \breve{\nu}_\fs \big\}_{\fs \in [0,\infty)}$ on $\O_0$
   such that $ \breve{\nu}_\fs  (\o_0)   \= \sI \big(  \ddot{\nu}_\fs (\o_0 )\big) $ for $d\fs \ti dP_0-$a.s. $(\fs,\o_0) \ins [0,\infty) \ti \O_0$. By    Fubini Theorem, it holds      for all $\o_0 \ins \O_0$ except on a $\cN_\nu \ins \sN_{P_0} \big(\cF^W_\infty\big)$ that
   $\breve{\nu}_\fs  (\o_0)   \= \sI \big(  \ddot{\nu}_\fs (\o_0 )\big)    $ for a.e. $\fs \ins [0,\infty)$.
  Define $  \nu^o_\fs(\o_0 ) \df \sI^{-1} \big(\breve{\nu}_\fs (\o_0)\big) \b1_{\{\breve{\nu}_\fs (\o_0) \in \fE\}}  \+ u_0 \b1_{\{\breve{\nu}_\fs (\o_0) \notin \fE\}}
  $, $ (\fs,\o_0) \ins [0,\infty) \ti \O_0 $, which is  a $\hU-$valued $ \bF^{W^t}-$predictable process. 
  Given $\oo \ins \big(\osW^t\big)^{-1}(\cN^c_\nu)$, it holds for a.e. $\fs \ins (0,\infty)$ that
  $\breve{\nu}_\fs  \big(\osW^t(\oo)\big)   \= \sI \big(  \ddot{\nu}_\fs \big(\osW^t(\oo) \big)\big)     $
  and thus   $ \nu^o_\fs \big(\osW^t(\oo)\big)  \= \ddot{\nu}_\fs \big(\osW^t(\oo) \big) $.
 As $\big(\osW^t\big)^{-1}(\cN_\nu) \ins \sN_\oP \big(\cF^{\oW^t}_\infty\big)$, \eqref{080622_27} leads to that
  $ \oP \big\{\oo \ins \oO \n : \osU^t_\fs (\oo) \= \nu^o_\fs (\osW^t(\oo)) \hb{ for a.e. } \fs \ins (0,\infty)\big\}  \= 1  $.
 To wit, $\oP$ satisfies (D1$'$) or equivalently (D1) of $\ocP_{t,\bx}$.

  \no {\bf 2c)}
 Let $ \wh{\nu} \= \{\wh{\nu}_s\}_{s \in [t,\infty)} $  be the $\hU-$valued, $ \bF^{W^t} -$predictable process 
 in (D1) such that the complement of $ \oO_\nu \df \big\{ \oo \ins \oO \n : \oU_s (\oo) \=  \ol{\nu}_s (\oo)   \hb{ for a.e. }  s \ins (t,\infty) \big\} $ is of $\sN_\oP \big(\sB(\oO)\big)$, where  $\ol{\nu}_s (\oo) \df  \wh{\nu}_s  \big(\oW (\oo)  \big)$, $\fa (s,\oo) \ins [t,\infty) \ti \oO $.

 Let $(\vf,n) \ins \fC(\hR^{d+l}) \ti \hN$.
 As $\oP  \big\{   \oX_s   \= \bx(s) , \fa s \ins [0,t]  \big\}   \=1$,
applying  Proposition \ref{prop_MPF1} with $\big(\O,\cF,P,B,X,\mu \big)\=\big(\oO,\sB(\oO),\oP ,\oW,\oX, \ol{\nu} \big) $ implies that $\big\{\oM^{t,\ol{\nu}}_{s \land \otau^t_n } (\vf) \big\}_{s \in [t,\infty)} $ is a bounded $\obF^t-$adapted continuous process under $ \oP  $.

  Let $ (\fs,\fr) \ins \hQ^{2,<}_+  $,  $ \big\{(t_i,\cO_i )\big\}^k_{i=1} \n \sb  \big( \hQ \cap [0,t] \big)  \ti \sO (\hR^l)  $ and    $   \big\{(s_j,\cO'_j )\big\}^m_{j=1} \n \sb \big(\hQ \cap (0,\fs]\big) \ti \sO (\hR^{d+l})   $.
  If $\bx(t_\fri) \notin \cO_\fri$ for some $\fri \ins \{1,  \cds \n , k\}$, then   $\oP\{\oX_{t_\fri} \ins \cO_\fri\} \= 0$ and thus
  $ E_\oP \Big[  \big( \oM^{t,\ol{\nu}}_{\otau^t_n \land (t+ \fr)} (\vf ) \- \oM^{t,\ol{\nu}}_{\otau^t_n \land (t+ \fs)} (\vf ) \big)   \underset{i=1}{\overset{k}{\prod}} \b1_{    \{  \oX_{t_i}  \in  \cO_i\}    } \underset{j=1}{\overset{m}{\prod}} \b1_{    \{ (\oW^t_{t+s_j},\oX_{t+s_j}) \in  \cO'_j\}    }   \Big]  \= 0 $.
  On the other hand, if $\bx(t_i) \ins \cO_i$ for each $i \ins \{1,  \cds \n , k\}$,
   since   $\oM^t_s(\vf) \= \oM^{t,\ol{\nu}}_s(\vf)$, $\fa s \ins [t,\infty)$ on $\oO_\nu$,
  then
    \beas
  && \hspace{-1.2cm}  E_\oP \Big[  \big( \oM^{t,\ol{\nu}}_{\otau^t_n \land (t+ \fr)} (\vf ) \- \oM^{t,\ol{\nu}}_{\otau^t_n \land (t+ \fs)} (\vf ) \big)     \underset{i=1}{\overset{k}{\prod}} \b1_{    \{  \oX_{t_i}  \in  \cO_i\}    }    \underset{j=1}{\overset{m}{\prod}}   \b1_{    \{ (\oW^t_{t+s_j},\oX_{t+s_j}) \in  \cO'_j\}    }  \Big] \\
  &&  \=  E_\oP \Big[  \big( \oM^t_{\otau^t_n \land (t+ \fr)} (\vf ) \- \oM^t_{\otau^t_n \land (t+ \fs)} (\vf ) \big)   \underset{j=1}{\overset{m}{\prod}}   \b1_{    \{ (\oW^t_{t+s_j},\oX_{t+s_j}) \in  \cO'_j\}    }     \Big] \= 0  .
    \eeas
    Since $\ocF^t_{t+\fs} $ is generated by the   Pi-system 
  $   \Big\{  \n \Big(   \underset{i=1}{\overset{k}{\cap}}    \oX_{t_i}^{\,-1}(\cO_i)   \Big) \Cp \Big(   \underset{j=1}{\overset{m}{\cap}}   (\oW^t_{t+s_j},\oX_{t+s_j})^{-1}(\cO'_j)   \Big)
       \n :    \big\{(t_i,\cO_i )\big\}^k_{i=1} \n \sb  \n  \big( \hQ \cap [0,t]  \big) \ti \sO (\hR^l) , \; \big\{(s_j,\cO'_j )\big\}^m_{j=1} \n \sb \n  \big(\hQ \cap (0,\fs]\big) \ti \sO (\hR^{d+l})  \Big\} $,
          Dynkin's Pi-Lambda Theorem renders that  \big(c.f. \eqref{010922_11}\big)
  \bea \label{010922_11b}
   E_\oP \Big[  \big( \oM^{t,\ol{\nu}}_{\otau^t_n \land (t+ \fr)} (\vf ) \- \oM^{t,\ol{\nu}}_{\otau^t_n \land (t+ \fs)} (\vf ) \big)  \b1_{\oA} \Big]
   \= 0 ,  \q \fa \oA \ins \ocF^t_{t+\fs}  .
  \eea

    \if{0}

 So the Lambda-system 
 $ \breve{\L}^{t,\ol{\nu}}_{\fs,\fr}  \df \Big\{\oA \ins \sB\big(\oO\big) \n : E_\oP \Big[  \big( \oM^{t,\ol{\nu}}_{\otau^t_n \land (t+ \fr)} (\vf ) \- \oM^{t,\ol{\nu}}_{\otau^t_n \land (t+ \fs)} (\vf ) \big)  \b1_{\oA} \Big]   \= 0 \Big\} $
  includes the   Pi-system 
  $   \Big\{  \n \Big(   \underset{i=1}{\overset{k}{\cap}}    \oX_{t_i}^{\,-1}(\cO_i)   \Big) \Cp \Big(   \underset{j=1}{\overset{m}{\cap}}   (\oW^t_{t+s_j},\oX_{t+s_j})^{-1}(\cO'_j)   \Big)
       \n :    \big\{(t_i,\cO_i )\big\}^k_{i=1} \n \sb  \n  \big( \hQ \cap [0,t]  \big) \ti \sO (\hR^l) , \; \big\{(s_j,\cO'_j )\big\}^m_{j=1} \n \sb \n  \big(\hQ \cap (0,\fs]\big) \ti \sO (\hR^{d+l})  \Big\} $,
  which generates $\ocF^t_{t+\fs} $.
  \if{0}

 As $\sO(\hR^{d+l})$ is closed under intersection, $ \ol{\sC}^t_s $ is a Pi-system  of $\oO$.
 It is clear that $ \si \big( \ol{\sC}^t_s \big) \sb \ocF^t_{t+s} $.
 The   continuity of process $ \{\oW^t_r\}_{r \in [t,\infty)} $ and process $ \{\oX_r\}_{r \in [0,\infty)} $ implies that
  \bea
   \ocF^t_{t+s}   \=    \si \big( \big(\oW^t_r,\oX_r\big) ; r \ins   (t,t\+s]\big) \ve \si \big(\oX_r; r \ins [0,t]\big)
   \=  \si \big( \big(\oW^t_{t+r},\oX_{t+r}\big) ; r \ins  \hQ  \Cp  (0,s]   \big) \ve \si \big(\oX_r; r \ins [0,t]\big) .   \label{Jan13_21}
  \eea

   Let $r \ins  \hQ  \Cp  (0,s]  $. Since
 $ \big(\oW^t_r,\oX_r\big)^{-1}(\cO)   \ins \ol{\sC}^t_s $, $ \fa \cO \ins \sO(\hR^{d+l}) $,
 the sigma-field 
 $ \L^t_s \df \big\{ \cE \sb \hR^{d+l} \n : \big(\oW^t_r,\oX_r\big)^{-1} (\cE) \ins \si \big(\ol{\sC}^t_s\big)\big\}$
 contains $ \sO(\hR^{d+l}) $. Then
 $ \sB(\hR^{d+l}) \= \si \big( \sO(\hR^{d+l}) \big) \sb \L^t_s $ or
 $ \big(\oW^t_r,\oX_r\big)^{-1} (\cE) \ins \si \big(\ol{\sC}^t_s\big) $ for any $   \cE \ins \sB(\hR^{d+l}) $.
 It follows from \eqref{Jan13_21} that
 $  \ocF^t_{t+s}  \=  \si \big( \big(\oW^t_{t+r},\oX_{t+r}\big) ; r \ins  \hQ  \Cp  (0,s]   \big) \ve \si \big(\oX_r; r \ins [0,t]\big)
 \sb \si \big(\ol{\sC}^t_s\big) \sb \ocF^t_{t+s} $.

  \fi
     Dynkin's Pi-Lambda Theorem renders that $ \ocF^t_{t+\fs} \sb \breve{\L}^{t,\ol{\nu}}_{\fs,\fr}  $, i.e.,
  \bea \label{010922_11b}
   E_\oP \Big[  \big( \oM^{t,\ol{\nu}}_{\otau^t_n \land (t+ \fr)} (\vf ) \- \oM^{t,\ol{\nu}}_{\otau^t_n \land (t+ \fs)} (\vf ) \big)  \b1_{\oA} \Big]
   \= 0 ,  \q \fa \oA \ins \ocF^t_{t+\fs}  .
  \eea

 \fi

   Let $t \ls s \< r \< \infty$ and $\oA \ins \ocF^t_s $.
   \if{0}

   Given $ k \ins \hN $  with $k \> -\log_2(r\-s) $,
   set $s_k \df \frac{\lceil s 2^k \rceil}{2^k} \ins \hQ_+ $ and $r_k \df \frac{\lceil r 2^k \rceil}{2^k} \ins \hQ_+ $.
   As $\oA \ins \ocF^t_{t+s_k}$,  taking $(s,r) \= (s_k,r_k)$ in \eqref{010922_11b} yields that
   $  E_\oP \Big[  \big( \oM^t_{\otau^t_n \land (t+ r_k)} (\vf ) \- \oM^t_{\otau^t_n \land (t+ s_k)} (\vf ) \big)   \b1_{\oA} \Big]   \= 0 $.
   Letting $k \nto \infty$, we can deduce from the continuity of bounded process $\big\{\oM^t_{s \land \otau^t_n } (\vf) \big\}_{s \in [t,\infty)} $     and the bounded convergence theorem that
   $  E_\oP \Big[ \big( \oM^t_{\otau^t_n \land (t+ r )} (\vf ) \- \oM^t_{\otau^t_n \land (t+ s )} (\vf ) \big) \b1_{\oA} \Big]
   \= \lmt{k \to \infty} E_\oP \Big[  \big( \oM^t_{\otau^t_n \land (t+ r_k)} (\vf ) \- \oM^t_{\otau^t_n \land (t+ s_k)} (\vf ) \big) \b1_{\oA} \Big]   \= 0 $.

   \fi
   Taking $  (\fs,\fr) \= \Big(\frac{\lceil (s -t) 2^k \rceil}{2^k}, \frac{1+\lceil (r -t) 2^k \rceil}{2^k}\Big)$, $ k \ins \hN $  in \eqref{010922_11b}
    and letting   $k \nto \infty$,  we can deduce from
      the continuity of bounded process $\big\{\oM^{t,\ol{\nu}}_{s \land \otau^t_n } (\vf) \big\}_{s \in [t,\infty)} $
   that
   \bea   \label{081722_71}
     E_\oP \Big[ \big( \oM^{t,\ol{\nu}}_{\otau^t_n \land r} (\vf ) \- \oM^{t,\ol{\nu}}_{\otau^t_n \land s} (\vf ) \big) \b1_{\oA} \Big]
   \n \= 0 ,
   \q \hb{i.e., $\big\{\oM^{t,\ol{\nu}}_{s \land \otau^t_n  } (\vf )  \big\}_{s \in [t,\infty)}$   is an $\big(\obF^t,\oP\big)-$martingale.}
   \eea
   By Remark \ref{rem_ocP} (ii), $\oP$ satisfies (D3) 
   of $\ocP_{t,\bx}$. \qed

\no {\bf Proof of Proposition \ref{prop_graph_ocP}:}
 According to Proposition \ref{prop_Ptx_char},
 $\big\lan\n\big\lan  \ocP  \big\ran\n\big\ran $ is the intersection of
 $ \big\lan\n\big\lan  \ocP  \big\ran\n\big\ran_1 \df \big\{\big(t,\bx,\oP \big)  \ins [0,\infty) \ti \OmX \ti \fP\big(\oO\big) \n :  \oP \ins \ocP^1_{t,\bx} \big\} $ and  $ \big\lan\n\big\lan  \ocP  \big\ran\n\big\ran_i \df \big\{\big(t,\bx,\oP \big)  \ins [0,\infty) \ti \OmX \ti \fP\big(\oO\big) \n :  \oP \ins \ocP^i_t \big\} $ for $i \= 2,3$.
 Similar to the proof of \cite[Proposition 4.2]{OSEC_stopping}, one can easily show that
  $  \big\lan\n\big\lan  \ocP  \big\ran\n\big\ran_1  $
    is   a Borel subset of $[0,\infty) \ti \OmX \ti \fP\big(\oO\big)$.

 \if{0}

   \no {\bf 1)} Since the  function $\fl_2 (t,\omX) \df \omX(t \ld \cd)  $ is continuous in $ (t,\omX) \ins  [0,\infty) \ti \OmX  $,
   the mapping  $  \psi_{\overset{}{X}}  (t, \bx,\oo) \df     \b1_{\{ \fl_2(t,\oX(\oo) )  - \fl_2(t,\bx )   = 0 \}}
   $, $   (t,\bx, \oo) \ins [0,\infty) \ti \OmX \ti  \oO  $  is $\sB[0,\infty) \oti   \sB(\OmX) \oti \sB(\oO)  -$measurable.
      Lemma \ref{lem_A1} implies that 
 $ \Psi_{\n X}  (t,\bx,\oP) \df \int_{\oo \in \oO} \psi_{\overset{}{X}}  (t, \bx,\oo)  \oP(d \, \oo)
   \= \oP  \big\{  \oX_s \= \bx(s),   \fa s \ins [0,t]   \big\}  $, $   (t,\bx, \oP) \ins [0,\infty) \ti \OmX  \ti   \fP\big(\oO\big) $
  is $ \sB[0,\infty) \oti \sB(\OmX)   \oti \sB\big(\fP\big(\oO\big)\big) 
  - $measurable. So    $  \big\lan\n\big\lan  \ocP  \big\ran\n\big\ran_1 \=
  \big\{ (t,\bx,\oP) \ins [0,\infty) \ti \OmX \ti \fP\big(\oO\big) \n : 
  \Psi_{\n X}  (t,\bx,\oP) \= 1  \big\} \ins  \sB[0,\infty) \oti \sB(\OmX)   \oti \sB\big(\fP\big(\oO\big)\big)$.

 \fi

 \ss \no {\bf 1)}  Since $\O_0$ is a Polish space and since $\hJ$ is a Borel space,
 we know from Proposition 7.13 and  Corollary 7.25.1  of \cite{Bertsekas_Shreve_1978}   that   $ \O_0 \ti \hJ \ti \hT$  with the product topology
 is   a Borel space and   $ \fP \big(\O_0 \ti \hJ \ti \hT\big) $ is also a Borel space.
 As Lemma \ref{lem_082020_15} and   Lemma \ref{lem_082020_17} show that
 $\Ga \n : \fU \ti \fS \ni (  \mu,\tau   ) \mto P_0 \nci  ( W,  \mu,\tau    )^{-1} \ins \fP \big( \O_0 \ti \hJ \ti \hT \big)$
   is a continuous injection from the Polish space $ \fU \ti \fS $ to 
 $   \fP \big(\O_0 \ti \hJ \ti \hT\big)$,
   the image $\Ga(\fU \ti \fS)$ is   a Lusin subset of $ \fP \big(\O_0 \ti \hJ \ti \hT\big) $.
    Theorem A.6 of \cite{Takesaki_1979} implies that  $\Ga(\fU \ti \fS)$ is even a Borel subset of
    the Borel space  $ \fP \big(\O_0 \ti \hJ \ti \hT\big) $.
    Then Lemma \ref{lem_082020_19} yields
 $ \big\lan\n\big\lan  \ocP  \big\ran\n\big\ran_2
  \= \big\{\big(t,\bx,\oP \big)  \ins [0,\infty) \ti \OmX \ti \fP\big(\oO\big) \n : \oQ_{t,\oP} \ins \Ga(\fU \ti \fS) \big\}   \ins \sB[0,\infty) \oti \sB(\OmX) \oti \sB\big(\fP\big(\oO\big)\big)  $.
 \if{0}

  By Lemma \ref{lem_082020_19}, $\breve{\Ga}(t,\bx, \oP) \df \oQ_{t,\oP} \ins \fP \big(\O_0 \ti \hJ \ti \hT\big)$, $\fa (t,\bx, \oP) \ins [0,\infty) \ti \OmX \ti  \fP\big(\oO\big)$ is a continuous mapping. So
  \beas
   \big\lan\n\big\lan  \ocP  \big\ran\n\big\ran_1 \= \big\{\big(t,\bx,\oP \big)
  \ins [0,\infty) \ti \OmX \ti \fP\big(\oO\big) \n :  \oP \ins \ocP^1_t \big\}
  \= \big\{\big(t,\bx,\oP \big)  \ins [0,\infty) \ti \OmX \ti \fP\big(\oO\big) \n : \breve{\Ga} (t, \oP) \ins \Ga(\fU \ti \fS) \big\}      \= \breve{\Ga}^{-1} \big(  \Ga(\fU \ti \fS)  \big)
  \eeas
  is a Borel subset of $[0,\infty) \ti \OmX \ti \fP\big(\oO\big)$.

 \fi

  \no {\bf 2)}
  Since 
  $W(s,\o_0) \df \o_0(s)$ is continuous in $ (s,\o_0) \ins  [0,\infty) \ti \O_0 $
  and 
    $W^X(s,\omX) \df \omX(s)$ is continuous in $ (s,\omX) \ins  [0,\infty) \ti \OmX $,
  the function $\Xi(t,\fs,\o_0,\omX) \df \big( W(t \+ \fs,\o_0) \- W(t,\o_0), W^X(t \+ \fs,\omX)  \big) $   
 is continuous in $   (t,\fs,\o_0,\omX) \ins [0,\infty) \ti [0,\infty) \ti \O_0 \ti \OmX $. 
 For any $n \ins \hN$,
 the mapping $\sT_n (t,\o_0,\omX) \df  \inf\big\{ \fs \ins [0,\infty) \n : |\Xi(t, \fs,\o_0,\omX)|    \gs n  \big\}$,
 $   (t, \o_0,\omX) \ins [0,\infty)   \ti \O_0 \ti \OmX $  is   Borel-measurable since for any $a \ins [0,\infty)$,
  \beas
   && \hspace{-1.2cm}   \big\{(t, \o_0,\omX) \ins [0,\infty)   \ti \O_0 \ti \OmX \n : \sT_n (t,\o_0,\omX) \> a \big\}
    \= \Big\{(t, \o_0,\omX) \ins [0,\infty)   \ti \O_0 \ti \OmX  \n : \Sup{ \fs' \in  [0,\fs]} |\Xi(t, \fs',\o_0,\omX)|   \< n   \Big\}   \\
   && \hspace{-0.5cm} 
    \= \Big(\ccup{k \in \hN}{} \, \ccap{q \in \hQ \cap [0,\fs] }{} \big\{(t, \o_0,\omX) \ins [0,\infty)   \ti \O_0 \ti \OmX  \n : |\Xi(t, q,\o_0,\omX)| \ls n  \- 1/k  \big\}\Big)
    \ins \sB[0,\infty) \oti \sB(\O_0) \oti \sB(\OmX)  .
  \eeas

 Let $ \vf   \ins   \sP(\hR^{d+l})$. Since the function
 $  H_\vf (r,\fx,\nxi, u) \df   \ol{b}  (r,  \fx,u )   \cd   D \vf (\nxi) \+
 \frac12 \ol{\si} \, \ol{\si}^T  (r,  \fx,u ) \n : \n D^2 \vf (\nxi)  $,  $\fa ( r, \fx , \nxi , u) \ins   (0,\infty) \\ \ti \OmX  \ti \hR^{d+l} \ti   \hU$ is Borel-measurable, Lemma \ref{lem_M29_01} (2) shows that the mapping
 \beas
 \cI_\vf (t,\fs,\o_0,\omX,\fu) \df   \int_t^{t+\fs} \n  H_\vf \big(r, \fl_2(r, \omX)  ,\Xi(t,(r\-t)^+,\o_0,\omX), \fu(r)\big)  dr ,
 ~ \fa (t,\fs,\o_0,   \omX,\fu) \ins [0,\infty) \ti [0,\infty) \ti \O_0 \ti \OmX \ti \hJ
 \eeas
  is   $ \sB[0,\infty) \oti \sB[0,\infty) \oti \sB(\O_0) \oti \sB(\OmX) \oti \sB(\hJ)     -$measurable.

  Given $n \ins \hN$ and $\fs \ins [0,\infty)$, since the random variables $(\oW,\oX,\oU)$ on
$\oO$ are $\sB(\O_0) \oti \sB(\OmX) \oti \sB(\hJ)-$measurable, we can derive from
    the Borel measurability of   $\cI_\vf$ and $ \sT_n $   that the mapping
    \bea
 \oM^\vf_{n,\fs} (t,\oo) & \tn \df & \tn  (\vf \nci \Xi)   \big( t, \sT_n \big(t,\oW(\oo),\oX(\oo)\big) \ld n \ld \fs ,\oW(\oo),\oX(\oo) \big)
 \- \cI_\vf  \big( t, \sT_n \big(t,\oW(\oo),\oX(\oo)\big) \ld n \ld \fs ,\oW(\oo),\oX(\oo),\oU(\oo) \big) \nonumber  \\
  & \tn  \= & \tn  \big( \oM^t (\vf) \big) \big(\otau^t_n (\oo) \ld (t\+\fs) ,\oo\big) , \q \fa (t,\oo) \ins [0,\infty) \ti \oO
    \label{Sep21_02}
 \eea
   is    $\sB[0,\infty) \oti \sB(\oO) -$measurable,
 where we used the fact $  \otau^t_n (\oo) \= t\+\sT_n \big(t,\oW(\oo),\oX(\oo)\big) \ld  n   $.

 Let    $\th \df \big(  \vf, n, (\fs,\fr) , \{(s_i,\cO_i )\}^k_{i=1} \big) \ins
   \fC(\hR^{d+l}) \ti \hN  \ti \hQ^{2,<}_+ \ti \wh{\sO} (\hR^{d+l}) $.
  Since
  $  \ff_\th (t,\oo) \df \big( \oM^\vf_{n,\fr} (t,\oo) \- \oM^\vf_{n,\fs} (t,\oo) \big) \ti   \underset{i=1}{\overset{k}{\prod}} \\  \b1_{\{ \Xi  ( t,s_i \land \fs,\oW(\oo),\oX(\oo)  )   \in \cO_i \}  }   $,
  $  (t,\oo) \ins [0,\infty) \ti   \oO$ is   $\sB[0,\infty) \oti \sB(\oO) -$measurable
   by \eqref{Sep21_02},
  an application of   Lemma A.3 of \cite{OSEC_stopping}   yields that the mapping
  $  (t,\oP) \mto \int_{\oo \in \oO} \, \ff_\th (t,\oo) \oP(d \, \oo)
    $
  is $ \sB[0,\infty)   \oti \sB\big(\fP\big(\oO\big)\big) 
  - $measurable and the set
  $ \Big\{ (t,\bx,\oP) \ins [0,\infty) \ti \OmX \ti \fP\big(\oO\big) \n :   E_\oP \Big[  \big( \oM^t_{\otau^t_n \land (t+\fr)} (\vf ) \- \oM^t_{\otau^t_n \land (t+\fs)} (\vf ) \big)  \underset{i=1}{\overset{k}{\prod}}   \b1_{ \{(\oW^t_{t+s_i \land \fs }, \oX_{t+s_i \land \fs }) \in \cO_i \}}   \Big]   \= 0 \Big\} $
  is thus Borel-measurable.
  Letting $\th$ run through the countable collection $\fC(\hR^{d+l}) \ti \hN  \ti \hQ^{2,<}_+ \ti \wh{\sO} (\hR^{d+l})$ shows
  $ \big\lan\n\big\lan  \ocP  \big\ran\n\big\ran_3 \n \ins   \sB[0,\infty) \otimes \sB(\OmX)   \otimes \sB\big(\fP\big(\oO\big)\big)$.

 Totally,  $ \big\lan\n\big\lan  \ocP  \big\ran\n\big\ran \=
 \big\lan\n\big\lan  \ocP  \big\ran\n\big\ran_1 \Cp \big\lan\n\big\lan  \ocP  \big\ran\n\big\ran_2 \Cp \big\lan\n\big\lan  \ocP  \big\ran\n\big\ran_3  $ is a Borel subset of $[0,\infty) \ti \OmX \ti \fP\big(\oO\big)$.   \qed

      \no {\bf Proof of Corollary \ref{cor_graph_ocP}: 1)}
   Let $\ff \n: (0,\infty) \ti \OmX \ti \hU \mto [0,\infty]$ be a Borel-measurable function.
 Taking $\psi(r,\fx,\nxi,u) \df \ff (r,\fx,u) $, $\fa (r,\fx,\nxi,u) \ins (0,\infty) \ti \OmX \ti \hR^{d+l} \ti \hU $
 in Lemma \ref{lem_M29_01} (2)
  shows that the mapping   $   \Psi_\ff  (t,\fs,\omX,\fu) \df  \Psi (t,\fs,\bz,\omX,\fu) \\
   \=  \int_t^{t+\fs} \psi \big(r, \fl_2(r,\omX), 0,\omX(r) ,   \fu(r)\big) dr
  \=    \int_t^{t+\fs} \ff \big(r, \fl_2(r,\omX),\fu(r)\big) dr  $,
  $  (t,\fs, \omX,\fu) \ins [0,\infty) \ti [0,\infty) \ti \OmX \ti \hJ$ is $\sB[0,\infty) \oti \sB[0,\infty) \oti \sB(\OmX) \oti  \sB(\hJ)-$measurable.

 Since the random variables $(\oX,\oU)$ on $\oO$ are $  \sB(\OmX) \oti \sB(\hJ)-$measurable, it follows that
 the mapping
 \bea
 \ol{\Psi}_\ff  (t,\oo) & \tn \df & \tn  
   \int_t^{\oT(\oo) \vee t } \ff \big(r, \oX_{r \land \cd}(\oo),\oU_r(\oo)\big) dr
 \= \lmt{n \to \infty}  \int_t^{(\oT(\oo) \land n) \vee  t}  \ff \big(r, \oX_{r \land \cd}(\oo),\oU_r(\oo)\big) dr  \nonumber \\
 & \tn \=  & \tn \lmt{n \to \infty}  \int_t^{t+(\oT(\oo) \land n -t )^+  } \ff \big(r, \oX_{r \land \cd}(\oo),\oU_r(\oo)\big) dr
  \= \lmt{n \to \infty} \Psi_\ff \big(t , (\oT(\oo) \ld n \-t)^+ ,  \oX (\oo) ,\oU (\oo)\big)  , \hspace{1.6cm}  \label{081222_11}
 \eea
 $\fa (t,\oo)  \ins [0,\infty) \ti \oO$  is $\sB[0,\infty) \oti \sB(\oO)-$measurable, and    Lemma A.3 of \cite{OSEC_stopping} implies that the mapping
 \bea  \label{081122_11}
 \wh{\Psi}_\ff (t,\oP) \df \int_{\oo \in \oO} \n  \ol{\Psi}_\ff  (t,\oo) \oP(d \, \oo) \=  E_\oP \Big[ \int_{\oT \land t}^\oT    \ff     \big(r,\oX_{r \land \cd} , \oU_r   \big) dr \Big], \q \fa (t,\oP) \ins  [0,\infty) \ti  \fP(\oO)
 \eea
 is $ \sB[0,\infty)    \oti \sB\big(\fP(\oO)\big)   - $measurable.

   Let $i \ins \hN$.  Taking $\ff \= g^\pm_i$ and $\ff \= h^\pm_i$ in \eqref{081122_11} yields that
both $ \wh{\Psi}_{g_i} (t,\oP) \df  E_\oP \big[ \int_{\oT \land t}^\oT    g_i \big(r,\oX_{r \land \cd} , \oU_r   \big) dr \big]$
and $ \wh{\Psi}_{h_i} (t,\oP) \df  E_\oP \big[ \int_{\oT \land t}^\oT    h_i \big(r,\oX_{r \land \cd} , \oU_r   \big) dr \big]$,
$ \fa  (t,\oP) \ins  [0,\infty) \ti  \fP(\oO)$     are $ \sB[0,\infty)   \oti   \sB\big(\fP(\oO)\big)   - $measurable.
   Then   the set
  \beas
   \sD  \df \big\{ (t,\bx,y,z,\oP) \ins [0,\infty) \ti \OmX \ti   \Re \ti \Re \ti  \fP(\oO) \n :  \wh{\Psi}_{g_i} (t,\oP)  \ls y_i ,\,
   \wh{\Psi}_{h_i} (t,\oP)  \= z_i ,\,\fa i \ins \hN \big\}
  \eeas
    is   $ \sB[0,\infty) \oti \sB(\OmX) \oti \sB(\Re) \oti \sB(\Re)   \oti \sB\big(\fP(\oO)\big) -$measurable.
      Since   $\big\lan\n\big\lan  \ocP  \big\ran\n\big\ran \ins  \sB[0,\infty) \oti \sB \big(\OmX\big) \oti \sB \big( \fP (\oO) \big)$
  by Proposition \ref{prop_Ptx_char},
  using  the   projection   $\ol{\Pi}_1 (t, \bx,y,z,\oP)   \df   \big(t,\bx ,\oP\big)$ yields that
 $\gP   \=   \big\{ (t,\bx, y,z, \oP ) \ins [0,\infty) \ti \OmX \ti \Re \ti \Re \ti \fP\big(\oO\big)  \n :    \oP \ins \ocP_{t,\bx} ; \;   E_\oP \big[ \int_{\oT \land t}^\oT   g_i\big(r,\oX_{r \land \cd} , \oU_r  \big) dr \big]    \ls y_i ,\, E_\oP \big[ \int_{\oT \land t}^\oT   h_i\big(r,   \oX_{r \land \cd} , \oU_r  \big) dr \big]    \= z_i , \, \fa i \ins \hN \big\}  \=     \ol{\Pi}^{-1}_1 \big( \big\lan\n\big\lan  \ocP  \big\ran\n\big\ran \big)     \Cp \sD$
 is a Borel subset of $D_\ocP \ti \fP\big(\oO\big)$.

   \no {\bf 2)}
   Since 
   $\fl_1 (t,\o_0) \df \o_0(t \ld \cd)  $ is continuous in $ (t,\o_0) \ins  [0,\infty) \ti \O_0  $,
   the mapping  $  \Phi_{\overset{}{W}}  (t, \bw,\oo) \df     \b1_{\{ \fl_1(t,\oW(\oo) )  - \fl_1(t,\bw )   = 0 \}}
   $, $   (t,\bw, \oo) \ins [0,\infty) \ti \O_0 \ti  \oO  $  is $\sB[0,\infty) \oti   \sB(\O_0) \oti \sB(\oO)  -$measurable.

Taking  $\psi(r,\fx,\nxi,u) \df e^{-r} \sI(u) $, $\fa (r,\fx,\nxi,u) \ins (0,\infty) \ti \OmX \ti \hR^{d+l} \ti \hU $
in Lemma \ref{lem_M29_01} (2) shows that the mapping   $   \Psi_{\dn \overset{}{\sI}}   (t,\fs,\fu) \df  \Psi (t,\fs,\bz,\bz,\fu)
   \=  \int_t^{t+\fs} \psi \big(r, \fl_2(r,\bz), 0, 0 ,   \fu(r)\big) dr
   \=  \int_t^{t+\fs} e^{-r} \sI\big(\fu(r)\big) dr    $,
  $  (t,\fs,  \fu) \ins [0,\infty) \ti [0,\infty) \ti   \hJ$ is $\sB[0,\infty) \oti \sB[0,\infty) \oti   \sB(\hJ)-$measurable.
 As the random variable $\oU$ on $\oO$ is $\sB(\oO)/\sB(\hJ)-$measurable,
 we can deduce that the mapping $ \ol{\Psi}_{\dn \overset{}{\sI}}   \big(t,\fu, \oo \big) \df \prod_{q \in \hQ_+} \b1_{\{ \Psi_{\dn \overset{}{\sI}}   (0,t \land q,\oU(\oo))
- \Psi_{\dn \overset{}{\sI}}   (0,t \land q,\fu) = 0 \}} $, $(t, \fu,\oo) \ins [0,\infty) \ti \hJ \ti \oO$
is $\sB[0,\infty) \ti \sB(\hJ) \ti \sB(\oO)-$measurable.
 Then an application of  Lemma A.3 of \cite{OSEC_stopping} again renders that
  the mapping
 \beas
  \ol{\fY} (t,\bw,\bu, \oP)  \df    \int_{\oo \in \oO} \Phi_{\overset{}{W}}  (t, \bw,\oo) \,  \ol{\Psi}_{\dn \overset{}{\sI}}   \big(t,\bu, \oo \big)  \oP(d \oo)
   \=    \oP \big\{  \oW_s  \=  \bw(s)  , \fa s \ins [0,t] ; \,
   \oU_s \= \bu(s) \hb{ for a.e. } s \ins (0,t) \big\}  ,
 \eeas
 $\fa (t,\bw,\bu, \oP) \ins [0,\infty) \ti \O_0 \ti \hJ \ti   \fP(\oO)$ is $\sB[0,\infty) \oti \sB(\O_0) \oti \sB(\hJ) \oti   \sB(\fP(\oO))
  -$measurable.
 Here we used the fact that  $ 
 \Psi_{\dn \overset{}{\sI}}   (0,t \ld q,\oU(\oo))  \= \Psi_{\dn \overset{}{\sI}}   (0,t \ld q,\bu)
   $, $\fa q \in  \hQ_+ $ iff
   $ \int_0^s e^{-r} \sI \big( \oU_r(\oo)\big) dr \= \int_0^s e^{-r} \sI \big(\bu(r)\big) dr $, $\fa s \ins [0,t]$
   iff $  \oU_s(\oo)   \= \bu(s) $ for a.e. $s \ins (0,t)$.
 By the projections
  $\ol{\Pi}_2 (t,\bw,\bu,\bx,y,z,\oP) \df   \big(t,\bx ,\oP\big)$,
$\ol{\Pi}_3 (t,\bw,\bu,\bx,y,z,\oP) \df   \big(t,\bx,y,z,\oP\big)$ and
$\ol{\Pi}_4 (t,\bw,\bu,\bx,y,z,\oP) \df (t,\bw,\bu,  \oP) $,
 \if{0}
 Since the projection $\ol{\Pi}_3 (t,\bw,\bu,\bx,y,z,\oP) \df   \big(t,\bx ,\oP\big)$, $\fa (t,\bw,\bu,\bx,y,z,\oP) \ins   [0,\infty) \ti \O_0 \ti \hJ \ti  \OmX   \ti  \Re \ti \Re  \ti \fP(\oO)$ is $  \sB[0,\infty) \oti \sB(\O_0) \oti \sB(\hJ) \oti  \sB(\OmX)   \oti \sB(\Re) \oti \sB(\Re)   \oti \sB\big(\fP(\oO)\big)\big/\sB[0,\infty) \oti \sB(\OmX) \oti \sB\big(\fP(\oO)\big) -$measurable,
 since the projection $\ol{\Pi}_4 (t,\bw,\bu,\bx,y,z,\oP) \df  \big(t,y,z,\oP\big)$, $\fa (t,\bw,\bu,\bx,y,z,\oP) \ins   [0,\infty) \ti \O_0 \ti \hJ \ti  \OmX   \ti  \Re \ti \Re  \ti \fP(\oO)$ is $  \sB[0,\infty) \oti \sB(\O_0) \oti \sB(\hJ) \oti  \sB(\OmX)   \oti \sB(\Re) \oti \sB(\Re)   \oti \sB\big(\fP(\oO)\big)\big/\sB[0,\infty)   \oti \sB(\Re) \oti \sB(\Re) \oti \sB\big(\fP(\oO)\big) -$measurable,
 since the projection $\ol{\Pi}_5 (t,\bw,\bu,\bx,y,z,\oP) \df (t,\bw,\bu,  \oP)$, $\fa (t,\bw,\bu,\bx,y,z,\oP) \ins   [0,\infty) \ti \O_0 \ti \hJ \ti  \OmX   \ti  \Re \ti \Re  \ti \fP(\oO)$ is $  \sB[0,\infty) \oti \sB(\O_0) \oti \sB(\hJ) \oti  \sB(\OmX)   \oti \sB(\Re) \oti \sB(\Re)   \oti \sB\big(\fP(\oO)\big)\big/\sB[0,\infty) \oti \sB(\O_0) \oti \sB(\hJ) \oti \sB\big(\fP(\oO)\big) -$measurable
 \fi
  we can derive that $\gcP  \= \ol{\Pi}^{-1}_2 \big( \big\lan\n\big\lan  \ocP  \big\ran\n\big\ran \big)
 \Cp \ol{\Pi}^{-1}_3 (\sD) \Cp \ol{\Pi}^{-1}_4 \big( \, \ol{\fY}^{-1} (1) \big)  $
 \if{0}
 we can derive that
 \beas
&& \hspace{-2.2cm} \gcP   \= \Big\{ (t, \bw,\bu,\bx,y,z, \oP ) \ins [0,\infty) \ti \O_0 \ti \hJ \ti \OmX \ti \Re \ti \Re \ti \fP\big(\oO\big)  \n :    \oP \ins \ocP_{t,\bw,\bu,\bx}(y,z)  \Big\} \\
&& \hspace{-1.5cm}  \q \= \Big\{ (t, \bw,\bu,\bx,y,z, \oP ) \ins [0,\infty) \ti \O_0 \ti \hJ \ti \OmX \ti \Re \ti \Re \ti \fP\big(\oO\big)  \n :    \oP \ins \ocP_{t,\bx}   \Big\} \\
&& \hspace{-1.5cm}  \q \cap \, \Big\{ (t, \bw,\bu,\bx,y,z, \oP ) \ins [0,\infty) \ti \O_0 \ti \hJ \ti \OmX \ti \Re \ti \Re \ti \fP\big(\oO\big)  \n :   \oP  \big\{  \oW_s  \=   \bw(s)   ,
    \fa s \ins [0,t] ; \, \oU_s \= \bu(s) \hb{ for a.e. } s \ins (0,t)  \big\} \= 1  \Big\} \\
&& \hspace{-1.5cm}  \q \cap \, \bigg\{ (t, \bw,\bu,\bx,y,z, \oP ) \ins [0,\infty) \ti \O_0 \ti \hJ \ti \OmX \ti \Re \ti \Re \ti \fP\big(\oO\big)  \n :     E_\oP \Big[ \int_{\oT \land t}^\oT   g_i\big(r,\oX_{r \land \cd}  \big) dr \Big]    \ls y_i ,\, E_\oP \Big[ \int_{\oT \land t}^\oT   h_i\big(r,\oX_{r \land \cd}  \big) dr \Big]    \= z_i , \, \fa i \ins \hN \bigg\} \\
&& \hspace{-1.5cm} \q  \= \ol{\Pi}^{-1}_2 \big( \big\lan\n\big\lan  \ocP  \big\ran\n\big\ran \big)
 \Cp \ol{\Pi}^{-1}_3 (\sD) \Cp \ol{\Pi}^{-1}_4 \big( \, \ol{\fY}^{-1} (1) \big)
\eeas
\fi
is a  
Borel subset of $ \cD_\ocP \ti \fP\big(\oO\big) $. \qed

\no {\bf Proof of Theorem \ref{thm_V_usa}:}
   Since   the measurability of functions   $\pi$ and $\fl_2$ shows that
    $\varpi(s,\omX) \df \b1_{\{s < \infty\}}\pi \big(s,\fl_2(s,\omX)\big)$, $(s,\omX) \ins   [0,\infty] \ti   \OmX$ is $\sB  [0,\infty] \oti   \sB(\OmX)-$measurable,
   \if{0}

  Let $\cE \ins \sB(-\infty,\infty]$.   Since $\pi \big(s,\fl_2(s,\omX)\big)$ is Borel-measurable in $(s,\omX) \ins   [0,\infty) \ti   \OmX$,
  one has   $\varpi^{-1} \big(\cE\big) \= \big\{ (s,\omX) \ins [0,\infty) \ti \OmX \n : \pi \big(s,\fl_2(s,\omX)\big) \ins \cE \big\} $ if $0 \notin \cE$
    and $\varpi^{-1} \big(\cE\big) \= \big\{ (s,\omX) \ins [0,\infty) \ti \OmX \n : \pi \big(s,\fl_2(s,\omX)\big) \ins \cE \big\} \cp \big(\{\infty\} \ti \OmX \big) $ if $ 0 \ins \cE$.

   \fi
  taking $\ff \= f^\pm$ in \eqref{081222_11} renders that the mapping
 \beas 
  \ol{\Psi}_{f,\pi}  (t,\oo)   \df     \int_t^{\oT(\oo) \vee t } f \big(r, \oX_{r \land \cd}(\oo),\oU_r(\oo)\big) dr
  \+   \varpi \big(\oT(\oo),  \oX(\oo)\big)     , \q \fa (t,\oo)  \ins [0,\infty) \ti \oO
  \eeas
    is $\sB[0,\infty) \oti \sB(\oO)-$measurable.
   Lemma A.3 of \cite{OSEC_stopping} implies that
   $   \ol{\sV}   (t,\oP) \df \int_{\oo \in \oO} \n \ol{\Psi}_{f,\pi}  (t,\oo)    \oP(d \, \oo)
   $,
 $ (t,\oP) \ins [0,\infty)  \ti \fP\big(\oO\big)$ is $ \sB[0,\infty)   \oti  \sB\big(\fP\big(\oO\big)\big) 
 - $measurable.
 \if{0}

    Lemma \ref{lem_A1} shows that
  $   \ol{\sV}   (t,\oP) \df \int_{\oo \in \oO}   \big( \ol{\cI}^+_f (t, \oo) \+ \ol{\phi}_\pi (\oo) \big) \oP(d \, \oo)
  \- \int_{\oo \in \oO}    \ol{\cI}^-_f (t, \oo)   \oP(d \, \oo)
   \= E_\oP \Big[ \int_t^{\oT  }  f  \big(r,\oX_{r \land \cd}  \big) dr
   \+   \b1_{\{\oT   < \infty\}}   \pi  \big(\oT  ,\oX_{\oT \land \cd}  \big) \Big] $,
 $ (t,\oP) \ins [0,\infty)  \ti \fP\big(\oO\big)$ is $ \sB[0,\infty)   \oti  \sB\big(\fP\big(\oO\big)\big) 
 - $measurable.

 \fi
 Then  Corollary  \ref{cor_graph_ocP} and  Proposition 7.47 of \cite{Bertsekas_Shreve_1978} yield that
 $  \oV (t,\bx,y,z) \=  \Sup{\oP \in \ocP_{t,\bx}(y,z)}  \ol{\sV}(t,\oP) \= \Sup{(t,\bx,y,z,\oP) \in [[\ocP]]}  \ol{\sV}(t,\oP) $
 is  \usa ~ on $  D_\ocP $ and
   $\oV(t,\bw,\bu,\bx,y,z)
   \=  \Sup{(t,\bw,\bu,\bx,y,z,\oP) \in \{\n\{\ocP\}\n\}}  \ol{\sV}(t,\oP)  $
   is  \usa ~ on $  \cD_\ocP$. \qed

\if{0}

  \no {\bf Proof of Corollary \ref{cor_070420_11}: 1)} Let $ \ofN   \ins  \sN_\oP\big(\ocF\big) $, so there is $\oA \ins \ocF \sb \sB(\oO)$
such that $\ofN \sb \oA$ and $\oP(\oA) \= 0 $. By (R2), we can find $\ocN_\oA   \ins \sN_\oP \big(\cF^{\oW^t}_\oga\big)$ such that
$\oP^t_{\oga,\oo} (\oA) \= E_{\oP^t_{\oga,\oo}} \big[ \b1_\oA \big] \= E_\oP \big[ \b1_\oA \big| \cF^{\oW^t}_\oga \big] (\oo) \= 0 $, $ \fa \oo \ins   \ocN^c_\oA   $.
 So for any $\oo \ins   \ocN^c_\oA  $,
  the subset $\ofN$ of $\oA$ is definitely a $\oP^t_{\oga,\oo}-$null set or $ \ofN \ins \sN_{\oP^t_{\oga,\oo}} \big(\ocF\big)$. We can further denote $\ocN_\oA  $ by $\ocN_\ofN  $.

\ss \no {\bf 2)} Let $\oxi$ be a  $ [0,\infty]-$valued, $\si \big( \ocF \cup \sN_\oP\big(\ocF\big) \big)-$measurable random variable  on $\oO$.
Applying Lemma \ref{lemm_cond_exp} with $(\O,\cF,P,\cG_1,\cG_2,\xi) \= \big(\oO,\sB(\oO),\oP,\ocF,\ocF,\oxi\big) $ shows that
$ \oxi \= E_\oP \Big[ \, \oxi \, \Big| \si \big( \ocF \cp \sN_\oP\big(\ocF\big) \big) \Big]
\= E_\oP \Big[ \, \oxi \, \big| \ocF \Big] \n : = \n \oxi ' $,
  $\oP-$a.s. or $\ofN \df \big\{ \oxi \nne \oxi' \big\} \ins \sN_\oP\big(\ocF\big)$.

According to Part 1), it holds  for all $\oo \ins \oO$ except on a $ \ocN^1_\oxi   \ins \sN_\oP\big(\cF^{\oW^t}_\oga\big)  $ that $\ofN  \ins  \sN_{\oP^t_{\oga,\oo}} \big(\ocF\big) $ and thus
$\oxi$ is $\si \big( \ocF \cp \sN_{\oP^t_{\oga,\oo}} \big(\ocF\big) \big)  -$measurable.
Also by (R2), there exists $ \ocN^2_\oxi   \ins \sN_\oP\big(\cF^{\oW^t}_\oga\big)  $ such that
\bea \label{070420_11}
E_{\oP^t_{\oga,\oo}} \big[ \, \oxi' \, \big] \= E_\oP \Big[ \, \oxi' \, \Big| \cF^{\oW^t}_\oga \Big] (\oo)
\= E_\oP \big[ \, \oxi \, \big| \cF^{\oW^t}_\oga \big] (\oo) , \q \fa \oo \ins \big( \ocN^2_\oxi   \big)^c .
\eea

Set $ \ocN_\oxi   \df \ocN^1_\oxi   \cp \ocN^2_\oxi   \ins \sN_\oP\big(\cF^{\oW^t}_\oga\big)  $.  For any $\oo \ins   \ocN^c_\oxi   $,   \eqref{070420_11} implies that
$ E_{\oP^t_{\oga,\oo}} \big[ \, \oxi \, \big] \= E_{\oP^t_{\oga,\oo}} \big[ \b1_{\ofN} \oxi \, \big]
\= E_{\oP^t_{\oga,\oo}} \Big[ \b1_{\ofN} \oxi' \, \Big]
\= E_{\oP^t_{\oga,\oo}} \big[ \, \oxi' \, \big] 
\= E_\oP \big[ \, \oxi \, \big| \cF^{\oW^t}_\oga \big] (\oo) $.  \qed

\fi

 \no {\bf Proof of Proposition \ref{prop_flow}: } We set $t_\oo \df \oga(\oo) \gs t $ for any   $  \oo \ins \oO$.

 \no {\bf  1)}  By (D1) of $ \ocP_{t,\bx} $,
  there is  a   $\hU-$valued, $ \bF^{W^t} -$predictable process
 $ \wh{\mu} \= \{\wh{\mu}_s\}_{s \in [t,\infty)} $   on $\O_0$ such that
 the $\oP-$measure of 
 $ \oO_\mu \df \big\{   \oU_r   \= \wh{\mu}_r (\oW  ) \hb{ for a.e. }  r \ins (t,\infty) \big\}  $
  is equal to $1$. And (R2) assures a $  \ocN_{\n \mu}   \ins \sN_\oP \big(\cF^{\oW^t}_\oga\big) $ such that
 \bea    \label{Jan16_11}
 \oP^t_{\oga,\oo} \big(\oO_\mu\big) \= E_{\oP^t_{\oga,\oo}} \big[\b1_{\oO_\mu}\big]
 \= E_\oP \big[ \b1_{\oO_\mu} \big| \cF^{\oW^t}_\oga\big] (\oo) \= 1 , \q \fa \oo \ins \ocN^c_{\n \mu}  .
 \eea
  According to Lemma \ref{lem_mu_oo},
  for any   $(s,\oo) \ins [t,\infty) \ti \oO   $ we can find   a $\hU-$valued, $ \bF^{W^s} -$predictable  process
   $   \wh{\mu}^{s,\oo} \= \big\{\wh{\mu}^{s,\oo}_r\big\}_{r \in [s,\infty)}$
 on $\O_0$ such that
 $ 
 \wh{\mu}^{s,\oo}_r \big( \oW (\oo')\big) \= \wh{\mu}_r \big( \oW  (\oo')\big)   $, $ \fa (r,\oo') \ins [s,\infty) \ti   \obW^t_{s,\oo}   $,
  where $\obW^t_{s,\oo}   \df \big\{\oo' \ins \oO \n : \oW^t_r (\oo') \= \oW^t_r (\oo), ~ \fa r \ins [t,s] \big\}$.

 Let $\oo \ins \big(   \ocN_0 \cp \ocN_{\n \mu}   \big)^c   $.
 We   set    $ \wh{\mu}^\oo_r ( \o_0) \df \wh{\mu}^{t_\oo,\oo}_r ( \o_0) $, $\fa (r,\o_0) \ins [t_\oo,\infty) \ti \O_0 $,
 which is a $\hU-$valued, $ \bF^{W^{t_\oo}} -$predictable  process  satisfying
 \bea \label{081622_31}
  \wh{\mu}^\oo_r \big( \oW  (\oo')\big)  \= \wh{\mu}_r \big( \oW  (\oo')\big)  ,
  \q \fa (r,\oo') \ins [t_\oo,\infty) \ti   \obW^t_{\oga,\oo}    .  
  \eea
 It follows that
 $ \Wtgo \Cp   \big\{\oo' \ins \oO \n : \oU_r  ( \oo')   \= \wh{\mu}^\oo_r   \big( \oW   (\oo') \big) \hb{ for a.e. }  r \ins (t_\oo,\infty) \big\}
    \n   \supset \n
  \Wtgo \Cp  \oO_\mu $,
  Then \eqref{Jan11_03} and \eqref{Jan16_11} render that
 $ \oP^t_{\oga,\oo} \big\{\oo' \ins \oO \n : \oU_r  ( \oo')   \= \wh{\mu}^\oo_r   \big( \oW   (\oo') \big)   \hb{ for }
 \hb{a.e. }  r \ins (t_\oo,\infty) \big\}
  \= 1  $. So $ \oP^t_{\oga,\oo} $ satisfies (D1) of $\ocP_{t_\oo,\oX_{\oga \land \cd}  (\oo)}$ with
 $ \wh{\mu}  \=   \wh{\mu}^\oo  $ for any $\oo \ins \big(   \ocN_0 \cp \ocN_{\n \mu}   \big)^c   $.

 \no {\bf 2a)} Set $\ol{\mu}_r \df \wh{\mu}_r(\oW)$, $r \ins [t,\infty)$.  From (D3) of $ \ocP_{t,\bx} $ we have
   $\ocN_{\n X} \df \big\{ \oo \ins \oO \n :  \oX_s (\oo)  \nne  \osX^{t,\bx,\ol{\mu}}_s  (\oo)    \hb{ for some } s \ins [0,\infty) \big\}
       \ins \sN_\oP \big(\ocF^t_{\n \infty}\big)$.  
 As $\big\{\osX^{t,\bx,\ol{\mu}}_s \big\}_{s \in [t,\infty)} $ is  an $  \bF^{\oW^t,\oP}  -$adapted continuous   process,
 an analogy to  Lemma 2.4 of \cite{STZ_2011a} allows us to construct   an  $\hR^l-$valued  $  \bF^{\oW^t }  -$predictable  process  $  \big\{ \oK^t_s \big\}_{s \in [t,\infty)}$
   such that $ \ocN_{\n K} \df  \big\{ \oo \ins \oO \n :   \oK^t_s  (\oo)   \nne \osX^{t,\bx,\ol{\mu}}_s  (\oo)   $ for some $s \ins [t,\infty)  \big\}  \ins  \sN_\oP \big(\cF^{\oW^t}_\infty\big)$.
  By (R2), it holds for all $\oo \ins \oO$ except on a   $ \wh{\cN}_{\n X,K} \ins   \sN_\oP\big(\cF^{\oW^t}_\oga\big)$ that
 \bea   \label{081722_33}
 \oP^t_{\oga,\oo} \big(\ocN_{\n X} \cp \ocN_{\n K}\big) \=  E_\oP  \big[ \b1_{\ocN_{\n X} \cup \ocN_{\n K}} \big|\cF^{\oW^t}_\oga  \big] (\oo) \= 0 .
 \eea
   Since   $\Ktgo  \df \ccap{r \in  \hQ \cap (t,\infty)}{} \big\{\oo' \ins \oO \n : \oK^t_{\oga \land r} (\oo') \= \oK^t_{\oga \land r} (\oo) \big\} $
    is an $ \cF^{\oW^t}_\oga -$measurable set   including $\oo$,   (R3) shows  that
  $ 
   \oP^t_{\oga,\oo} \big(\Ktgo\big)   \= 1 $, $ \fa \oo \ins \ocN^c_0$. 
   For any $\oo \ins \big(\ocN_{\n X}  \cp \ocN_{\n K}\big)^c$,    we can deduce from \eqref{090520_11} that
   \bea
   \q &    & \hspace{-1.6cm}   \Wtgo  \Cp    \big(\ocN_{\n X} \n \cp \ocN_{\n K}\big)^c  \n  \Cp \Ktgo
   \=   \Wtgo  \Cp   \big(\ocN_{\n X}  \n  \cp \ocN_{\n K}\big)^c  \n  \Cp   \big\{\oo'  \n \ins \oO \n: \oX_{\n s}(\oo') \= \bx(s),   \fa s \ins [0,t]; \,
   \oK^t_{\n \oga(\oo) \land r}   (\oo') \= \oK^t_{\n \oga(\oo) \land r} (\oo),   \fa r \ins  \hQ \Cp (t,\infty)   \big\} \nonumber \\
    & & \=    \Wtgo  \Cp \big(\ocN_{\n X}  \cp \ocN_{\n K}\big)^c \n \Cp   \big\{\oo'  \n \ins \oO \n: \oX_{\n s}(\oo') 
    \= \oX_{\n s}(\oo) ,   \fa s \ins [0,t]; \,
   \oX_{\n \oga(\oo) \land r}   (\oo') \= \oX_{\n \oga(\oo) \land r} (\oo),   \fa r \ins  \hQ \Cp (t,\infty)   \big\} \nonumber \\
  & & \=   \Wtgo  \Cp \big(\ocN_{\n X}  \cp \ocN_{\n K}\big)^c  \n \Cp \big\{\oo'  \n \ins \oO \n:   \oX_r  ( \oo') \= \oX_{\oga \land r} (  \oo), \, \fa r \ins    [0, \oga(\oo) ]   \big\} . \label{091620_21}
\eea

 Set $\ocN_1 \df    \ocN_{\n X} \cp    \ocN_{\n K} \cp \wh{\cN}_{\n X,K} \ins \sN_\oP \big( \ocF^t_{\n \infty} \big) $.
 Given  $\oo \ins \big(\ocN_0 \cp \ocN_1 \big)^c   $,
  taking $\oP(\cd)$ in \eqref{091620_21} and using \eqref{Jan11_03} yield  that
$ \oP^t_{\oga,\oo} \big\{\oo' \ins \oO \n:   \oX_r  ( \oo') \= \oX_{\oga \land r} (  \oo),    \fa r \ins    [0,t_\oo ]   \big\} \= 1 $.

\ss  \no {\bf 2b)}   For any $\vf \ins \fC(\hR^{d+l})$ and $q \ins \hQ^d$,   define a function $\vf_q (w,x) \df \vf (w\-q,x)$, $(w,x) \ins \hR^{d+l}  $.
We set $\sC \df \{ \vf_q  : \vf \ins \fC(\hR^{d+l}), q \ins \hQ^l\}$, which is a countable sub-collection of $ C^2(\hR^{d+l})$.
 For any $n \ins \hN $, define an  $\obF^t-$stopping time by $ \oz_n (\oo) \df \inf\big\{r \ins [\oga(\oo),\infty) \n : | \oW^t_r (\oo) \- \oW^t_\oga (\oo) |^2  
 \+ |\oX_r (\oo)|^2 \gs n^2  \big\} \ld \big(\oga (\oo)\+n \big)$, $\oo \ins \oO$.
    \if{0}

 In general,  let $\tau$ be an $\bF-$stopping time and let $X$ be an $\bF-$adapted continuous process. Define $\z (\o) \df \inf\big\{r \ins [\tau(\o),\infty) \n: | X_r (\o)|  \gs n \big\}$, $\o \ins \O$. We can show that for any $s \ins [0,\infty)$
 \beas
 \{\z \> s\} \= \ccup{k \in \hN}{} \ccap{r \in \hQ \cap [0,s]}{} \Big(\{r \ls \tau\} \cp \big(\{r \> \tau\} \Cp \{ |X_r| \ls n \- 1/k \}\big)  \Big) \ins \cF_s .
 \eeas
 So $\z $ is also an $\bF-$stopping time.

    \fi

Let  $\th \df \big(  \phi, n,j, (\fs,\fr) , \{(s_i,\cO_i )\}^k_{i=1} \big) \ins
   \sC \ti \hN \ti \hN  \ti \hQ^{2,<}_+ \ti \wh{\sO} (\hR^{d+l})    $.
 Since    $ \big\{ \oM^{t,\ol{\mu}}_{s \land \otau^t_j } (\phi) \big\}_{ s \in  [t,\infty) } $
 is a bounded $(\obF^t,  \oP)-$martingale by 
 applying Proposition \ref{prop_MPF1} with $\big(\O,\cF,P,B,X,\mu \big)\=\big(\oO,\sB(\oO),\oP ,\oW,\oX, \ol{\mu} \big) $,
 the  optional sampling theorem implies  that
 $  E_\oP  \Big[   \oM^{t,\ol{\mu}}_{  ( \oga  +\fr) \land \oz_n  \land \otau^t_j} (\phi  )    \Big|\ocF^t_{\n \oga+\fs}  \Big]
  \=  \oM^{t,\ol{\mu}}_{  ( \oga  +\fs) \land \oz_n \land \otau^t_j } ( \phi  )  $,  $\oP-$a.s.
 Set $ \oxi_\th \df \oM^{t,\ol{\mu}}_{  ( \oga  +\fr) \land \oz_n  \land \otau^t_j} (\phi  )   \- \oM^{t,\ol{\mu}}_{  ( \oga  +\fs) \land \oz_n \land \otau^t_j } ( \phi  ) \= \b1_{\{\otau^t_j > \oga \}} \Big( \oM^{t,\ol{\mu}}_{  ( \oga  +\fr) \land \oz_n  \land \otau^t_j} (\phi  )   \- \oM^{t,\ol{\mu}}_{  ( \oga  +\fs) \land \oz_n \land \otau^t_j } ( \phi  ) \Big) $
 and set $\oeta_\th \df
 \underset{i=1}{\overset{k}{\prod}}   \b1_{    \{(\oW^t_{ \oga +s_i \land \fs } - \oW^t_\oga,\oX_{ \oga+s_i \land \fs  }) \in \cO_i  \}    } \ins \ocF^t_{\n \oga+\fs} $.
 As $ \cF^{\oW^t}_\oga \sb \ocF^t_{\n \oga} \sb \ocF^t_{\n \oga+\fs} $, the tower property   renders that
 $ E_\oP  \big[  \,  \oxi_\th \oeta_\th \big|\cF^{\oW^t}_\oga  \big]
 \=  E_\oP  \Big[    \oeta_\th  E_\oP  \big[  \,  \oxi_\th   \big|\ocF^t_{\n \oga+\fs}  \big]  \Big|\cF^{\oW^t}_\oga  \Big] \= 0 $, $ \oP-$a.s.
 By (R2) again, there exists an $\ocN_\th \ins \sN_\oP \big(\cF^{\oW^t}_\oga\big)$ such that
 \bea \label{012822_17}
  E_{\oP^t_{\oga,\oo}} \big[ \, \oxi_\th \oeta_\th \big]
   \= E_\oP  \big[  \,  \oxi_\th \oeta_\th \big|\cF^{\oW^t}_\oga  \big] (\oo) \= 0  , \q \fa \oo \ins \ocN^c_\th  .
 \eea

 Define $\ocN_2 \df     \bigcup \big\{\ocN_\th \n :  \th \ins \sC \ti \hN \ti \hN  \ti \hQ^{2,<}_+ \ti \wh{\sO} (\hR^{d+l}) \big\}   \ins \sN_\oP \big(\ocF^t_{\n \infty}\big)$. We   fix   $\oo \ins  \big( \ocN_0 \cp \ocN_1 \cp \ocN_2 \cp \ocN_{\n \mu}  \big)^c $
 and set $\ol{\mu}^\oo_r (\oo') \df \wh{\mu}^\oo_r (\oo')  $, $ (r,\oo') \ins [t_\oo,\infty) \ti \oO$.
 Let    $\big(\vf,n,(\fs,\fr),   \{(s_i,   \cO_i )\}^k_{i=1} \big) \ins \fC(\hR^{d+l}) \ti \hN \ti \hQ^{2,<}_+  \ti  \wh{\sO} (\hR^{d+l}) $
 and let $j \ins \hN$.

 There exists a sequence  $\{q_m\=q_m(\oo)\}_{m \in   \hN} $ of $ \hQ^d$ that converges to $  \oW^t_\oga(\oo)$.
 Let $ m \ins \hN$. We set $\th_m \df \big(  \vf_{q_m},   n,j,   (\fs, \\ \fr) , \{(s_i,\cO_i )\}^k_{i=1} \big)$ and define $\d^{j,m}_\oo    \df \Sup{|(w,x)| \le j} \Big( \sum^2_{i=0} \big|D^i\vf_{q_m}(w,x) \- D^i \vf \big(w\-\oW^t_\oga(\oo),x\big)\big| \Big)   \= \Sup{|(w,x)| \le j} \Big( \sum^2_{i=0}\big|D^i\vf (w\-q_m,x) \- D^i\vf \big(w\-\oW^t_\oga(\oo),x\big)\big| \Big)$.

  Given $\oo' \ins \Wtgo \Cp \ocN^c_{\n X} \Cp \big\{  \otau^t_j   \> \oga \big\}   $, \eqref{090520_11} implies that
 $\otau^t_j (\oo') \> \oga(\oo') 
 \= t_\oo$ and
 $ 
  \oz_n(\oo') \= \inf\big\{r \ins [t_\oo,\infty) \n : | \oW_r (\oo') \- \oW_{t_\oo} (\oo') |^2
 \+ |\oX_r (\oo')|^2 \gs n^2  \big\} \ld \big(t_\oo\+n \big) \= \otau^{t_\oo}_n (\oo')$. 
 Since $ \oW^{t_\oo}_{\n r}(\oo') \= \oW^t_{\n r}(\oo') \- \oW^t_{\n t_\oo}(\oo')  
 \= \oW^t_{\n r}(\oo') \- \oW^t_\oga(\oo) $, $\fa r \ins [t_\oo,\infty)$,
   \eqref{081622_31} shows  that   for any   $t_\oo \ls s_1  \ls s_2 \< \infty  $
  \beas
\; && \hspace{-2cm} \big(\oM^{t_\oo,\ol{\mu}^\oo}_{s_2} (\vf) \- \oM^{t_\oo,\ol{\mu}^\oo}_{s_1} (\vf)\big)(\oo')
\= \vf \big(\oW^{t_\oo}_{\n s_2} (\oo')  , \oX_{s_2} (\oo') \big) \- \vf \big(\oW^{t_\oo}_{\n s_1} (\oo')  , \oX_{s_1} (\oo') \big)
  \- \n \int_{s_1}^{s_2}  \n  \ol{b}  \big( r, \oX_{r \land \cd}(\oo'),\wh{\mu}^\oo_r(\oW(\oo')) \big) \n \cd \n D \vf \big( \oW^{t_\oo}_{\n r} (\oo') , \oX_r(\oo') \big) dr \\
  && \hspace{-0.7cm}  -   \frac12 \n \int_{s_1}^{s_2}  \n  \ol{\si} \, \ol{\si}^T  \big( r,  \oX_{r \land \cd}(\oo'),\wh{\mu}^\oo_r(\oW(\oo'))  \big) \n : \n D^2 \vf  ( \oW^{t_\oo}_{\n r}(\oo')  , \oX_r(\oo')  )   dr \\
 && \hspace{-1cm} \=   \vf \big(\oW^t_{\n s_2} (\oo') \- \oW^t_{\n \oga} (\oo) , \oX_{\n s_2} (\oo') \big) \- \vf \big( \oW^t_{\n s_1} (\oo') \- \oW^t_{\n \oga} (\oo) , \oX_{\n s_1} (\oo') \big)
   \- \n \int_{s_1}^{s_2}  \n  \ol{b}  \big( r, \oX_{r \land \cd} (\oo'), \wh{\mu}_r(\oW(\oo'))  \big) \n \cd \n D \vf \big( \oW^t_{\n r} (\oo') \- \oW^t_\oga (\oo) , \oX_r (\oo') \big) dr \\
  && \hspace{-0.7cm}  -   \frac12 \n \int_{s_1}^{s_2} \n  \ol{\si} \, \ol{\si}^T  \big( r,  \oX_{r \land \cd}(\oo'), \wh{\mu}_r(\oW(\oo')) \big) \n : \n D^2 \vf  ( \oW^t_{\n r} (\oo') \- \oW^t_\oga (\oo) , \oX_r  (\oo') )   dr .
\eeas
 As $   \big| \big(\oW^t_r (\oo')  , \oX_r (\oo') \big) \big|  \ls j $ for any $r \ins [t_\oo, \otau^t_j(\oo') ]$,
 we can deduce from   \eqref{coeff_cond1}, \eqref{coeff_cond2}  and Cauchy-Schwarz inequality   that
   for any   $t_\oo \ls s_1  \ls s_2 \ls \otau^t_j(\oo')  $
 \beas
&& \hspace{-1.4cm}
  \Big| \big( \oM^{t_\oo,\ol{\mu}^\oo}_{s_2} (\vf) \- \oM^{t_\oo,\ol{\mu}^\oo}_{s_1} (\vf) \- \oM^{t,\ol{\mu}}_{s_2} (\vf_{q_m} )  \+ \oM^{t,\ol{\mu}}_{s_1} ( \vf_{q_m} ) \big) (\oo') \Big|   \\
&&  \hspace{-0.9cm} \ls 2 \d^{j,m}_\oo  \+ \d^{j,m}_\oo \n \int_t^{\otau^t_j(\oo')}  \n \Big(   \big| \ol{b}  \big( r, \oX_{r \land \cd} (\oo'), \wh{\mu}_r(\oW(\oo'))  \big) \big|   \+ \frac12   \big| \ol{\si}  \big( r, \oX_{r \land \cd} (\oo'), \wh{\mu}_r(\oW(\oo'))  \big) \big|^2 \Big) dr \\
&&  \hspace{-0.9cm} \ls 2 \d^{j,m}_\oo  \+ \d^{j,m}_\oo \n \int_t^{\otau^t_j(\oo')} \n \Big(     \k( r)   \big\|\oX_{r \land \cd}(\oo')\big\|_r
    \+   \big| b \big( r,\bz,\wh{\mu}_r(\oW(\oo')) \big) \big|    \+      \frac{d}{2}   \+    \k^2( r)   \big\|\oX_{r \land \cd}(\oo')\big\|^2_r
    \+   \big|\si \big( r,\bz,\wh{\mu}_r(\oW(\oo')) \big) \big|^2 \Big) dr \ls \d^{j,m}_\oo (2   \+ c^j_{t,\bx}) . 
 \eeas
 where $  c^j_{t,\bx} \df \big[ d / 2 \+  \k( t\+j)  ( \|\bx\|_t  \+ j    )  \+     \k^2( t\+j)  ( \|\bx\|_t  \+  j     )^2      \big] j
 \+ \int_t^{t+j}  \Sup{u \in \hU} \big(  |b( r,\bz,u)| \+ |\si( r,\bz,u)|^2 \big)  dr \< \infty $.
 Taking $ s_1  \=   \big(( \oga   \+\fs) \ld \oz_n   \ld \otau^t_j\big) (\oo') \=  ( t_\oo  \+\fs) \ld \otau^{t_\oo}_n (\oo') \ld \otau^t_j (\oo')  $ and $ s_2  
 \=   ( t_\oo  \+\fr) \ld \otau^{t_\oo}_n (\oo') \ld \otau^t_j (\oo') $ yields that
 $ \Big| \big( \oM^{t_\oo,\ol{\mu}^\oo}_{( t_\oo  +\fr) \land \otau^{t_\oo}_n \land \otau^t_j} (\vf)  - \oM^{t_\oo,\ol{\mu}^\oo}_{( t_\oo  +\fs) \land \otau^{t_\oo}_n \land \otau^t_j} (\vf)\big) (\oo') \- \oxi_{\th_m}   (\oo') \Big|
 \ls \d^{j,m}_\oo (2   + c^j_{t,\bx}) $.
 As $ \oeta_{\th_m}   (\oo') \=  \underset{i=1}{\overset{k}{\prod}}   \b1_{    \{(\oW^{t_\oo}_{ t_\oo+s_i \land \fs } (\oo'),\oX_{ t_\oo+s_i \land \fs  } (\oo')) \in \cO_i  \}    }    $ by \eqref{090520_11},
  we see from \eqref{Jan11_03} and \eqref{081722_33}  that
   \beas
   \q E_{\oP^t_{\oga,\oo}}  \Big[ \b1_{\{  \otau^t_j   > \oga \}} \Big| \Big( \oM^{t_\oo,\ol{\mu}^\oo}_{  ( t_\oo +\fr) \land \otau^{t_\oo}_n \land \otau^t_j } (\vf )  \- \oM^{t_\oo,\ol{\mu}^\oo}_{  ( t_\oo +\fs) \land \otau^{t_\oo}_n \land \otau^t_j } (\vf ) \Big) \underset{i=1}{\overset{k}{\prod}} \,  \b1_{    \{(\oW^{t_\oo}_{ t_\oo+s_i \land \fs },\oX_{ t_\oo+s_i \land \fs  }) \in \cO_i  \}    } \- \oxi_{\th_m} \oeta_{\th_m} \Big|   \Big] \ls \d^{j,m}_\oo (2   \+ c^j_{t,\bx})  .
    \eeas

 The uniform  continuity of $D^i \vf  $'s    over compact sets 
   implies   $ \lmtd{m \to \infty} \d^{j,m}_\oo \= 0$, and one can  then  obtain  from  \eqref{012822_17} that
  \bea \label{012922_17}
  \hspace{-0.3cm}
  E_{\oP^t_{\oga,\oo}} \n \Big[ \b1_{\{  \otau^t_j   > \oga \}} \Big( \oM^{t_\oo,\ol{\mu}^\oo}_{  ( t_\oo +\fr) \land \otau^{t_\oo}_n \land \otau^t_j } (\vf )  \- \oM^{t_\oo,\ol{\mu}^\oo}_{  ( t_\oo +\fs) \land \otau^{t_\oo}_n \land \otau^t_j } (\vf ) \Big) \underset{i=1}{\overset{k}{\prod}} \,  \b1_{    \{(\oW^{t_\oo}_{ t_\oo+s_i \land \fs },\oX_{ t_\oo+s_i \land \fs  }) \in \cO_i  \}    }    \Big]
   \n \=  \n  \lmt{m \to \infty} E_{\oP^t_{\oga,\oo}}  \big[ \,  \oxi_{\th_m} \oeta_{\th_m}  \big]  \n \= 0. \q
  \eea
 Since $  \oP^t_{\oga,\oo} \big\{\oo' \ins \oO \n:   \oX_r  ( \oo') \= \oX_{\oga \land r}  (  \oo),    \fa r \ins    [0,t_\oo ]   \big\} \= 1 $
 by Part (2a),
  applying Proposition \ref{prop_MPF1} with $\big(\O,\cF,P,B,X \big)\=\big(\oO,\sB(\oO),\oP^t_{\oga,\oo},\oW,\oX \big) $
  and   $(t,\bx,\mu) \= \big(t_\oo,\oX_{\oga \land \cd} (\oo),   \ol{\mu}^\oo\big)$ renders that
 $  \big\{\oM^{t_\oo,\ol{\mu}^\oo}_{  s \land \otau^{t_\oo}_n  } (\vf )\big\}_{s \in [t_\oo,\infty)}  $ is a bounded   process under $\oP^t_{\oga,\oo}$.
 \if{0}
 applying Proposition \ref{prop_MPF1} with $(t,\bx) \= \big(t_\oo,\oX_{\oga \land \cd}(\oo)\big)$ and
 $(\O,\cF,P,B,X) \= \big( \oO,\sB(\oO),\oP^t_{\oga,\oo},\oW,\oX\big)$ renders  that
 $  \big\{\oM^{t_\oo,\ol{\mu}^\oo}_{  s \land \otau^{t_\oo}_n  } (\vf )\big\}_{s \in [t_\oo,\infty)}  $ is a bounded $\obF^{t_\oo}-$adapted continuous process under $\oP^t_{\oga,\oo}$.
 \fi
   As $\lmtu{j \to \infty} \otau^t_j (\oo') \= \infty$ for any $\oo' \ins \oO$,
   letting $j \nto \infty$ in \eqref{012922_17}
 and using the  bounded convergence theorem reach that
 $   E_{\oP^t_{\oga,\oo}} \n \Big[   \Big( \oM^{t_\oo,\ol{\mu}^\oo}_{  ( t_\oo +\fr) \land \otau^{t_\oo}_n   } (\vf )  \- \oM^{t_\oo,\ol{\mu}^\oo}_{  ( t_\oo +\fs) \land \otau^{t_\oo}_n   } (\vf ) \Big) \underset{i=1}{\overset{k}{\prod}}   \b1_{    \{(\oW^{t_\oo}_{ t_\oo+s_i \land \fs },\oX_{ t_\oo+s_i \land \fs  }) \in \cO_i  \}    }    \Big]    \= 0  $.
Following similar arguments to those that lead to \eqref{081722_71}, we can   derive that $  \big\{\oM^{t_\oo,\ol{\mu}^\oo}_{  s \land \otau^{t_\oo}_n  } (\vf )\big\}_{s \in [t_\oo,\infty)}  $ is a   $\big(\obF^{t_\oo},\oP^t_{\oga,\oo}\big)-$martingale.
 Then Remark \ref{rem_ocP} (ii) shows that $ \oP^t_{\oga,\oo}$ satisfies (D2)+(D3) of $\ocP_{t_\oo,\oX_{\oga \land \cd}  (\oo)}$ for
 any  $\oo \ins  \big( \ocN_0 \cp \ocN_1 \cp \ocN_2 \cp \ocN_{\n \mu}  \big)^c $.

\no {\bf 3)} By (D4) of $\ocP_{t,\bx}$,   there is   a $[t,\infty]-$valued   $ \bF^{W^t,P_0} -$stopping time $\wh{\tau}$    such that
 $ \oP \big\{  \oT  \=   \wh{\tau} (\oW)   \big\} \= 1$.
 Analogous  to Part (2) in the proof of \cite[Proposition 5.1]{OSEC_stopping}, we can find  $ \ocA_*  \ins \cF^{\oW^t}_\oga $ satisfying
    \bea \label{020322_11}
 \big\{ \wh{\tau}(\oW) \gs \oga \big\} \D \ocA_* \ins \sN_\oP \big(\cF^{\oW^t}_\infty\big) ,
 \eea
  and there exists $  \ocN_3   \ins \sN_\oP \big(\cF^{\oW^t}_\oga\big) $   such that
 for any  $ \oo \ins \ocN^c_0 \Cp \ocN^c_3 \Cp \ocA_* $,  $\oP^t_{\oga,\oo}   \big\{ \oT \=   \wh{\tau}^\oo (\oW  ) \big\}   \= 1$
 for some   $[t_\oo,\infty]-$valued   $\bF^{W^{t_\oo},P_0}-$stopping time $ \wh{\tau}^\oo $, namely, $\oP^t_{\oga,\oo} $ satisfies
  (D3)  of $\ocP_{t_\oo,\oX_{\oga \land \cd}  (\oo)}$.

Let $i \ins \hN$.
   \if{0}

 Since the $\bF^{\oW^t}-$stopping time $\oga$ is $\cF^{\oW^t}_\oga-$measurable
 and since  $\cF^{\oW^t}_\oga \sb \cF^{\oW^t}_\infty \sb \sB(\oO) $,
  we can deduce from the measurability of function  $ \ol{\cI}_{g_i}  $   defined in \eqref{Sep21_04}  that
  $ \oO \n \ni \n  \oo \mto \ol{\cI}_{g_i} \big(\oga(\oo),\oo\big) \= \int_{\oT(\oo) \land \oga(\oo)}^{\oT(\oo)}  g_i  \big(r,\oX (r \ld \cd,\oo) \big) dr  $
  is a $\sB(\oO) -$measurable random variable on $\oO$.

   \fi
 According to (R2), it holds for all $\oo \ins \oO$ except on a $\ocN^{\,i}_{\n g,h}   \ins \sN_\oP \big(\cF^{\oW^t}_\oga\big) $   that
 $ E_{\oP^t_{\oga,\oo}} \big[  \int_{\oT \land  \oga }^\oT   g_i  (r, \oX_{r \land \cd},\oU_r) dr    \big] \\
 \= \big(\oY^i_{\n \oP}  (\oga) \big) (\oo) $
 \if{0}

 \beas
 E_{\oP^t_{\oga,\oo}} \Big[  \int_{\oT \land  \oga }^\oT   g_i  (r, \oX_{r \land \cd}) dr    \Big]
 & \tn \= & \tn   E_{\oP^t_{\oga,\oo}} \Big[  \int_{\oT \land  \oga }^\oT   g^+_i  (r, \oX_{r \land \cd}) dr    \Big]
 \-  E_{\oP^t_{\oga,\oo}} \Big[  \int_{\oT \land  \oga }^\oT   g^-_i  (r, \oX_{r \land \cd}) dr    \Big] \\
  & \tn \= & \tn   E_\oP \Big[   \int_{\oT \land  \oga }^\oT    g^+_i (r,\oX_{r \land \cd} ) dr \Big| \cF^{\oW^t}_\oga \Big] (\oo)
 \-  E_\oP \Big[   \int_{\oT \land  \oga }^\oT    g^-_i (r,\oX_{r \land \cd} ) dr \Big| \cF^{\oW^t}_\oga \Big] (\oo)
 \= \big(\oY^i_{\n \oP}  (\oga) \big) (\oo)
 \eeas

 \fi
 and $  E_{\oP^t_{\oga,\oo}} \big[  \int_{\oT \land  \oga }^\oT   h_i  (r, \oX_{r \land \cd},\oU_r) dr    \big]  \= \big(\oZ^i_\oP  (\oga) \big) (\oo) $.
   Given $ \oo \ins   \Big( \ocN_0  \cup \ocN_3 \cup \ocN^{\,i}_{\n g,h}  \Big)^c \Cp \ocA_* $,
 \eqref{090520_11}, \eqref{Jan11_03}
 and  $\oP^t_{\oga,\oo}   \big\{ \oT \=   \wh{\tau}^\oo (\oW  ) \gs t_\oo \big\}   \= 1$
  imply
  $ \big( \oY^i_{\n \oP}  (\oga) \big)  (\oo)
    \=    E_{\oP^t_{\oga,\oo}} \big[  \int_{\oT \land  \oga }^\oT   g_i  (r, \oX_{r \land \cd},\oU_r) dr    \big]
  \= E_{\oP^t_{\oga,\oo}} \big[ \b1_{\Wtgo} \int_{\oT \land  t_\oo }^\oT   g_i (r, \oX_{r \land \cd},   \oU_r) dr    \big] \\
    \=     E_{\oP^t_{\oga,\oo}} \big[  \int_{  t_\oo }^{\oT  }  g_i (r, \oX_{r \land \cd},   \oU_r) dr    \big]
   $
 and similarly   $       E_{\oP^t_{\oga,\oo}} \big[  \int_{  t_\oo }^{\oT  } h_i (r, \oX_{r \land \cd},\oU_r) dr    \big] \= \big( \oZ^i_\oP  (\oga) \big) (\oo) $. Therefore,
   \bea   \label{062821_11b}
  \oP^t_{\oga,\oo} \ins \ocP_{ \oga(\oo),\oX_{\oga \land \cd}  (\oo)} \Big(  \big(\oY_{\n \oP}  (\oga )\big)    (\oo), \big(\oZ_\oP  (\oga )\big)    (\oo) \Big), \q \fa \oo \ins \ocA_*   \Cp  \ocN^c_* ,
  \eea
  where $ \ocN_* \df \ocN_0 \cup \ocN_1 \cup \ocN_2 \cup \ocN_3  \cup \ocN_{\n \mu} \cup \Big( \ccup{i \in \hN}{} \ocN^{\,i}_{\n g,h}  \Big) \ins \sN_\oP \big(\ocF^t_{\n \infty}\big) $.
 In particular, \eqref{062821_11}  holds for $\oP-$null set $\ocN \df \ocN_*  \cp   \big\{\oT \nne   \wh{\tau}(\oW )\big\} \cp \big( \{ \wh{\tau}(\oW) \gs \oga  \} \D \ocA_* \big) \ins \sN_\oP  \big(\sB(\oO) \big) $.  \qed

   \no {\bf Proof of Theorem \ref{thm_DPP1}:}
   \if{0}

 We first show  the  measurability of $ \b1_{\{\oT \ge \ogaP\}} \oV \big(  \ogaP ,\oX_{ \ogaP  \land \cd} ,    \oY_{\n \oP}  \big( \ogaP \big) , \oZ_\oP  \big( \ogaP \big)   \big) $ for each $\oP \ins \ocP_{t,\bx}(y,z)$ so that the right hand side of \eqref{020422_14}  is well-defined.

   Let $\oP \ins \ocP_{t,\bx}(y,z)$ and  simply denote $\ogaP$ by $\oga$.
   By (D1) of $ \ocP_{t,\bx} $,
  there is  a   $\hU-$valued, $ \bF^{W^t} -$predictable process
 $ \wh{\mu} \= \{\wh{\mu}_s\}_{s \in [t,\infty)} $   on $\O_0$ such that
   $\oP  \big\{   \oU_s   \= \wh{\mu}_s (\oW  ) \hb{ for a.e. }  s \ins (t,\infty) \big\} \= 1 $,
  where $\ol{\mu}_s \df \wh{\mu}_s(\oW)$, $\fa s \ins [t,\infty)$.
  Like in Part (2a) of the proof of Proposition \ref{prop_flow},
 we still set   $\ocN_{\n X} \df \big\{ \oo \ins \oO \n :   \oX_s  (\oo)   \nne  \osX^{t,\bx,\ol{\mu}}_s  (\oo)    \hb{ for some } s \ins [0,\infty) \big\}
       \ins \sN_\oP \big(\ocF^t_{\n \infty}\big)$
    and let $  \big\{ \oK^t_s \big\}_{s \in [t,\infty)}$ be  the  $\bF^{\oW^t}-$predictable  process
   such that $ \ocN_{\n K} \df  \big\{ \oo \ins \oO \n :   \oK^t_s  (\oo)   \nne \osX^{t,\bx,\ol{\mu}}_s  (\oo)   $ for some $s \ins [t,\infty)  \big\}  \ins  \sN_\oP \big(\cF^{\oW^t}_\infty\big)$.
  Since $ \oX|_{[0,t)} \= \bx|_{[0,t)}$   and
   $ \oX|_{[t,\infty)} \= \oK^t  $ on $ \big(\ocN_{\n X} \cp \ocN_{\n K}\big)^c $, one can deduce that
    the random variable $ \oX_{\oga \land \cd} \n : \oO \mto \OmX $ is   $ \si \big( \cF^{\oW^t}_\oga \cp  \sN_\oP  (\ocF^t_{\n \infty} ) \big) \big/   \sB(\OmX)-$measurable.

 Set $\breve{\O} \df [0,\infty) \ti \OmX \times \Re \times \Re \n \supset \n D_\ocP$.
  Let $\wh{\tau}$ be the $[t,\infty]-$valued   $ \bF^{W^t,P_0} -$stopping time with
 $ \oP \big\{  \oT  \=   \wh{\tau} (\oW)   \big\} \= 1$.
  By \eqref{020322_11} and \eqref{062821_11b}, there exist
   $ \ocA_*   \ins \cF^{\oW^t}_\oga$ and $\ocN_* \ins \sN_\oP \big(\ocF^t_{\n \infty}\big)   $ such that   $ \Big( \oga (\oo) ,\oX_{\oga \land \cd}  (\oo) ,   \big(\oY_\oP  (\oga )\big)  (\oo),    \big(\oZ_\oP  (\oga )\big)    (\oo) \Big) \ins D_\ocP $ for any $ \oo \in \ocA_* \cap   \ocN^c_* $. By the measurability of $ \oX_{\oga \land \cd}$,
\bea  \label{020622_23}
\breve{\Psi} (\oo) \df \b1_{\{\oo \in \ocA^c_* \cup   \ocN_*\}} (t,\bx,y,z) \+ \b1_{\{\oo \in \ocA_* \cap   \ocN^c_*\}} \Big( \oga (\oo) ,\oX_{\oga \land \cd}  (\oo) ,   \big(\oY_\oP  (\oga )\big)    (\oo), \big(\oZ_\oP  (\oga )\big)    (\oo) \Big) \ins D_\ocP , \q \fa \oo \ins \oO
\eea
 is  a $ \si \big( \cF^{\oW^t}_\oga \cp   \sN_\oP \big(\ocF^t_{\n \infty}\big) \big) \big/ \sB  (D_\ocP )  -$measurable random variable,
 which induces a probability measure
  $ \breve{P}   \df   \oP  \circ  \breve{\Psi}^{-1}  $ on $ \big( \breve{\O}, \sB  (\breve{\O} ) \big)$.
    Then $\breve{\Psi}$ is further $ \si \big(  \cF^{\oW^t}_\oga \cp   \sN_\oP  (\ocF^t_{\n \infty}   ) \big) \big/  \si \big( \sB  ( D_\ocP  ) \cp \sN_{\breve{P}}  \big( \sB   ( D_\ocP   )  \big)  \big) -$measurable.
  As the universally measurable function $ (t',\bx',   y',   z') \mto \oV (t',\bx',   y',   z')$  is $\si \big( \sB  ( D_\ocP ) \cp \sN_{\breve{P}} \big( \sB  ( D_\ocP ) \big)  \big) \big/ \sB[-\infty,\infty] -$measurable,  
  \bea  \label{020622_21}
  \breve{V}(\oo) \df \b1_{\{\oo \in \ocA_* \cap   \ocN^c_*\}} \oV \big( \breve{\Psi} (\oo) \big)
  \= \b1_{\{\oo \in \ocA_* \cap   \ocN^c_*\}} \oV \Big(  \oga (\oo) ,\oX_{\oga \land \cd}  (\oo) ,   \big(\oY_\oP  (\oga )\big)    (\oo), \big(\oZ_\oP  (\oga )\big)    (\oo) \Big) 
  , \q \fa \oo \ins \oO
  \eea
   is $ \si \big( \cF^{\oW^t}_\oga \cp   \sN_\oP \big(\ocF^t_{\n \infty}\big) \big) \big/ \sB[-\infty,\infty] -$measurable.
   We see from \eqref{020322_11}   that  $ (\ocA_* \cap   \ocN^c_*) \D \{\oT \gs \oga\} \sb \big(\ocA_* \D \{\oT \gs \oga\}\big) \cp \ocN_*
   \sb \big(\ocA_* \D \{\wh{\tau}(\oW) \gs \oga\}\big) \cup \{\oT \nne \wh{\tau}(\oW)\} \cup \ocN_* \n \ins \n \sN_\oP \big(\sB(\oO)\big)  $.
   It follows that   $  \b1_{\{\oT (\oo) \ge \oga (\oo)\}} \oV \big(  \oga (\oo) ,\oX_{\oga \land \cd}  (\oo) ,   \big(\oY_\oP  (\oga )\big)    (\oo),   \big(\oZ_\oP  (\oga )\big)    (\oo) \big) $, $  \oo \ins \oO$ is $ \si \big(   \cF^{\oW^t}_\oga \cp  \sN_\oP \big(\sB(\oO)\big) \big) \big/ \sB[-\infty,\infty] -$measurable.
  As $\oo \mto \big(\oR(t)\big)(\oo)$ is $\sB(\oO)-$measurable by  \eqref{082122_17},
 $ \b1_{\{\oT (\oo) < \oga (\oo) \}} \\ \big(\oR(t)\big) (\oo) $, $  \oo \ins \oO$ is also $ \sB(\oO)-$measurable.
 Hence, the right hand side of \eqref{020422_14} is   well-defined.

   \fi
 For any $[t,\infty)-$valued $\bF^{\oW^t}-$stopping time $\oz$, we denote  $\oR  (\oz) \df \int_{\oT \land \oz}^\oT   f   (r, \oX_{r \land \cd}, \oU_r  ) dr  \+ \b1_{\{\oT < \infty\}} \pi   \big(\oT, \oX_{\oT \land \cd}\big) $.


  \no {\bf (I) (sub-solution side)} Fix  $\oP \ins \ocP_{t,\bx}(y,z)$   and  simply denote $\ogaP$ by $\oga$.

  Let $\wh{\tau}$ be the $[t,\infty]-$valued   $ \bF^{W^t,P_0} -$stopping time with
 $ \oP \big\{  \oT  \=   \wh{\tau} (\oW)   \big\} \= 1$ and let  $ \ocA_*   \ins \cF^{\oW^t}_\oga$, $\ocN_* \ins \sN_\oP \big(\ocF^t_{\n \infty}\big)   $ be  as in \eqref{020322_11} and \eqref{062821_11b}.
  By (R2), 
  there is a $\ocN_{f,\pi}  \ins \sN_\oP \big(\cF^{\oW^t}_\oga\big)$ such that
 $  E_{\oP^t_{\oga,\oo}} \big[  \oR (\oga)      \big]
  \=  E_\oP \big[   \oR (\oga)    \big| \cF^{\oW^t}_\oga \big] (\oo) $  for any $\oo \ins   \ocN^c_{f,\pi}    $.
  For any $\oo \ins \ocA_* \Cp \big( \ocN_*  \cp \ocN_{f,\pi}   \big)^c$,
 as $\ocN_0 \sb   \ocN_*$,    \eqref{090520_11}, \eqref{Jan11_03} and \eqref{062821_11b} imply  that
 \beas 
 \hspace{-0.5cm}
  E_\oP \big[ \, \oR (\oga)   \big| \cF^{\oW^t}_\oga \big] (\oo)    \=  E_{\oP^t_{\oga,\oo}} \big[ \, \oR (\oga)   \big]
 \= E_{\oP^t_{\oga,\oo}} \big[ \b1_{\Wtzo} \oR (\oga(\oo))  \big]   \=  E_{\oP^t_{\oga,\oo}} \big[ \, \oR (\oga(\oo))   \big]
 \ls \oV \big(  \oga(\oo),\oX_{\oga \land \cd}(\oo) ,  \big(  \oY_{\n \oP}  (\oga ) \big)   (\oo), \big( \oZ_\oP  (\oga ) \big)   (\oo) \big) .
 \eeas

  Since $\b1_{\{\oT    \ge   \oga \}} \= \b1_{\{ \wh{\tau}(\oW)   \ge   \oga \}} \= \b1_{\ocA_*}$, $\oP-$a.s.  by \eqref{020322_11}
 and since $ \ocA_*    \ins \cF^{\oW^t}_\oga$,  
 the tower property renders that
   \beas
  && \hspace{-1.2cm}  E_\oP \Big[  \b1_{\{\oT    \ge   \oga \}} \oV \big( \oga , \oX_{\oga \land \cd} ,     \oY_{\n \oP}  (\oga ) , \oZ_\oP  (\oga )         \big) \Big]
 \= E_\oP \Big[ \b1_{\ocA_*} \oV \big( \oga , \oX_{\oga \land \cd} ,     \oY_{\n \oP}  (\oga ) , \oZ_\oP  (\oga )         \big) \Big]
    \gs      E_\oP \Big[ \b1_{\ocA_*} E_\oP \big[ \, \oR (\oga)   \big| \cF^{\oW^t}_\oga \big] \Big] \\
  &&  \=    E_\oP \Big[  E_\oP \big[ \b1_{\ocA_*} \oR (\oga)   \big| \cF^{\oW^t}_\oga \big] \Big]
     \=    E_\oP \big[  \b1_{\ocA_*} \oR (\oga)   \big]  \= E_\oP \big[   \b1_{\{\oT    \ge   \oga \}} \oR (\oga)  \big]   .
     \eeas
    \if{0}

     This implies that
     $ \b1_{\ocA_*}  \oV^- \big(  \oga ,\oX_{\oga \land \cd}  ,    \oY_{\n \oP}  (\oga )    ,   \oZ_\oP  (\oga )     \big)
     \ls \b1_{\ocA_*}  E_\oP \big[ \, \oR^- (\oga)   \big| \cF^{\oW^t}_\oga \big]  $
      and thus
     \beas
      && \hspace{-1.2cm}  E_\oP \Big[  \b1_{\{\oT    \ge   \oga \}}  \oV^- \big(  \oga ,\oX_{\oga \land \cd}  ,   \oY_{\n \oP}  (\oga )    ,   \oZ_\oP  (\oga )     \big) \Big]
    \=  E_\oP \Big[ \b1_{\ocA_*}  \oV^- \big(  \oga ,\oX_{\oga \land \cd}  ,    \oY_{\n \oP}  (\oga )     ,   \oZ_\oP  (\oga )     \big) \Big] \\
    && \ls   E_\oP \Big[  \b1_{\ocA_*}  E_\oP \big[ \, \oR^- (\oga)   \big| \cF^{\oW^t}_\oga \big]   \Big]
     \ls  E_\oP \big[ \, \b1_{\ocA_*}   \oR^- (\oga)     \big]
     \ls  E_\oP \big[ \, \oR^- (\oga)     \big] \ls E_{P_0} \Big[ \int_0^\infty \n  f^- (r,   X^{t,\bx,\mu^o}_{r \land \cd}, \mu^o_r )  dr \Big] \+ c_\pi \< \infty  .
     \eeas

    \fi
 It follows that
   $ E_\oP \big[ \, \oR(t)   \big]
    \ls  E_\oP \Big[ \b1_{\{\oT <  \oga \}} \oR(t) \+ \b1_{\{\oT  \ge \oga \}}  \Big( \n \int_t^{\oga } f  (r,\oX_{r \land \cd}, \oU_r   ) dr
   \+     \oV \big( \oga , \oX_{\oga \land \cd} ,     \oY_{\n \oP}  (\oga ) , \oZ_\oP  (\oga )         \big) \Big) \Big]  $.
  Letting $\oP$ vary over $\ocP_{t,\bx}(y,z)$ yields that
  $ \oV(t,\bx,y,z)   \=  \n   \Sup{\oP \in \ocP_{t,\bx}(y,z)} E_\oP \big[ \, \oR(t) \big]  \ls \n  \Sup{\oP \in \ocP_{t,\bx}(y,z)}  E_\oP \Big[
  \b1_{\{\oT < \ogaP \}} \Big( \n \int_t^\oT \n f(r,\oX_{r \land \cd}, \oU_r) dr   \+ \pi \big(\oT, \oX_{\oT \land \cd}\big) \Big)
      \+ \b1_{\{\oT \ge \ogaP \}} \Big( \n \int_t^{\ogaP} \n f(r,\oX_{r \land \cd}, \oU_r) dr  \+ \oV \big(   \ogaP ,\oX_{\ogaP \land \cd} ,    \oY_{\n \oP}   (\ogaP ),  \oZ_\oP   (\ogaP )   \big) \Big) \Big] $.

  \no {\bf (II) (super-solution side)}
  Let   $\oP \ins \ocP_{t,\bx}(y,z)$ and   simply denote $\oga_\oP$ by $\oga$.
   We shall  show that
  \bea \label{081720_15}
   E_\oP \bigg[
  \b1_{\{\oT < \oga  \}} \oR(t) \+ \b1_{\{\oT \ge \oga  \}} \Big( \n \int_t^\oga  \n    f(r,\oX_{r \land \cd},\oU_r) dr  \+ \oV \big( \oga  ,\oX_{ \oga   \land \cd} ,    \oY_{\n \oP}   (\oga  ) , \oZ_\oP   (\oga  )  \big) \Big) \bigg] \ls  \oV (t,\bx,y,z) .
  \eea

  As $\cF^{\oW^t}_t \= \{\es,\oO\}$,   the $[t,\infty)-$valued  $\bF^{\oW^t} - $stopping time $\oga$ satisfies either $\{\oga \= t\} \= \oO $ or
 $ \{\oga \> t\} \= \oO $.

     Suppose first  that $\{\oga \= t\} \= \oO $:    for any $i \ins \hN$,
 $ 
 \oY^i_{\n \oP} (t)  \= E_\oP \big[ \int_{\oT \land t}^\oT    g_i (r,\oX_{r \land \cd},\oU_r ) dr \big| \cF^{\oW^t}_t \big] \= E_\oP \big[ \int_t^\oT    g_i(r,\oX_{r \land \cd},\oU_r ) dr   \big] \ls y_i$ and 
 $ 
 \oZ^i_\oP (t)
 \= E_\oP \big[ \int_t^\oT    h_i(r,\oX_{r \land \cd},\oU_r ) dr   \big] \= z_i$.
 \if{0}
 For any $(t,\bx,y,z) \ins D_\ocP$ and for any $y' \= \{y'_i\}_{i \in \hN} \ins \Re$ such that $ (t,\bx,y',z) \ins D_\ocP $ and that
 $y_i \ls y'_i$ for any $i \ins \hN$, it is clear that  
 $\oV (t,\bx,y,z) \ls \oV (t,\bx,y',z)$.
 \fi
  Then 
 $$  E_\oP \bigg[  \b1_{\{\oT < \oga \}} \oR(t)
  \+ \b1_{\{\oT \ge \oga \}} \Big(   \int_t^{\oga} \n f(r,\oX_{r \land \cd},\oU_r) dr  \+  \oV \big(  \oga ,\oX_{\oga \land \cd} ,    \oY_{\n \oP}  (\oga) , \oZ_\oP  (\oga)  \big)  \Big) \bigg]
  \n \= \n E_\oP \big[  \,   \oV \big( t   ,\oX_{t \land \cd} ,    \oY_{\n \oP}  ( t ) , \oZ_\oP  ( t )   \big)   \big]
    \ls \oV \big( t   , \bx  ,   y ,z  \big)    . $$  

 Let us assume  $ \big\{\oga \> t \big\} \= \oO $ in the rest of this proof.
  By (D1) of $ \ocP_{t,\bx} $, there exists   a   $\hU-$valued, $ \bF^{W^t} -$predictable process
 $ \wh{\mu} \= \{\wh{\mu}_s\}_{s \in [t,\infty)} $   on $\O_0$ such that
  the $\oP-$measure of 
 $ \oO_\mu \df \big\{   \oU_r   \= \wh{\mu}_r (\oW  ) \hb{ for a.e. }  r \ins (t,\infty) \big\}  $
  is equal to $1$.
    Since  $\ol{\mu}_s \df  \wh{\mu}_s (\oW) $, $ \fa s \ins [t,\infty)  $
  is a $\hU-$valued  $\bF^{\oW^t }-$predictable process on $\oO$,
  the $[0,1)-$valued $\bF^{\oW^t} -$adapted continuous process
  $ \oJ^{\, t}_s \df \int_t^s e^{-r}   \sI  ( \ol{\mu}_r   ) dr $, $ \fa s \ins [t,\infty)$ satisfies that
  for any  $\oo \ins \oO_\mu   $
  \bea \label{052320_25}
  \ol{\U}^t_s (\oo) \= \int_t^s e^{-r} \sI\big(\oU_r(\oo)\big) dr
  \= \int_t^s e^{-r} \sI\big(\ol{\mu}_r(\oo)\big) dr
  \= \oJ^{\, t}_s (\oo), \q \fa s \ins [t,\infty) .
  \eea
  Set $\ocN_{\n X} \df  \big\{ \oo \ins \oO \n :  \oX_s (\oo)  \nne  \osX^{t,\bx,\ol{\mu}}_s (\oo)   \hb{ for some } s \ins [0,\infty) \big\}  \ins \sN_\oP \big(\ocF^t_{\n \infty}\big) $.
   We also let $\wh{\tau}$ be the $[t,\infty]-$valued   $ \bF^{W^t,P_0} -$stopping time with
 $ \oP \big\{  \oT  \=   \wh{\tau} (\oW)   \big\} \= 1$ and let  $ \ocA_*   \ins \cF^{\oW^t}_\oga$, $\ocN_* \ins \sN_\oP \big(\ocF^t_{\n \infty}\big)   $ be  as in \eqref{020322_11} and \eqref{062821_11b}.

 \no {\bf II.a.1)}
  Define  mappings  $    \big( \oW^{t,\oga},\oU^{t,\oga}  \big) \n : \oO \mto \O_0 \ti \hJ  $   by
   \beas
  \big( \oW^{t,\oga}_r,\oU^{t,\oga}_r \big)  (\oo)  \df  \Big( \oW^t  \big(  (r \ve t) \ld   \oga(\oo)  , \oo \big)  , \b1_{\{r \in [0,t) \cup ( \oga(\oo),\infty)\}} u_0  \+ \b1_{\{r \in [t, \oga(\oo)]\}} \oU_r (   \oo  )  \Big)  , \q  \fa (r,\oo) \ins [0,\infty) \ti \oO .
 \eeas
  Clearly, $\oW^{t,\oga}$ is  $   \cF^{\oW^t}_\oga   \big/ \sB(\O_0) -$measurable.
  To show the measurability of $ \oU^{t,\oga} $,
  we let   $\vf \ins L^0 \big((0,\infty) \ti \hU;\hR\big)$.
     The $ \bF^{\oW^t } -$predictability of $\big\{\ol{\mu}_s\big\}_{s \in [t,\infty)}$ implies that
      $ \big\{ \vf \big(s,\ol{\mu}_s  \big) \big\}_{s \in [t,\infty)}$  is also an $ \bF^{\oW^t } -$predictable process and
  $\int_t^\oga   \vf  (s,\ol{\mu}_s) ds $ is  thus  $\cF^{\oW^t }_\oga -$measurable.
  Then
    $\oxi_\mu \df \int_0^t \vf(s,u_0)ds \+ \int_t^\oga   \vf  (s,\ol{\mu}_s) ds \+ \int_\oga^\infty \vf(s,u_0)ds $
  is an $\cF^{\oW^t }_\oga -$measurable random variable such that
  $ \oxi_\mu (\oo) \= \int_0^\infty \vf \big(s, \oU^{t,\oga}_s (\oo) \big) ds \= I_\vf \big(\oU^{t,\oga}  (\oo)\big) $
  for any $\oo \ins \oO_\mu$.
  Since  $ \oO_\mu \= \big\{\oo \ins \oO \n : \ol{\U}^t_s(\oo) \= \int_t^s e^{-r} \sI (\oU_r(\oo)) dr  \= \int_t^s e^{-r} \sI \big(\ol{\mu}_r (\oo)\big)  dr   , \, \fa s \ins (t,\infty) \big\}
   \= \underset{s \in \hQ \in (t,\infty)}{\cap} \big\{\oo \ins \oO \n :  \ol{\U}^t_s(\oo) \= \int_t^s e^{-r} \sI \big(\ol{\mu}_r (\oo)\big)  dr   \big\}
  \ins \cF^{\oW^t,\ol{\U}^t,\oP}_\infty $,
  it holds  for any $  \cE \ins \sB(\hR)$  that
 $ \big(\oU^{t,\oga} \big)^{-1} \big( (I_\vf)^{-1}(\cE) \big)
    \=     \big\{ \oo \ins \oO  \n : I_\vf \big(\oU^{t,\oga}  (\oo)\big) \ins  \cE  \big\}
     \=     \big\{ \oo \ins \oO_\mu  \n : \oxi_\mu (\oo) \ins  \cE  \big\} \cp   \big\{ \oo \ins \oO^c_\mu  \n : I_\vf \big(\oU^{t,\oga}  (\oo)\big) \ins  \cE  \big\}  \ins \si \Big( \cF^{\oW^t}_\oga \cp \sN_\oP \big(  \cF^{\oW^t,\ol{\U}^t}_\infty\big) \Big)  $.
   By Lemma \ref{lem_M29_01}  (1),  the sigma-field 
 $    \Big\{ A \sb \hJ \n : \big(\oU^{t,\oga} \big)^{-1}  (A) \ins \si \Big( \cF^{\oW^t}_\oga \cp \sN_\oP \big(  \cF^{\oW^t,\ol{\U}^t}_\infty\big) \Big) \Big\}$
 includes all generating sets of $ \sB(\hJ)$ and thus contains $ \sB(\hJ)$.
 Hence,   $ \oU^{t,\oga} $ is $\si \Big( \cF^{\oW^t}_\oga \cp \sN_\oP \big(  \cF^{\oW^t,\ol{\U}^t}_\infty\big) \Big)-$measurable.

 Since $\big\{\osX^{t,\bx,\ol{\mu}}_s \big\}_{s \in [t,\infty)} $ is  an $  \bF^{\oW^t,\oP}  -$adapted continuous   process,
 we can emulate   Lemma 2.4 of \cite{STZ_2011a}
 to construct an $\hR^l-$valued $\bF^{\oW^t}-$predictable  process $  \big\{ \oK^t_s \big\}_{s \in [t,\infty)}$
   such that $ \ocN_{\n K} \df  \big\{ \oo \ins \oO \n :   \oK^t_s  (\oo)   \nne \osX^{t,\bx,\ol{\mu}}_s  (\oo)   $ for some $s \ins [t,\infty)  \big\}  \ins  \sN_\oP \big(\cF^{\oW^t}_\infty\big)$.
  Since $ \oX|_{[0,t)} \= \bx|_{[0,t)}$   and
   $ \oX|_{[t,\infty)} \= \oK^t  $ on $ \big(\ocN_{\n X} \cp \ocN_{\n K}\big)^c $, one can deduce that
    the random variable $ \oX_{\oga \land \cd} \n : \oO \mto \OmX $ is   $ \si \big( \cF^{\oW^t}_\oga \cp  \sN_\oP  (\ocF^t_{\n \infty} ) \big) \big/   \sB(\OmX)-$measurable.

 \no {\bf II.a.2)}
 Set  $\ocG^t_s \df \cF^{\oW^t}_s \ve \cF^{\ol{\U}^t}_s \ve \cF^\oX_s \= \si \big(\oW^t_{\n r}; r \ins [t,s]\big) \ve \si \big(\ol{\U}^t_{\n r}; r \ins [t,s]\big) \ve \si \big(\oX_r;r \ins [0,s]\big)$, $\fa s \in [t,\infty)$ and set $\ddot{\O} \df [0,\infty) \ti   \O_0   \ti \hJ \ti \OmX   \ti \Re \ti \Re \n \supset \n \cD_\ocP $.
 We   arbitrarily pick     $(\bw,\bu)$ from $\O_0 \ti \hJ$.

  Since $(t,\bx,y,z) \ins D_\ocP$, Theorem \ref{thm_V=oV} and \eqref{062821_11b} show that
 $(t,\bw,\bu,\bx,y,z) \ins \cD_\ocP$ and  that  $   \Big(  \oga (\oo) ,  \oW^{t,\oga} (\oo), \oU^{t,\oga} (\oo), \\  \oX_{\oga \land \cd}  (\oo) ,  \big(\oY_\oP  (\oga )\big)    (\oo),   \big(\oZ_\oP  (\oga )\big)    (\oo) \Big)   \ins   \cD_\ocP $ for any $\oo \ins \ocA_* \Cp \ocN^c_*  $.
 By the measurability of random variables $   \oW^{t,\oga}   $, $ \oU^{t,\oga}   $ and $  \oX_{\oga \land \cd}   $ in Step (II.a.1),
 \bea
 \hspace{-0.5cm}  \ddot{\Psi} (\oo) \df \b1_{\{\oo \in \ocA^c_* \cup \ocN_*\}} (t,\bw,\bu,\bx,y,z) \+ \b1_{\{\oo \in \ocA_* \cap \ocN^c_*\}} \Big( \oga (\oo) , \oW^{t,\oga} (\oo) , \oU^{t,\oga} (\oo) ,   \oX_{\oga \land \cd} (\oo) ,
\big(\oY_{\n \oP}  (\oga ) \big)   (\oo), \big(\oZ_\oP  (\oga ) \big)   (\oo) \Big)
 \ins   \cD_\ocP  ,   \q     \label{020822_11}
\eea
 $\fa \oo \ins \oO$ is    $\si \Big( \cF^{\oW^t}_\oga \cp \sN_\oP \big(  \ocG^t_\infty\big) \Big)  \Big/ \sB  ( \cD_\ocP ) -$measurable,
 which induces a probability measure   $ \ddot{P}    \df   \oP  \circ \ddot{\Psi}^{-1}  $ on $ \big( \ddot{\O},\sB  ( \ddot{\O}  )\big)$.
  Then $\ddot{\Psi} $ is further $ \si \Big( \cF^{\oW^t}_\oga \cp \sN_\oP \big( \ocG^t_\infty\big) \Big)  \Big/ \si \big( \sB  ( \cD_\ocP  ) \cp \sN_{\ddot{P} }  ( \sB   ( \cD_\ocP  )  )  \big) -$measurable.

 \no {\bf II.b)} Fix   $\e \ins (0,1)$   through Part (II.f). For any $\oo \ins \oO$, we may denote   $ t_\oo \df \oga(\oo)$.

 \no {\bf II.b.1)}   According to Jankov-von Neumann  Theorem (Proposition 7.50 of \cite{Bertsekas_Shreve_1978}),  Corollary  \ref{cor_graph_ocP} and Theorem \ref{thm_V_usa},
  there exists an analytically measurable function $ \ol{\bQ}_\e \n : \cD_\ocP \mto \fP\big(\oO\big) $ such that for any  $(\ft,\fw,\fu,\fx,\fy,\fz) \ins \cD_\ocP $, $  \ol{\bQ}_\e (\ft,\fw,\fu,\fx,\fy,\fz) $ belongs to $   \ocP_{\ft,\fw,\fu,\fx}(\fy,\fz) $ and satisfies
\bea \label{081620_19}
\hspace{-1cm}
  E_{\ol{\bQ}_\e (\ft,\fw,\fu,\fx,\fy,\fz)} \big[ \, \oR(\ft)   \big] \gs
\left\{
\ba{ll}
\oV(\ft,\fw,\fu,\fx,\fy,\fz) \- \e , & \hb{ if } \oV(\ft,\fw,\fu,\fx,\fy,\fz) \< \infty ; \\
 1/\e ,  & \hb{ if } \oV(\ft,\fw,\fu,\fx,\fy,\fz) \= \infty .
 \ea
 \right.
\eea
   As $ \ol{\bQ}_\e$ is    universally measurable,
   it is also  $  \si \big( \sB  ( \cD_\ocP  ) \cup \sN_{\ddot{P}}  ( \sB  ( \cD_\ocP )  )  \big) \big/  \sB \big( \fP\big(\oO\big) \big) -$measurable and
  $ \oQ^\oo_\e   \df      \b1_{\{\oo \in \ocA^c_*  \cup   \ocN_* \}} \oP  \\ + \n \b1_{\{\oo \in \ocA_* \cap   \ocN^c_* \}}  \ol{\bQ}_\e \big( \ddot{\Psi} (\oo) \big) $, $ \fa  \oo \ins \oO $
    is thus $  \si \Big(  \cF^{\oW^t}_\oga \cp   \sN_\oP  \big( \ocG^t_\infty \big) \Big) \Big/ \sB\big(\fP\big(\oO\big)\big) -$measurable.

   Given a $[0,\infty]-$valued $\sB(\oO)-$measurable random variable $\ol{\phi}$,
   Proposition 7.25 of  \cite{Bertsekas_Shreve_1978} implies that
the mapping
$ \fP\big(\oO\big) \ni \oQ \mto E_\oQ  \big[ \, \ol{\phi} \, \big] $ is $\sB\big(\fP\big(\oO\big)\big)   -$measurable.
 The  
 measurability of $\big\{\oQ^\oo_\e\big\}_{\oo \in \oO} $ renders that
\bea \label{022022_23}
 \hb{the random variable  $ \oO \ni \oo \mto E_{\oQ^\oo_\e} \big[ \, \ol{\phi} \, \big] $
 is $ \si \Big( \cF^{\oW^t}_\oga \cp \sN_\oP \big(\ocG^t_\infty\big) \Big) -$measurable.}
 \eea

  \no {\bf II.b.2)}     Let $\oo \ins \ocA_* \Cp \ocN^c_*$.
    We know from \eqref{020822_11} that
    \bea
   \oQ^\oo_\e    \= \ol{\bQ}_\e \big( \ddot{\Psi} (\oo) \big)   \ins   \ocP_{\oga(\oo) , \oW^{t,\oga} (\oo) , \oU^{t,\oga} (\oo) ,  \oX_{\oga \land \cd} (\oo)} \Big( \big(\oY_{\n \oP}  (\oga ) \big)   (\oo), \big(\oZ_\oP  (\oga ) \big)   (\oo) \Big)   .    \label{April07_11}
  \eea

  Set $ \oO^t_{\oga,\oo}   \df    \big\{ \oo' \ins \oO \n: (\oW_{\n s},\oX_s) (\oo') \=  ( \oW^{t,\oga}_s, \oX_s    ) (\oo)   ,    \fa s \ins  [0,\oga(\oo) ]  ; \; \oU_s (\oo') \= \oU^{t,\oga}_s (\oo) \hb{ for a.e. } s \ins (0,\oga(\oo)) \big\} $ and $ \oTh^t_{\oga,\oo}
 \df \big\{ \oo' \ins \oO \n:  (\oW^t_s,\ol{\U}^t_s)  (\oo')  \=  (\oW^t_s,\ol{\U}^t_s)  (\oo)  , \fa s \ins \big[t, \oga(\oo)\big] ; \,
   \oX_s(\oo') \= \oX_s(\oo), \fa s \ins [0, \oga(\oo)]   \big\} $.
 Since   $   \oO^t_{\oga,\oo}    \sb       \big\{ \oo' \ins \oO \n:  \oW_{\n s}  (\oo') \=  0   ,
    \fa s \ins [0,t] ; \oW_{\n s}  (\oo') \=   \oW^t_{\n s}  (\oo)    ,
    \fa s \ins (t,\oga(\oo)] ; \oX_s(\oo') \= \oX_s(\oo) , \fa s \ins [0,\oga(\oo)]; \, \oU_s(\oo') \= \oU_s(\oo) \hb{ for a.e. } s \ins (t, \oga(\oo)) \big\}
      \sb  
     \oTh^t_{\oga,\oo}      \sb \Wtgo  $,
 we see from \eqref{April07_11}   that
    \bea \label{072820_15}
  \oQ^\oo_\e  \big( \oO^t_{\oga,\oo}  \big) \=    1 , \q \hb{and thus} \q  \oQ^\oo_\e \big( \Wtgo \big) \= \oQ^\oo_\e \big( \oTh^t_{\oga,\oo} \big) \=  1 .
   \eea

By \eqref{April07_11},  there is   a   $\hU-$valued, $ \bF^{W^{t_\oo}} -$predictable process
 $ \wh{\mu}^\oo \= \big\{\wh{\mu}^\oo_s 
\big\}_{s \in [t_\oo,\infty)} $   on $\O_0$ 
   such that the $\oQ^\oo_\e -$measure of
   $ \oO^\oo_\mu  \df   \big\{\oo' \ins \oO \n : \oU_s (\oo') \= \ol{\mu}^\oo_s  ( \oo')   \hb{ for a.e. } s \ins (t_\oo,\infty) \big\} $
    is equal to  $1$ with $  \ol{\mu}^\oo_s
 \df   \wh{\mu}^\oo_s  (\oW  ) $, $\fa s \ins [t_\oo,\infty)$ and that   $ \ocN^{\,\oo}_{\n X} \df   \big\{ \oo' \ins \oO \n :    \oX_s (\oo')  \nne \osX^\oo_s  (\oo')   \; \hb{ for some } s \ins [0,\infty) \big\} \ins \sN_{\oQ^\oo_\e} \big(\ocF^{t_\oo}_\infty \big)$, where
    $\Big\{ \osX^\oo_s \=  \osX^{t_\oo,\oX_{\oga \land \cd} (\oo),\ol{\mu}^\oo}_s \Big\}_{s \in  [0,\infty)}$
    is an $\Big\{\cF^{\oW^{t_\oo},\oQ^\oo_\e}_{s \vee t_\oo}\Big\}_{s \in [0,\infty)}-$adapted continuous   process
    that uniquely solves the following SDE with the open-loop control $\ol{\mu}^\oo$ on   $\big(\oO , \sB\big(\oO\big) , \oQ^\oo_\e\big)   $:
\beas
\osX_s   \=   \oX  \big(t_\oo,\oo\big) + \int_{t_\oo}^s b \big(r, \osX_{r \land \cd}, \ol{\mu}^\oo_r \big)dr
      \+ \int_{t_\oo}^s \si \big(r, \osX_{r \land \cd}, \ol{\mu}^\oo_r\big) d \oW_r, \q \fa s \ins [t_\oo,\infty)
\eeas
with initial condition $\osX_s \= \oX_s  ( \oo ) $, $\fa s \ins [0,t_\oo]$.

    Since  $ \big\{\ol{\mu}^\oo_s \= \wh{\mu}^\oo_s (\oW) \big\}_{s \in [t_\oo,\infty)} $
  is a $\hU-$valued, $\bF^{\oW^{t_\oo} }-$predictable process on $\oO$,
  the $[0,1)-$valued $\bF^{\oW^{t_\oo}} -$adapted continuous process
  $ \oJ^{\,\oo}_s \df \int_{t_\oo}^s e^{-r}   \sI  ( \ol{\mu}^\oo_r   ) dr $, $ \fa s \ins [t_\oo,\infty)$ satisfies that
  for any  $\oo' \ins \oO^\oo_\mu   $
  \bea    \label{090522_31}
  \ol{\U}^\oo_s (\oo') \= \int_{t_\oo}^s e^{-r} \sI\big(\oU_r(\oo')\big) dr
  \= \int_{t_\oo}^s e^{-r} \sI\big(\ol{\mu}^\oo_r(\oo')\big) dr
  \= \oJ^{\,\oo}_s (\oo'), \q \fa s \ins [t_\oo,\infty) .
  \eea
  Like Lemma 2.4 of \cite{STZ_2011a}, one can   construct   an  $\hR^l-$valued  $  \bF^{\oW^{t_\oo} }  -$predictable  process  $  \big\{ \oK^\oo_s \big\}_{s \in [t_\oo,\infty)}$
   such that $ \ocN^{\,\oo}_{\n K} \df  \big\{ \oo' \ins \oO \n :   \oK^\oo_s  (\oo')   \nne \osX^\oo_s  (\oo')   $ for some $s \ins [t_\oo,\infty)  \big\}  \ins  \sN_{\oQ^\oo_\e} \big(\cF^{\oW^{t_\oo}}_\infty\big)$.

  \if{0}

  By  (D4) of Definition \ref{def_ocP},
 there is   a $[t_\oo,\infty]-$valued   $ \bF^{W^{t_\oo},P_0} -$stopping time $\wh{\tau}_\loo$ with  
 \bea \label{022122_11}
  \oQ^\oo_\e \big( \big\{  \oT  \=  \wh{\tau}_\loo (\oW  )   \big\} \big) \= 1 .
  \eea

  \fi

  \no {\bf II.b.3)}    Let $\oA \ins \sB (\oO)$.  We claim that
  \bea \label{May06_25}
 \oQ^\oo_\e (\ocA \Cp \oA) \= \b1_{ \{\oo \in \ocA  \}}  \oQ^\oo_\e (\oA) , \q \fa \ocA \ins \ocG^t_\oga ,
   ~ \fa \oo \ins \ocA_* \Cp \ocN^c_* .
 \eea
 To see this, we  take  $    \ocA  \ins \ocG^t_{\n \oga} $.
 Let $\oo_1 \ins \ocA  \Cp \ocA_* \Cp   \ocN^c_*  $ and   set    $s_1 \df \oga(\oo_1)$.
 Since $  \ocA  \Cp \{\oga \ls s_1\} $ is an $ \ocG^t_{s_1}-$measurable set including $\oo_1 $,
 one can deduce that
 \bea
   \oTh^t_{\oga,\oo_1} & \tn \= & \tn  \big\{ \oo' \ins \oO \n:  (\oW^t_r,\ol{\U}^t_r)  (\oo')  \= (\oW^t_r,\ol{\U}^t_r) (\oo_1)  , \fa r \ins  [t, s_1] ;
   \nonumber \\
    & \tn  & \qq    \oX_r(\oo') \= \oX_r(\oo_1), \fa r \ins [0, s_1]  \big\}
    \hb{ is also contained in $ \ocA  \Cp \{\oga \ls s_1\} $. } \label{090920_15}
 \eea
 By \eqref{072820_15},  $ \oQ^{\oo_1}_\e \big(\ocA \big)   \= 1 $ and thus  $\oQ^{\oo_1}_\e \big( \ocA  \Cp \oA \big) \= \oQ^{\oo_1}_\e \big(   \oA \big) \= \b1_{ \{\oo_1 \in \ocA  \}}  \oQ^{\oo_1}_\e (\oA) $. We next let $\oo_2 \ins \ocA^c \Cp \ocA_* \Cp   \ocN^c_*   $ and   set $s_2  \df \oga(\oo_2)$.
As $  \ocA^c  \Cp \{\oga \ls s_2 \} $ is an $ \ocG^t_{s_2}-$measurable set including $\oo_2 $,
 $ \oTh^t_{\oga,\oo_2} \= \big\{ \oo' \ins \oO \n: (\oW^t_r,\ol{\U}^t_r) (\oo')  \= (\oW^t_r,\ol{\U}^t_r) (\oo_2)  , \fa r \ins  [t, s_2] ; \,
   \oX_r(\oo') \= \oX_r(\oo_2), \fa r \ins [0, s_2]  \big\}  $
 is also included in $ \ocA^c  \Cp \{\oga \ls s_2\} $.
 We correspondingly have  $ \oQ^{\oo_2}_\e (\ocA^c )   \= 1 $
 and thus $\oQ^{\oo_2}_\e \big( \ocA  \Cp \oA \big) \= 0 \= \b1_{ \{\oo_2 \in \ocA  \}}  \oQ^{\oo_2}_\e (\oA) $. 

 Consider a pasted probability measure $\oP_\e \ins \fP\big(\oO\big)$:
    \bea   \label{092520_11}
  \oP_\e (\oA)  \df  \oP \big(\ocA^c_* \Cp \oA \big)   \+
  \int_{\oo \in \ocA_*}   \oQ^\oo_\e (\oA ) \oP(d\oo)    ,
  \q \fa \oA \ins \sB (\oO) .
 \eea
 \if{0}

    Let $\oA \ins \sB \big(\oO\big)$. By Proposition 7.25 of  \cite{Bertsekas_Shreve_1978},
 the function  $\phi_\oA \big(\oQ\big) \df E_\oQ  [\b1_\oA ] \= \oQ (\oA) $, $ \fa  \oQ  \ins   \fP\big(\oO\big) $ is
 $ \sB\big(\fP\big(\oO\big)\big) \big/ \sB[0,1] -$ measurable.
 As  $ \big\{ \oQ^\oo_\e   \big\}_{ \oo \in  \oO } $ is  $  \si \big(  \cF^{\oW^t}_\oga \cp   \sN_\oP  ( \ocG^t_\infty ) \big) \big/ \sB\big(\fP\big(\oO\big)\big) -$measurable,     the random variable  $ \oO \ni \oo \mto \oQ^\oo_\e(\oA) \= \phi_\oA \big(\oQ^\oo_\e\big)   $
  is $ \si \big( \cF^{\oW^t}_\oga \cp \sN_\oP  ( \ocG^t_\infty ) \big) \big/ \sB[0,1] -$measurable.

   It is clear that
  $ \oP_\e (\es) \= \oP \big(\ocA^c_* \Cp \es  \big) \+ \int_{\oo \in \ocA_*} \oQ^\oo_\e (\es ) \oP(d\oo)
 \= \oP (\es) \+ E_\oP \big[ \b1_{\ocA_*} 0 \big] \= 0 $ and $ \oP_\e \big(\oO\big) \= \oP \big(\ocA^c_* \Cp \oO  \big) \+
 \int_{\oo \in \ocA_*} \oQ^\oo_\e \big(\oO\big) \oP(d\oo)
 \= \oP \big( \ocA^c_* \big)     \+ E_\oP \big[ \b1_{\ocA_*} \big]
 \= 1  $.
 Given   a disjoint sequence $\{\oA_i\}_{i \in \hN}$ in $ \sB \big(\oO\big) $, the monotone convergence theorem implies that
 \beas
    \oP_\e \Big( \underset{i \in \hN}{\cup}\oA_i \Big)
   & \tn  \=  & \tn   \oP \Big(\ocA^c_* \Cp \Big( \underset{i \in \hN}{\cup}\oA_i \Big)  \Big)
  \+ \int_{\oo \in \ocA_*}
  \oQ^\oo_\e \Big(\underset{i \in \hN}{\cup}\oA_i  \Big) \oP(d\oo)
    \=     \oP \Big(  \underset{i \in \hN}{\cup} \big( \ocA^c_* \Cp \oA_i \big)  \Big)
  \+ \int_{\oo \in \ocA_*} \underset{i \in \hN}{\sum}  \oQ^\oo_\e  \big( \oA_i  \big) \oP(d\oo) \\
    & \tn  \=  & \tn   \underset{i \in \hN}{\sum} \oP   \big( \ocA^c_* \Cp \oA_i \big)
   \+ \underset{i \in \hN}{\sum} \int_{\oo \in \ocA_*} \n \oQ^\oo_\e  \big( \oA_i  \big) \oP(d\oo)
    \=    \underset{i \in \hN}{\sum} \Big( \oP   \big( \ocA^c_* \Cp \oA_i \big)
   \+   \int_{\oo \in \ocA_*} \n \oQ^\oo_\e  \big( \oA_i  \big) \oP(d\oo) \Big)
   \= \underset{i \in \hN}{\sum} \, \oP_\e \big(\oA_i\big) .
 \eeas
 So  $\oP_\e$ is a probability measure  of   $\fP\big(\oO\big)$.

 \fi
 In particular, taking $\oA \= \oO$ in  \eqref{May06_25} renders that
 \bea
  \oP_\e (\ocA)   \=     \oP \big(\ocA^c_* \Cp \ocA \big)   \+
  \int_{\oo \in \ocA_*} \b1_{
  \{\oo \in \ocA\}  }   \oP(d\oo)   \= \oP  (\ocA)  , \q \fa   \ocA \ins \ocG^t_\oga  .   \label{052420_21}
  \eea

 In the next four parts, we demonstrate that   $\oP_\e $ also belongs to  $\ocP_{t,\bx}(y,z)$,
 i.e., the probability class $\ocP_{t,\bx}(y,z)$ is stable under the pasting \eqref{092520_11}.

 \ss \no {\bf II.c)}
   In this part, we demonstrate that    $ \oW^t $  is a 
 Brownian motion
   with respect to the filtration $\obG^t \df \big\{ \ocG^t_s \big\}_{s \in [t,\infty)}  $
 under $\oP_\e$.
  More precisely,   we need to verify that  for any $ t \ls s \< r \< \infty $ and $ \cE \ins  \sB(\hR^d) $
 \beas  
  \oP_\e  \big\{ \big(\oW^t_r   \- \oW^t_s \big)^{-1} (\cE) \Cp \oA \big\}
  \=  \oP_\e  \big\{ \big(\oW^t_r   \- \oW^t_s \big)^{-1} (\cE)   \big\}
  \oP_\e \big( \oA \big) \= \phi(r\-s,\cE) \oP_\e \big( \oA \big) , \q \fa \oA \ins \ocG^t_s    ,
 \eeas
 where $\phi(a,\cE) \df \big(2 \pi a\big)^{-d/2} \int_{z \in \cE} e^{ - \frac{z^2}{2 a}} dz $, $\fa  a  \ins (0,\infty) $.

 \ss \no {\bf II.c.1)}  Let $ t \ls s \< r \< \infty $, $ \cE \ins  \sB(\hR^d) $ and $ \oo \ins \ocA_* \Cp \ocN^c_* $.
  Since   $ \{\oga \gs r\} \ins \cF^{\oW^t}_{\oga \land r} \sb \ocG^t_\oga$ and
  since $ \{\oga \gs r\} \Cp (\oW^t_r   \- \oW^t_s  )^{-1} (\cE) \=
  \{\oga \gs r\} \Cp (\oW^t_{\oga \land r}   \- \oW^t_{\oga \land s}  )^{-1} (\cE)  \ins \cF^{\oW^t}_{\oga \land r} \sb \ocG^t_\oga $, \eqref{May06_25} implies that
  \bea
  \b1_{\{\oga (\oo) \ge r\}} \oQ^\oo_\e  \Big(  (\oW^t_r   \- \oW^t_s  )^{-1} (\cE)  \Big)
  \= \oQ^\oo_\e  \Big( \{\oga \gs r\} \Cp  (\oW^t_r   \- \oW^t_s  )^{-1} (\cE)  \Big)
  \= \b1_{\{\oga (\oo) \ge r\} \cap \big\{\oo \in  (\oW^t_r  - \oW^t_s  )^{-1} (\cE) \big\}} .  \label{063020_25}
  \eea

 \bul    If $  \oga(\oo) \ls s $, since $(\oW^t_r \- \oW^t_s)  (\oo')
 \= \oW_r ( \oo') \- \oW_s ( \oo') \= \oW^{t_\oo}_r ( \oo') \- \oW^{t_\oo}_s ( \oo')$, $\fa \oo' \ins \oO$
 and since $ \oW^{t_\oo} $  is a 
 Brownian motion with respect to the filtration $\bF^{\oW^{t_\oo}}$ under $\oQ^\oo_\e$
 by \eqref{April07_11},
  \bea \label{063020_23}
  \oQ^\oo_\e  \big(  (\oW^t_r   \- \oW^t_s  )^{-1} (\cE)  \big) \=
    \oQ^\oo_\e  \Big\{   \Big( \oW^{t_\oo}_r \- \oW^{t_\oo}_s \Big)^{-1} (\cE) \Big\}
     \=  \phi(r\-s,\cE) . 
 \eea

 \bul  We next suppose that $s \< \oga(\oo) \< r$ and
  set  $\cE_\oo \df \big\{\fx \- \oW^t_\oga  ( \oo ) \+ \oW^t_s ( \oo ) \n : \fx \ins \cE \big\} \ins \sB(\hR^d)$.
  For any $\oo' \ins \Wtgo $, $ (\oW^t_r   \- \oW^t_s  ) (\oo')   \ins \cE $ if and only if $   \oW^{t_\oo}_r  (  \oo' ) \= \oW_r( \oo') \- \oW( \oga(\oo),\oo')
   \=  \oW^t_r (\oo') \- \oW^t \big(\oga(\oo) ,\oo'\big) \=  \oW^t_r (\oo') \- \oW^t_s  (\oo') \- \oW^t_\oga  ( \oo ) \+ \oW^t_s  (\oo) \ins \cE_\oo $, which shows
  \bea \label{063020_17}
   (\oW^t_r   \- \oW^t_s  )^{-1} (\cE) \Cp \Wtgo \= \Big( \oW^{t_\oo}_r \Big)^{-1} \big( \cE_\oo \big) \Cp \Wtgo \, .
  \eea
  So \eqref{072820_15} gives that
 \bea \label{063020_21}
  \oQ^\oo_\e \big(  (\oW^t_r   \- \oW^t_s  )^{-1} (\cE)   \big)
 \=      \oQ^\oo_\e \Big\{ \Big(\oW^{t_\oo}_r \Big)^{-1} \big( \cE_\oo \big) \Big\}
 \=    \phi \big( r\- t_\oo ,\cE_\oo \big) . 
 \eea

  By  (R2), there exists $  \ocN_{\n s,r,\cE}   \ins \n \sN_\oP \big(\cF^{\oW^t}_\oga\big) $ such  that
 $  E_\oP \Big[ \b1_{(\oW^t_r   - \oW^t_s  )^{-1} (\cE)  } \Big| \cF^{\oW^t}_\oga \Big] (\oo)
    \=   E_{\oP^t_{\oga,\oo}} \Big[ \b1_{(\oW^t_r   - \oW^t_s  )^{-1} (\cE)  } \Big] $
    for any $\oo \ins \ocN^{\, c}_{\n s,r,\cE} $.   
  Given $ \oo \ins \{s \< \oga \< r\} \Cp \ocN^c_{s,r,\cE} \Cp \ocA_* \Cp \ocN^c_*   $,
  since  $ \oW^{t_\oo} $  is a 
    Brownian motion   with respect to the filration $ \bF^{\oW^{t_\oo}} $  under $\oP^t_{\oga,\oo}$ by \eqref{062821_11b},
 we can deduce from 
 \eqref{Jan11_03}, \eqref{063020_17} and \eqref{063020_21} that
 \bea
 \hspace{-0.7cm}
 E_\oP \Big[ \b1_{(\oW^t_r   - \oW^t_s  )^{-1} (\cE)  } \Big| \cF^{\oW^t}_\oga \Big] (\oo)
 & \tn \= & \tn \oP^t_{\oga,\oo} \big( (\oW^t_r \- \oW^t_s  )^{-1} (\cE)  \big)
   \=   \oP^t_{\oga,\oo} \big( (\oW^t_r \- \oW^t_s  )^{-1} (\cE) \Cp \Wtgo \big)
 \= \oP^t_{\oga,\oo} \Big\{ \Big(\oW^{t_\oo}_r\Big)^{-1} \big( \cE_\oo \big) \Cp \Wtgo \Big\} \nonumber \\
 & \tn  \= & \tn  \oP^t_{\oga,\oo} \Big\{ \Big(\oW^{t_\oo}_r \Big)^{-1} \big( \cE_\oo \big)   \Big\}
 \= \phi \big( r\- t_\oo ,\cE_\oo \big)
 \=  \oQ^\oo_\e \big(  (\oW^t_r   \- \oW^t_s  )^{-1} (\cE)   \big) .   \label{063020_27}
 \eea
 Since $\{s \< \oga  \< r \} \ins \cF^{\oW^t}_{\oga \land r}   \sb \cF^{\oW^t}_\oga$
 and since  $   \ocA_* \ins \cF^{\oW^t}_\oga  $, it follows that
 $ \oP \big( \{s \< \oga  \< r \} \Cp   (\oW^t_r   \- \oW^t_s  )^{-1} (\cE)   \Cp \ocA_* \big)
   \=   E_\oP \Big[ \b1_{\{s < \oga < r \} \cap \ocA_*}  E_\oP \big[ \b1_{(\oW^t_r   - \oW^t_s  )^{-1} (\cE)  } \big| \cF^{\oW^t}_\oga \big]  \Big]
   \=    \int_{\oo \in \oO} \b1_{\{s < \oga(\oo)< r \} \cap  \{\oo \in \ocA_* \}} \oQ^\oo_\e \big(  (\oW^t_r   \- \oW^t_s  )^{-1} (\cE)   \big) \oP(d\oo) $.
  Then we see from \eqref{063020_23} and \eqref{063020_25}     that
\beas
 && \hspace{-1cm} \oP_\e  \big( \big(\oW^t_r   \- \oW^t_s \big)^{-1} (\cE)    \big)
   \=   \oP  \big( \ocA^c_* \Cp \big(\oW^t_r   \- \oW^t_s \big)^{-1} (\cE)   \big)
 \+  \int_{\oo \in \ocA_*}   \oQ^\oo_\e \big(  (\oW^t_r   \- \oW^t_s)^{-1} (\cE)   \big)   \oP(d\oo) \\
  & &   \=   \oP  \big( \ocA^c_* \Cp \big(\oW^t_r   \- \oW^t_s \big)^{-1} (\cE)   \big)
  \+ \phi(r\-s,\cE)  \oP \big( \{\oga \ls s  \} \Cp \ocA_* \big)
  \+ \oP \big( \{ s \< \oga \< r   \}  \Cp   (\oW^t_r   \- \oW^t_s  )^{-1} (\cE)   \Cp \ocA_* \big) \\
 && \q  + \oP \big( \{  \oga \gs r  \}  \Cp   (\oW^t_r   \- \oW^t_s  )^{-1} (\cE)   \Cp \ocA_* \big) \\
  & &   \=   \oP  \big( \ocA^c_* \Cp \{  \oga \ls s  \} \Cp \big(\oW^t_r   \- \oW^t_s \big)^{-1} (\cE)   \big)
  \+ \phi(r\-s,\cE)  \oP \big( \{\oga \ls s  \} \Cp \ocA_* \big)
  \+ \oP \big( \{  \oga \> s  \}  \Cp   (\oW^t_r   \- \oW^t_s  )^{-1} (\cE)    \big)   .
\eeas
  Since $\oW^t$ is a 
  Brownian motion   under $\oP$
  and since $ \ocA^c_* \Cp \{  \oga \ls s  \}   \ins \cF^{\oW^t}_s $,
  \beas 
 \oP_\e  \big( \big(\oW^t_r   \- \oW^t_s \big)^{-1} (\cE)    \big)
 \= \phi(r\-s,\cE) \Big( \oP  \big( \ocA^c_* \Cp \{  \oga \ls s  \}    \big)
 \+   \oP \big( \{\oga \ls s  \} \Cp \ocA_* \big)  \+   \oP   \{  \oga \> s  \} \Big)
 = \phi(r\-s,\cE) .
  \eeas

\ss \no {\bf II.c.2)}  Let $ t \ls s \< r \< \infty $ and $ \cE \ins  \sB(\hR^d) $.
 We  now show that
 \bea \label{082922_11}
 \oP_\e  \Big( \big(\oW^t_r   \- \oW^t_s \big)^{-1} (\cE) \Cp \oA   \Big)
 \=  \phi(r\-s,\cE) \oP_\e  (    \oA    ) , \q \fa \oA \ins \ocG^t_s \, .
 \eea

 Let $\oA \ins \ocG^t_s$.  Since $   \oA \Cp \{\oga \> s\} \ins \ocG^t_\oga $,
  one can derive from \eqref{May06_25}, \eqref{063020_25}, \eqref{063020_27}, the tower property and \eqref{052420_21}   that
 \bea
&& \hspace{-1.5cm} \oP_\e  \Big( \big(\oW^t_r   \- \oW^t_s \big)^{-1} (\cE) \Cp \oA \Cp   \{\oga \> s\} \Big)
 \= \oP  \Big( \ocA^c_* \Cp \big(\oW^t_r   \- \oW^t_s \big)^{-1} (\cE) \Cp \oA \Cp   \{\oga \> s\} \Big)
 \+  \int_{\oo \in \ocA_*} \b1_{\{\oo \in   \oA \cap \{\oga > s\} \}} \oQ^\oo_\e \big(  (\oW^t_r   \- \oW^t_s)^{-1} (\cE)   \big)   \oP(d\oo) \nonumber \\
&& \= \oP  \Big( \ocA^c_* \Cp \big(\oW^t_r   \- \oW^t_s \big)^{-1} (\cE) \Cp \oA \Cp   \{\oga \> s\} \Big)
 \+ \oP  \Big( \ocA_* \Cp \big(\oW^t_r   \- \oW^t_s \big)^{-1} (\cE) \Cp \oA \Cp   \{\oga \ge r\} \Big) \nonumber \\
&& \q +  \int_{\oo \in \ocA_*} \b1_{\{\oo \in   \oA \cap \{s < \oga < r  \} \}} E_\oP \Big[ \b1_{(\oW^t_r   - \oW^t_s  )^{-1} (\cE)  } \Big| \cF^{\oW^t}_\oga \Big] (\oo)  \oP(d\oo)  \nonumber \\
&& \= \oP  \Big( \ocA^c_* \Cp \big(\oW^t_r   \- \oW^t_s \big)^{-1} (\cE) \Cp \oA \Cp   \{\oga \> s\} \Big)
 \+ \oP  \Big( \ocA_* \Cp \big(\oW^t_r   \- \oW^t_s \big)^{-1} (\cE) \Cp \oA \Cp   \{\oga \ge r\} \Big) \nonumber \\
&& \q +  E_\oP \bigg[  E_\oP \Big[ \b1_{\{  \ocA_* \cap \oA \cap \{s < \oga < r  \} \}} \b1_{(\oW^t_r   - \oW^t_s  )^{-1} (\cE)  } \Big| \cF^{\oW^t}_\oga \Big] \bigg]  \nonumber \\
& & \= \oP  \Big(   \big(\oW^t_r   \- \oW^t_s \big)^{-1} (\cE) \Cp \oA \Cp   \{\oga \> s\} \Big)  \=   \oP  \big(   \big(\oW^t_r   \- \oW^t_s \big)^{-1} (\cE)   \big) \ti \oP  \big(   \oA \Cp   \{\oga \> s\} \big)
\= \phi(r\-s,\cE) \oP_\e \big(    \oA \Cp   \{\oga \> s\} \big) , \label{063020_29}
 \eea
 where we used the   independence of
 $ \big(\oW^t_r   \- \oW^t_s \big)^{-1} (\cE) $ from $\cF^{\oW^t}_s $ under $\oP$ in the fifth equality above.
  This equality   directly verifies \eqref{082922_11} for the case $s \= t$ as we assume $\{\oga \> t\} \= \oO$.
 \if{0}

 Assume first that $ s \= t $. We let $0 \ls t_1 \ls \cds \ls t_n = t $, let $\{\cE_i\}^n_{i=1} \sb \sB(\hR^l)$
 and set   $ \oA \df \ccap{i=1}{n} \oX^{-1}_{t_i} (\cE_i)  \ins \cF^\oX_t \= \ocG^t_t \sb \ocG^t_\oga $.
 If $ \bx(t_i) \n \notin \n \cE_i $ for some $i \ins \{1,\cds,n\}$,
 since $\oP \{\oX_{t'} \= \bx(t'), \fa t' \ins [0,t]\} \= 1$ by (D2') of Remark \ref{rem_ocP},
 \eqref{052420_21} shows that $\oP_\e \big( \oA \big) \= \oP \big( \oA \big) \= 0$
 and that
 $  \oP_\e  \big( \big(\oW^t_r     \big)^{-1} (\cE) \Cp \oA   \big) \= \oP \big( \oA^c_* \Cp \big(\oW^t_r    \big)^{-1} (\cE) \Cp \oA \big)
\+ \int_{\oo \in \oA_*} \b1_{\{\oo \in  \oA \}} \oQ^\oo_\e \big( (\oW^t_r   )^{-1} (\cE)   \big)   \oP(d\oo) \= 0
\= \phi(r\-t,\cE)   \oP_\e (\oA)  $.
 On the other hand, if $ \bx(t_i) \ins \cE_i $ for all $i \ins \{1,\cds,n\}$,
 we see from \eqref{052420_21} and \eqref{070120_25} that $\oP_\e \big( \oA \big) \= \oP \big( \oA \big) \= 1$
 and that  $  \oP_\e  \big( \big(\oW^t_r     \big)^{-1} (\cE) \Cp \oA   \big) \= \oP \big( \oA^c_* \Cp \big(\oW^t_r    \big)^{-1} (\cE) \Cp \oA \big)
\+ \int_{\oo \in \oA_*} \b1_{\{\oo \in  \oA \}} \oQ^\oo_\e \big( (\oW^t_r   )^{-1} (\cE)   \big)   \oP(d\oo)
\= \oP \big( \oA^c_* \Cp \big(\oW^t_r    \big)^{-1} (\cE)   \big)
\+ \int_{\oo \in \oA_*}   \oQ^\oo_\e \big( (\oW^t_r   )^{-1} (\cE)   \big)   \oP(d\oo)
\= \oP_\e  \big( \big(\oW^t_r    \big)^{-1} (\cE)     \big)
\= \phi(r\-t,\cE) \oP_\e (\oA)  $.
So the Lambda-system $\big\{ \oA \ins \sB(\oO) \n :  \oP_\e  \big( \big(\oW^t_r     \big)^{-1} (\cE) \Cp \oA   \big) \= \phi(r\-t,\cE) \oP_\e (\oA) \big\}$ contains
the generating Pi-system $\Big\{ \ccap{i=1}{n} \oX^{-1}_{t_i} (\cE_i)  \n :  0 \ls t_1 \ls \cds \ls t_n = t , \, \{\cE_i\}^n_{i=1} \sb \sB(\hR^l)  \Big\} $ of $\cF^\oX_t   $.   Dynkin's Pi-Lambda Theorem implies that \eqref{082922_11}  holds for $s \= t$.

\fi

Next suppose  that $s \> t$. We let $ 0 \ls t_1 \ls \cds \ls t_n \ls t  $,
   $\{\cE^o_i\}^n_{i=1} \sb \sB(\hR^l)$ and set $\oA_X \df \ccap{i=1}{n} \oX^{-1}_{t_i} (\cE^o_i) \ins \cF^\oX_t \sb \ocG^t_\oga \,$.
 We also let $   t \= s_1 \< s_2 \< \cds \< s_{m-1} \< s_m \= s$ with $k \gs 2 $,
   let $\{ \cE_j \}^m_{j=1} \sb \sB(\hR^{d+1+l})$ and set $\oA_m \df \ccap{j=1}{m}    \big( \oW^t_{   s_j } , \ol{\U}^t_{   s_j }, \oX_{ s_j }       \big)^{-1} (\cE_j) \ins \ocG^t_s  $.
 Taking $\oA \= \oA_X \Cp \oA_m $ in \eqref{063020_29} renders that
 \bea \label{090522_37}
 \oP_\e  \Big( \big(\oW^t_r   \- \oW^t_s \big)^{-1} (\cE) \Cp \oA_X \Cp \oA_m \Cp   \{\oga \> s\} \Big)
 \= \phi(r\-s,\cE) \oP_\e \big(  \oA_X \Cp \oA_m \Cp   \{\oga \> s\} \big) .
 \eea

By \eqref{052320_25}, the set $\ocA \df \Big( \ccap{i=1}{n} \{\bx(t_i) \ins \cE^o_i\}\Big) \Cp \Big( \ccap{j=1}{m}   \big( \oW^t_{ s_j  } , \oJ^{\,t}_{ s_j  }, \oK^t_{ s_j  }  \big)^{-1} (\cE_j) \Big)
 \ins \cF^{\oW^t}_s$ satisfies $\oO_\mu \Cp \big(\ocN_{\n X} \cp \ocN_{\n K}\big)^c \Cp \oA_X \Cp \oA_m
    \=   \  \oO_\mu \Cp \big(\ocN_{\n X} \cp \ocN_{\n K}\big)^c \Cp \ocA$.
    Since $\oW^t$ is a 
 Brownian motion   under $\oP$
 and since $ \ocA^c_* \Cp   \{\oga \ls s\} \ins \cF^{\oW^t}_s $,
 \bea
 && \hspace{-1cm}\oP  \big( \ocA^c_* \Cp \big(\oW^t_r   \- \oW^t_s \big)^{-1} (\cE) \Cp  \oA_X \Cp \oA_m \Cp   \{\oga \ls s\} \big)
 \=  \oP  \big( \ocA^c_* \Cp \big(\oW^t_r   \- \oW^t_s \big)^{-1} (\cE) \Cp \ocA \Cp   \{\oga \ls s\} \big) \nonumber \\
 &&\= \oP  \big(   \big(\oW^t_r   \- \oW^t_s \big)^{-1} (\cE)  \big) \ti \oP  \big( \ocA^c_*   \Cp \ocA \Cp   \{\oga \ls s\} \big)
 \= \phi(r\-s,\cE) \oP  \big( \ocA^c_*   \Cp  \oA_X \Cp \oA_m \Cp   \{\oga \ls s\} \big) . \label{090522_35}
 \eea

 Fix $k \= 1, \cds, m \- 1$. We set $\oA^\oga_k \df \ccap{j=1}{k}    \big( \oW^t_{ \oga \land s_j } , \ol{\U}^t_{ \oga \land s_j }, \oX_{ \oga \land s_j }       \big)^{-1} (\cE_j) \ins \ocG^t_{\oga \land s_{m-1}} \sb \ocG^t_\oga  $ and let  $\oo \ins \{s_k \< \oga   \ls s_{k+1}\} \Cp \ocA_* \Cp \ocN^c_* $. By    \eqref{090520_11},
 \bea  \label{082922_43}
  \hspace{-1cm} \Wtgo   \Cp \Big(  \ccap{j=1}{k}   \big( \oW^t_{ s_j  } , \ol{\U}^t_{ s_j  }, \oX_{ s_j  }  \big)^{-1} (\cE_j)  \Big)
      \= \Wtgo   \Cp \Big(  \ccap{j=1}{k}   \big( \oW^t_{ \oga \land s_j  } , \ol{\U}^t_{ \oga \land s_j }, \oX_{ \oga \land s_j }    \big)^{-1} (\cE_j)  \Big)    \= \Wtgo   \Cp \oA^\oga_k .
 \eea
  We   set   $\fra_\loo \df \big( \oW^t_\oga ( \oo ), \ol{\U}^t_\oga ( \oo ) , \bz \big) \ins \hR^{d+1+l} $  and define
 $ \ocA^\oo_k  \df   \ccap{j=k+1}{m}    \Big(  \oW^{t_\oo}_{s_j  }, \oJ^{\,\oo}_{   s_j }, \oK^\oo_{   s_j }  \Big)^{-1} \big( \cE_{j,\oo} \big)
\ins \cF^{\oW^{t_\oo}}_s  $,
  where  $\cE_{j,\oo} \df \big\{\fx \- \fra_\loo \n : \fx \ins \cE_j \big\} \ins \sB(\hR^{d+1+l})$.
 For   $j \= k\+1,\cds,m$ and $\oo' \ins \oTh^t_{\oga,\oo} \Cp \oO^\oo_\mu \Cp \big(\ocN^{\,\oo}_{\n X} \cp \ocN^{\,\oo}_{\n K}\big)^c $,
 one can deduce from \eqref{090522_31} that     $ \big( \oW^t_{ s_j }  , \ol{\U}^t_{ s_j  }   , \oX_{ s_j  }   \big) (\oo') \ins  \cE_j $ if and only if
  $  \big( \oW^{t_\oo}_{s_j  }, \oJ^{\,\oo}_{   s_j }, \oK^\oo_{   s_j } \big) (\oo')
   \= \big( \oW^{t_\oo}_{s_j  }, \ol{\U}^\oo_{   s_j }, \oX_{   s_j } \big) (\oo')
   \= \big( \oW^t_{s_j  } (\oo') \- \oW^t_{t_\oo} (\oo'), \ol{\U}^t_{s_j} (\oo') \- \ol{\U}^t_{t_\oo} (\oo'), \oX_{   s_j } (\oo') \big)
   \= \big( \oW^t_{s_j  } , \ol{\U}^t_{s_j}   , \oX_{   s_j }  \big) (\oo')
   \- \fra_\loo \ins \cE_{j,\oo} $.
  So $  \Big(  \ccap{j=k+1}{m}   \big( \oW^t_{ s_j  }  , \ol{\U}^t_{ s_j  }   , \oX_{ s_j  } \big)^{-1} (\cE_j)  \Big) \Cp \oTh^t_{\oga,\oo} \Cp \oO^\oo_\mu \Cp \big(\ocN^{\,\oo}_{\n X} \cp \ocN^{\,\oo}_{\n K}\big)^c
   \=  \\  \ocA^\oo_k  \Cp \oTh^t_{\oga,\oo} \Cp \oO^\oo_\mu \Cp \big(\ocN^{\,\oo}_{\n X} \cp \ocN^{\,\oo}_{\n K}\big)^c $,
   which together with \eqref{082922_43} shows that  $ \oA_m \Cp \oTh^t_{\oga,\oo} \Cp \oO^\oo_\mu \Cp \big(\ocN^{\,\oo}_{\n X} \cp \ocN^{\,\oo}_{\n K}\big)^c \= \oA^\oga_k \Cp  \ocA^\oo_k \Cp \oTh^t_{\oga,\oo}  \Cp \oO^\oo_\mu \Cp \big(\ocN^{\,\oo}_{\n X} \cp \ocN^{\,\oo}_{\n K}\big)^c $.

   As $ \oW^{t_\oo} $  is a 
   Brownian motion with respect to the filtration  $\bF^{\oW^{t_\oo}}$ under $\oQ^\oo_\e$ by \eqref{April07_11},
   we can deduce from \eqref{072820_15}   and \eqref{May06_25}   that
  \beas
   && \hspace{-1cm} \oQ^\oo_\e  \Big(   (\oW^t_r   \- \oW^t_s  )^{-1} (\cE) \Cp \oA_X \Cp \oA_m \Big) \=
 \oQ^\oo_\e  \Big\{ \oA_X \Cp \oA^\oga_k \Cp \ocA^\oo_k \Cp (\oW^t_r   \- \oW^t_s  )^{-1} (\cE) \Big\}
 \= \b1_{\{\oo \in \oA_X \cap \oA^\oga_k\}}  \oQ^\oo_\e  \Big\{   \ocA^\oo_k  \Cp   \Big( \oW^{t_\oo}_r \- \oW^{t_\oo}_s \Big)^{-1} (\cE) \Big\} \nonumber \\
 && \hspace{-0.4cm}  \= \b1_{\{\oo \in \oA_X \cap \oA^\oga_k\}} \oQ^\oo_\e  \big(  \ocA^\oo_k \big) \ti  \oQ^\oo_\e  \Big\{  \Big( \oW^{t_\oo}_r \- \oW^{t_\oo}_s \Big)^{-1} (\cE) \Big\}
 \=  \oQ^\oo_\e  \big( \oA_X \Cp \oA^\oga_k \Cp   \ocA^\oo_k    \big)  \phi(r\-s,\cE)
 \= \oQ^\oo_\e  \big(  \oA_X \Cp \oA_m \big) \phi(r\-s,\cE) ,    
 \eeas
 and thus $ \int_{\oo \in \ocA_*} \b1_{\{\oo \in   \{s_k < \oga \le s_{k+1}\} \}} \oQ^\oo_\e  \big(  (\oW^t_r   \- \oW^t_s  )^{-1} (\cE) \Cp \oA_X \Cp \oA_m   \big) \oP(d\oo)
   \= \phi(r\-s,\cE) \ti \int_{\oo \in \ocA_*} \b1_{\{\oo \in   \{s_k < \oga \le s_{k+1}\} \}} \oQ^\oo_\e  \big(  \oA_X \Cp \oA_m \big) \oP(d\oo)   $.
  Since $\{\oga \> 0\} \= \oO$ and since  $\{  \oga   \ls s \} \ins \cF^{\oW^t}_{\oga \land s}
  \sb \ocG^t_\oga  $,
 taking summation   from $k\=1$ through $k\=m \- 1$ and using \eqref{May06_25}  yield that
 \beas
&& \hspace{-1.5cm} \int_{\oo \in \ocA_*}   \oQ^\oo_\e  \big( \{  \oga   \ls s\} \Cp (\oW^t_r   \- \oW^t_s  )^{-1} (\cE) \Cp \oA_X \Cp \oA_m   \big) \oP(d\oo)
\= \int_{\oo \in \ocA_*} \b1_{\{\oo \in   \{  \oga \le s  \} \}} \oQ^\oo_\e  \big(   (\oW^t_r   \- \oW^t_s  )^{-1} (\cE) \Cp \oA_X \Cp \oA_m   \big) \oP(d\oo) \nonumber \\
&& \= \phi(r\-s,\cE) \ti \int_{\oo \in \ocA_*} \b1_{\{\oo \in   \{  \oga \le s  \} \}} \oQ^\oo_\e  \big(  \oA_X \Cp \oA_m   \big) \oP(d\oo)
 \= \phi(r\-s,\cE) \ti \int_{\oo \in \ocA_*}   \oQ^\oo_\e  \big( \{  \oga \ls s\} \Cp \oA_X \Cp \oA_m \big) \oP(d\oo) . 
 \eeas

 Adding it to \eqref{090522_37} and \eqref{090522_35} reaches   $ \oP_\e  \big( \big(\oW^t_r   \- \oW^t_s \big)^{-1} (\cE) \Cp  \oA_X \Cp \oA_m     \big) \= \phi(r\-s,\cE) \oP_\e  \big(   \oA_X \Cp \oA_m     \big) $.
  So   the Lambda-system 
 $ \ol{\L}^t_{s,r} \df \Big\{\oA \ins \sB\big(\oO\big) \n : \oP_\e  \big( \big(\oW^t_r   \- \oW^t_s \big)^{-1} (\cE) \Cp \oA   \big)
 \= \phi(r\-s,\cE) \oP_\e  \big(   \oA   \big) \Big\}$
  contains the Pi-system
 $ \Big\{   \Big(\ccap{i=1}{n} \oX^{-1}_{t_i} (\cE^o_i) \Big) \Cp \Big( \ccap{j=1}{m}    \big( \oW^t_{   s_j } , \ol{\U}^t_{   s_j }, \oX_{ s_j }       \big)^{-1} (\cE_j) \Big)    \n :  0 \ls t_1 \ls \cds \ls t_n \ls t \= s_1 \< s_2 \< \cds \< s_{m-1} \< s_m \= s  , \,
    \{\cE^o_i\}^n_{i=1} \sb \sB(\hR^l)  , \, \{ \cE_j \}^m_{j=1} \sb \sB(\hR^{d+1+l})   \Big\} $,
 which generates $\ocG^t_s$.
 In light of  Dynkin's Pi-Lambda Theorem, we obtain $\ocG^t_s \sb \ol{\L}^t_{s,r}$, proving \eqref{082922_11}.

  Hence,  $ \oW^t $  is a 
  Brownian motion with respect to the filtration $\obG^t$  under $\oP_\e$.
  Namely,  $ \oP_\e $ satisfies (D2) of $\ocP_{t,\bx}$.

 \ss \no {\bf II.d)}   In this part, we demonstrate that  $ \oP_\e $ satisfies (D1) of $\ocP_{t,\bx}$.

 \ss \no {\bf II.d.1)}  For any $s \ins [t,\infty)$, there is a $[0,1]-$valued $\cF^{W^t}_s-$measurable random variable $\U^\e_s$ on $\O_0$
 such that
 \beas 
\U^\e_s \big( \oW(\oo) \big) \=  E_{\oP_\e}  \Big[  \ol{\U}^t_s \big| \cF^{\oW^t}_s \Big] ( \oo ),
\q \fa   \oo  \ins   \oO .
\eeas
   Since $ \oW^t $  is a 
 Brownian motion with respect to the filtration $\obG^t$   under $\oP_\e$  by Part (II.c),
 applying  Lemma \ref{lem_122921_11} with $t_0 \= t$, $(\O_1, \cF_1, P_1,B^1)   \= \big(\oO ,  \sB(\oO ),    \oP_\e, \oW \big) $, $(\O_2,  \cF_2, P_2,B^2) \= \big(\O_0,  \sB(\O_0),  P_0,  W \big) $  and $\Phi \= \oW$   yields that
  $\big\{ \U^\e_s (\oW) \big\}_{s \in [t,\infty)}$ is an $\bF^{\oW^t}-$adapted process   and that
    $E_{P_0} \big[\U^\e_s\big] \= E_{\oP_\e} \big[ \U^\e_s (\oW) \big]
 \= E_{\oP_\e} \big[ \, \ol{\U}^t_s  \big] $ is right-continuous in $s \ins [t,\infty)$.
 As  $\bF^{W^t,P_0}$ is a right-continuous complete filtration,
the process $ \{\U^\e_s\}_{s \in [t,\infty)}$ admits a  \cad modification $ \big\{\wh{\U}^\e_s\big\}_{s \in [t,\infty)} $,
which  is a $[0,1]-$valued $\bF^{W^t,P_0}-$adapted process.

   As $ \oW^t $  is      a 
    Brownian motion both under      $\oP$
and under $\oP_\e$,  
  taking   $t_0 \= t$, $(\O_1, \cF_1, P_1,B^1)   \= \big(\oO ,  \sB(\oO ),   \oP \hb{ or } \oP_\e, \oW \big) $, $(\O_2, \cF_2, P_2,B^2) \= \big(\O_0,  \sB(\O_0),  P_0,  W \big) $  and $\Phi \= \oW$ in  Lemma \ref{lem_122921_11} renders that
\bea 
&&\hb{the process $\big\{\wh{\U}^\e_s (\oW )\big\}_{s \in [t,\infty)}$ is both
$\bF^{\oW^t,\oP}-$adapted and $\bF^{\oW^t,\oP_\e}-$adapted} \nonumber \\
 &&\hspace{-3.4cm} \hb{as well as that } \;
   \oW^{-1} (\cN) \ins \sN_\oP (\cF^{\oW^t}_\infty) \Cp \sN_{\oP_\e} (\cF^{\oW^t}_\infty) , ~ \; \fa \cN \ins \sN_{P_0} \big(\cF^{W^t}_\infty\big) .
  \label{070320_17}
 \eea

Given $s \ins [t,\infty)$,
since  $ \big\{ \wh{\U}^\e_s \nne \U^\e_s \big\} \ins \sN_{P_0}   \big(\cF^{W^t}_\infty\big)$,
  \eqref{070320_17} shows   that
 $ \wh{\cN}^\e_s \df  \big\{ \wh{\U}^\e_s (\oW ) \nne \U^\e_s (\oW ) \big\}
 \=     \oW^{-1}  \Big( \big\{ \wh{\U}^\e_s \nne \U^\e_s \big\} \Big)  \ins \sN_{\oP} \big(\cF^{\oW^t}_\infty\big) \Cp \sN_{\oP_\e} (\cF^{\oW^t}_\infty) $.
 Applying Lemma \ref{lem_Markov}    with $\big(\oP,\ocW,\fF_\cd ,\oxi \big) \= \big(\oP_\e , \{\oW^t_s\}_{s \in [t,\infty)}, \obG^t , \ol{\U}^t_s \big)$, we obtain
 \bea  \label{090522_51}
 \wh{\U}^\e_s (\oW )   \=   \U^\e_s (\oW ) \=  E_{\oP_\e}  \Big[  \ol{\U}^t_s \Big| \cF^{\oW^t}_s \Big]
 \=  E_{\oP_\e}  \Big[  \ol{\U}^t_s \Big| \cF^{\oW^t}_\infty \Big] , \q  \hb{$\oP_\e-$a.s.}
\eea

Let $\oz$ be a    $[t,\infty)-$valued $\bF^{\oW^t}-$stopping time, let $\oA \ins \cF^{\oW^t}_\oz$ and let $n \ins \hN$.
We set $s^n_i \df t \+ i 2^{-n} $, $\fa i \ins \hN \cp \{0\}$
and define $\oz_n \df \sum_{i \in \hN} s^n_i \b1_{\{\oz \in [s^n_{i-1},s^n_i)\}} $.
For any $i \ins \hN$, one  can deduce from \eqref{090522_51} and the monotone convergence theorem   that
 $ \wh{\U}^\e_{\oz_n} \n ( \oW  )   \=    \sum_{i \in \hN} \b1_{\{\oz \in [s^n_{i-1},s^n_i)\}} \wh{\U}^\e_{s^n_i} \big( \oW  \big)
\= \sum_{i \in \hN} \b1_{\{\oz \in [s^n_{i-1},s^n_i)\}}   E_{\oP_\e}  \Big[  \ol{\U}^t_{s^n_i} \big| \cF^{\oW^t}_\infty \Big]
\=  E_{\oP_\e}  \Big[ \sum_{i \in \hN} \b1_{\{\oz \in [s^n_{i-1},s^n_i)\}}  \ol{\U}^t_{s^n_i} \big| \cF^{\oW^t}_\infty \Big]
   \=    E_{\oP_\e}  \Big[   \ol{\U}^t_{\oz_n} \big| \cF^{\oW^t}_\infty \Big]  $,  $\oP_\e -$a.s.
  As  $\oz \= \lmtd{n \to \infty} \oz_n $,   the right-continuity of process $\wh{\U}^\e$, $\ol{\U}^t$
  and 
  the   bounded convergence theorem imply that
  \bea    \label{090422_37}
  \wh{\U}^\e_\oz  \big( \oW  \big) \= \lmt{n \to \infty} \wh{\U}^\e_{\oz_n} \n ( \oW  )
  \=  \lmt{n \to \infty}  E_{\oP_\e}  \Big[   \ol{\U}^t_{\oz_n} \big| \cF^{\oW^t}_\infty \Big]
  \=  E_{\oP_\e}  \Big[   \ol{\U}^t_\oz  \big| \cF^{\oW^t}_\infty \Big] , \q    \hb{$\oP_\e -$a.s.}
  \eea

As the right-continuous $\bF^{W^t,P_0}-$adapted  process $ \big\{\wh{\U}^\e_s\big\}_{s \in [t,\infty)} $ is $\bF^{W^t,P_0}-$optional,
  \beas
\wh{\sU}^\e_s (\o_0) \df e^s \linf{\d \to 0+} \frac{1}{\d} \Big( \wh{\U}^\e_s (\o_0)  \- \wh{\U}^\e_{(s-\d) \vee t} (\o_0) \Big) , \q \fa (s,\o_0) \ins [t,\infty) \ti \O_0
\eeas
is  a $[-\infty,\infty]-$valued, $ \bF^{W^t,P_0}  -$optional process.
 Similar  to   Lemma 2.4 of \cite{STZ_2011a},
 we can construct a  $[-\infty,\infty]-$valued $\bF^{W^t} -$predictable process
  $ \sU^\e \= \{\sU^\e_s\}_{s \in [t,\infty)} $ with
  $ \sU^\e_s(\o_0) \= \wh{\sU}^\e_s (\o_0) $ for $ds \ti d P_0-$a.s. $(s,\o_0) \ins [t,\infty) \ti \O_0$.
  \if{0}
  For any $n \ins \hN$, there exists an  $\hR  -$valued $\bF^{W^t} -$predictable process
  $ \sU^{\e,n} \= \{\sU^{\e,n}_s\}_{s \in [t,\infty)} $  such that
  $ \sU^{\e,n}_s(\o_0) \= \wh{\sU}^\e_s (\o_0) \ld n $ for $ds \ti d P_0-$a.s. $(s,\o_0) \ins [t,\infty) \ti \O_0$.
  Then the $\hR \cp \{\infty\}-$valued $\bF^W -$predictable process
  $ \sU^\e_s(\o_0) \df \linf{n \to \infty} \sU^{\e,n}_s(\o_0) $, $ (s,\o_0) \ins [t,\infty) \ti \O_0 $ satisfies that
  \beas
  \q \wh{\sU}^\e_s (\o_0) \= \lmtu{n \to \infty} \wh{\sU}^\e_s (\o_0) \ld n
  \= \lmtu{n \to \infty} \sU^{\e,n}_s(\o_0) 
  \=   \sU^\e_s(\o_0)   \q \hb{for $ds \ti d P_0-$a.s. } (s,\o_0) \ins [t,\infty) \ti \O_0 .
  \eeas
  \fi
 Using Fubini Theorem, one can find a $  \cN^\e_\sU \ins \sN_{P_0} \big(\cF^{W^t}_\infty\big)$ such that
  for any  $\o_0 \ins \big(\cN^\e_\sU\big)^c $,
   \bea \label{052420_11}
    \sU^\e_s(\o_0) \= \wh{\sU}^\e_s (\o_0) \q \hb{for a.e. } s \ins [t,\infty) .
    \eea
 By  \eqref{070320_17},  $\ocN^{\,\e}_{\n \sU} \df  \oW^{-1} \big(\cN^\e_\sU\big) \ins \sN_\oP  \big( \cF^{\oW^t}_\infty \big) \Cp \sN_{\oP_\e} \big( \cF^{\oW^t}_\infty \big)$.  
  We define a $\hU-$valued $ \bF^{W^t}   -$predictable process by
\beas
  \mu^\e_s (\o_0)   \df    \sI^{-1} \big( \sU^\e_s (\o_0)   \big) \b1_{\{\sU^\e_s (\o_0) \in \fE\}}
\+ u_0 \b1_{\{\sU^\e_s (\o_0) \notin \fE\}}  \ins \hU , \q \fa (s,\o_0) \ins [t,\infty) \ti \O_0 .       
\eeas

 \ss \no {\bf II.d.2)}
Let $s \ins [t,\infty)$ and let $\oA \ins \cF^{\oW^t}_\infty   $.
Since $ \ol{\U}^t_{\oga \land s} \ins  
\ocG^t_\oga   $
and    $ \oJ^{\,t}_{\oga \land s}  \ins \cF^{\oW^t}_{\oga \land s} \sb \ocG^t_\oga   $,
 \eqref{May06_25}  and  \eqref{052320_25} show that
 $ E_{\oP_\e} \big[ \b1_\oA \ol{\U}^t_{\oga \land s} \big]
   \=    E_\oP \Big[ \b1_{\ocA^c_*   \cap \oA} \ol{\U}^t_{\oga \land s}  \Big]   \+
   \int_{\oo \in \ocA_*}   \ol{\U}^t_{\oga \land s} (\oo) \oQ^\oo_\e \big(  \oA \, \big) \oP(d\oo)
    \=   E_\oP \big[ \b1_{\ocA^c_*  \cap \oA} \oJ^{\,t}_{\oga \land s}  \big]
    \+   \int_{\oo \in \ocA_*}  \oJ^{\,t}_{\oga \land s}  (\oo) \oQ^\oo_\e \big(  \oA \, \big) \oP(d\oo) \= E_{\oP_\e} \big[ \b1_\oA \oJ^{\,t}_{\oga \land s}  \big] $.
  Letting $ \oA $ varies over $  \cF^{\oW^t}_\infty $ and taking $\oz \= \oga \ld s $ in \eqref{090422_37}, we obtain
$ 
 \oJ^{\,t}_{\oga \land s} \= E_{\oP_\e} \big[\ol{\U}^t_{\oga \land s}\big| \cF^{\oW^t}_\infty \big]
\= \wh{\U}^\e_{\oga \land s} (\oW )  $,  $\oP_\e - $a.s.
   By the right-continuity of processes $\oJ^{\,t}$ and $\wh{\U}^\e$,
  it holds for all $ \oo \ins \oO$ except on
 a $ \ocN^{\,\e}_1 \ins  \sN_{\oP_\e} \big(\cF^{\oW^t}_\infty\big)$  that
 $ \wh{\U}^\e_s  \big(  \oW (\oo)\big) \= \oJ^{\,t}_s  (  \oo  ) $, $ \fa s \ins \big[t,\oga (\oo)\big] $.

 Let  $\oo \ins   \big( \ocN^{\,\e}_{\n \sU}  \cp  \ocN^{\,\e}_1\big)^c
 $. Since the  Lebesgue differentiation theorem yields that
 \beas
  \lmt{\d \to 0+} \frac{1 }{\d} \Big(\wh{\U}^\e_s  \big(  \oW (\oo)\big)  \- \wh{\U}^\e_{(s-\d) \vee t}  \big(  \oW (\oo)\big)\Big)
  \= \lmt{\d \to 0+} \frac{1 }{\d}\Big( \oJ^{\,t}_s   (  \oo )  \- \oJ^{\,t}_{(s-\d) \vee t}   (  \oo )\Big)
  \= \lmt{\d \to 0+} \, \frac{1 }{\d} \int_{(s-\d) \vee t}^s e^{-r} \sI \big( \ol{\mu}_r (\oo)   \big) dr
  \= e^{-s} \sI \big( \ol{\mu}_s  ( \oo  )   \big)
   \eeas
 for a.e. $ s \ins \big(t,\oga(\oo)\big) $, we see from \eqref{052420_11} that
 $ \sU^\e_s  \big( \oW (\oo) \big) \= \wh{\sU}^\e_s  \big( \oW (\oo) \big)
   \=    e^s \lmt{\d \to 0+} \, \frac{1}{\d}
 \Big(\wh{\U}^\e_s  \big(  \oW (\oo) \big)  \- \wh{\U}^\e_{(s-\d) \vee t} \big(  \oW (\oo)\big) \Big)
   \=   \sI \big( \ol{\mu}_s  ( \oo  )   \big) $ for a.e. $ s \ins \big(t,\oga(\oo)\big) $.  
  Thus,
  \bea \label{090422_41}
    \hAe_1    \df    \big\{\oo \ins \oO \n :  \mu^\e_s  \big(   \oW (\oo) \big)
  \=   \ol{\mu}_s  (  \oo  )   \hb{ for a.e. } s \ins (t,\oga (\oo)) \big\} \n \supset \n  \big( \ocN^{\,\e}_{\n \sU}  \cp  \ocN^{\,\e}_1\big)^c  .
  \eea

As   $\big\{\int_t^s  e^{-r} \sI \big( \mu^\e_r  \big)  dr\big\}_{s \in  [t,\infty)}$
is an  $\bF^{W^t}-$adapted process,
    Lemma \ref{lem_122921_11} (1) shows that
  $ \big\{ \int_t^s  e^{-r} \sI \big( \mu^\e_r (\oW   )\big)  dr \big\}_{s \in  [t,\infty)}$
is an $\bF^{\oW^t}-$adapted continuous process, which together with  the $\bF^{\oW^t}-$adaptedness of continuous process  $\oJ^{\,t}$ implies
 \bea
  \hAe_1  
     \=    \Big\{   \int_t^{\oga   \land s }  e^{-r}  \sI \big( \mu^\e_r \big(\oW  \big)\big)  dr
   \= \oJ^{\,t}_{\oga   \land s}   , \; \fa s \ins [t,\infty) \Big\}
     \=    \ccap{s \in \hQ \cap [t,\infty)}{} \Big\{ \int_t^{\oga   \land s }  e^{-r}  \sI \big( \mu^\e_r \big(\oW  \big)\big)  dr
   \= \oJ^{\,t}_{\oga   \land s}  \Big\} \ins \cF^{\oW^t }_\oga \n \sb \ocG^t_\oga \, . \q  \label{090422_43}
  \eea
 Then we can derive from \eqref{052420_21} and \eqref{090422_41}    that
$ 
 1   \=   \oP_\e \big(  ( \ocN^{\,\e}_{\n \sU}   \cp  \ocN^{\,\e}_1\big)^c  \big)
 \ls \oP_\e \big( \hAe_1 \big) \= \oP  \big( \hAe_1 \big) $.
 Similar to \eqref{090422_43},
  \beas
  \ocA^\e_1 \df  \big\{\oo \ins \oO \n :
   \mu^\e_s \big(  \oW (\oo) \big) \=  \oU_{\n s}  (  \oo  )      \hb{ for a.e. } s \ins (t,\oga(\oo)) \big\} \= \ccap{s \in \hQ \cap [t,\infty)}{} \Big\{   \ol{\U}^t_{\oga \land s }   \=   \int_t^{\oga  \land s }  e^{-r} \sI \big( \mu^\e_r  (\oW  )\big) dr   \Big\} \ins \ocG^t_\oga   \,   .
   \eeas
 \if{0}

Since
\bea
&& \hspace{-1.2cm} \ocA^\e_1 \df  \big\{\oo \ins \oO \n :     \oU \big(\oga (\oo) \ld s, \oo \big)
\=   \mu^\e \big(\oga (\oo) \ld s, \oW (\oo) \big)    ,  \hb{ for a.e. } s \ins \big[t,\infty) \big\}   \\
&& \=  \bigg\{ \oo \ins \oO \n : \int_t^{\oga (\oo) \land s }  e^{-r} \sI \big( \oU _r(\oo)\big) dr \=  \int_t^{\oga (\oo) \land s }  e^{-r} \sI \big( \mu^\e_r \big(\oW (\oo)\big)\big) dr , ~ \fa s \ins [t,\infty) \bigg\}  \nonumber \\
&& \= \bigg\{ \oo \ins \oO \n : \ol{\U}^t_{\oga \land s } (\oo) \=  \int_t^{\oga (\oo) \land s }  e^{-r} \sI \big( \mu^\e_r \big(\oW(\oo)\big)\big) dr , ~ \fa s \ins \hQ \Cp [t,\infty) \bigg\} \nonumber \\
&& \= \ccap{s \in \hQ \cap [t,\infty)}{} \bigg\{ \oo \ins \oO \n : \ol{\U}^t_{\oga \land s } (\oo) \=  \int_t^{\oga (\oo) \land s }  e^{-r} \sI \big( \mu^\e_r \big(\oW(\oo)\big)\big) dr  \bigg\} \ins \ocG^t_\oga  \,  .
\label{070720_21}
\eea

 \fi
 Applying \eqref{052420_21} again 
 renders  that
 \bea \label{090822_19}
 \oP_\e \big( \ocA^\e_1 \big) \= \oP \big( \ocA^\e_1 \big)
 \= \oP \big( \oO_\mu \Cp \ocA^\e_1 \big) \= \oP \big( \oO_\mu \Cp \hAe_1 \big)
 \= \oP \big( \hAe_1 \big) \= 1.
 \eea

 \ss \no {\bf II.d.3)}
Let $s \ins [t,\infty)$ and let $\oA \ins \cF^{\oW^t}_\infty   $.
 Since
 $ E_{\oP_\e} \big[ \b1_{\ocA^c_*   \cap \oA} \big( \ol{\U}^t_s \- \ol{\U}^t_{\oga \land s} \big) \big]
   \=    E_\oP \Big[ \b1_{\ocA^c_*   \cap \oA} \big( \ol{\U}^t_s \- \ol{\U}^t_{\oga \land s} \big)  \Big]
    \=   E_\oP \big[ \b1_{\ocA^c_*  \cap \oA} \big( \oJ^{\,t}_s \- \oJ^{\,t}_{\oga \land s} \big)   \big]
    \= E_{\oP_\e} \big[ \b1_{\ocA^c_*   \cap \oA} \big( \oJ^{\,t}_s \- \oJ^{\,t}_{\oga \land s} \big) \big] $
    by  \eqref{052320_25}, 
  letting $ \oA $ varies over $  \cF^{\oW^t}_\infty $ and using \eqref{090422_37} again yield  that
\beas
 \q \b1_{\ocA^c_*   } \big( \oJ^{\,t}_s \- \oJ^{\,t}_{\oga \land s} \big)
 \= E_{\oP_\e} \big[\b1_{\ocA^c_*   } \big( \ol{\U}^t_s \- \ol{\U}^t_{\oga \land s} \big)\big| \cF^{\oW^t}_\infty \big]
 \= \b1_{\ocA^c_*   } E_{\oP_\e} \big[ \big( \ol{\U}^t_s \- \ol{\U}^t_{\oga \land s} \big)\big| \cF^{\oW^t}_\infty \big]
\= \b1_{\ocA^c_*   } \big( \wh{\U}^\e_s (\oW ) \- \wh{\U}^\e_{\oga \land s} (\oW ) \big)  , \q \hb{$\oP_\e - $a.s.} \q
\eeas
  The right-continuity of processes $\oJ^{\,t}$ and $\wh{\U}^\e$ assures an  $ \ocN^{\,\e}_2 \ins  \sN_{\oP_\e} \big(\cF^{\oW^t}_\infty\big)$
  such that for any $\oo \ins \big(  \ocA_* \cp     \ocN^{\,\e}_2\big)^c$ and any $  s \ins [\oga(\oo),\infty)$,
  one has $ \wh{\U}^\e_s  \big(  \oW (\oo)\big) \- \wh{\U}^\e  \big( \oga (\oo)   ,\oW (\oo)\big)
 \= \oJ^{\,t}_s (\oo)\- \oJ^{\,t} \big( \oga (\oo)   , \oo \big) $.

 Let  $\oo \ins     \big( \ocA_* \cp  \ocN^{\,\e}_{\n \sU}  \cp  \ocN^{\,\e}_2\big)^c  $.
 Since the  Lebesgue differentiation theorem implies that
 $ \lmt{\d \to 0+} \frac{1 }{\d} \Big(\wh{\U}^\e_s  \big(  \oW (\oo)\big)  \- \wh{\U}^\e  \big( (  s \- \d) \ve \oga (\oo) ,\oW (\oo)\big)\Big)
 \= \lmt{\d \to 0+} \frac{1 }{\d}\Big( \oJ^{\,t}_s   (  \oo )  \- \oJ^{\,t}  \big( (  s \- \d) \ve \oga (\oo) ,\oo\big)\Big)
   \= \lmt{\d \to 0+} \, \frac{1 }{\d} \int_{(  s  - \d) \vee \oga (\oo)}^s  e^{-r} \sI \big( \ol{\mu}_r (\oo)   \big) dr
   \= e^{-s} \sI \big( \ol{\mu}_s  ( \oo  )   \big) $ for a.e. $ s \ins \big(\oga(\oo),\infty\big) $,
    \eqref{052420_11} shows that
 $ \sU^\e_s  \big(   \oW (\oo) \big) \= \wh{\sU}^\e_s  \big(  \oW (\oo) \big)
   \=    e^s  \lmt{\d \to 0+} \, \frac{1}{\d}
 \Big(\wh{\U}^\e_s  \big(  \oW (\oo) \big)  \- \wh{\U}^\e \big( (   s \- \d) \ve \oga (\oo) ,\oW (\oo)\big) \Big)
   \=  \sI \big( \ol{\mu}_s  ( \oo  )   \big) $ for a.e. $ s \ins \big(\oga(\oo),\infty\big) $.  
  So   $ \hAe_2    \df    \big\{\oo \ins \oO \n :  \mu^\e_s  \big(  \oW (\oo) \big)
  \=   \ol{\mu}_s  (  \oo  )   \hb{ for a.e. } s \ins \big(\oga(\oo),\infty\big) \big\} \n \supset \n \ocA^c_* \Cp  \big( \ocN^{\,\e}_{\n \sU}  \cp  \ocN^{\,\e}_2\big)^c  $. Set
 \bea   \label{090822_17}
 \ocA^\e_2 \df  \big\{\oo \ins \oO \n :   \mu^\e_s \big(  \oW (\oo) \big)  \=   \oU_{\n s}  (  \oo  )
        \hb{ for a.e. } s \ins \big(\oga(\oo) ,\infty\big) \big\} .
 \eea

 \if{0}

By an analogy to \eqref{090422_43},
 $ \hAe_2   \=    \Big\{   \int_{\oga   \land s }^s  e^{-r}  \sI \big( \mu^\e_r  (\oW  )\big)  dr
  \= \oJ^{\,t}_s    \-  \oJ^{\,t}_{\oga   \land s}    , ~ \fa s \ins [t,\infty) \Big\}
    \=    \ccap{s \in \hQ \cap [t,\infty)}{} \Big\{  \int_{\oga   \land s }^s  e^{-r}  \sI \big( \mu^\e_r  (\oW  )\big)  dr
  \= \oJ^{\,t}_s    \-  \oJ^{\,t}_{\oga   \land s}  \Big\} \ins \cF^{\oW^t }_\infty $
   and  $\ocA^\e_2 \df  \big\{\oo \ins \oO \n :   \mu^\e_s \big(  \oW (\oo) \big)  \=   \oU_{\n s}  (  \oo  )
       \hb{ for a.e. } s \ins \big(\oga(\oo) ,\infty\big) \big\} \= \ccap{s \in \hQ \cap [t,\infty)}{} \big\{  \ol{\U}^t_s \- \ol{\U}^t_{\oga \land s }   \=  \int_{\oga  \land s }^s  e^{-r} \sI \big( \mu^\e_r  (\oW )\big) dr  \big\} \ins \ocG^t_\infty $.

  \fi
 As $  \oP_\e \big( \ocA^c_* \big) \= \oP_\e \big( \ocA^c_* \Cp ( \ocN^{\,\e}_{\n \sU}   \cp  \ocN^{\,\e}_2\big)^c  \big)
 \ls \oP_\e \big( \ocA^c_* \Cp  \hAe_2  \big) \ls \oP_\e \big( \ocA^c_* \big)  $, we obtain that
 \bea  \label{090422_45}
 \oP \big( \ocA^c_* \big)   \=    \oP_\e \big( \ocA^c_* \big) \=   \oP_\e \big( \ocA^c_* \Cp  \hAe_2  \big) \= \oP  \big( \ocA^c_* \Cp \hAe_2  \big)  \= \oP \big( \oO_\mu \Cp \ocA^c_* \Cp \hAe_2 \big) \= \oP \big( \oO_\mu \Cp \ocA^c_* \Cp \ocA^\e_2 \big)
 \= \oP \big( \ocA^c_* \Cp \ocA^\e_2 \big)   .
 \eea

 \ss \no {\bf II.d.4)}
 We denote $\hQ_t \df \{t+q \n : q \ins \hQ_+\} $ and set   $\ocN^\e_{\n \U} \df \ccup{s \in  \hQ_t}{} \big\{ \wh{\U}^\e_s (\oW) \nne \U^\e_s (\oW)  \big\} \ins \sN_{\oP} \big(\cF^{\oW^t}_\infty\big) \Cp \sN_{\oP_\e} \big(\cF^{\oW^t}_\infty\big)$.

 Fix $s \ins \hQ_t$. Given $k \ins \hN$, we set $s^k_i \df t \+ i 2^{-k} (s\-t)$ for $i \= 0, 1, \cds, 2^k$.
 Since $\big\{\U^\e_r (\oW )\big\}_{r \in [t,\infty)}$ is a $[0,1]-$valued $\bF^{\oW^t}-$adapted process,
 the random variable
 $\oxi^{\e,k}_s  \df \b1_{\{\oga \ge s\}} \U^\e_s (\oW ) \+\sum^{2^k}_{i=1} \b1_{\{ s^k_{i-1}  \le \oga    < s^k_i \}} \U^{\phantom{\rule[-0.8mm]{0.45pt}{2.6mm}} \e}_{s^k_i} (\oW ) \ins [0,1] $ is $ \cF^{\oW^t}_s  -$measurable.
  Then   $ \oxi^\e_s \df \linf{k \to \infty} \oxi^{\e,k}_s$ is also a $[0,1]-$valued $ \cF^{\oW^t}_s  -$measurable random variable.

For any $\oo \ins \big( \ocN^\e_\U \big)^c  $, the right-continuity of process $\wh{\U}^\e$ shows that
 $ \wh{\U}^\e \big(\oga(\oo) \ld s ,\oW(\oo)\big)
\= \lmt{k \to \infty} \Big(  \b1_{\{\oga (\oo) \ge s\}} \wh{\U}^\e_s \big(\oW (\oo) \big) \+\sum^{2^k}_{i=1} \b1_{\{ s^k_{i-1}  \le \oga (\oo)   < s^k_i \}} \wh{\U}^\e_{s^k_i} \big(\oW (\oo) \big) \Big)
\= \lmt{k \to \infty} \oxi^{\e,k}_s  (\oo) 
\= \oxi^\e_s  (\oo)  $      
 and thus
 \bea
 \wh{\U}^\e_s \big(\oW(\oo)\big) \- \wh{\U}^\e \big(\oga(\oo) \ld s ,\oW(\oo)\big)
 \=  \U^\e_s  \big(\oW (\oo)\big) \- \oxi^\e_s (\oo)  .   \label{070320_23}
 \eea
 As $\ocN^{\,\e}_{\n \sU} \cp \ocN^\e_{\n \U}  \ins \n   \sN_{\oP_\e} \big(\cF^{\oW^t}_\infty\big) $,   there exists $\oA^\e_U \ins \cF^{\oW^t}_\infty  $ such that
$ \ocN^{\,\e}_{\n \sU} \cp \ocN^\e_{\n \U} \sb \oA^\e_U $ and $ \oP_\e \big(\oA^\e_U\big) \= 0 $, which implies
$0 \ls  \int_{\oo \in \ocA_*}   \oQ^\oo_\e \big(\oA^\e_U\big)   \oP(d\oo)   \ls \oP_\e \big(\oA^\e_U\big) \= 0 $.
Since the random variable $ \oO \ni \oo \mto \oQ^\oo_\e \big(\oA^\e_U\big)$ is $ \si \Big( \cF^{\oW^t}_\oga \cp \sN_\oP \big(\ocG^t_\infty\big) \Big)  -$measurable by  \eqref{022022_23},
  there exists  $\ofN^\e_U \ins \sN_\oP \big(\ocG^t_\infty\big)$
  such that for any $\oo \ins \ocA_* \Cp \big( \ocN_* \cp \ofN^\e_U \big)^c $
\bea \label{070620_21}
\oQ^\oo_\e \big(\oA^\e_U\big) \= 0 \q \hb{and thus} \q  \ocN^{\,\e}_{\n \sU} \cp \ocN^\e_{\n \U} \ins \sN_{\oQ^\oo_\e} \big(\cF^{\oW^t}_\infty\big) .
\eea

 Let $ \big\{\ol{O}_j\big\}_{j \in \hN} $ be a countable Pi-system 
 that generates $ \cF^{\oW^t}_s $. 
  Let $j \ins \hN$ and $\oA  \ins \cF^{\oW^t}_\oga$.
 One  can deduce from  \eqref{May06_25},  \eqref{070320_23} and \eqref{090422_37}  that
\beas
 && \hspace{-1.2cm} \int_{\oo \in \oO}  \b1_{\{\oo \in \ocA_* \cap \oA   \}}    E_{\oQ^\oo_\e} \Big[  \b1_{ \ol{O}_j}    \Big(   \U^\e_s   (\oW ) \- \oxi^\e_s \Big) \Big]  \oP(d\oo)
 \= E_{\oP_\e} \Big[ \b1_{\ocA_* \cap \oA \cap \ol{O}_j}   \Big(  \U^\e_s   (\oW ) \- \oxi^\e_s \Big) \Big]
  \=  E_{\oP_\e} \Big[   \b1_{\ocA_* \cap \oA \cap \ol{O}_j} \Big(\wh{\U}^\e_s (\oW ) \- \wh{\U}^\e_{\oga \land s} (\oW ) \Big) \Big]   \nonumber \\  && \=    E_{\oP_\e} \Big[  \b1_{\ocA_* \cap \oA \cap \ol{O}_j}    E_{\oP_\e} \Big[   \ol{\U}^t_s \- \ol{\U}^t_{\oga \land s}  \Big| \cF^{\oW^t}_\infty \Big]  \Big] \= E_{\oP_\e} \Big[  E_{\oP_\e} \Big[   \b1_{\ocA_* \cap \oA \cap \ol{O}_j} \big(\ol{\U}^t_s \- \ol{\U}^t_{\oga \land s}\big)   \Big| \cF^{\oW^t}_\infty \Big]  \Big]  \nonumber \\
 &&     \= E_{\oP_\e} \Big[  \b1_{\ocA_* \cap \oA \cap \ol{O}_j} \big(\ol{\U}^t_s \- \ol{\U}^t_{\oga \land s}\big)  \Big]
 \= \int_{\oo \in \oO}  \b1_{\{\oo \in \ocA_* \cap \oA   \}}    E_{\oQ^\oo_\e} \big[  \b1_{ \ol{O}_j}    \big(\ol{\U}^t_s \- \ol{\U}^t_{\oga \land s}\big) \big]  \oP(d\oo)    \, .
\eeas
 So the Lambda-system 
$  \ol{\L}^\e_{s,j}   \df   \Big\{\oA \ins   \sB_\oP(\oO) 
\n :  \int_{\oo \in \oO}  \b1_{\{\oo \in \ocA_* \cap \oA   \}}    E_{\oQ^\oo_\e} \big[  \b1_{ \ol{O}_j}    \big(  \U^\e_s   (\oW ) \- \oxi^\e_s \big) \big]  \oP(d\oo)
      \=   \int_{\oo \in \oO}  \b1_{\{\oo \in \ocA_* \cap \oA   \}}    E_{\oQ^\oo_\e} \big[  \b1_{ \ol{O}_j}  \\  \big(\ol{\U}^t_s \- \ol{\U}^t_{\oga \land s}\big) \big]  \oP(d\oo)  \Big\} $
contains $ \cF^{\oW^t}_\oga $ and $\sN_\oP\big(\sB (\oO )\big)$.
As $  \cF^{\oW^t}_\oga \cp \sN_\oP\big(\sB\big(\oO\big)\big)   $ is closed under intersection,
we know from Dynkin's Pi-Lambda Theorem that
$\si \big( \cF^{\oW^t}_\oga \cp \sN_\oP \big(\sB (\oO )\big) \big) \sb \ol{\L}^\e_{s,j}$ or
\bea
&& \hspace{-1.5cm} \int_{\oo \in \oO}  \b1_{\{\oo \in \ocA_* \cap \oA   \}}    E_{\oQ^\oo_\e} \Big[  \b1_{ \ol{O}_j}    \Big(  \U^\e_s   (\oW ) \- \oxi^\e_s \Big) \Big]  \oP(d\oo) \nonumber \\
 &&  \=   \int_{\oo \in \oO}  \b1_{\{\oo \in \ocA_* \cap \oA   \}}    E_{\oQ^\oo_\e} \big[  \b1_{ \ol{O}_j}    \big(\ol{\U}^t_s \- \ol{\U}^t_{\oga \land s}\big) \big]  \oP(d\oo), \q \fa \oA \ins \si \Big( \cF^{\oW^t}_\oga \cp \sN_\oP \big(\sB\big(\oO\big)\big) \Big) .   \label{070520_29}
\eea

  \if{0}

 By \eqref{070520_27}, applying Lemma \ref{lem_A1} (1) with $(\fY,\xi) \= \Big(\oO,  \b1_{ \ol{O}_j}    \big(  \U^\e_s   (\oW ) \- \oxi^\e_s \big)\Big) $ and $(\fY,\xi) \= \Big(\oO,  \b1_{ \ol{O}_j}  \big(\ol{\U}^t_s \- \ol{\U}^t_{\oga \land s}\big) \Big) $ respectively shows that the mapping
 \beas
 \ol{\phi}^{\l,\e,s}_1\big(\oQ\big) \df E_\oQ  \Big[  \b1_{ \ol{O}_j}    \Big(  \U^\e_s  (\oW ) \- \oxi^\e_s \Big) \Big] , \q \fa \oQ \ins \fP\big(\oO\big)
\eeas
 is $\sB\big(\fP\big(\oO\big)\big) \big/\sB[-\infty,\infty]-$measurable and the mapping
 \beas
 \ol{\phi}^{\l,\e,s}_2\big(\oQ\big) \df
     E_\oQ  \Big[  \b1_{ \ol{O}_j}    \big( \ol{\U}^t_s \- \ol{\U}^t_{\oga \land s} \big) \Big] , \q \fa \oQ \ins \fP\big(\oO\big)
\eeas
is $\sB\big(\fP\big(\oO\big)\big) \big/\sB[-1,1]-$measurable.

Then \eqref{April07_11} renders that
the random variable
\beas
 \oxi^{\e,s}_1 (\oo) \df   \ol{\phi}^{\l,\e,s}_1 \big( \oQ^\oo_\e  \big)
\= E_{\oQ^\oo_\e}  \Big[  \b1_{ \ol{O}_j}    \Big(  \U^\e_s   (\oW ) \- \oxi^\e_s \Big) \Big] ,
\q \fa \oo \ins \oO
\eeas
is $ \si \Big( \cF^{\oW^t}_\oga \cp \sN_\oP \big(\ocG^t_\infty\big) \Big) \Big/ \sB[-\infty,\infty]-$measurable
and the random variable
\beas
 \oxi^{\e,s}_2 (\oo) \df   \ol{\phi}^{\l,\e,s}_2 \big( \oQ^\oo_\e  \big)
\= E_{\oQ^\oo_\e}  \Big[  \b1_{ \ol{O}_j}    \big( \ol{\U}^t_s \- \ol{\U}^t_{\oga \land s} \big) \Big] ,
\q \fa \oo \ins \oO
\eeas
 is $ \si \Big( \cF^{\oW^t}_\oga \cp \sN_\oP \big(\ocG^t_\infty\big) \Big) \Big/ \sB[-1,1]-$measurable.

  \fi
 Taking $\phi \= \b1_{ \ol{O}_j}    \big(  \U^\e_s  (\oW ) \- \oxi^\e_s \big) $
and $\phi \= \b1_{ \ol{O}_j}    \big(\ol{\U}^t_s \- \ol{\U}^t_{\oga \land s}\big) $
respectively in \eqref{022022_23} shows that
$\oO \ni \oo \mto E_{\oQ^\oo_\e}  \big[  \b1_{ \ol{O}_j}    \big(  \U^\e_s   (\oW ) \- \oxi^\e_s \big) \big]$
and $\oO \ni \oo \mto E_{\oQ^\oo_\e}  \big[  \b1_{ \ol{O}_j}    \big( \ol{\U}^t_s \- \ol{\U}^t_{\oga \land s} \big) \big]$
are two $[-1,1]-$valued $ \si \Big( \cF^{\oW^t}_\oga \cp \sN_\oP \big(\ocG^t_\infty\big) \Big) -$measurable random variables.
  Letting $\oA$ run through $\si \Big( \cF^{\oW^t}_\oga \cp \sN_\oP \big(\ocG^t_\infty\big) \Big)$ in \eqref{070520_29},
  we can find an $ \ofN^\e_{s,j} \ins \sN_\oP \big(\ocG^t_\infty\big)$ such that for any $\oo \ins   \big(   \ofN^\e_{s,j} \big)^c$
 \bea
\b1_{ \{\oo \in \ocA_* \cap \ocN^c_* \}}    E_{\oQ^\oo_\e} \Big[  \b1_{ \ol{O}_j}    \Big(  \U^\e_s  (\oW ) \- \oxi^\e_s \Big) \Big]
\=  \b1_{ \{\oo \in \ocA_* \cap \ocN^c_* \}}    E_{\oQ^\oo_\e} \Big[  \b1_{ \ol{O}_j}    \big(\ol{\U}^t_s \- \ol{\U}^t_{\oga \land s}\big) \Big]   .   \label{070620_23}
\eea

 Set $  \ofN^\e_s   \df   \ccup{j \in \hN}{ }  \ofN^\e_{s,j}    \ins \sN_\oP \big(\ocG^t_\infty\big) $ and
 let $ \oo  \ins  \ocA_* \Cp \big( \ocN_* \cp \ofN^\e_U \cp \ofN^\e_s   \big)^c   $.
 Since it holds for any  $\oo' \ins \Wtgo \Cp \oO^\oo_\mu  $ that
      $ \ol{\U}^t_s(\oo') \-\ol{\U}^t \big(\oga(\oo') \ld s ,\oo'\big)
 \n \= \ol{\U}^t_s(\oo') \-\ol{\U}^t \big(\oga(\oo) \ld s ,\oo'\big)
 \n \= \n   \int_{\oga(\oo) \land s}^s \n e^{-r} \sI  \big( \oU_r  (   \oo' ) \big) dr
  \=   \int_{\oga(\oo)  }^{\oga(\oo) \vee s}  e^{-r} \sI  \big( \ol{\mu}^\oo_r  (   \oo' ) \big) dr
  \= \ol{\U}^\oo_{t_\oo \vee s}(\oo')   $ by \eqref{090520_11},       
    we can deduce from \eqref{070320_23}, \eqref{070620_21}, \eqref{070620_23}, \eqref{072820_15} and \eqref{090522_31}  that
    for any $j \ins \hN$
 \beas
    E_{\oQ^\oo_\e} \Big[  \b1_{ \ol{O}_j}    \Big( \wh{\U}^\e_s (\oW ) \- \wh{\U}^\e_{\oga \land s} (\oW ) \Big) \Big]
\n \=    E_{\oQ^\oo_\e} \Big[  \b1_{ \ol{O}_j}    \Big(  \U^\e_s  (\oW) \- \oxi^\e_s \Big) \Big]
\n \=     E_{\oQ^\oo_\e} \Big[  \b1_{ \ol{O}_j}    \big(\ol{\U}^t_s \- \ol{\U}^t_{\oga \land s}\big) \Big]
\n \=     E_{\oQ^\oo_\e} \Big[ \b1_{ \ol{O}_j} \ol{\U}^\oo_{ t_\oo \vee s} \Big]  \n \=   E_{\oQ^\oo_\e} \Big[ \b1_{ \ol{O}_j} \oJ^{\,\oo}_{ t_\oo \vee s} \Big] .
\eeas
 Then Dynkin's Pi-Lambda Theorem implies that  the Lambda-system 
 $   \Big\{ \cE \ins \sB_{\oQ^\oo_\e}\big(\oO\big)   \n : E_{\oQ^\oo_\e} \Big[  \b1_\cE    \Big( \wh{\U}^\e_s (\oW) \- \wh{\U}^\e_{\oga \land s} (\oW) \Big) \Big]
\=  E_{\oQ^\oo_\e} \Big[ \b1_\cE \oJ^{\,\oo}_{t_\oo \vee s} \Big] \Big\} $
  includes $ \cF^{\oW^t}_s $ and thus contains
  $ \cF^{\oW^t, \oQ^\oo_\e}_s \= \si \Big( \cF^{\oW^t}_s \cp \sN_{\oQ^\oo_\e} \big(\cF^{\oW^t}_\infty \big) \Big) $\,:
 \bea \label{070620_27}
E_{\oQ^\oo_\e}  \Big[ \b1_{\cE}   \Big(\wh{\U}^\e_s (\oW) \- \wh{\U}^\e_{\oga \land s} (\oW) \Big) \Big]
\=   E_{\oQ^\oo_\e} \Big[  \b1_{\cE} \oJ^{\,\oo}_{t_\oo \vee s} \Big] , \q \fa \cE \ins \cF^{\oW^t, \oQ^\oo_\e}_s .
 \eea

If $\oga(\oo) \gs s   $, then $\cF^{\oW^{t_\oo}}_{t_\oo \vee s} \= \cF^{\oW^{t_\oo}}_{t_\oo} \= \{\es, \oO\} \sb \cF^{\oW^t}_s $;  if $\oga(\oo) \< s   $, then $\cF^{\oW^{t_\oo}}_{t_\oo \vee s}
 \= \cF^{\oW^{t_\oo}}_s \=  \si \Big( \big(\oW^{t_\oo}_r\big)^{-1}(\cE) \n :  r \ins \big[\oga(\oo), s   \big], \cE \ins \sB(\hR^d) \Big)
  \= \si \Big(  \big(\oW^t_r \- \oW^t_{t_\oo}\big)^{-1}(\cE) \n :  r \ins \big[\oga(\oo),s   \big], \cE \ins \sB(\hR^d) \Big) \sb \cF^{\oW^t}_s $.
  In both cases, we see that $ \oJ^{\,\oo}_{t_\oo \vee s} \ins \cF^{\oW^{t_\oo}}_{t_\oo \vee s} \sb \cF^{\oW^t}_s $.
  Since $ \wh{\U}^\e_s  (\oW ) \- \wh{\U}^\e_{\oga \land s} (\oW )
 $ is $ \cF^{\oW^t,\oQ^\oo_\e}_s -$measurable by   \eqref{070320_23}  and \eqref{070620_21},
 letting $\cE$ run through $\cF^{\oW^t, \oQ^\oo_\e}_s$ in \eqref{070620_27},
 we can find some  $ \ofN^\e_{s,\oo}   \ins  \sN_{\oQ^\oo_\e}\big( \cF^{\oW^t}_\infty \big) $ such that
  \bea \label{070620_35}
   \wh{\U}^\e_s  \big( \oW(\oo')\big) \- \wh{\U}^\e  \big(\oga(\oo') \ld s,\oW(\oo')\big)
   \= \oJ^{\,\oo}_{t_\oo \vee s}  \big( \oo' \big)   , \q \fa \oo' \ins \big( \ofN^\e_{s,\oo} \big)^c .
 \eea

 Now, set $ \ofN^\e_\sharp \df \ccup{s \in \hQ_t}{}  \ofN^\e_s \= \ccup{s \in \hQ_t}{}  \ccup{j \in \hN}{ }  \ofN^\e_{s,j}  \ins \sN_\oP \big(\ocG^t_\infty\big) $ and fix  $ \oo  \ins  \ocA_* \Cp \big( \ocN_* \cp \ofN^\e_U \cp \ofN^\e_\sharp   \big)^c    $.
 We also set $ \ofN^\e_{\sharp,\oo} \df \ccup{s \in \hQ_t}{}  \ofN^\e_{s,\oo} \ins  \sN_{\oQ^\oo_\e}\big( \cF^{\oW^t}_\infty \big) $
 and let $\oo' \ins \Wtgo \Cp  \big(\ocN^{\,\e}_{\n \sU} \cp \ofN^\e_{\sharp,\oo} \big)^c   $.
  The right-continuity of    process $\wh{\U}^\e$, \eqref{090520_11} and \eqref{070620_35} render that
 $ \wh{\U}^\e_s  \big( \oW(\oo')\big) \- \wh{\U}^\e  \big(t_\oo  ,\oW(\oo')\big)
 \= \int_{t_\oo}^s  e^{-r} \sI \big( \ol{\mu}^\oo_r   \big(  \oo' \big)\big)  dr $, $ \fa s \ins  [t_\oo,\infty ) $.
  By the Lebesgue differentiation theorem,  we have
$  \lmt{\d \to 0+} \frac{1 }{\d}\Big(\wh{\U}^\e_s  \big(   \oW (\oo')\big)  \- \wh{\U}^\e  \big( (  s \- \d)  \ve t_\oo ,\oW (\oo')\big)\Big)
 \= \lmt{\d \to 0+} \,  \frac{1 }{\d} \int_{(  s  - \d)  \vee t_\oo}^s e^{-r} \sI \big( \ol{\mu}^\oo_r   \big(  \oo' \big)\big)   dr
   \= e^{-  s} \sI \big( \ol{\mu}^\oo_s   \big(  \oo' \big)\big)  $ for a.e. $ s \ins  (t_\oo,\infty ) $
 and thus
$ \sU^\e_s  \big( \oW (\oo') \big) \= \wh{\sU}^\e_s  \big( \oW (\oo') \big)
  \=   e^s \lmt{\d \to 0+} \, \frac{1}{\d}
 \Big(\wh{\U}^\e_s  \big(  \oW(\oo') \big)  \- \wh{\U}^\e \big( (s \- \d)  \ve t_\oo ,\oW(\oo')\big) \Big)
   \=   \sI \big( \ol{\mu}^\oo_s  \big(  \oo' \big)\big) $  for a.e. $ s \ins  (t_\oo,\infty ) $.
 It follows that $ \hAe_{2,\oo}  \df \big\{\oo' \ins \oO \n :    \mu^\e_s   \big(  \oW(\oo') \big)
   \=    \ol{\mu}^\oo_s   (  \oo'  )  \hb{ for a.e. } s \ins  (t_\oo ,   \infty ) \big\} \n \supset \n \Wtgo \Cp  \big(\ocN^{\,\e}_{\n \sU} \cp \ofN^\e_{\sharp,\oo}\big)^c   $.
  \if{0}

Since $\cF^{\oW^{t_\oo}}_s  \sb \cF^{\oW^t}_s$ for any $s \gs  t_\oo  $,
an analogy to \eqref{090422_43} shows that
 $ \hAe_{2,\oo}    \=    \Big\{ \oo' \ins \oO \n :   \int_{t_\oo}^s  e^{-r}  \sI \big( \mu^\e_r  (\oW (\oo') )\big)  dr
  \= \oJ^{\,\oo}_s(\oo')  , ~ \fa s \ins [t_\oo,\infty) \Big\}
     \=    \ccap{s \in \hQ \cap [t_\oo,\infty)}{} \Big\{  \oo' \ins \oO \n :   \int_{t_\oo}^s  e^{-r}  \sI \big( \mu^\e_r  (\oW (\oo') )\big)  dr
   \= \oJ^{\,\oo}_s(\oo') \Big\} \ins \cF^{\oW^t }_\infty $.

  \fi
 Then \eqref{072820_15} and \eqref{070620_21} show that
$  1 \= \oQ^\oo_\e \Big\{ \Wtgo \Cp \big(\ocN^{\,\e}_{\n \sU}\cp \ofN^\e_{\sharp,\oo}\big)^c   \Big\} \ls \oQ^\oo_\e \big(\hAe_{2,\oo}\big)$.
 Since $ \Wtgo \Cp \oO^\oo_\mu \Cp \ocA^\e_2
\= \Wtgo \Cp \oO^\oo_\mu \Cp   \big\{\oo' \ins \oO \n :   \mu^\e_s \big(  \oW (\oo') \big)  \=   \oU_{\n s}  (  \oo'  )
        \hb{ for a.e. } s \ins \big(t_\oo ,\infty\big) \big\}
        \= \Wtgo \Cp \oO^\oo_\mu \Cp \hAe_{2,\oo}  $ by  \eqref{090822_17}, we further see that
\bea \label{091022_17}
  \oQ^\oo_\e \big( \ocA^\e_2  \big) \=   \oQ^\oo_\e \big(\hAe_{2,\oo}\big)  \= 1 ,
  \q \fa  \oo  \ins  \ocA_* \Cp \big( \ocN_* \cp \ofN^\e_U \cp \ofN^\e_\sharp   \big)^c   ,
\eea
which together with \eqref{090822_19} and \eqref{090422_45} yields that
\beas
\oP_\e \big\{     \oU_s  \=   \mu^\e_s  (\oW  )     ,  \hb{ for a.e. } s \ins  (t,\infty) \big\}
\= \oP_\e \big(  \ocA^\e_1 \Cp \ocA^\e_2 \big) \= \oP_\e \big( \ocA^\e_2 \big)
\= \oP  \big( \ocA^c_* \Cp \ocA^\e_2 \big) \+ \int_{\oo \in \ocA_*} \oQ^\oo_\e \big( \ocA^\e_2 \big) \oP(d \oo)
\= \oP  \big( \ocA^c_* \big) \+ \oP \big( \ocA_* \big) \= 1 .
\eeas

 Hence, $ \oP_\e $ satisfies (D1) of $\ocP_{t,\bx}$. We set $\ol{\mu}^\e_s \df  \mu^\e_s  (\oW  )$, $\fa s \ins [t,\infty)$.

 \no {\bf II.e)}
 We next show that  $ \oP_\e $   satisfies (D2$'$) of $\ocP_{t,\bx}$ and thus satisfies  \(D3\)   of $\ocP_{t,\bx}$.

 \no {\bf II.e.1)} Set $ \oO_X \df \big\{    \oX_s    \= \bx(s) ,\,\fa s \ins [0,t ] \big\}$.
 We know from the proof of Proposition \ref{prop_flow}   that  $ \oO^c_X \sb \ocN_{\n X} \=   \big\{ \oo \ins \oO \n :  \oX_s (\oo)  \nne  \osX^{t,\bx,\ol{\mu}}_s (\oo)   \hb{ for some } s \ins [0,\infty) \big\} \sb \ocN_1 \sb   \ocN_* $.
 Given $\oo \ins \ocA_* \Cp \ocN^c_* \sb \oO_X$, 
 one has   $\oO^t_{\oga,\oo} \sb   \big\{ \oo' \ins \oO \n:   \oX_s(\oo') \= \oX_s(\oo) , \fa s \ins [0,t ] \big\}
\= \oO_X$ and \eqref{072820_15} implies that $ \oQ^\oo_\e \big(   \oO_X  \big) \= 1 $.
 As $\oP \big(   \oO_X  \big) \= 1$ by (D2$'$) of Remark \ref{rem_ocP}, one can deduce that
  $  \oP_\e \big( \oO_X \big)  \=  \oP \big(\ocA^c_* \Cp \oO_X \big)   \+
  \int_{\oo \in \ocA_* \cap \ocN^c_*}  1 \n \cd \n \oP(d\oo)  \= \oP \big(\ocA^c_* \big) \+ \oP(\ocA_*) \= 1   $.
  Applying Proposition \ref{prop_MPF1} with $(\O,\cF,P,B,X,\mu) \= \big(\oO,\sB(\oO),\oP_\e,\oW,\oX,\ol{\mu}^\e \big) $ renders that
   $ \big\{\oM^{t,\ol{\mu}^\e}_{s \land \otau^t_n}(\vf)\big\}_{s \in [t,\infty)} $
   is a bounded $\obF^t-$adapted process under $  \oP_\e  $.

\no {\bf II.e.2)}
  Fix $(\vf,n) \ins \fC (\hR^{d+l}) \ti \hN$ and let    $ \oo  \ins  \ocA_* \Cp \big( \ocN_* \cp \ofN^\e_U \cp \ofN^\e_\sharp   \big)^c    $.
  We  define a $C^2 (\hR^{d+l})  $ function
 $\vf_\loo(w,x) \df \vf \big(w \+ \oW^t_\oga  ( \oo ), x \big) $, $(w,x) \ins \hR^{d+l}$
 and   define an   $\obF^{t_\oo}-$stopping time $ \oz^n_\loo (\oo') \df \inf\big\{s  \ins [t_\oo,\infty) \n : \big|(\oW^{t_\oo}_s,\oX_s)  (\oo') \- \fra_{\overset{}{\oo}}   \big|   \gs n  \big\}  $, $  \oo' \ins \oO$
 with $\fra_{\overset{}{\oo}} \df   \big( \dn - \n \oW^t_\oga( \oo ) ,   \bz \big) \ins \hR^{d+l}$.
   For $i \= 0, 1,2 $ and $ \oo' \ins  \Wtgo \Cp \hAe_{2,\oo}$, since $D^i\vf  \big(  \oW^t_r (\oo'),\oX_r (\oo')\big)   \=   D^i\vf  \big( \oW^t_r(\oo') \- \oW^t  (\oga(\oo),\oo' )
  \+ \oW^t   (\oga(\oo), \oo ) , \oX_r (\oo') \big)
 \= D^i \vf_\loo \big( \oW^{t_\oo}_r ( \oo' ) ,\oX_r (\oo')\big) $, $\fa r \ins  [t_\oo,\infty)$,
 we obtain
 \bea
 && \hspace{-1.2cm} \big( \,\oM^{t,\ol{\mu}^\e}_{r_2} (\vf)\big)  (\oo'  ) \- \big(\,\oM^{t,\ol{\mu}^\e}_{r_1} (\vf)\big)  (\oo'  ) \nonumber \\
  && \hspace{-0.5cm} \=     - \n \int_{r_1}^{r_2}  \n  \ol{b}  \big( r, \oX_{r \land \cd}(\oo'),\ol{\mu}^\e_r (\oo') \big) \n \cd \n D \vf \big( \oW^t_{\n r}(\oo')  , \oX_r(\oo') \big) dr
    \-   \frac12 \n \int_{r_1}^{r_2} \n  \ol{\si} \, \ol{\si}^T  \big( r,  \oX_{r \land \cd}(\oo'),\ol{\mu}^\e_r (\oo') \big) \n : \n
    D^2 \vf  \big( \oW^t_{\n r}(\oo')  , \oX_r(\oo')  \big)   dr  \nonumber \\
  &&  \hspace{-0.5cm} \=    - \n \int_{r_1}^{r_2}  \n  \ol{b}  \big( r, \oX_{r \land \cd}(\oo'),\ol{\mu}^\oo_r (\oo') \big) \n \cd \n D \vf_\loo \big( \oW^{t_\oo}_{\n r}(\oo')  , \oX_r(\oo') \big) dr
    \-   \frac12 \n \int_{r_1}^{r_2} \n  \ol{\si} \, \ol{\si}^T  \big( r,  \oX_{r \land \cd}(\oo'),\ol{\mu}^\oo_r (\oo') \big) \n : \n
    D^2 \vf_\loo \big( \oW^{t_\oo}_{\n r}(\oo')  , \oX_r(\oo')  \big)   dr  \nonumber \\
  &&  \hspace{-0.5cm} \=  \Big( \,\oM^{t_\oo,\ol{\mu}^\oo}_{r_2} (\vf_\oo)\Big)  (\oo'  ) \- \Big( \,\oM^{t_\oo,\ol{\mu}^\oo}_{r_1} (\vf_\oo) \Big)  (\oo'  ) , \q \fa t_\oo \ls r_1 \ls r_2 \< \infty.   \label{091920_17}
 \eea

  Let   $t \ls s \< r \< \infty$,  $   \{(s_i,\cE_i )\}^k_{i=1} \sb   [t,s]  \ti \sB (\hR^{d+l}) $
 and set $\oA \df \ccap{i=1}{k}  (\oW^t_{  s_i  },\oX_{  s_i   })^{-1} (\cE_i) \ins \ocF^t_s $.
  If $t\+n  \ls s$, one directly has
  $   E_{\oP_\e} \big[   \big( \oM^{t,\ol{\mu}^\e}_{\otau^t_n \land r} (\vf )  \- \oM^{t,\ol{\mu}^\e}_{\otau^t_n \land s} (\vf ) \big)   \b1_{\oA}    \big] \=  E_{\oP_\e} \big[   \big( \oM^{t,\ol{\mu}^\e}_{\otau^t_n  } (\vf )  \- \oM^{t,\ol{\mu}^\e}_{\otau^t_n  } (\vf ) \big)   \b1_{\oA}    \big] \= 0 $ since $\otau^t_n \ls t \+n \ls s $.
 So we only need to verify $E_{\oP_\e} \big[   \big( \oM^{t,\ol{\mu}^\e}_{\otau^t_n \land r} (\vf )  \- \oM^{t,\ol{\mu}^\e}_{\otau^t_n \land s} (\vf ) \big)   \b1_{\oA}    \big]   \= 0 $ for the case ``$t\+n \> s $".

  \no {\bf (i)} Let $t\+n \> s$.   We first show that
   \bea \label{071720_31}
   E_{\oP_\e} \Big[ \b1_{\{\oga >  s\}}  \big( \oM^{t,\ol{\mu}^\e}_{\otau^t_n \land r} (\vf )  \- \oM^{t,\ol{\mu}^\e}_{\otau^t_n \land s} (\vf ) \big)   \b1_{\oA}    \Big]   \= 0 .
 \eea

 For any $ \oo  \ins  \hAe_1  $ and any $ t \ls r_1 \ls r_2 \ls \oga(\oo)$, we have
 $  \big( \,\oM^{t,\ol{\mu}^\e}_{r_2} (\vf)\big)  (\oo   ) \- \big(\,\oM^{t,\ol{\mu}^\e}_{r_1} (\vf)\big)  (\oo   ) 
   \= \big( \,\oM^{t,\ol{\mu} }_{r_2} (\vf)\big)  (\oo   ) \- \big(\,\oM^{t,\ol{\mu} }_{r_1} (\vf)\big)  (\oo   )   $.  
 So $ \b1_{ \wh{A}^\e_1   }  \big( \oM^{t,\ol{\mu}^\e}_{\otau^t_n \land \oga \land r} (\vf ) \- \oM^{t,\ol{\mu}^\e}_{\otau^t_n \land \oga \land s} (\vf ) \big)    \=   \b1_{ \wh{A}^\e_1   }  \big( \oM^{t,\ol{\mu} }_{\otau^t_n \land \oga \land r} (\vf ) \- \oM^{t,\ol{\mu} }_{\otau^t_n \land \oga \land s} (\vf ) \big) $
 and  \eqref{090822_19} renders that
 $ E_{\oP_\e} \Big[   \big( \oM^{t,\ol{\mu}^\e}_{\otau^t_n \land \oga \land r} (\vf ) \- \oM^{t,\ol{\mu}^\e}_{\otau^t_n \land s} (\vf ) \big)   \b1_{\{\oga  > s\}  \cap \oA}    \Big]
    \=    E_{\oP_\e} \Big[   \big( \oM^{t,\ol{\mu}^\e}_{\otau^t_n \land \oga \land r} (\vf ) \- \oM^{t,\ol{\mu}^\e}_{\otau^t_n \land \oga \land s} (\vf ) \big)   \b1_{\{\oga  > s\}  \cap \oA}     \Big]
     \=  E_{\oP_\e} \Big[    \big( \oM^{t,\ol{\mu} }_{\otau^t_n \land \oga \land r} (\vf ) \- \oM^{t,\ol{\mu} }_{\otau^t_n \land \oga \land s} (\vf ) \big)   \b1_{\{\oga  > s\}  \cap \oA}     \Big]$.
 Since $ \{\oga \> s\}  \Cp \oA \= \{\oga \> s\}  \Cp \Big( \ccap{i=1}{k}  (\oW^t_{ \oga \land  s_i   },\oX_{ \oga \land  s_i  })^{-1} (\cE_i) \Big) \n \ins \ocF^t_{\n \oga \land s} $ and
$     \oM^{t,\ol{\mu} }_{\otau^t_n \land \oga \land r} (\vf ) \- \oM^{t,\ol{\mu} }_{\otau^t_n \land \oga \land s} (\vf )   \ins
 \ocF^t_{\n \oga \land r}$,  
using \eqref{052420_21} and applying \eqref{020522_17} with
$( \fra,\oz_1,\oz_2)   \= \big(  \bz, \oga \ld s, \oga \ld r \big)$
 yield that
   \bea
   E_{\oP_\e} \Big[   \big( \oM^{t,\ol{\mu}^\e}_{\otau^t_n \land \oga \land r} (\vf ) \- \oM^{t,\ol{\mu}^\e}_{\otau^t_n \land s} (\vf ) \big)   \b1_{\{\oga  > s\}  \cap \oA}    \Big]
   \=    E_\oP  \Big[   \big( \oM^{t,\ol{\mu} }_{\otau^t_n \land \oga \land r} (\vf ) \- \oM^{t,\ol{\mu} }_{\otau^t_n \land \oga \land s} (\vf ) \big)   \b1_{\{\oga  > s\}  \cap \oA}     \Big]  \= 0  . \qq  \label{071420_21}
  \eea

 It also holds for any $ \oo  \ins  \hAe_2  $ and any $ \oga(\oo) \ls r_1 \ls r_2 \< \infty $ that
 \bea \label{091022_31}
  \big( \,\oM^{t,\ol{\mu}^\e}_{r_2} (\vf)\big)  (\oo   ) \- \big(\,\oM^{t,\ol{\mu}^\e}_{r_1} (\vf)\big)  (\oo   )  
    \= \big( \,\oM^{t,\ol{\mu} }_{r_2} (\vf)\big)  (\oo   ) \- \big(\,\oM^{t,\ol{\mu} }_{r_1} (\vf)\big)  (\oo   )  .
    \eea
    As $\oP \Big(\ocA^c_* \cap \big(\hAe_2\big)^c\Big) \= 0 $ by \eqref{090422_45},  taking $ \big(\fra,\oz_1,\oz_2 \big) \= \big(\bz, \oga, \oga \ve r \big)$ in \eqref{020522_17}, we obtain
   \beas
  \q  && \hspace{-1cm}  E_\oP \Big[ \b1_{\ocA^c_*  } \big( \, \oM^{t,\ol{\mu}^\e}_{\otau^t_n \land (\oga \vee r)} (\vf ) \- \oM^{t,\ol{\mu}^\e}_{\otau^t_n \land \oga} (\vf ) \big)   \b1_{\{\oga  > s\}  \cap \oA}     \Big]
  \=  E_\oP \Big[ \b1_{\ocA^c_* \cap \hAe_2 }    \big( \, \oM^{t,\ol{\mu}^\e}_{\otau^t_n \land (\oga \vee r)} (\vf ) \- \oM^{t,\ol{\mu}^\e}_{\otau^t_n \land \oga} (\vf ) \big)   \b1_{\{\oga  > s\}  \cap \oA}     \Big]   \\
  && \=  E_\oP \Big[ \b1_{\ocA^c_* \cap \hAe_2 } \big( \, \oM^{t,\ol{\mu} }_{\otau^t_n \land (\oga \vee r)} (\vf )
   \- \oM^{t,\ol{\mu} }_{\otau^t_n \land \oga} (\vf ) \big)   \b1_{\{\oga  > s\}  \cap \oA}     \Big]
   \=  E_\oP \Big[ \b1_{\ocA^c_*   }   \big( \, \oM^{t,\ol{\mu} }_{\otau^t_n \land (\oga \vee r)} (\vf )
   \- \oM^{t,\ol{\mu} }_{\otau^t_n \land \oga} (\vf ) \big)   \b1_{\{\oga  > s\}  \cap \oA}     \Big]     \= 0  .   \nonumber
   \eeas
 It follows  from  \eqref{May06_25}   that
  \bea
  \q   E_{\oP_\e} \Big[   \big( \, \oM^{t,\ol{\mu}^\e }_{\otau^t_n \land r} (\vf ) \- \oM^{t,\ol{\mu}^\e}_{\otau^t_n \land \oga \land r} (\vf ) \big)   \b1_{\{\oga  > s\}  \cap \oA}     \Big]
    \=    \int_{\oo \in \ocA_*} \dn  \b1_{ \{\oga (\oo) > s \}} \b1_{ \{\oo \in \oA\,\}}  E_{\oQ^\oo_\e} \n \big[ \,  \oM^{t,\ol{\mu}^\e }_{ \otau^t_n \land (\oga \vee r) } (\vf ) \- \oM^{t,\ol{\mu}^\e }_{\otau^t_n \land \oga} (\vf ) \big]  \oP(d\oo) . \q \;\; \label{071720_27}
\eea

 Fix $ \oo  \ins \big\{ \otau^t_n   \> \oga   \big\} \Cp \ocA_* \Cp \big( \ocN_* \cp \ofN^\e_U \cp \ofN^\e_\sharp   \big)^c  $ and set $\fr_\loo \df t_\oo \ve r $.
 Since $ t_\oo \= \oga (\oo) \< \otau^t_n (\oo) \ls t\+n $,
 applying \eqref{020522_17} with $(t,\bx,\oP,\ol{\mu},\vf, \fra,\oz_1,\oz_2) \= \big(t_\oo ,  \oX_{\oga \land \cd} (\oo), \oQ^\oo_\e ,\ol{\mu}^\oo,\vf_\loo,   \fra_{\overset{}{\oo}},  t_\oo, \fr_\loo \ld (t\+n)\big)$ yields that
  \bea \label{071720_23}
  0 \= E_{\oQ^\oo_\e} \Big[   \oMoo_{\oz^n_\loo \land ( t_\oo + n  ) \land \fr_\loo \land (t + n)} (\vf_\loo) \- \oMoo_{\oz^n_\loo \land ( t_\oo + n  )   \land t_\oo} (\vf_\loo) \Big] \=   E_{\oQ^\oo_\e} \Big[   \oMoo_{\oz^n_\loo   \land \fr_\loo \land (t+n) } (\vf_\loo) \- \oM^{t_\oo,\ol{\mu}^\oo}_{t_\oo} (\vf_\loo) \Big] .
 \eea

 Because $\oo \ins \big\{ \oo' \ins \oO \n :  \otau^t_n(\oo') \> \oga (\oo) \big\} \ins \ocF^t_{\oga (\oo)} \sb \ocG^t_{\oga (\oo)} $,
 an analogy to \eqref{090920_15} shows that
 $\oTh^t_{\oga,\oo} \sb \big\{ \oo' \ins \oO \n :  \otau^t_n(\oo') \> \oga (\oo) \big\}$.
 Let $\oo' \ins \oTh^t_{\oga,\oo} \Cp \hAe_{2,\oo} $. Since $\inf \n \big\{s \ins [t,\infty) \n : \big| (\oW^t_s,\oX_s) (\oo') \big|   \gs n  \big\} \gs \otau^t_n(\oo') \> \oga(\oo)$, one has $\big| (\oW^t_s,\oX_s) (\oo') \big|\< n $, $\fa s \ins [t,t_\oo] $  and thus
 \beas
   \;\;\;   \inf \n \big\{s \ins [t,\infty) \n : |(\oW^t_s,\oX_s) (\oo')|   \gs n  \big\}
    \= \inf\big\{s \ins [t_\oo,\infty) \n : \big|\big(\oW^t_s (\oo') \- \oW^t (\oga(\oo),\oo'),\oX_s(\oo')\big) \+ \big(\oW^t_\oga (\oo),\bz\big) \big|   \gs n  \big\}
 \= \oz^n_\loo (\oo')   .
 \eeas
  It follows that $  \otau^t_n(\oo')   \=   \oz^n_\loo (\oo')   \ld (t\+n)   $.
  Taking $(r_1,r_2) \= \big(t_\oo, \otau^t_n(\oo') \ld \fr_\loo \big)$ in \eqref{091920_17},
   we can   deduce from \eqref{090520_11}   that
 \beas
&& \hspace{-1.2cm} \big(\,\oM^{t,\ol{\mu}^\e} \n (\vf)\big) \big(   \otau^t_n(\oo')    \ld  (\oga (\oo') \ve r ) ,\oo' \big) \- \big(\,\oM^{t,\ol{\mu}^\e} \n (\vf ) \big) \big( \oga (\oo'), \oo'  \big)
  \= \big(\,\oM^{t,\ol{\mu}^\e} \n (\vf)\big) \big(  \otau^t_n(\oo')    \ld  (\oga (\oo) \ve r ) ,\oo' \big) \- \big(\oM^{t,\ol{\mu}^\e} \n (\vf ) \big) \big( \oga (\oo), \oo'  \big)  \nonumber \\
 &&  \q  
 \= \big( \oM^{t_\oo,\ol{\mu}^\oo}  (\vf_\loo) \big) \big( \oz^n_\loo (\oo') \ld (t \+ n)  \ld \fr_\loo, \oo' \big) \- \big( \oM^{t_\oo,\ol{\mu}^\oo}  (\vf_\loo) \big)  ( t_\oo , \oo'  ) . \q  
\eeas
 As $\{ \otau^t_n   \> \oga   \}  \ins \ocF^t_{\otau^t_n \land \oga} \sb \ocG^t_\oga $,
   \eqref{May06_25},  \eqref{072820_15}, \eqref{091022_17}
  and \eqref{071720_23} then imply that
 \beas
  && \hspace{-1cm} E_{\oQ^\oo_\e} \big[ \,  \oM^{t,\ol{\mu}^\e}_{\otau^t_n \land (\oga \vee r)} (\vf ) \- \oM^{t,\ol{\mu}^\e }_{\otau^t_n \land \oga} (\vf )  \big]
 \= E_{\oQ^\oo_\e} \big[ \b1_{ \{ \otau^t_n   > \oga   \}}  \big( \,   \oM^{t,\ol{\mu}^\e}_{ \otau^t_n   \land (\oga \vee r)} (\vf ) \- \oM^{t,\ol{\mu}^\e}_\oga (\vf ) \big)   \big]
\= \b1_{ \{ \otau^t_n (\oo)  > \oga (\oo)  \}}  E_{\oQ^\oo_\e} \big[  \,  \oM^{t,\ol{\mu}^\e}_{ \otau^t_n   \land (\oga \vee r)} (\vf ) \- \oM^{t,\ol{\mu}^\e}_\oga (\vf )   \big] \\
  &&    \=   \b1_{ \{ \otau^t_n(\oo)   > \oga (\oo)  \} }  E_{\oQ^\oo_\e} \Big[   \oMoo_{\oz^n_\loo   \land  (t+n)  \land \fr_\loo }  (\vf_\loo)    \-   \oM^{t_\oo,\ol{\mu}^\oo}_{t_\oo}  (\vf_\loo)    \Big] \= 0  , ~ \;  \fa \oo  \ins \ocA_* \Cp \big( \ocN_* \cp \ofN^\e_U \cp \ofN^\e_\sharp   \big)^c  .
 \eeas
 Thus  $  \int_{\oo \in \ocA_*}  \n  \b1_{ \{\oga (\oo) > s \}}  \b1_{\{\oo \in \oA\,\} }   E_{\oQ^\oo_\e} \big[  \, \oM^{t,\ol{\mu}^\e }_{\otau^t_n   \land (\oga \vee r)} (\vf )  \-   \oM^{t,\ol{\mu}^\e }_{\otau^t_n   \land \oga}   (\vf )  \big]  \oP(d\oo)   \= 0 $,
 which together with \eqref{071720_27} and \eqref{071420_21}   leads to \eqref{071720_31}.

  \no {\bf (ii)} If $t\+n \> s \= t $,    as $\{\oga \> t\} \= \oO$,    \eqref{071720_31} directly becomes
 $  E_{\oP_\e} \big[   \big( \, \oM^{t,\ol{\mu}^\e}_{\otau^t_n \land r} (\vf ) \- \oM^{t,\ol{\mu}^\e}_t (\vf ) \big)   \b1_{\oA}  \big]  \= 0  $.

 Next, let $t\+ n \> s \> t $. In this case, we can assume with loss of generality that  $ t \=   s_1 \< \cds   \< s_k \= s $ with $k \gs 2 $.
 Since  \eqref{091022_31} renders that   $ \b1_{\wh{A}^\e_2} \b1_{\{\oga \le s\}} \big( \, \oM^{t,\ol{\mu}^\e}_{\otau^t_n \land r} (\vf ) \- \oM^{t,\ol{\mu}^\e}_{\otau^t_n \land s} (\vf ) \big) \= \b1_{\wh{A}^\e_2} \b1_{\{\oga \le s\}} \big( \, \oM^{t,\ol{\mu} }_{\otau^t_n \land r} (\vf ) \- \oM^{t,\ol{\mu} }_{\otau^t_n \land s} (\vf ) \big) $
  and since   $ \ocA^c_* \Cp   \{\oga \ls s  \} \ins \cF^{\oW^t}_s \sb \ocF^t_s $,
  using $\oP \Big(\ocA^c_* \cap \big(\hAe_2\big)^c\Big) \= 0 $
  and taking $(\fra,\oz_1,\oz_2 ) \= (\bz,s,r) $  in  \eqref{020522_17}     yield that
 \bea  \label{071820_27}
  E_\oP \Big[ \b1_{\ocA^c_* \cap \{\oga \le s\}} \big( \, \oM^{t,\ol{\mu}^\e}_{\otau^t_n \land r} (\vf ) \- \oM^{t,\ol{\mu}^\e}_{\otau^t_n \land s} (\vf ) \big)   \b1_{ \oA}     \Big]
  \=  E_\oP \Big[ \b1_{\ocA^c_* \cap \{\oga \le s\}} \big( \, \oM^{t,\ol{\mu} }_{\otau^t_n \land r} (\vf ) \- \oM^{t,\ol{\mu} }_{\otau^t_n \land s} (\vf ) \big)   \b1_{ \oA}     \Big] \= 0 .
  \eea

 Fix $i \ins \{ 1, \cds \n , k \- 1\}$   and fix $\oo \ins \big\{ \otau^t_n   \> \oga   \big\} \Cp \big\{s_i \< \oga  \ls   s_{i+1} \big\} \Cp \ocA_* \Cp  \big( \ocN_* \cp \ofN^\e_U \cp \ofN^\e_\sharp   \big)^c  $.
     Since $ \Wtgo \sb \{\oga \> s_i  \} $ by  \eqref{090520_11},   $\oA_i \df \ccap{j=1}{i}    \big( \oW^t_{ \oga \land s_j }  , \oX_{ \oga \land s_j }      \big)^{-1} (\cE_j) \ins 
 \ocF^t_{\n \oga} $ satisfies that
\bea
    \Wtgo \Cp \Big(  \ccap{j=1}{i}   \big( \oW^t_{s_j  }, \oX_{s_j  }  \big)^{-1} (\cE_j)  \Big)
 \= \Wtgo \Cp \oA_i .   \label{070120_11}
\eea
 Also,  \eqref{090520_11} shows that $ \Wtgo \sb \{\oga \ls s   \} $ and thus $\Wtgo \Cp \{\otau^t_n \ls \oga\} 
  \sb \{\otau^t_n \ls s\}$. We see from  \eqref{072820_15} that
  \bea  \label{021322_11}
    E_{\oQ^\oo_\e} \big[ \b1_{\{\otau^t_n \le \oga\}} \big( \oM^{t,\ol{\mu}^\e}_{\otau^t_n \land r} (\vf ) \- \oM^{t,\ol{\mu}^\e}_{\otau^t_n \land s} (\vf ) \big)  \b1_{ \oA }  \big]
   \=  E_{\oQ^\oo_\e} \big[ \b1_{\Wtgo \cap \{\otau^t_n \le \oga\}} \big( \oM^{t,\ol{\mu}^\e}_{\otau^t_n \land r} (\vf ) \- \oM^{t,\ol{\mu}^\e}_{\otau^t_n \land s} (\vf ) \big)  \b1_{ \oA }  \big] \= 0  .
   \eea

  Define  $\oA^\oo_i   \df   \ccap{j=i+1}{k}   \big( \oW^{t_\oo}_{s_j   }, \oK^\oo_{s_j } \big)^{-1}  ( \cE_{j,\oo}  )   \ins \cF^{\oW^{t_\oo}}_{s   } $
 with   $\cE_{j,\oo} \df \big\{\fx \+ \fra_{\overset{}{\oo}} \n : \fx \ins \cE_j \big\} \ins \sB(\hR^{d+l})$.
Using \eqref{020522_17} with $(t,\bx,\oP,\ol{\mu},\vf, \fra,\oz_1,\oz_2) \= \big(t_\oo ,   \oX_{\oga \land \cd} (\oo), \oQ^\oo_\e ,\ol{\mu}^\oo,\vf_\loo,   \fra_{\overset{}{\oo}} ,   s  ,    (t\+n) \ld r \big)$ renders that
\bea
  0 & \tn \= & \tn  E_{\oQ^\oo_\e} \Big[ \Big(    \oMoo_{\oz^n_\loo   \land (t_\oo+n )   \land (t +n) \land r)}(\vf_\loo)    \- \oMoo_{\oz^n_\loo   \land (t_\oo+n )   \land s}(\vf_\loo) \Big) \b1_{\oA^\oo_i } \Big] \nonumber \\
 & \tn \= & \tn  E_{\oQ^\oo_\e} \Big[ \Big(    \oMoo_{\oz^n_\loo   \land (t +n) \land r)}(\vf_\loo)
 \- \oMoo_{\oz^n_\loo     \land s}(\vf_\loo) \Big) \b1_{\oA^\oo_i } \Big] .   \label{071820_25}
\eea

 Let $\oo' \ins \oTh^t_{\oga,\oo} \Cp \hAe_{2,\oo} $.
 Like in  Step (i), we still have $\oTh^t_{\oga,\oo} \sb \big\{ \oo' \ins \oO \n :  \otau^t_n(\oo') \> \oga (\oo) \big\}$
 and $  \otau^t_n(\oo')   \=   \oz^n_\loo (\oo')   \ld (t\+n)   $. 
  As $ \otau^t_n(\oo') \ld s \gs \oga(\oo) $,
  taking $(r_1,r_2) \= \big(\otau^t_n(\oo') \ld s, \otau^t_n(\oo') \ld r \big)$ in \eqref{091920_17} shows that
  $ \big(\oM^{t,\ol{\mu}^\e} \n (\vf )\big) \big(\otau^t_n(\oo') \ld r,\oo'\big) \- \big(\oM^{t,\ol{\mu}^\e} \n (\vf )\big) \big(\otau^t_n(\oo') \ld s,\oo'\big) \= \big( \oM^{t_\oo,\ol{\mu}^\oo}  (\vf_\loo) \big) \big( \oz^n_\loo (\oo')   \ld  (t\+n) \ld r , \oo' \big) \- \big( \oM^{t_\oo,\ol{\mu}^\oo}  (\vf_\loo) \big) \big( \oz^n_\loo (\oo')   \ld s, \oo' \big)  $. It follows from \eqref{072820_15}, \eqref{091022_17} and \eqref{071820_25} that
  \bea \label{070622_19}
  E_{\oQ^\oo_\e} \big[    \big( \oM^{t,\ol{\mu}^\e}_{\otau^t_n \land r} (\vf )
 \- \oM^{t,\ol{\mu}^\e}_{\otau^t_n \land s } (\vf ) \big) \b1_{   \oA^\oo_i} \big]
  \=    E_{\oQ^\oo_\e} \Big[  \Big(     \oMoo_{\oz^n_\loo \land (t +n) \land r    }(\vf_\loo)    \- \oMoo_{\oz^n_\loo  \land s    }(\vf_\loo) \Big) \b1_{  \oA^\oo_i} \Big] \= 0 .
  \eea

 For any $j \ins \{ i\+1,\cds \n , k \} $ and $\oo' \ins  \Wtgo  \Cp \big( \ocN^{\,\oo}_{\n X} \cp \ocN^{\,\oo}_K \big)^c$,
   \eqref{090520_11} implies that  $ \big(\oW^t_{s_j}  , \oX_{s_j}  \big) (   \oo')  \ins  \cE_j $ if and only if
 $ \big( \oW^{t_\oo}_{s_j} ,   \oK^\oo_{s_j} \big)  (   \oo' )
     \=   \big( \oW^t_{  s_j} ( \oo') \- \oW^t \big(\oga(\oo) ,\oo'\big),
\osX^\oo_{  s_j}   ( \oo' )  \big)
\=  \big(\oW^t_{  s_j}  , \oX_{  s_j}  \big) (   \oo') \+ \fra_{\overset{}{\oo}} \ins  \cE_{j,\oo}  $.
 By \eqref{070120_11}, one has
 $ \oA   \Cp \Wtgo  \Cp \big( \ocN^{\,\oo}_{\n X} \cp  \ocN^{\,\oo}_{\n K} \big)^c
     \=  \oA_i \Cp  \oA^\oo_i   \Cp  \Wtgo \Cp \big( \ocN^{\,\oo}_{\n X} \cp  \ocN^{\,\oo}_{\n K} \big)^c   $.
    Then we can deduce from \eqref{021322_11},    \eqref{072820_15},    \eqref{May06_25} and  \eqref{070622_19}   that
 \beas
&& \hspace{-1.2cm}  E_{\oQ^\oo_\e} \big[ \big( \oM^{t,\ol{\mu}^\e}_{\otau^t_n \land r} (\vf ) \- \oM^{t,\ol{\mu}^\e}_{\otau^t_n \land s} (\vf ) \big)  \b1_{ \oA }  \big]
\= E_{\oQ^\oo_\e} \big[ \b1_{\{\otau^t_n > \oga\}} \big( \oM^{t,\ol{\mu}^\e}_{\otau^t_n \land r} (\vf ) \- \oM^{t,\ol{\mu}^\e}_{\otau^t_n \land s} (\vf ) \big)  \b1_{ \oA }  \big]
 \= E_{\oQ^\oo_\e} \big[  \b1_{\{ \otau^t_n  > \oga \}} \big( \oM^{t,\ol{\mu}^\e}_{\otau^t_n \land r} (\vf ) \- \oM^{t,\ol{\mu}^\e}_{\otau^t_n \land s} (\vf ) \big) \b1_{\oA_i \cap  \oA^\oo_i} \big] \\
&&  \hspace{-0.3cm}  \= \b1_{\{\oo \in \oA_i\}} \b1_{\{ \otau^t_n (\oo)  >  \oga (\oo) \}}  E_{\oQ^\oo_\e} \Big[  \big( \oM^{t,\ol{\mu}^\e}_{\otau^t_n \land r} (\vf )
 \- \oM^{t,\ol{\mu}^\e}_{\otau^t_n \land s } (\vf ) \big) \b1_{  \oA^\oo_i} \Big] \= 0 , \q   \fa \oo \ins \{s_i \< \oga   \ls s_{i+1}\} \Cp \ocA_* \Cp  \big( \ocN_* \cp \ofN^\e_U \cp \ofN^\e_\sharp   \big)^c  ,
\eeas
 and thus
 $ \int_{\oo \in \ocA_*} \n  \b1_{ \{ s_i < \oga(\oo) \le s_{i+1}  \}}
       E_{\oQ^\oo_\e} \big[ \big( \oM^{t,\ol{\mu}^\e}_{\otau^t_n \land r} (\vf ) \- \oM^{t,\ol{\mu}^\e}_{\otau^t_n \land s} (\vf ) \big)  \b1_{ \oA } \big]  \oP(d\oo) \= 0 $.
 Taking summation    from $i\=1$ through $i\=k \- 1$, we obtain from \eqref{071820_27}    that
$ E_{\oP_\e} \big[ \b1_{\{  \oga \le s \}} \big( \oM^{t,\ol{\mu}^\e}_{\otau^t_n \land r} (\vf ) \- \oM^{t,\ol{\mu}^\e}_{\otau^t_n \land s} (\vf ) \big)  \b1_{ \oA }    \big]
    \=  E_\oP \big[ \b1_{\ocA^c_* \cap \{  \oga \le s\}  }  \big( \oM^{t,\ol{\mu}^\e}_{\otau^t_n \land r} (\vf ) \- \oM^{t,\ol{\mu}^\e}_{\otau^t_n \land s} (\vf ) \big) \b1_\oA  \big]
\+  \int_{\oo \in \ocA_*}  \n   \b1_{ \{   \oga(\oo) \le s  \}}
       E_{\oQ^\oo_\e} \big[ \big( \oM^{t,\ol{\mu}^\e}_{\otau^t_n \land r} (\vf ) \- \oM^{t,\ol{\mu}^\e}_{\otau^t_n \land s} (\vf ) \big) \b1_{ \oA } \big]  \oP(d\oo)
\= 0 $.
 Adding it to \eqref{071720_31} yields that
 the Lambda-system $\Big\{\oA \ins \sB(\oO) \n : E_{\oP_\e} \big[ \big( \oM^{t,\ol{\mu}^\e}_{\otau^t_n \land r} (\vf ) \- \oM^{t,\ol{\mu}^\e}_{\otau^t_n \land s} (\vf ) \big)  \b1_{\oA} \big] \= 0 \Big\}$
 $   $ contains the Pi-system $\Big\{ \ccap{i=1}{k}  (\oW^t_{  s_i  },\oX_{  s_i   })^{-1} (\cE_i) \n : \{(s_i,\cE_i )\}^k_{i=1} \sb   [t,s]  \ti \sB (\hR^{d+l}) \Big\}$ and thus includes $  \ocF^t_s $ thanks to Dynkin's Pi-Lambda Theorem.
 Hence, $\big\{ \oM^{t,\ol{\mu}^\e}_{ s \land \otau^t_n  } (\vf) \big\}_{s \in  [t,\infty) } $  is a bounded $ \big( \obF^t , \oP_\e \big)   -$martingale.
 According to Remark \ref{rem_ocP} (ii), $ \oP_\e $  satisfies  \(D3\)   of $\ocP_{t,\bx}$.

 Similar to Part (II.d) in the proof of \cite[Theorem 5.1]{OSEC_stopping}, we can construct a $[t,\infty]-$valued   $ \bF^{W^t,P_0} -$stopping time $\wh{\tau}_\e$ with  $ \oP_\e \big\{  \oT  \=  \wh{\tau}_\e (\oW)   \big\} \= 1$, So $ \oP_\e $ also satisfies (D4) of $   \ocP_{t,\bx} $.

 \no {\bf II.f)} Fix $i \ins \hN$.  Since $ \big\{  \int_t^s g_i  \big(r, \osX^{t,\bx,\ol{\mu}}_{r \land \cd},\ol{\mu}_r  \big) dr  \big\}_{s \in [t,\infty)} $
 and $ \big\{  \int_t^s h_i  \big(r, \osX^{t,\bx,\ol{\mu}}_{r \land \cd},\ol{\mu}_r   \big) dr  \big\}_{s \in [t,\infty)} $ are two
 $ \bF^{\oW^t,\oP}  - $adapted continuous processes,
  Lemma 2.4 of \cite{STZ_2011a} assures two   $  \bF^{\oW^t }  -$predictable  processes  $  \big\{ \ol{\Phi}^i_s  \big\}_{s \in [t,\infty)}$ and
  $  \big\{ \ol{\Psi}^i_s  \big\}_{s \in [t,\infty)}$
   such that   $ \ocN^{\, i,1}_{g,h} \df \Big\{\oo \ins \oO \n : \ol{\Phi}^i_s (\oo) \nne \int_t^s g_i  \big(r, \osX^{t,\bx,\ol{\mu}}_{r \land \cd} (\oo),\ol{\mu}_r(\oo) \big) dr \hb{ or } \ol{\Psi}^i_s (\oo) \nne \int_t^s h_i  \big(r, \osX^{t,\bx,\ol{\mu}}_{r \land \cd} (\oo),\ol{\mu}_r(\oo) \big) dr \hb{ for some }   s \ins [t,\infty) \Big\} \ins  \sN_\oP \big(\cF^{\oW^t}_\infty \big)$.
   \if{0}

 Since $ \big\{  \int_t^s g^\pm_i  \big(r, \osX^{t,\bx,\ol{\mu}}_{r \land \cd},\ol{\mu}_r  \big) dr  \big\}_{s \in [t,\infty)} $
 and $ \big\{  \int_t^s h^\pm_i  \big(r, \osX^{t,\bx,\ol{\mu}}_{r \land \cd},\ol{\mu}_r  \big) dr  \big\}_{s \in [t,\infty)} $ are   $[0,\infty]-$valued
 $ \bF^{\oW^t,\oP}  - $adapted continuous processes,  Lemma 2.4 of \cite{STZ_2011a} assures   $[0,\infty]-$valued  $  \bF^{\oW^t }  -$predictable  processes  $  \big\{ \ol{\Phi}^{i,\pm}_s  \big\}_{s \in [t,\infty)}$ and
  $  \big\{ \ol{\Psi}^{i,\pm}_s  \big\}_{s \in [t,\infty)}$
   such that all $\oo \ins \oO$ except on some $ \ocN^\pm_{\n g_i,h_i} \ins  \sN_\oP \big(\cF^{\oW^t}_\infty \big)$,
   $\ol{\Phi}^{i,\pm}_s (\oo) \= \int_t^s g^\pm_i  \big(r, \osX^{t,\bx,\ol{\mu}}_{r \land \cd} (\oo),\ol{\mu}_r (\oo) \big) dr$
   and $\ol{\Psi}^{i,\pm}_s (\oo) \= \int_t^s h^\pm_i  \big(r, \osX^{t,\bx,\ol{\mu}}_{r \land \cd} (\oo),\ol{\mu}_r (\oo) \big) dr$ for any $ s \ins [t,\infty) $.
   So for any $\oo \ins \big( \ocN^+_{\n g_i,h_i} \cp \ocN^-_{\n g_i,h_i}\big)^c$,
    \beas
    \ol{\Phi}^i_s (\oo) \df \ol{\Phi}^{i,+}_s (\oo) \- \ol{\Phi}^{i,-}_s (\oo) \= \int_t^s g_i  \big(r, \osX^{t,\bx,\ol{\mu}}_{r \land \cd} (\oo),\ol{\mu}_r (\oo) \big) dr
    \aand
   \ol{\Psi}^i_s (\oo) \df \ol{\Psi}^{i,+}_s (\oo) \- \ol{\Psi}^{i,-}_s (\oo) \= \int_t^s h_i  \big(r, \osX^{t,\bx,\ol{\mu}}_{r \land \cd} (\oo),\ol{\mu}_r (\oo) \big) dr , \q \fa s \ins [t,\infty) .
    \eeas

      \fi
 By Remark \ref{rem_ocP2} (1),
 $E_\oP \big[ \int_t^\infty   \n g^-_i  \big( r,  \oX_{r \land \cd},\oU_r   \big) \ve h^-_i  \big( r,  \oX_{r \land \cd},\oU_r  \big)   dr \big] \< \infty$.
 So  it holds for any $\oo \ins \oO$ except on some $ \ocN^{\, i,2}_{g,h} \ins \sN_\oP \big(\sB(\oO)\big) $ 
  that $ \int_t^\infty   \n g^-_i  \big( r,  \oX_{r \land \cd} (\oo),\oU_r(\oo)   \big) \ve h^-_i  \big( r,  \oX_{r \land \cd} (\oo),\oU_r(\oo)  \big)   dr \< \infty$.

  For any $\oo \ins \ocA_*   \Cp  \ocN^c_* \Cp \ocN^c_{\n X} \Cp \oO_\mu \Cp \big(  \ocN^{\, i,1}_{g,h} \cp \ocN^{\, i,2}_{g,h} \big)^c   $,
 as  $ \oO^t_{\oga,\oo} \sb \big\{ \oo' \ins \oO \n:  \oX_s  (\oo') \=  \oX_s   (\oo)   ,
    \fa s \ins \big[0,\oga(\oo)\big] ;  \oU_s  (\oo') \=  \oU_s   (\oo) \hb{ for a.e. } s \ins (t,\oga(\oo)) \big\} $,
     \eqref{072820_15} and \eqref{April07_11} show that
  \bea
  && \hspace{-1.2cm}   E_{\oQ^\oo_\e} \Big[ \int_t^\oT   \n g_i \big( r,  \oX_{r \land \cd},\oU_r   \big)   dr \Big]
      \=      \int_t^{\oga(\oo)}   \n g_i \big(   r,  \oX_{r \land \cd} (\oo),\oU_r (\oo) \big)   dr
    \+ E_{\oQ^\oo_\e} \Big[  \int_{\oga(\oo)}^\oT   \n g_i \big( r,  \oX_{r \land \cd},\oU_r   \big)   dr   \Big] \label{091020_15} \\
    && \ls    \int_t^{\oga(\oo)}   \n g_i \big(   r,  \osX^{t,\bx,\ol{\mu}}_{r \land \cd} (\oo), \ol{\mu}_r (\oo)  \big)   dr
    \+   \big( \oY^i_{\n \oP}  ( \oga ) \big) (\oo)
    \= \ol{\Phi}^i_\oga (\oo) \+  E_\oP \Big[ \int_{\oT \land  \oga }^\oT    g_i (r,\oX_{r \land \cd},\oU_r ) dr \Big| \cF^{\oW^t}_\oga \Big] (\oo)
    \nonumber
    \eea
    and similarly that
     $  E_{\oQ^\oo_\e} \Big[ \int_t^\oT   \n h_i \big( r,  \oX_{r \land \cd},\oU_r    \big)   dr \Big]
    \= \ol{\Psi}^i_\oga (\oo) \+  E_\oP \Big[ \int_{\oT \land  \oga }^\oT    h_i (r,\oX_{r \land \cd},\oU_r  ) dr \Big| \cF^{\oW^t}_\oga \Big] (\oo) $.
     Since $ \ocA_* ,  \ol{\Phi}^i_\oga  \ins \cF^{\oW^t}_\oga$ and
  since $\b1_{\ocA_*}  \= \b1_{\{ \wh{\tau}(\oW)   \ge   \oga \}} \=  \b1_{\{\oT    \ge   \oga \}}   $, $\oP-$a.s.  by \eqref{020322_11},
  we can deduce from the tower property that
    \beas
 && \hspace{-1.5cm} \int_{\oo \in \ocA_*}   E_{\oQ^\oo_\e} \Big[  \int_t^\oT   g_i \big(r,\oX_{r \land \cd},\oU_r  \big) dr \Big]     \oP(d\oo)
     \ls  
     E_\oP \bigg[    E_\oP \Big[ \b1_{\ocA_*} \Big( \ol{\Phi}^i_\oga
    \+ \int_{\oT \land \oga}^\oT    g_i(r,\oX_{r \land \cd},\oU_r ) dr \Big) \Big| \cF^{\oW^t}_\oga \Big]   \bigg] \\
   &&  \=    E_\oP \Big[ \b1_{\ocA_*} \Big( \int_t^\oga  g_i  \big(r, \osX^{t,\bx,\ol{\mu}}_{r \land \cd}, \ol{\mu}_r    \big) dr
    \+ \int_\oga^\oT    g_i(r,\oX_{r \land \cd},\oU_r ) dr \Big)    \Big]
     \=  E_\oP \Big[   \b1_{\ocA_*}  \int_t^\oT   \n g_i \big( r,  \oX_{r \land \cd},\oU_r    \big)   dr \Big]
 \eeas
 and thus $ E_{\oP_\e} \big[    \int_t^\oT   g_i \big(r,\oX_{r \land \cd},\oU_r  \big) dr \big]
 \ls   E_\oP  \big[    \int_t^\oT   g_i \big(r,\oX_{r \land \cd} ,\oU_r \big) dr \big] \ls y_i $.
 Analogously, we have  $ E_{\oP_\e} \big[    \int_t^\oT  h_i \big(r,\oX_{r \land \cd},\oU_r  \big) dr \big] \\
 \=  E_\oP  \big[    \int_t^\oT  h_i \big(r,\oX_{r \land \cd},\oU_r  \big) dr \big] \= z_i $.  Hence,   $\oP_\e $ belongs to $  \ocP_{t,\bx}(y,z)$.

 \no {\bf II.g)} Since the  mapping $\ddot{\Psi} $ is   $ \si \Big( \cF^{\oW^t}_\oga \cp \sN_\oP \big( \ocG^t_\infty\big) \Big)  \Big/ \si \big( \sB  ( \cD_\ocP  ) \cp \sN_{\ddot{P} }  ( \sB   ( \cD_\ocP  )  )  \big) -$measurable by Step (II.a.2),
 Theorem \ref{thm_V_usa} renders that 
 $D^V_\infty \df  \big\{\oo \ins \ocA_* \cap   \ocN^c_* \n :    \oV \big( \ddot{\Psi} (\oo) \big) \= \infty  \big\} $ is $\si \Big(  \cF^{\oW^t}_\oga \cp   \sN_\oP \big(\ocF^t_{\n \infty}  \big) \Big) -$measurable.
  As $E_\oP \big[ \int_t^\infty   \n f^-   \big( r,  \oX_{r \land \cd},\oU_r    \big)    dr \big]   \< \infty$,
  there is a $ \ocN_{\n f} \ins \sN_\oP \big(\sB(\oO)\big) $ such that
 $ \int_t^\infty   \n f^-   \big( r,  \oX_{r \land \cd} (\oo), \oU_r (\oo)   \big)    dr \< \infty$
 for any $\oo \ins \ocN^c_{\n f}$.

  Let $\e \ins (0,1)$.
    \if{0}

    since \eqref{081620_19} implies that
\beas 
    E_{\oQ^\oo_\e} \Big[ \int_{\oga(\oo)}^\oT   \n f  \big( r,  \oX_{r \land \cd}   \big)   dr   \+    \b1_{\{\oT < \infty\}} \pi  \big( \oT,\oX_{\oT \land \cd} \big) \Big] \gs
\left\{
\ba{ll}
1/\e , & \fa \oo \ins   D^V_\infty ;\\
\oV \big( \ddot{\Psi} (\oo) \big)    \- \e , & \fa \oo \ins   \big(D^V_\infty\big)^c,
 \ea
 \right.
 \eeas

   \fi
  For any  $\oo \ins \ocA_* \Cp \ocN^c_* \Cp \ocN^c_{\n f} $,
   an analogy to \eqref{091020_15}, \eqref{081620_19} and Theorem \ref{thm_V=oV}  imply    that
   \beas
    E_{\oQ^\oo_\e} \big[ \, \oR(t) \big]
   & \tn \=  & \tn    \int_t^{\oga(\oo)}   \n f  \big(  r,  \oX_{r \land \cd} (\oo),\oU_r(\oo)  \big)   dr \+ E_{\ol{\bQ}_\e  ( \ddot{\Psi} (\oo) )} \big[   \, \oR(\oga(\oo)) \big] \\
    & \tn  \gs  & \tn   \int_t^{\oga(\oo)}   \n f  \big(  r,  \oX_{r \land \cd} (\oo),\oU_r(\oo)  \big)   dr \+  \b1_{ \{\oo \in  (D^V_\infty )^c  \}} \Big( \oV \big( \oga(\oo)   ,  \oX_{\oga \land \cd}(\oo)  ,   \big(\oY_{\n \oP}  ( \oga )\big) (\oo),    \big( \oZ_\oP  ( \oga ) \big) (\oo)   \big)   \- \e \Big)
   \+   \frac{1}{\e} \b1_{ \{\oo \in D^V_\infty  \}}    .
 \eeas
  Since $\oP_\e \ins \ocP_{t,\bx}(y,z)$ and since $\b1_{\ocA_*} 
  \=  \b1_{\{\oT    \ge   \oga \}}   $, $\oP-$a.s. 
  \bea
  && \hspace{-1cm}\oV (t,\bx,y,z)    \gs    E_{\oP_\e} \n  \big[ \, \oR(t) \big]
  \gs    E_\oP   \bigg[ \b1_{\ocA^c_*} \oR(t) \+ \b1_{\ocA_*} \Big(  \int_t^{\oga }   \n f  \big(  r,  \oX_{r \land \cd},\oU_r   \big)   dr \+
   \b1_{  (D^V_\infty )^c  } \big[ \, \oV \big( \oga   ,  \oX_{\oga \land \cd}  ,   \oY_{\n \oP}  ( \oga ) ,     \oZ_\oP  ( \oga )    \big)   \- \e \big]  \n \+   \frac{1}{\e} \b1_{  D^V_\infty  }  \Big)     \bigg] \nonumber \\
   & & \hspace{-0.5cm} \=    E_\oP   \bigg[ \b1_{\{\oT   <  \oga \}} \oR(t) \+ \b1_{\{\oT    \ge   \oga \}} \Big( \int_t^{\oga }   \n f  \big(  r,  \oX_{r \land \cd},\oU_r    \big)   dr \+
   \b1_{  (D^V_\infty )^c  } \big[ \, \oV \big( \oga   ,  \oX_{\oga \land \cd}  ,   \oY_{\n \oP}  ( \oga ) ,     \oZ_\oP  ( \oga )    \big)   \- \e \big]    \+   \frac{1}{\e} \b1_{  D^V_\infty  }  \Big)     \bigg] \nonumber \\
    && \hspace{-0.5cm}  \gs    E_\oP   \bigg[ \b1_{\{\oT   <  \oga \}} \oR(t) \+ \b1_{\{\oT    \ge   \oga \}} \Big( \int_t^{\oga }   \n f  \big(  r,  \oX_{r \land \cd},\oU_r    \big)   dr \+
   \b1_{  (D^V_\infty )^c  }   \oV \big( \oga   ,  \oX_{\oga \land \cd}  ,   \oY_{\n \oP}  ( \oga ) ,     \oZ_\oP  ( \oga )    \big)       \Big)     \bigg] \- \e \+ \frac{1}{\e} \oP \big(   \{\oT    \gs   \oga \} \Cp D^V_\infty \big)
   .  \qq \qq
   \label{081720_17}
 \eea

 To verify \eqref{081720_15}, we set $\ocI^t_\oP \df \b1_{\{\oT < \oga  \}} \oR(t) \+ \b1_{\{\oT \ge \oga  \}} \Big( \n \int_t^\oga  \n    f(r,\oX_{r \land \cd},\oU_r ) dr  \+ \oV \big( \oga  ,\oX_{ \oga   \land \cd} ,    \oY_{\n \oP}   (\oga  ) , \oZ_\oP   (\oga  )  \big) \Big) $.

\bul  If $ \oP \big( \{\oT    \gs   \oga \} \cap D^V_\infty \big) \= 0  $, then
   $ \oV (t,\bx,y,z)   \gs  E_\oP \Big[ \b1_{\{\oT   <  \oga \}} \oR(t) \+ \b1_{\{\oT    \ge   \oga \}} \Big( \int_t^{\oga }   \n f  \big(  r,  \oX_{r \land \cd},\oU_r    \big)   dr \+
       \oV \big( \oga   ,  \oX_{\oga \land \cd}  ,   \oY_{\n \oP}  ( \oga ) ,     \oZ_\oP  ( \oga )    \big) \Big) \Big] \\ \- \e    $
 holds for any $\e \ins (0,1)$.    Letting $\e \nto 0$ gives   \eqref{081720_15}.

 \bul   If $ \oP \big( \{\oT    \gs   \oga \} \Cp D^V_\infty \big) \> 0  $  and
 $  E_\oP \big[ \big(\ocI^t_\oP\big)^- \big] \= \infty$, then $ E_\oP \big[  \ocI^t_\oP  \big] \= -\infty \ls  \oV (t,\bx,y,z) $,
 so \eqref{081720_15} holds automatically.

 \bul If $ \oP \big( \{\oT    \gs   \oga \} \cap D^V_\infty \big) \> 0  $  and $ E_\oP \big[ \big(\ocI^t_\oP\big)^- \big] \< \infty $,
 since Remark \ref{rem_ocP2} (1) shows that  $E_\oP   \Big[ \-\b1_{\{\oT   <  \oga \}} \oR(t) \- \b1_{\{\oT    \ge   \oga \}} \Big(  \n \int_t^{\oga }   \n f  \big(  r,  \oX_{r \land \cd}, \\ \oU_r    \big)   dr   \+
   \b1_{  (D^V_\infty )^c  }   \oV \big( \oga   ,  \oX_{\oga \land \cd}  ,   \oY_{\n \oP}  ( \oga ) ,     \oZ_\oP  ( \oga )      \big) \Big)        \Big]
      \= E_\oP   \Big[ \- \b1_{  (D^V_\infty )^c  }  \ocI^t_\oP
     \- \b1_{  D^V_\infty  }    \int_t^{\oga \land \oT}   \n f  \big(  r,  \oX_{r \land \cd} ,\oU_r   \big)   dr
     \- \b1_{\{\oT   <  \oga \} \cap  D^V_\infty   }   \pi   \big(\oT, \oX_{\oT \land \cd}\big)   \Big]
     \ls E_\oP   \big[   \big( \ocI^t_\oP \big)^-
     \+    \int_t^\infty   \n f^-  \big(  r,  \oX_{r \land \cd},\oU_r    \big)   dr \big] \- c_\pi \< \infty  $,
    we can deduce from   \eqref{081720_17}  that
 \beas
  \oV (t,\bx,y,z)
  \gs - E_\oP   \Big[   \big( \ocI^t_\oP \big)^-
     \+    \int_t^\infty   \n f^-  \big(  r,  \oX_{r \land \cd},\oU_r    \big)   dr \Big] \+ c_\pi \- \e \+ \frac{1}{\e} \oP \big(   \{\oT    \gs   \oga \} \Cp D^V_\infty \big) , \q  \fa \e \ins (0,1)  .
  \eeas
  Sending $\e \nto 0$   yields $ \oV (t,\bx,y,z) \= \infty $, so \eqref{081720_15} still holds. This completes the proof of Theorem \ref{thm_DPP1}. \qed

\appendix
\renewcommand{\thesection}{A}
\refstepcounter{section}
\makeatletter
\renewcommand{\theequation}{\thesection.\@arabic\c@equation}
\makeatother

\section{Appendix}

\begin{lemm} \label{lem_122921_11}
Let  $t_0 \ins [0,\infty)$. For $i \= 1,2$,  let $(\O_i, \cF_i, P_i)$ be   a    probability space
and  let $B^i\= \{B^i_s\}_{s \in [0,\infty)}$ be an  $\hR^d-$valued continuous process on $\O$ with $B^i_0 \= \bz$ such that
 $\fB^i_s \df B^i_s \- B^i_{t_0}$, $  s \ins [t_0,\infty)$ is a  Brownian motion on $(\O_i, \cF_i, P_i)$.
 Let $\Phi \n: \O_1 \mto \O_2$ be a mapping such that  $ \fB^2_s (\Phi(\o))\= \fB^1_s(\o)  $  for any $(s,\o) \ins [t_0,\infty) \ti \O_1 $, then
 \(i\) $\Phi^{-1}\big(\cF^{\fB^2}_s\big) \= \cF^{\fB^1}_s  $, $\fa s \ins [t_0,\infty]$;
 \(ii\) $\Phi^{-1}\big(\sN_{P_2}(\cF^{\fB^2}_\infty)\big) \sb \sN_{P_1}(\cF^{\fB^1}_\infty) $;
 \(iii\) $\Phi^{-1}\big(\cF^{\fB^2,P_2}_s\big) \sb \cF^{\fB^1,P_1}_s  $, $\fa s \ins [t_0,\infty]$
 and \(iv\) $P_1 \nci \Phi^{-1} (A) \= P_2(A)$ for any $ A \ins \cF^{\fB^2,P_2}_\infty$.

\end{lemm}

The proof of Lemma \ref{lem_122921_11} is basic. We refer interested readers to the ArXiv version of \cite{OSEC_stopping} for it.
We also recall the following result from \cite{OSEC_stopping}.

 \begin{lemm} \label{lem_M31_01}
 Let $(\O, \cF, P)$ be   a    probability space and let $\,t \ins [0,\infty)$.
 Let $B\= \{B_s\}_{s \in [0,\infty)}$ be an $\hR^d-$valued continuous process on $\O$ with $B_0 \= \bz$ such that
 $B^t_s \df B_s \- B_t$, $  s \ins [t,\infty)$ is a  Brownian motion on $(\O, \cF, P)$.

\no \(1\) For any $[t,\infty]-$valued $ \bF^{W^t,P_0}-$stopping time $\wh{\tau}$ on $\O_0$,
$\wh{\tau}(B)$ is an    $  \bF^{B^t,P}-$stopping time on $\O$.

\no \(2\)  Let $\Phi \n: \O  \mto \O_0$ be a mapping such that  $ W^t_s (\Phi(\o))\= B^t_s (\o)  $  for any $(s,\o) \ins [t,\infty) \ti \O  $.
 For any   $[t,\infty]-$valued $  \bF^{B^t,P}-$stopping time $\tau$ on $\O$,  there exists a $[t,\infty]-$valued $ \bF^{W^t,P_0}-$stopping time $\wh{\tau} $ on $\O_0$  such that $\tau \= \wh{\tau}  (\Phi )$, $P-$a.s.

 \end{lemm}

   \if{0}

\begin{lemm} \label{lem_predict0}
  Let $(\O, \cF, P)$ be   a    probability space. Given $t \ins [0,\infty)$,
 let $\wh{\mu} \= \{\wh{\mu}_s\}_{s \in [t,\infty)} $ be a $\hU-$valued, $\bF^{W^t} -$predictable  process on $\O_0$
  and let $ B \= \{ B_s\}_{s \in [0,\infty)}$ is an  $\hR^d-$valued continuous process with $B_0 \= \bz$  on   $ \O $.
 Then $\mu_s (\o) \df \wh{\mu}_s \big( B (\o)\big)$, $(s,\o) \ins [t,\infty) \ti \O $ is an $\bF^{B^t}-$predictable process.

\end{lemm}

\ss \no {\bf Proof:}  Let  $\sP^{W^t}$ be the predictable sigma-field of $[t,\infty) \ti \O_0$ with respect to the filtration $\bF^{W^t}$   \big(i.e., the sigma-field of $[t,\infty) \ti \O_0$ generated by all $\big\{\cF^{W^t}_{s-}\big\}_{s \in [t,\infty)}-$adapted, c\`agl\`ad processes\big),
 and let $\sP^{B^t}$ be the predictable sigma-field of $[t,\infty) \ti \O$ with respect to the filtration $\bF^{B^t}$.

  Let $s \ins (t,\infty)$.
  Since it holds for any $ r \ins [t,s)   $ and $\cE \ins \sB(\hR^d)$ that
 $ 
 B^{-1} \big( (W^t_r)^{-1} (\cE)\big) \= \big\{\o \ins \O\n : W^t_r \big(B  (\o)\big) \ins \cE \big\}
  \= \big\{\o \ins \O\n : B^t_r ( \o) \ins \cE \big\} \= (B^t_r)^{-1} (\cE) \ins \cF^{B^t}_r \sb \cF^{B^t}_{s-} $,
      the sigma-field of $\O_0$,
  $   \{A \sb \O_0 \n : B^{-1} (A) \ins \cF^{B^t}_{s-} \}$, contains $\cF^{W^t}_r$,
  $\fa r \ins [t,s)$ and thus includes $\cF^{W^t}_{s-} \= \si\Big(\underset{r \in  [t,s)}{\cup} \cF^{W^t}_r \Big) $. Namely,
  \bea \label{Sep28_01c}
  B^{-1}(A) \ins \cF^{B^t}_{s-} , \q \fa A \ins \cF^{W^t}_{s-} , ~ \fa s \ins (t,\infty) .
  \eea
 It is also clear that $B^{-1}(\cF^{W^t}_t)   \= B^{-1}(\{\es,\O_0\}) \= \{\es,\O\} \= \cF^{B^t}_t $.

Define $ \Phi_B (r,\o) \df \big(r,B(\o)\big) $, $ \fa (r,\o) \ins [t,\infty) \ti \O $ and set
$ \L_t  \df \big\{\{t\} \ti A_o \n : A_o \ins \cF^{W^t}_t \big\} \cp \big\{ (s,\infty) \ti A \n :    s \ins [t,\infty) \Cp \hQ, \,
   A  \ins \cF^{W^t}_{s-}  \big\} $,
 which generates $\sP^{W^t}$.
For any $ A_o \ins \cF^{W^t}_t $, since $B^{-1}(\cF^{W^t}_t)   \= \cF^{B^t}_t $,
\bea \label{Sep28_02}
 \Phi^{-1}_B \big( \{t\} \ti A_o \big)
 \= \big\{ (s,\o) \ins [t,\infty) \ti \O \n : \Phi_B (s,\o) \ins \{t\} \ti A_o \big\}
 \= \{t\} \ti B^{-1}(A_o) \ins \sP^{B^t} .
\eea
On the other hand, let $ s \ins [t,\infty) \Cp \hQ $ and $  A  \ins \cF^{W^t}_{s-} $.
We see from   \eqref{Sep28_01c}   that $B^{-1}(A) \ins \cF^{B^t}_{s-}$ and thus
\beas
 \Phi^{-1}_B \big( (s,\infty) \ti A \big)
 \= \big\{ (r,\o) \ins [t,\infty) \ti \O \n : \Phi_B (r,\o) \ins  (s,\infty) \ti A  \big\}
 \= (s,\infty) \ti B^{-1}(A) \ins \sP^{B^t} .
\eeas
This together with \eqref{Sep28_02} shows that
$ 
\L_t \sb \sS_t  \df \{ D \sb [t,\infty) \ti \O_0  \n : \Phi^{-1}_B (D) \ins  \sP^{B^t}  \} $,
 which is a sigma-field of   $ [t,\infty) \ti \O_0 $.
 So
$  \sP^{W^t} 
\sb \sS_t $ or  $ \Phi^{-1}_B (D) \ins  \sP^{B^t}  $ for any $ D \ins \sP^{W^t} $.

   For any $\cE \ins \sB(\hU)$,  since $\wh{\mu}^{-1} (\cE) \=  \big\{ (s,\o_0) \ins [t,\infty) \ti \O_0 \n : \wh{\mu}(s,\o_0) \ins \cE \big\} \ins \sP^{W^t}$, it follows that
\bea  \label{Oct01_04}
\big\{ (s,\o) \ins [t,\infty) \ti \O \n : \mu  (s,\o) \ins \cE \big\}
\= \big\{ (s,\o) \ins [t,\infty) \ti \O \n : \wh{\mu} \big( \Phi_B (s,\o) \big)   \ins \cE \big\}
\= \Phi^{-1}_B \big( \wh{\mu}^{-1} (\cE) \big) \ins  \sP^{B^t}  .
\eea
 Hence $\{\mu_s\}_{s \in [t,\infty)}$ is an  $\bF^{B^t}-$predictable process. \qed

\begin{lemm} \label{lem_predict}
  Let $(\O, \cF, P)$ be   a    probability space.
Given $t \ins [0,\infty)$, let $\wh{\mu} \= \{\wh{\mu}_s\}_{s \in [t,\infty)} $ be a $\hU-$valued, $\bF^{W^t,P_0} -$ predictable  process on $\O_0$ and let   $ B \= \{B_s\}_{s \in [0,\infty)}$ be an $\hR^d-$valued continuous process on $\O$ with $B_0 \= \bz$ such that
    $B^t$ 
 is a  Brownian motion on $(\O, \cF, P)$.
 Then  $\mu_s (\o) \df \wh{\mu}_s \big(B(\o)\big)$, $(s,\o) \ins [t,\infty) \ti \O$ is an $\bF^{B^t,P}-$predictable process.

\end{lemm}

\ss \no {\bf Proof:}  Let  $\sP^{W^t,P_0}$ be the predictable sigma-field of $[t,\infty) \ti \O_0$ with respect to the filtration $\bF^{W^t,P_0}$   \big(i.e., the sigma-field of $[t,\infty) \ti \O_0$ generated by all $\big\{\cF^{W^t,P_0}_{s-}\big\}_{s \in [t,\infty)}-$adapted, c\`agl\`ad processes\big), and let $\sP^{B^t,P}$ be the predictable sigma-field of $[t,\infty) \ti \O$ with respect to the filtration $\bF^{B^t,P}$.

      Applying Lemma \ref{lem_122921_11} (iii) with $t_0 \= t$,  $(\O_1, \cF_1, P_1,B^1)   \= \big(\O, \cF, P , B \big) $, $(\O_2, \cF_2, P_2,B^2) \= \big(\O_0, \sB(\O_0) , P_0 , W\big) $ and $\Phi \= B$  implies that
  $B^{-1} \big(\cF^{W^t,P_0}_t\big) \sb \cF^{B^t,P}_t $,
  and that for any $s \ins (t,\infty)$   the sigma-field of $\O_0$,
  $   \{A \sb \O_0 \n : B^{-1}(A) \ins \cF^{B^t,P}_{s-} \}$, contains $\cF^{W^t,P_0}_r$,
  $\fa r \ins [t,s)$ and thus includes $\cF^{W^t,P_0}_{s-} \= \si\Big(\underset{r \in  [t,s)}{\cup} \cF^{W^t,P_0}_r \Big) $. Namely,
  \bea \label{Sep28_01}
  B^{-1}(A) \ins \cF^{B^t,P}_{s-} , \q \fa A \ins \cF^{W^t,P_0}_{s-} , ~ \fa s \ins (t,\infty)  .
  \eea


Define $ \Phi_B (r,\o) \df \big(r,B(\o)\big) $, $ \fa (r,\o) \ins [t,\infty) \ti \O $ and set
$\L_t \df \big\{\{t\} \ti A_o \n : A_o \ins \cF^{W^t,P_0}_t \big\} \cp \big\{ (s,\infty) \ti A \n :    s \ins [t,\infty) \Cp \hQ, \,
   A  \ins \cF^{W^t,P_0}_{s-}  \big\} $,
 which generates $\sP^{W^t,P_0}$.
For any $ A_o \ins \cF^{W^t,P_0}_t $, since $B^{-1} \big(\cF^{W^t,P_0}_t\big) \sb \cF^{B^t,P}_t $,
\bea \label{Sep28_02}
 \Phi^{-1}_B \big( \{t\} \ti A_o \big)
 \= \big\{ (s,\o) \ins [t,\infty) \ti \O \n : \Phi_B (s,\o) \ins \{t\} \ti A_o \big\}
 \= \{t\} \ti B^{-1}(A_o) \ins \sP^{B^t,P} .
\eea
On the other hand, let $ s \ins [t,\infty) \Cp \hQ $ and $  A  \ins \cF^{W^t,P_0}_{s-} $.
We see from   \eqref{Sep28_01}   that $B^{-1}(A) \ins \cF^{B^t,P}_{s-}$ and thus
\beas
 \Phi^{-1}_B \big( (s,\infty) \ti A \big)
 \= \big\{ (r,\o) \ins [t,\infty) \ti \O \n : \Phi_B (r,\o) \ins  (s,\infty) \ti A  \big\}
 \= (s,\infty) \ti B^{-1}(A) \ins \sP^{B^t,P} .
\eeas
This together with \eqref{Sep28_02} shows that
$ 
\L_t \sb \sS_t \df \{ D \sb [t,\infty) \ti \O_0  \n : \Phi^{-1}_B (D) \ins  \sP^{B^t,P}  \} $,
 which is a sigma-field of $ [t,\infty) \ti \O_0 $.
 So $  \sP^{W^t,P_0} 
 \sb \sS_t $ or  $ \Phi^{-1}_B (D) \ins  \sP^{B^t,P}  $ for any $ D \ins \sP^{W^t,P_0} $.

   For any $\cE \ins \sB(\hU)$,  since $\wh{\mu}^{-1} (\cE) \=  \big\{ (s,\o_0) \ins [t,\infty) \ti \O_0 \n : \wh{\mu}(s,\o_0) \ins \cE \big\} \ins \sP^{W^t,P_0}$, it follows that
\beas
\big\{ (s,\o) \ins [t,\infty) \ti \O \n : \mu  (s,\o) \ins \cE \big\}
\= \big\{ (s,\o) \ins [t,\infty) \ti \O \n : \wh{\mu} \big( \Phi_B (s,\o) \big)   \ins \cE \big\}
\= \Phi^{-1}_B \big( \wh{\mu}^{-1} (\cE) \big) \ins  \sP^{B^t,P}  .
\eeas
 Hence $\{\mu_s\}_{s \in [t,\infty)}$ is an  $\bF^{B^t,P}-$predictable process. \qed

\begin{lemm} \label{lem_061222}
Given $t \ins [0,\infty)$, $\{\mu_s\}_{s \in [t,\infty)}$ is a $\hU-$valued, $\bF^{W^t}-$predictable process on $\O_0$
if and only if there exists a $\hU-$valued, $\bF^W-$predictable process $\{\nu_\fs\}_{\fs \in [0,\infty)}$ on $\O_0$ such that
$ \nu_\fs \big(\sW^t(\o_0)\big) \= \mu_{t+\fs} (\o_0)    $ for any $  (\fs,\o_0) \ins [0,\infty) \ti \O_0$, where
$\sW^t_\fs  \df W_{t+\fs}  \- W_t   \= W^t_{t+\fs}  $, $  \fa \fs  \ins [0,\infty)  $.
\end{lemm}

\no {\bf Proof:}    Define two processes on $\O_0$ by:
\beas
    \beta_s (\o_0) \df W_{s \vee t} (\o_0) \- W_t (\o_0) \aand
  \ddot{\beta}_s(\o_0) \df W_{(s-t)^+}(\o_0)   ,  \q \fa (s,\o_0) \ins [0,\infty) \ti \O_0 .
\eeas

\no {\bf 1)}  Let $  \{\mu_s\}_{s \in [t,\infty)}$ be a   $\hU-$valued, $ \bF^{W^t } - $predictable process  on $\O_0$.

\no {\bf 1a)}  Since it holds for any $(r,\cE) \ins [t,\infty) \ti \sB(\hR^d)$ that    $\beta^{-1} \big(\{W^t_r \ins \cE\}\big) \= \big\{W^t_r(\beta) \ins \cE \big\} \=  \{ W_r(\beta) \- W_t(\beta) \ins \cE \big\}
 \=  \{  \beta_r \-  \beta_t \ins \cE \big\}
\= \{W_r  \- W_t  \ins \cE\} \= \{W^t_r   \ins \cE\} $,   the sigma-field $\{A_0 \sb \O_0 : \beta^{-1}(A_0) \= A_0 \} $ contains all generating sets of $\cF^{W^t}_\infty$ and thus includes $\cF^{W^t}_\infty$.

Define a mapping $\Phi_t(s,\o_o) \df \big(s,  \beta (\o_0)\big)$, $\fa (s,\o_0) \ins [t,\infty) \ti \O_0$ and set $ \L_t  \df \big\{\{t\} \ti A_o \n : A_o \ins \cF^{W^t}_t \big\} \cp \big\{ (s,\infty) \ti A \n :    s \ins [t,\infty) \Cp   \hQ, \,
   A  \ins \cF^{W^t}_{s-}  \big\} $,
 which generates the $\bF^{W^t}-$predictable sigma$-$field $\sP^{W^t}$. 
  For any $A_o \ins \cF^{W^t}_t  $,   one has
  $ \Phi^{-1}_t \big(\{t\} \ti A_o \big) \= \big\{(s,\o_0) \ins [t,\infty) \ti \O_0 \n : (s,\beta (\o_0)) \ins \{t\} \ti A_o \big\}
 \= \{t\} \ti  \beta^{-1}(A_o) \= \{t\} \ti A_o   $.
 For any $s \ins [t,\infty) \cap \hQ$  and $ A  \ins \cF^{W^t}_{s-}   $,
 we get   $\Phi^{-1}_t \big((s,\infty) \ti  A\big) \= \big\{(s',\o_0) \ins [t,\infty) \ti \O_0 \n : (s', \beta  (\o_0)) \ins (s,\infty) \ti  A \big\}
 \= (s , \infty) \ti  \beta^{-1}(A) \= (s , \infty) \ti A     $.
 So  the sigma-field $\big\{ D \sb [t,\infty) \ti \O_0 \n : \Phi^{-1}_t (D) \= D \big\}$
  contains $\L_t$    and thus includes  $\sP^{W^t}$.
  By a standard approximation,
   \bea \label{061222_31}
   \mu(s,\o_0) \= \mu \big(\Phi_t (s,\o_0) \big) \= \mu (s, \beta(\o_0)), \q \fa (s,\o_0) \ins [t,\infty) \ti \O_0 .
   \eea

\if{0}

 Let $\{u_i\}_{i \in \hN}$ be a countable dense subset of $\big(\hU, \Rho{\hU}\big)$.
 Given $i,n \ins \hN$, we set $o^n_i \df \big\{u \ins \hU: \Rho{\hU} (u,u_i) \< 2^{-n} \big\}$
 as  the open ball centered at $u_i$ with radius $2^{-n}$.
   We also set   $\wt{o}^n_1 \df o^n_1 $ and $\wt{o}^n_i \df o^n_i \big\backslash \Big( \underset{j < i }{\cup} o^n_j \Big) $ for $i \gs 2$.
   For any $(s,\o_0) \ins [t,\infty) \ti \O_0$, define  $\mu^n (s,\o_0) \df \sum_{i \in \hN} u_i \b1_{\{ (s,\o_o) \in D^n_i \}} $, where
    $ D^n_i \df \big\{(s',\o'_0) \ins [t,\infty) \ti \O_0 \n :   \mu  (s',\o'_0) \ins \wt{o}^n_i \big\} \ins \sP^{W^t} $.

  Given $(s,\o_0) \ins [t,\infty) \ti \O_0$ and $n \ins \hN$, since
  $  \Rho{\hU} \big(\mu^n (s,\o_0), \mu  (s,\o_0) \big)
    \=   \sum_{i \in \hN} \b1_{\{  (s,\o_0)   \in D^n_i \}}  \Rho{\hU} \big( u_i , \mu  (s,\o_0)  \big)
   \< 2^{-n} $ \big(similarly $\Rho{\hU} \big(\mu^n \big( \Phi_t (s,\o_0)\big), \mu \big( \Phi_t (s,\o_0)\big) \big) \< 2^{-n} $\big)
   and since
   \beas
   \mu^n \big( \Phi_t (s,\o_0)\big) \= \sum_{i \in \hN} u_i \b1_{\{ (s,\o_o) \in \Phi^{-1}_t (D^n_i) \}}
   \= \sum_{i \in \hN} u_i \b1_{\{ (s,\o_o) \in  D^n_i  \}} \= \mu^n (s,\o_0)  ,
   \eeas
   we see that both $\mu  (s,\o_0)$ and $\mu \big( \Phi_t (s,\o_0)\big)$ are the limits of $\big\{ \mu^n (s,\o_0) \big\}_{n \in \hN}$ under
   $\Rho{\hU}$ and thus $\mu  (s,\o_0) \= \mu \big( \Phi_t (s,\o_0)\big)$.

 \fi

\no {\bf 1b)}   As a shifted canonical process on $\O_0$,   $\big\{\sW^t_\fs \big\}_{\fs \in  [0,\infty)}  $  is also a Brownian motion under $P_0$.  
   Since
   $ \sW^t_\fs \big( \ddot{\beta} (\o_0)\big) \= W_{t+\fs} \big( \ddot{\beta} (\o_0)\big) \- W_t \big( \ddot{\beta} (\o_0)\big)
    \= \ddot{\beta}_{t+\fs} (\o_0)  \-   \ddot{\beta}_t (\o_0) 
    \=  W_\fs(\o_0) $ for any $(\fs,\o_0) \ins [0,\infty) \ti \O_0$,
 applying Lemma \ref{lem_122921_11}   with $t_0 \= 0$, $(\O_1, \cF_1, P_1,B^1)   \= \big(\O_0,  \sB(\O_0),  P_0,W\big) $, $(\O_2, \cF_2, P_2,B^2) \= \big(\O_0,  \sB(\O_0),  P_0, \sW^t\big) $ and $\Phi \= \ddot{\beta}$  yields  that
  \bea  \label{010622_11b}
  \ddot{\beta}^{-1} \big(\cF^{\sW^t}_\fs\big) \= \cF^W_\fs  , \q \fa \fs \ins [0,\infty] ,
  \eea
 where  $\cF^{\sW^t}_\fs \= \si\big(\sW^t_\fr; \fr \ins [0,\fs] \Cp \hR \big) 
 \= \si\big( W^t_r; r \ins [t,t\+\fs]  \Cp \hR \big) \= \cF^{W^t}_{t+\fs} $. 

 Define a mapping $\Psi_t(\fs,\o_o) \df \big(t\+\fs, \ddot{\beta}(\o_0)\big)$, $\fa (\fs,\o_0) \ins [0,\infty) \ti \O_0$.
   For any $A_o \ins \cF^{W^t}_t \= \cF^{\sW^t}_0$, we see from \eqref{010622_11b}   that $\ddot{\beta}^{-1}(A_o) \ins \cF^W_0$,
  so    $\Psi^{-1}_t \big(\{t\} \ti A_o \big) \= \big\{(\fs,\o_0) \ins [0,\infty) \ti \O_0 \n : (t\+\fs,\ddot{\beta} (\o_0)) \ins \{t\} \ti A_o \big\}
 \= \{0\} \ti \ddot{\beta}^{-1}(A_o) $ belongs to $\bF^W-$predictable sigma$-$field $   \sP^W  $
   and  $\Psi^{-1}_t \big((t,\infty) \ti A_o \big) \= \big\{(\fs,\o_0) \ins [0,\infty) \ti \O_0 \n : (t\+\fs,\ddot{\beta} (\o_0)) \ins (t,\infty) \ti A_o \big\} \= (0,\infty) \ti \ddot{\beta}^{-1}(A_o)  \ins \sP^W  $.
 Given $s \ins (t,\infty) \cap \hQ$  and $ A  \ins \cF^{W^t}_{s-}   $,
  as \eqref{010622_11b} shows  $\ddot{\beta}^{-1} \big(\cF^{W^t}_r\big) \= \ddot{\beta}^{-1} \big(\cF^{\sW^t}_{r-t}\big) \= \cF^W_{r-t} \sb \cF^W_{(s-t)-}$,     $\fa r \ins [t,s)$, one can deduce that
  $ \ddot{\beta}^{-1} (A)  \ins \ddot{\beta}^{-1} \big(\cF^{W^t}_{s-} \big) \sb \cF^W_{(s-t)-}  $
  and thus     $\Psi^{-1}_t \big((s,\infty) \ti  A\big) \= \big\{(\fs,\o_0) \ins [0,\infty) \ti \O_0 \n : (t\+\fs,\ddot{\beta} (\o_0)) \ins (s,\infty) \ti  A \big\}
 \= (s\-t,\infty) \ti \ddot{\beta}^{-1}(A)  \ins \sP^W  $.
  It follows that   the sigma-field $\big\{ D \sb [t,\infty) \ti \O_0 \n : \Psi^{-1}_t (D) \ins \sP^W \big\}$
  contains $\L_t$    and thus includes  $\sP^{W^t}$.
  Namely,
  \bea \label{061222_21}
   \Psi^{-1}_t (D) \ins \sP^W , \q \fa D \ins \sP^{W^t}  .
   \eea

  Define $ \nu (\fs,\o_0) \df \mu \big(\Psi_t(\fs,\o_0)\big) \= \mu_{t+\fs} \big(\ddot{\beta} (\o_0)\big)   $, $\fa (\fs,\o_0) \ins [0,\infty)$.
   For any $A \ins \sB(\hU)$, as $D_A \df \big\{(s,\o_0) \ins [t,\infty) \ti \O_0 \n : \mu(s,\o_0) \ins A \big\} \ins \sP^{W^t}$,
  one has $\big\{(\fs,\o_0) \ins [0,\infty) \ti \O_0 \n : \nu(\fs,\o_0) \= \mu \big(\Psi_t(\fs,\o_0)\big) \ins A \big\} \= \Psi^{-1}_t \big(D_A\big) \ins \sP^W $. So $\{\nu_\fs\}_{\fs \in [0,\infty)}$ is  a $\hU-$valued, $\bF^W-$predictable process  on $\O_0$.
  Since
 $ \ddot{\beta}_s (\sW^t)   \= \sW^t_{(s-t)^+} \= W_{t + (s-t)^+} \- W_t  \= W_{s \vee t} \- W_t \=    \beta_s     $
 for any $ s \ins [0,\infty) $,
     \eqref{061222_31} renders that
  $ \nu_\fs \big(\sW^t(\o_0)\big) \= \mu_{t+\fs} \big( \ddot{\beta}(\sW^t(\o_0)) \big) \=  \mu_{t+\fs} (\beta(\o_0)) \=  \mu_{t+\fs} ( \o_0 ) $
  for any $(\fs,\o_0) \ins [0,\infty) \ti \O_0$.

\no {\bf 2)} Let $\{\nu_\fs\}_{\fs \in [0,\infty)}$ is  a $\hU-$valued, $\bF^W-$predictable process  on $\O_0$.

 Since   ${\ddot{\beta}}^t_s \df \ddot{\beta}_s \- \ddot{\beta}_t \= W_{s-t}$, $s \ins [t,\infty)$ is   a Brownian motion on $\big(\O_0,  \sB(\O_0),  P_0\big)$
 and since $  {\ddot{\beta}}^t_s\big(\sW^t(\o_0)\big) \= \sW^t_{s-t}(\o_0) \= W_s(\o_0) \- W_t(\o_0) \= W^t_s (\o_0)$ for any $(s,\o_0) \ins [t,\infty) \ti \O_0 $,
 applying Lemma \ref{lem_122921_11} with $t_0 \= t$,  $(\O_1, \cF_1, P_1,B^1)   \= \big(\O_0,  \sB(\O_0),  P_0,W\big) $, $(\O_2, \cF_2, P_2,B^2) \= \big(\O_0,  \sB(\O_0),  P_0, \ddot{\beta}\big) $ and $\Phi \= \sW^t$   shows that
  $  
   (\sW^t)^{-1} ( \cF^{{\ddot{\beta}}^t }_s) \= \cF^{W^t }_s $, $ \fa s \ins [t,\infty) $.
  As $\cF^{{\ddot{\beta}}^t}_s \= \si\big({\ddot{\beta}}^t_r; r \ins [t,s] \Cp \hR \big)   \= \si\big( W_{r-t}; r \ins [t,s]  \Cp \hR \big)
 \= \si\big( W_{r'}; r' \ins [0,s\-t]  \Cp \hR  \big) \= \cF^W_{s-t} $ for any $s \ins [t,\infty]$,  we see that
  $(\sW^t)^{-1} \big( \cF^W_{s-t} \big) \= \cF^{W^t }_s $ and $(\sW^t)^{-1} \big( \cF^W_{(s-t)-} \big) \= \cF^{W^t }_{s-} $ for any $   s \ins [t,\infty) $.

Define a mapping $\U_t(s,\o_o) \df \big(s \- t,  \sW^t (\o_0)\big)$, $\fa (s,\o_0) \ins [t,\infty) \ti \O_0$ and set $ \L   \df \big\{\{0\} \ti A_o  \n : A_o  \ins \cF^W_0 \big\} \cp \big\{ (s,\infty) \ti A \n :    s \ins [0,\infty) \Cp   \hQ, \,
   A  \ins \cF^W_{s-}  \big\} $,
 which generates $\sP^W$.  
  For any $A_o \ins \cF^W_0  $,  as $ (\sW^t)^{-1}(A_o) \ins \cF^{W^t}_t $, one has
  $ \U^{-1}_t \big(\{0\} \ti A_o \big) \= \big\{(s,\o_0) \ins [t,\infty) \ti \O_0 \n : (s\-t,\sW^t (\o_0)) \ins \{0\} \ti A_o \big\}
 \= \{t\} \ti  (\sW^t)^{-1}(A_o) \ins \sP^{W^t}  $.
 For any $\fs \ins [0,\infty) \cap \hQ$  and $ A  \ins \cF^W_{\fs-}   $,
 we obtain that $(\sW^t)^{-1}(A ) \ins \cF^{W^t}_{(t+\fs)-}$ and thus   $\U^{-1}_t \big((\fs,\infty) \ti  A\big) \= \big\{(s,\o_0) \ins [t,\infty) \ti \O_0 \n : (s\-t, \sW^t  (\o_0)) \ins (\fs,\infty) \ti  A \big\}
 \= (t\+\fs , \infty) \ti  (\sW^t)^{-1}(A) \ins \sP^{W^t}   $.
 So  the sigma-field $\big\{ D \sb [0,\infty) \ti \O_0 \n : \U^{-1}_t (D) \ins \sP^{W^t} \big\}$
  contains $\L $    and thus includes  $\sP^W $. Namely,
  $ \U^{-1}_t (D) \ins \sP^{W^t} $ for any $D \ins \sP^W $.

 Define $\mu_s (\o_0) \df \nu \big(\U_t(s,\o_0)\big) \=  \nu_{s-t} \big(\sW^t (\o_0)\big)   $, $\fa (s,\o_0) \ins [t,\infty) \ti \O_0$.
 For any $\wt{A} \ins \sB(\hU)$, as $D_{\wt{A}} \df \big\{(s,\o_0) \ins [0,\infty) \ti \O_0 \n : \nu(s,\o_0) \ins \wt{A} \big\} \ins \sP^W $,
  one can deduce that $\big\{(s,\o_0) \ins [t,\infty) \ti \O_0 \n : \mu(s,\o_0)   \= \nu \big(\U_t(s,\o_0)\big)
  \ins \wt{A} \big\} \= \U^{-1}_t \big(D_{\wt{A}}\big) \ins \sP^{W^t} $. So $\{\mu_s\}_{s \in [t,\infty)}$ is  a $\hU-$valued, $\bF^{W^t}-$predictable process  on $\O_0$. \qed

   \fi

 \begin{lemm}     \label{lem_061122}
Let  $\mu \=\{\mu_s\}_{s \in [0,\infty)} $ be a $\hU-$valued, $\bF^{W,P_0} -$predictable  process on $\O_0$ with all paths in $\hJ$.
Then $\mu^{-1} (A) \df \{\o_0 \ins \O_0 \n : \mu_\cd (\o_0) \ins A\} \ins \cF^{W,P_0}_\infty $ for any $A \ins \sB(\hJ)$.
\end{lemm}

\no {\bf Proof:}
  Let   $\vf \ins L^0 \big((0,\infty) \ti \hU;\hR\big)$.
     The $ \bF^{W,P_0} -$predictability of $\mu$ implies that
     $ \nu_s(\o_0) \df  \vf \big(s,\mu_s (\o_0)\big)$, $(s,\o_0) \ins (0,\infty) \times  \O_0  $
 is also an $ \bF^{W,P_0} -$predictable process and   $\int_0^\infty \nu_s ds \= \int_0^\infty \vf  (s,\mu_s) ds $
 is thus  $\cF^{W,P_0}_\infty-$measurable. Then it holds  for any $  \cE \ins \sB(\hR)$  that
 $ \mu^{-1} \big( (I_\vf)^{-1}(\cE) \big)
    \=     \big\{ \o_0 \ins \O_0 \n : I_\vf \big(\mu_\cd (\o_0)\big) \ins  \cE  \big\}
 \=   \big\{ \o_0 \ins \O_0 \n : \int_0^\infty   \vf  \big(s, \mu_s (\o_0)\big) ds  \ins  \cE  \big\} \ins \cF^{W,P_0}_\infty   $,
   which together with Lemma \ref{lem_M29_01}  (1) shows that the sigma-field 
 $    \big\{ A \sb \hJ \n : \mu^{-1}  (A) \ins \cF^{W,P_0}_\infty \big\}$
 includes all generating sets of $ \sB(\hJ)$ and thus contains $ \sB(\hJ)$. \qed

 \if{0}

\no {\bf Method 2:}    
   By \eqref{J04_01},  the topology $\fT_\sharp(\hJ)$ of $\hJ$ has a countable subbase
  $\big\{ \fri^{-1}_\hJ \big(O_{\frac{1}{n}} (\fm_k,\phi_j)\big) \big\}_{n,k,j \in \hN}$.
  Given $n,k,j \ins \hN$, the continuity of function $\phi_j$ and the $ \bF^{W,P_0} -$predictability of process $\mu$
  imply that the process $ \mu^j_s(\o_0) \df \phi_j \big(s, \mu_s (\o_0)\big) $, $(s,\o_0) \ins [0,\infty) \ti \O_0$
  is also $ \bF^{W,P_0} -$predictable.
  Then the random variable $\xi_j(\o_0) \df \int_0^\infty e^{-s}   \phi_j (s, \mu_s (\o_0))   ds$, $\o_0  \ins  \O_0$ is $ \cF^{W,P_0}_\infty -$measurable and thus
 \beas
\hspace{-0.5cm} \mu^{-1} \big( \fri^{-1}_\hJ \big(O_{\frac{1}{n}} (\fm_k,\phi_j)\big) \big)
 & \tn  \=  & \tn  \bigg\{\o_0 \ins \O_0 \n :  \Big| \int_0^\infty e^{-s}   \phi_j (s, \mu_s (\o_0)) ds
 \- \int_0^\infty \n \int_\hU
  \phi_j (s,u) \fm_k (ds,du)      \Big| \< 1/n  \bigg\} \\
 & \tn \= & \tn  \big\{\o_0 \ins \O_0 \n :    \xi_j(\o_0) \ins  ( c_{k,j} \- 1/n,  c_{k,j} \+ 1/n ) \big\} \ins \cF^{W,P_0}_\infty ,
 \eeas
 where $ c_{k,j} \df  \int_0^\infty \n \int_\hU
 \phi_j (s,u) \fm_k (ds,du) \ins \hR   $ is a constant independent of $\o_0 \ins \O_0$.
 It follows that the sigma-field
 $ \L \df  \{A \ins \sB(\hJ) \n : \mu^{-1} (A) \ins \cF^{W,P_0}_\infty  \} $
 of $\hJ$ contains all open sets of $\fT_\sharp(\hJ)$ and thus equals $\sB(\hJ)$.

\fi

\begin{lemm} \label{lem_080422}

For any  $t \ins [0,\infty)$ and $(\o_0,\fu ) \ins \O_0 \ti \hJ$,   define
$ \sW^t_\fs(\o_0) \df \o_0(t\+\fs) \- \o_0(t) $ and $  \sU^t_\fs(\fu) \df \fu(t\+\fs) $,
   $\fa \fs \ins [0,\infty) $.
   Then $(t,\o_0) \mto \sW^t (\o_0)$ is a continuous mapping from $[0,\infty) \ti \O_0 $ to $\O_0$
   and $(t,\fu) \mto \sU^t (\fu)$ is a continuous mapping from $[0,\infty) \ti \hJ $ to $\hJ$.

\end{lemm}

\no {\bf Proof: 1)} Let $(t,\o_0) \ins [0,\infty) \ti \O_0$ and let $\e \ins (0,1)$.
Set $N \df \lceil 2 \- \log_2 \e \rceil$ and $T \df \lceil t\+1 \rceil $. Since $\o_0(s) $ is uniformly continuous in $s \ins [0,N\+T]$ ,
there exists $\dis \d \= \d ( t,\o_0,\e ) \ins \Big(0,\frac{\e}{2^{T+3}}\Big)$ such that $\dis |\o_0(s_1) \- \o_0(s_2)| \ls \frac{\e}{4N}$ for any
$s_1,s_2 \ins [0,N\+T]$ with $|s_2\-s_1| \< \d$.

 For any $(t',\o'_0) \ins [0,\infty) \ti \O_0$ with $ |t'\-t| \ve \Rho{\O_0}(\o_0,\o'_0) \< \d  $,
 we can deduce that
 \beas
  \Rho{\O_0} (\sW^t(\o_0),\sW^{t'}(\o'_0) \big)
 & \tn  \= & \tn  \sum_{n \in \hN} \Big( 2^{-n} \ld \Sup{\fs \in [0,n]}  \big|\o_0(t\+\fs)  \- \o_0(t) \- \o'_0(t'\+\fs)   \+ \o'_0(t') \big| \Big)  \\
   & \tn \ls & \tn  \sum^N_{n = 1} \Big( 2^{-n} \ld 2 \Sup{\fs \in [0,n]}  \big|\o_0(t\+\fs)   \- \o'_0(t'\+\fs)    \big| \Big) \+ \sum^\infty_{n=N+1} 2^{-n} \\
      & \tn  \ls & \tn  2 \sum^N_{n = 1} \Big( 2^{-n} \ld \Sup{\fs \in [0,n]}   \big|\o_0(t\+\fs)   \- \o_0(t'\+\fs)   \big|   \Big)
      \+ 2 \sum^N_{n = 1} \Big( 2^{-n} \ld \Sup{\fs \in [0,n]}   \big|\o_0(t'\+\fs)   \- \o'_0(t'\+\fs)   \big|    \Big)   \+ 2^{-N}  \\
      & \tn  \ls & \tn  2 N \Sup{\substack{s_1,s_2 \in [0,N+T]\\|s_2 - s_1| < \d} }    \big|\o_0(s_2)   \- \o_0(s_1)   \big|
      \+ 2 \sum^N_{n = 1} \Big( 2^{-n} \ld \Sup{s \in [0,n+T]}    \big|\o_0(s)   \- \o'_0(s)   \big|    \Big)       \+ \e/4 \\
      & \tn  \ls  & \tn  3\e/4 \+   2 \sum^N_{n = 1} 2^T \Big( 2^{-n-T} \ld \Sup{s \in [0,n+T]}    \big|\o_0(s)   \- \o'_0(s)   \big|    \Big)
     \ls 3\e/4 \+   2^{1+T} \Rho{\O_0} \big(\o_0, \o'_0\big) \< \e .
\eeas
So $(t,\o_0) \mto \sW^t (\o_0)$ is a continuous mapping from $[0,\infty) \ti \O_0 $ to $\O_0$.

 \no {\bf 2)} We next discuss  the continuity of   mapping $(t,\fu) \mto \sU^t (\fu)$   from $[0,\infty) \ti \hJ $ to $\hJ$.

 Denote by $\fT[0,\infty)$   the Euclidean topology on $[0,\infty)$.
 As  $\fT_\sharp(\hJ)$ is   generated by the subbase $\big\{\fri^{-1}_\hJ \big( O_{\frac{1}{n}} (\fm_k,\phi_j ) \big)\big\}_{ n,k,j \in  \hN }$,
 it suffices to show that
 $ A_{n,k,j} \df \big\{(t,\fu) \ins [0,\infty) \ti \hJ \n : \sU^t (\fu) \ins \fri^{-1}_\hJ \big( O_{\frac{1}{n}} (\fm_k,\phi_j ) \big)\big\}$
 belongs to the product topology $\fT[0,\infty) \oti \fT_\sharp(\hJ)$ for any $n,k,j \ins  \hN$.

 Fix $n,k,j \ins  \hN $ and set $\|\phi_j\|_\infty \df \Sup{(s,u) \in [0,\infty) \times \hU} \big| \phi_j (s,u) \big|   $.
 We pick    $(t,\fu) \ins A_{n,k,j} $ and set $\dis \fc \df \frac{1}{5 e^{t+1} (\|\phi_j\|_\infty\+1)}\Big(1/n  \- \hb{$\big| \int_0^\infty e^{-\fs}   \phi_j (\fs,\fu(t\+\fs) ) d\fs \- \int_0^\infty \n \int_\hU   \phi_j (t,u) \fm_k (dt,du)    \big|$} \Big)\> 0 $.
  By the uniform continuity of   $\phi_j$, 
 there exists $\l  \ins (0,\fc)$ such that $\big| \phi_j (s_1,u_1) \- \phi_j (s_2,u_2) \big| \<   \fc $ for any
 $ (s_1,u_1), (s_2,u_2) \ins [0,\infty)\ti \hU$ with $|s_1\-s_2| \ve \Rho{\hU}(u_1,u_2) \< \l $.
  We define another function $ \phi^\l_j   $ of $\wh{C}_b\big([0,\infty) \ti \hU\big)$ by
 $ \phi^\l_j   (s,u) \df    \big( \big( 1 \+  (s\- t)/ \l  \big)^+\ld 1 \big) \phi_j \big(  (s\- t)^+  ,u\big)  $, $ \fa (s,u) \ins [0,\infty) \ti \hU $.
 Clearly, $\cO_\l(\fu) \df \fri^{-1}_\hJ \big( O_\l (\fri_\hJ(\fu),\phi^\l_j )\big) $ is a member of $  \fT_\sharp(\hJ)$.

 Let $(t',\fu') \ins \big((t\-\l)^+, t\+\l\big) \ti \cO_\l(\fu) $.
   Since   $ \l \> \big| \int_0^\infty e^{-s}  \big[ \phi^\l_j (s, \fu(s) )   \-   \phi^\l_j   (s,\fu'(s) ) \big] ds    \big|
 \= \Big| \int_{(t-\l)^+}^t  e^{-s}    \big( 1 \+  (s\- t)/ \l  \big) \big[ \phi_j  ( 0  ,\fu(s)  ) \- \phi_j  ( 0  ,\fu'(s)  ) \big] ds
 \+ \int_t^\infty  e^{-s}    \big[  \phi_j  ( s\-t  ,\fu(s)  ) \- \phi_j  ( s\-t  ,\fu'(s)  ) \big] ds \Big| $,
 one has
 \beas
  \cI_1 & \tn \df & \tn   \Big| \int_t^\infty  e^{-s+t}    \big[  \phi_j  ( s\-t  ,\fu(s)  ) \- \phi_j  ( s\-t  ,\fu'(s)  ) \big] ds \Big|
 \ls e^t \Big| \int_{(t-\l)^+}^t  e^{-s}    \big( 1 \+  (s\- t)/ \l  \big) \big( \phi_j \big( 0  ,\fu(s) \big) \- \phi_j \big( 0  ,\fu'(s) \big) \big) ds \Big| \+ e^t \l \\
  & \tn  \ls  & \tn   e^t \|\phi_j\|_\infty e^{-(t-\l)^+} \l \+ e^t \l  \ls e \l \|\phi_j\|_\infty  \+ e^t \l .
 \eeas
 We can also    estimate:

 \no \bul
$ \cI_2 \df \big| \int_t^\infty e^{-s+t}  \big[ \phi_j (s\-t, \fu'(s) )   \-   \phi_j ( (s\-t')^+ , \fu'(s) ) \big]  d s    \big|
\ls e^t \int_t^\infty e^{-s} \big|\phi_j (s\-t, \fu'(s) ) \- \phi_j ( (s\-t')^+ , \fu'(s) )     \big|  d s
\ls   \fc $.

\no \bul
$ \cI_3 \df
\big| \int_t^\infty \n  e^{-s} (e^t \- e^{t'} )   \phi_j ( (s\-t')^+ , \fu'(s) ) d s   \big|
\ls \|\phi_j\|_\infty \big|e^{t'}\-e^t\big| e^{-t} 
\ls   \|\phi_j\|_\infty  e^{t \vee t'-t} |t \- t'|   \ls  e \l  \|\phi_j\|_\infty    $.

\no \bul
$ \, \cI_4 \df \big| \int_t^\infty  \n   e^{-s+t'}   \phi_j ((s\-t')^+, \fu'(s) ) d s \- \int_{t'}^\infty \n  e^{-s+t'}   \phi_j ( s\-t' , \fu'(s) ) d s   \big|
\ls \int_{t \land t'}^{t \vee t'}  \n  e^{-s+t'}  \big| \phi_j ((s\-t')^+, \fu'(s) ) \big| d s
\ls \|\phi_j\|_\infty     e^{t'- t \land t'}    |t\-t'|  
\ls e  \l \|\phi_j\|_\infty   $.

 Putting them together leads to that
 \beas
\qq && \hspace{-1cm} \Big| \int_0^\infty e^{-\fs}   \phi_j (\fs,\sU^{t'}_\fs (\fu') ) d\fs \- \int_0^\infty \n \int_\hU   \phi_j (t,u) \fm_k (dt,du)   \Big| \\
&& \ls  \Big| \int_0^\infty e^{-\fs}   \phi_j (\fs, \fu'(t'\+\fs) ) d\fs \- \int_0^\infty e^{-\fs}   \phi_j (\fs, \fu (t\+\fs) ) d\fs   \Big|
\+  \Big| \int_0^\infty e^{-\fs}   \phi_j (\fs, \fu (t\+\fs) ) d\fs \- \int_0^\infty \n \int_\hU   \phi_j (t,u) \fm_k (dt,du)   \Big| \\
&& \=  \Big| \int_{t'}^\infty e^{-s+t'}   \phi_j (s\-t', \fu'(s) ) ds \- \int_t^\infty e^{-s+t}   \phi_j (s\-t, \fu (s) ) ds   \Big| \+ (1/n \- 4 e^{t+1} (\|\phi_j\|_\infty\+1) \fc ) \\
&& \ls \sum^4_{i=1} \cI_i \+ \big[1/n \- 5 e^{t+1} (\|\phi_j\|_\infty\+1) \fc \big] \< 1/n .
 \eeas
 This shows  $ \sU^{t'}  (\fu') \ins \fri^{-1}_\hJ \big( O_{\frac{1}{n}} (\fm_k,\phi_j ) \big)$
 and thus   $ (t,\fu) \ins \big((t\-\l)^+, t\+\l\big) \ti \cO_\l(\fu) \sb A_{n,k,j}$.
 Then $A_{n,k,j}$ is an open set of $\fT[0,\infty) \oti \fT_\sharp(\hJ)$, proving the lemma.   \qed

 \begin{lemm}

 \label{lem_mu_oo}
 Given $t \ins [0,\infty)$,
 let $\mu \= \{\mu_r\}_{r \in [t,\infty)}$ be a   $\hU-$valued, $ \bF^{W^t} -$predictable process on $\O_0$.
 For any   $(s,\oo) \ins [t,\infty) \ti \oO   $,  there exists
  a $\hU-$valued, $ \bF^{W^s} -$predictable  process $   \big\{\mu^{s,\oo}_r\big\}_{r \in [s,\infty)}$
 on $\O_0$ such that
 $ \mu^{s,\oo}_r \big( \oW (\oo')\big) \= \mu_r \big( \oW  (\oo')\big)   $, $ \fa (r,\oo') \ins [s,\infty) \ti   \obW^t_{s,\oo}   $,
  where $\obW^t_{s,\oo}   \df \big\{\oo' \ins \oO \n : \oW^t_r (\oo') \= \oW^t_r (\oo), ~ \fa r \ins [t,s] \big\}$.

 \end{lemm}

 \no {\bf Proof: 1)}
    Define      $ \L    \df  \Big\{   D \sb [t,\infty) \ti \O_0   \n:   $  for any $(s,\oo) \ins [t,\infty) \ti \oO$ there exists  $ D^{s,\oo} \ins \sP^{W^s} $     such that  $ \b1_{ \{  (r,\oW  (\oo')  ) \in  D  \}}
     \= \b1_{\big\{  ( r,\oW  (\oo')  ) \in  D^{s,\oo}\big\}}, \;
     \fa (r,\oo') \ins [s,\infty) \ti   \obW^t_{s,\oo}      \Big\} $.
     Clearly, $\es \ins \L $ with $ D^{s,\oo} \= \es $ for any $(s,\oo) \ins [t,\infty) \ti \oO$.
     Given $D \ins \L $, by taking the complement of each $D^{s,\oo}$, we see that  $D^c \ins \L $.
  \if{0}

 Let $r \ins [t,\infty)$ and  define
    $  \L_r     \df     \big\{  A \sb   \O_0   \n: $  for any $(s,\oo) \ins [t,r] \ti \oO$   there exists   $A^{s,\oo}  \ins \cF^{W^s}_r$   satisfying    $ \b1_{\{ \oW  (\oo')  \in  A \}}     \= \b1_{\{ \oW   (\oo')  \in  A^{s,\oo} \}},
     \fa  \oo'  \ins  \obW^t_{s,\oo}   \big\} $.

  Clearly,  $\es \ins \L_r $ with $A^{s,\oo} \= \es$ for any  $(s,\oo) \ins [t,r] \ti \oO$.
  Given $A \ins \L_r  $, by taking the complement of each $A^{s,\oo}$, we see that  $A^c \ins \L_r $.
  \if{0}

  As $ \b1_{\{ \oW (\oo')  \in  \es \}} \= 0 \= \b1_{\{ \oW (\oo')  \in  \es\}}  $
   and $ \b1_{\{ \oW (\oo')  \in  \O_0 \}} \= 1 \= \b1_{\{ \oW (\oo')  \in  \O_0\}} $
 for any  $   \oo'   \ins   \oO $,
 it is clear that   both $\es$ and $\O_0  $ belong to $\L_r $.

 When $A \ins \L_r  $, one can find an $\ocN_r\ins \sN_\oP \big(\cF^{\oW^t}_\infty\big) $  such that
 for each  $(s,\oo) \ins [t,r] \ti \oO$, there exists   $A^{s,\oo}  \ins \cF^{W^s}_r$   satisfying
    $ \b1_{\{ \oW (\oo')  \in  A \}}
     \= \b1_{\{ \oW (\oo')  \in  A^{s,\oo} \}}$,
     $\fa  \oo'  \ins  \obW^t_{s,\oo} \Cp \ocN^c_r $. Then for any $(s,\oo) \ins [t,r] \ti \oO$,
   $ (A^{s,\oo} )^c \ins \cF^{W^s}_r $  satisfies
   $ \b1_{\{ \oW (\oo')  \in  A^c \}} \= 1\- \b1_{\{ \oW (\oo')  \in  A \}}
     \= 1 \-\b1_{\{ \oW (\oo')  \in A^{s,\oo} \}} \= \b1_{\{ \oW (\oo')  \in  (A^{s,\oo} )^c \}} $, $ \fa \oo'  \ins  \obW^t_{s,\oo} \Cp \ocN^c_r $,
   so $A^c \ins \L_r $.

 \fi
   Let $\{A_n\}_{n \in \hN} \sb \L_r $. For any $n \ins \hN$ and $(s,\oo) \ins [t,r] \ti \oO$,
        there exists   $A^{s,\oo}_n \ins \cF^{W^s}_r$   satisfying
    $  \b1_{ \{\oW  (\oo')\in A_n \} } \= \b1_{\{ \oW (\oo')  \in A^{s,\oo}_n \} }    $ for any $\oo'  \ins  \obW^t_{s,\oo}   $.
    Then for any $(s,\oo) \ins [t,r] \ti \oO$,
      $ \underset{n \in \hN}{\cap} A^{s,\oo}_n \ins \cF^{W^s}_r $ satisfies that
    $  \b1_{\big\{\oW  (\oo')\in \underset{n \in \hN}{\cap}  A_n \big\} }
    \= \underset{n \in \hN}{\prod}  \b1_{ \{\oW  (\oo')\in A_n \} }
    \= \underset{n \in \hN}{\prod} \b1_{\{ \oW (\oo')  \in A^{s,\oo}_n \} }
    \= \b1_{\big\{ \oW  (\oo')  \in \underset{n \in \hN}{\cap} A^{s,\oo}_n \big\} }  $, $ \fa \oo'  \ins  \obW^t_{s,\oo}   $,
  which shows   $ \underset{n \in \hN}{\cap}  A_n \ins \L_r  $.
  So $\L_r $ is a sigma$-$field of $\O_0$.

  Let $(a,\cE) \ins [t,r] \ti \sB(\hR^d)$.
  We verify that  $ (W^t_a)^{-1} (\cE) \ins \L_r $ by three cases: Let $(s,\oo) \ins [t,r] \ti \oO $.

 \no (i) If $s \gs a$ and $ \oW^t_a   (\oo) \ins \cE  $, then
 $ \b1_{\{ \oW (\oo')  \in (W^t_a)^{-1} (\cE) \}}
  \=  \b1_{\{ \oW^t_a (\oo')  \in  \cE  \}}
  \=  \b1_{\{ \oW^t_a (\oo)  \in  \cE  \}} \= 1
     \= \b1_{\{ \oW  (\oo')  \in  \O_0\}} $,  $  \fa \oo'  \ins  \obW^t_{s,\oo} $.

  \no (ii) If $ s \gs  a   $ but $ \oW^t_a   (\oo) \n  \notin \n \cE  $, then
 $  \b1_{\{ \oW (\oo')  \in (W^t_a)^{-1} (\cE) \}} 
  \=  \b1_{\{ \oW^t_a (\oo)  \in  \cE  \}} \= 0 \= \b1_{\{ \oW  (\oo')  \in  \es\}} $, $ \fa \oo'  \ins  \obW^t_{s,\oo} $.

  \no (iii) If $ s \<  a   $,   set $ \cE_{s,\oo} \df   \big\{\fx \- \oW^t_s   (\oo) \n : \fa \fx \ins \cE \} \ins \sB(\hR^d) $, then $ (W^s_a)^{-1} (\cE_{s,\oo}) \ins 
  \cF^{W^s}_r $ satisfies that
  $ \b1_{\{ \oW (\oo')  \in (W^t_a)^{-1} (\cE) \}} 
  \=  \b1_{\{ \oW^t_a (\oo')  - \oW^t_s (\oo')  \in  \cE_{s,\oo}  \}} 
  \= \b1_{\{ \oW  (\oo')   \in (W^s_a)^{-1}  (\cE_{s,\oo}) \}} $, $ \fa \oo'  \ins  \obW^t_{s,\oo} $.
  So $ (W^t_a)^{-1} (\cE) \ins \L_r $,
  it follows that
     \bea \label{Jan16_03}
    \cF^{W^t}_r \sb   \L_r .
     \eea

     Next, we define
     $ \L    \df  \Big\{   D \sb [t,\infty) \ti \O_0   \n:   $  for any $(s,\oo) \ins [t,\infty) \ti \oO$ there exists an $\bF^{W^s}-$predictable set $ D^{s,\oo} $
       satisfying $ \b1_{\big\{  (r,\oW  (\oo')  ) \in  D \big\}}
     \= \b1_{\big\{  ( r,\oW  (\oo')  ) \in  D^{s,\oo}\big\}}, \;
     \fa (r,\oo') \ins [s,\infty) \ti   \obW^t_{s,\oo}      \Big\} $.
     It is clear that for any $s \ins  [t,\infty)    $ and $(r,\oo') \ins [s,\infty) \ti \oO$ that
  \beas
   \b1_{ \big\{  (r,\oW  (\oo') )  \in \es \big\} } \= 0 \= \b1_{\big\{  (r,\oW  (\oo') )  \in \es \big\} }   \q \hb{and} \q
   \b1_{\big\{  (r,\oW  (\oo') )  \in [t,\infty) \times \O_0 \big\} } \= 1 \=   \b1_{ \big\{  (r,\oW  (\oo') ) \in [s,\infty) \times \O_0 \big\} } ,
   \eeas
   so both $\es$ and $[t,\infty) \ti \O_0$ belong to $\L$.

 Let $D \ins \L$. For any $(s,\oo) \ins [t,\infty) \ti \oO$ there exists an
    $\bF^{W^s}-$predictable set $D^{s,\oo}   $ satisfying
     $  \b1_{\big\{ (  r,\oW  (\oo') ) \in  D \big\}}
     \= \b1_{\big\{ ( r,\oW  (\oo') ) \in  D^{s,\oo}\big\}}$
     for any $ (r,\oo') \ins [s,\infty) \ti   \obW^t_{s,\oo}   $.
    Then for any  $(s,\oo) \ins [t,\infty) \ti \oO$,
    the $\bF^{W^s}-$predictable set $(D^{s,\oo})^c   $  satisfies that
   \beas
   \hspace{-0.3cm}
     \b1_{\big\{ (  r,\oW  (\oo') ) \in  D^c \big\}}
     \=  1 \- \b1_{\big\{ (  r,\oW  (\oo') ) \in  D \big\}}
     \= 1 \- \b1_{\big\{ ( r,\oW  (\oo') ) \in  D^{s,\oo}\big\}}
     \= \b1_{\big\{ ( r,\oW  (\oo') ) \in  (D^{s,\oo})^c\big\}} ,
   \eeas
   $\fa (r,\oo') \ins [s,\infty) \ti   \obW^t_{s,\oo}  $, which shows  $D^c \ins \L$.

   \fi
   Let $\{D_n\}_{n \in \hN} \sb \L$. For any $n \ins \hN$
    and $(s,\oo) \ins [t,\infty) \ti \oO$, there exists $D^{s,\oo}_n \ins \sP^{W^s} $ satisfying
     $  \b1_{\big\{ (  r,\oW  (\oo') ) \in  D_n \big\}}
     \= \b1_{\big\{ ( r,\oW  (\oo') ) \in  D^{s,\oo}_n\big\}}$
     for any $ (r,\oo') \ins [s,\infty) \ti   \obW^t_{s,\oo}   $.
     Then  for any $(s,\oo) \ins [t,\infty) \ti \oO$,
     the $\bF^{W^s}-$predictable set $ \underset{n \in \hN}{\cap} D^{s,\oo}_n   $ satisfies that
    \beas
    \b1_{\big\{ (  r,\oW  (\oo')) \in \underset{n \in \hN}{\cap}  D_n \big\} }
    \= \underset{n \in \hN}{\prod}  \b1_{ \big\{  (  r,\oW  (\oo') ) \in D_n \big\} }
    \= \underset{n \in \hN}{\prod} \b1_{\big\{  (r,\oW  (\oo') )  \in D^{s,\oo}_n \big\} }
    \= \b1_{\big\{  (r,\oW  (\oo') )  \in \underset{n \in \hN}{\cap} D^{s,\oo}_n \big\} } \, ,
     \eeas
     $\fa (r,\oo') \ins [s,\infty) \ti   \obW^t_{s,\oo}  $.
  Namely,  $ \underset{n \in \hN}{\cap} D_n \ins \L $.
  Hence,  $\L$ is a sigma$-$field of $[t,\infty) \ti \O_0$.

  \no {\bf 2)} We next demonstrate  that $ \L$ contains all generating sets of
  the $\bF^{W^t }-$predictable sigma$-$field $\sP^{W^t }$ of $[t,\infty) \ti \O_0$:
  \beas
  \{t\} \ti A \hb{ for }A \ins \cF^{W^t }_t  \aand (q,\infty) \ti A \hb{ for }  q \ins  [t,\infty) \Cp \hQ   , ~
   A  \ins \cF^{W^t }_{q-} .
  \eeas

 \no {\bf (2a)} For any $s \ins [t,\infty)  $ and $ (r,\oo') \ins [s,\infty) \ti  \oO  $, one has
 $   \b1_{\big\{ (  r,\oW  (\oo') ) \in \{t\} \times \O_0 \big\}}
      \= \b1_{\{   r = t  \}} $ and  $ \b1_{\big\{ (  r,\oW  (\oo') ) \in (t,\infty) \times \O_0 \big\}}
      \= \b1_{\{   r = (t,\infty)  \}}$. They show  that  $ \{t\} \ti \O_0 \ins \L $ with $ D^{s,\oo} \= \b1_{\{s=t\}} (\{s\} \ti \O_0)
     \+ \b1_{\{s>t\}} \es$, $\fa (s,\oo) \ins [t,\infty) \ti \oO $
     and that  $ (t,\infty) \ti \O_0 \ins \L $ with $ D^{s,\oo} \= \b1_{\{s=t\}} ( (s,\infty) \ti \O_0)
     \+ \b1_{\{s>t\}} ( [s,\infty) \ti \O_0) $, $\fa (s,\oo) \ins [t,\infty) \ti \oO $.
     Since $\cF^{W^t}_t \= \cF^{W^t}_{t-} \= \{\es,\O_0\}$ and since $ \es \ins \L $, we see that
     $ \{t\} \ti A \ins \L $ for any $A \ins \cF^{W^t}_t$ and that  $ (t,\infty) \ti A \ins \L $ for any $A \ins \cF^{W^t}_{t-}$.

 \no {\bf (2b)} Fix $q \ins   (t,\infty) \Cp \hQ  $ and set
   $ \wh{\L}_q \df \big\{ A \sb \O_0 \n : (q,\infty) \ti A \ins \L   \big\}$.
   Clearly, $\es \ins \wh{\L}_q$.
   Since it holds for any $s \ins [t,\infty)$ and $(r,\oo') \ins [s,\infty) \ti \oO$
   that $\b1_{\big\{ (  r,\oW  (\oo') ) \in (q,\infty) \times \O_0 \big\}}
        \= \b1_{\{ r \in (q,\infty)\}}$, we obtain that
     $ (q,\infty) \ti  \O_0 \ins \L $ with   $ D^{s,\oo} \= \b1_{\{s \le q \}} ((q,\infty)   \ti  \O_0)
     \+ \b1_{\{s > q \}} ([s,\infty)   \ti  \O_0)$, $\fa (s,\oo) \ins [t,\infty) \ti \oO $.
     So   $\O_0$ belongs to $ \wh{\L}_q $.
     Given $A \ins \wh{\L}_q $, since $ (q,\infty) \ti A \ins \L $ and $(q,\infty) \ti  \O_0 \ins \L$,
     we can deduce that  $(q,\infty) \ti A^c \= \big( (q,\infty) \ti  \O_0 \big) \Cp \big( (q,\infty) \ti A \big)^c \ins \L $ and thus
     $A^c \ins \wh{\L}_q$. If $\{A_n\}_{n \in \hN} \ins  \wh{\L}_q$, then
     $ (q,\infty) \ti \Big( \ccup{n \in \hN}{} A_n \Big) \= \ccup{n \in \hN}{} \big( (q,\infty) \ti A_n \big) \ins \L  $, i.e.,
     $ \ccup{n \in \hN}{} A_n \ins \wh{\L}_q $. Thus $\wh{\L}_q$ is a sigma-field of $\O_0$.

 \ss  Let $q' \ins  [t,q)$ and $\cE \ins \sB(\hR^d)$. We show by three cases that
  $ (q,\infty) \ti (W^t_{q'})^{-1} (\cE) \ins \L $:

 \no (i) If   $(s,\oo) \ins [t,q') \ti \oO   $, then  $\cE_{s,\oo} \df   \big\{\fx \- \oW^t_s   (\oo) : \fa \fx \ins \cE \} \ins \sB(\hR^d) $ satisfies that
$ \b1_{\big\{ (  r,\oW  (\oo') ) \in (q,\infty) \times  (W^t_{q'})^{-1}(\cE) \big\}}
   \= \b1_{\{   r  \in (q,\infty)   \}} \b1_{\{ \oW^t_{q'} (\oo')  - \oW^t_s (\oo')  \in \cE_{s,\oo} \}}
   \= \b1_{\big\{ (  r,\oW  (\oo') ) \in (q,\infty) \times ( W^s_{q'})^{-1}(\cE_{s,\oo}) \big\}}   $,
   $ \fa (r,\oo') \ins [s,\infty) \ti \obW^t_{s,\oo} $.

 \no (ii) If $(s,\oo) \ins [q',\infty) \ti (\oW^t_{q'})^{-1}(\cE)  $,
  it   holds for any $  (r,\oo') \ins [s,\infty) \ti \obW^t_{s,\oo} $ that
 $ \b1_{\big\{ (  r,\oW  (\oo') ) \in (q,\infty) \times  (W^t_{q'})^{-1}(\cE) \big\}}
    \=   \b1_{\{   r  \in (q,\infty)   \}} \b1_{\{  \oW^t_{q'}( \oo' )   \in  \cE  \}}
    \=   \b1_{\{   r  \in (q,\infty)   \}} \b1_{\{  \oW^t_{q'}( \oo )   \in  \cE  \}}
    \= \b1_{\{   r  \in (q,\infty)   \}}
    \= \b1_{\big\{ (  r,\oW  (\oo') ) \in (q,\infty) \times  \O_0 \big\}}  $.

 \no (iii) If $(s,\oo) \ins [q',\infty) \ti (\oW^t_{q'})^{-1}(\cE^c)  $,
  it   holds for any   $  (r,\oo') \ins [s,\infty) \ti \obW^t_{s,\oo} $ that
$  \b1_{\big\{ (  r,\oW  (\oo') ) \in (q,\infty) \times  (W^t_{q'})^{-1}(\cE) \big\}}
    \=   \b1_{\{   r  \in (q,\infty)   \}} \b1_{\{  \oW^t_{q'}( \oo )   \in  \cE  \}}
    \= 0   \= \b1_{\big\{ (  r,\oW  (\oo') ) \in   \es  \big\}}  $.

 So $(q,\infty) \ti (W^t_{q'})^{-1}(\cE) \ins \L $  with
$   D^{s,\oo}   \=    \b1_{\{(s,\oo) \in [t,q') \times \oO\}} \big((q,\infty) \ti  (W^s_{q'})^{-1}(\cE_{s,\oo}) \big)
     \+ \b1_{\big\{(s,\oo) \in  [q',\infty) \times (\oW^t_{q'})^{-1}(\cE^c) \big\}} \es
     \+  \b1_{\big\{(s,\oo) \in  [q',q] \times (\oW^t_{q'})^{-1}(\cE) \big\}} \big( (q,\infty) \ti  \O_0 \big)
     \+ \b1_{\big\{(s,\oo) \in  (q,\infty) \times (\oW^t_{q'})^{-1}(\cE) \big\}} \big( [s,\infty) \ti  \O_0 \big) $
  for any   $  (s,\oo) \ins [t,\infty) \ti \oO$. It follows that
     $   \cF^{W^t}_{q-} 
     \= \si \big(W^t_{q'}; q' \ins [t,q)\big)  \sb \wh{\L}_q  $.
  Then
  $\L$ contains all generating sets of  $\sP^{W^t}$ and thus includes $\sP^{W^t}$.

  \no {\bf 3)}
 Let $  \{\mu_s\}_{s \in [t,\infty)}$ be a general $\hU-$valued, $ \bF^{W^t } - $predictable  process  on $\O_0$.

 Let $n \ins \hN$, we set $a^n_i \df i 2^{-n}$, $\fa i \ins \{0,1,\cds,1\+2^n\}$
 and $D^n_i \df \big\{(r,\o_0) \ins [t,\infty) \ti \O_0 \n :  \sI \big(\mu_r(\o_0)\big) \ins [a^n_i,a^n_{i+1}) \big\} \ins \sP^{W^t} \sb \L $, $\fa i \ins \{0,1,\cds, 2^n\}$. So for $i \= 0,1,\cds,2^n$ and   $(s,\oo) \ins [t,\infty) \ti \oO$,  there exists   $ D^{s,\oo,n}_i \ins \sP^{W^s} $    satisfying $ \b1_{\big\{  (r,\oW  (\oo')  ) \in  D^n_i \big\}}
     \= \b1_{\big\{  ( r,\oW  (\oo')  ) \in  D^{s,\oo,n}_i \big\}}, \;
     \fa (r,\oo') \ins [s,\infty) \ti   \obW^t_{s,\oo}  $.

   Fix $(s,\oo ) \ins [t,\infty) \ti \oO  $. For any $n \ins \hN$,  define an    $\bF^{W^s}-$predictable process $\nu^{s,\oo,n}$ 
   by  $ 
 \dis  \nu^{s,\oo,n}_r(\o_0) \df   \sum_{i = 1}^{2^n}\b1_{\big\{(r,\o_0) \in \wt{D}^{s,\oo,n}_i\big\}} a^n_i $,
 $ \fa (r,\o_0) \ins [s,\infty) \ti \O_0 $,
   where  $\wt{D}^{s,\oo,n}_1 \= D^{s,\oo,n}_1 $ and $\wt{D}^{s,\oo,n}_i \df D^{s,\oo,n}_i \big\backslash \Big( \underset{j = 1}{\overset{i-1}{\cup}} D^{s,\oo,n}_j\Big) \ins \sP^{W^s} $ for $i \= 2, \cds \n , 2^n$.
   Then  $ \ul{\nu}^{s,\oo}_r (\o_0) \df \linf{n \to \infty} \nu^{s,\oo,n}_r (\o_0)$, $ \fa (r,\o_0) \ins [s,\infty) \ti \O_0 $
   is a  $[0,1]-$valued, $\bF^{W^s}-$predictable process 
   and
   \beas
   \mu^{s,\oo}_r (\o_0) \df \sI^{-1} \big(\ul{\nu}^{s,\oo}_r (\o_0)\big)  \b1_{  \big\{\ul{\nu}^{s,\oo}_r (\o_0) \in \fE \big\}}
\+ u_0  \b1_{  \big\{\ul{\nu}^{s,\oo}_r (\o_0) \notin \fE\big\}} , \q \fa (r,\o_0) \ins [s,\infty) \ti \O_0
  \eeas
    defines a  $\hU-$valued, $\bF^{W^s}-$predictable process. 

 Let  $(r,\oo') \ins [s,\infty) \ti \obW^t_{s,\oo}  $ and let $n \ins \hN$. For any $i \= 2, \cds, 2^n $,
 since $    0      \ls \b1_{ \big\{(r,\oW (\oo')) \in  D^{s,\oo,n}_i \cap (\underset{j = 1}{\overset{i-1}{\cup}} D^{s,\oo,n}_j )  \big\} }
    \ls  \sum^{i-1}_{j=1} \b1_{ \big\{(r,\oW (\oo')) \in  D^{s,\oo,n}_i \cap   D^{s,\oo,n}_j   \big\}}
    \= \sum^{i-1}_{j=1}  \b1_{\big\{(r,\oW (\oo')) \in D^n_i \cap D^n_j\big\}} \= 0  $,
   we obtain
    $      \b1_{ \big\{(r,\oW (\oo')) \in  \wt{D}^{s,\oo,n}_i\big\}} \=  \b1_{ \big\{(r,\oW (\oo')) \in  D^{s,\oo,n}_i\big\}} \= \b1_{\big\{(r,\oW (\oo')) \in D^n_i  \big\}}  $.
  It follows that
  $ \nu^{s,\oo,n}_r ( \oW (\oo')) \= \sum_{i = 1}^{2^n}  \b1_{ \big\{(r,\oW (\oo')) \in  \wt{D}^{s,\oo,n}_i\big\}} a^n_i \= \sum_{i = 1}^{2^n}  \b1_{\big\{(r,\oW (\oo')) \in D^n_i\big\}} a^n_i  $, $\fa n \ins \hN$ and thus
  $  \sI \big(  \mu_r(\oW (\oo'))\big)\= \lmtu{n \to \infty} \sum_{i = 1}^{2^n}\b1_{\big\{(r,\oW (\oo')) \in  D^n_i \big\}} a^n_i
  \= \lmtu{n \to \infty} \nu^{s,\oo,n}_r \big( \oW (\oo')\big) \= \ul{\nu}^{s,\oo}_r \big( \oW (\oo')\big)$.
  Then    $  \mu^{s,\oo}_r \big(\oW (\oo')\big) \= \mu_r \big(\oW (\oo')\big)  $, $\fa (r,\oo') \ins [s,\infty) \ti   \obW^t_{s,\oo}$. \qed

\begin{lemm} 
\label{lem_Markov}
Given $t \ins [0,\infty)$ and $\oP \ins \fP\big(\oO\big)$, let $\ocW \= \big\{\ocW_s\big\}_{s \in [t,\infty)}$ be a $d-$ dimensional   Brownian motion
  with respect to some filtration $ \fF_\cd \= \{\fF_s\}_{s \in [t,\infty)} $ on $\big(\oO,\sB(\oO), \oP \big)$.
 For any $s \ins [t,\infty)$, if $\oxi$ is a real-valued, $\fF_s-$measurable random variable that is $\oP-$integrable,
 then  $ 
  E_\oP  \big[ \, \oxi \big| \cF^{\ocW}_\infty \big]
\=  E_\oP  \big[ \, \oxi \big| \cF^{\ocW}_s \big]$, $\oP -$a.s.

\end{lemm}

\no {\bf Proof:} Fix $s \ins [t,\infty)$. Let $t \ls t_1 \< \cds \< t_n \ls s \= s_0 \< s_1 \< \cds \< s_k \< \infty $
 and let   $\{\cE_i\}^n_{i=1} \cp \{\cE'_j\}^k_{j=1} \sb \sB(\hR^d)$.

 Set $ \psi_k ( x) \df 1 $, $\fa  x \ins \hR^d$.
 For $j\=k,\cds,1$ recursively,  the Markov property of the Brownian motion $\ocW$
 shows that there exists  another Borel-measurable function  $ \psi_{j-1} \n : \hR^d \mto \hR$ satisfying
 $  E_\oP \Big[ \b1_{ \ocW_{s_j}^{-1}(\cE'_j) } \psi_j \big(\ocW_{s_j}\big) \big| \fF_{s_{j-1}} \Big] \= \psi_{j-1} \big(\ocW_{s_{j-1}} \big) $, $\oP-$a.s. 
 Then we can deduce that
  \beas
 \q  && \hspace{-1.2cm} E_\oP \bigg[ \underset{j=1}{\overset{k}{\prod}} \b1_{\ocW_{s_j}^{-1}(\cE'_j)}  \Big| \fF_s \bigg]
 \= E_\oP \Bigg[ E_\oP \bigg[  \underset{j=1}{\overset{k}{\prod}} \b1_{\ocW_{s_j}^{-1}(\cE'_j)}  \Big| \fF_{s_{k-1}} \bigg] \bigg| \fF_s \Bigg]
 \= E_\oP \bigg[ \, \underset{j=1}{\overset{k-1}{\prod}} \b1_{\ocW_{s_j}^{-1}(\cE'_j)} E_\oP \Big[ \b1_{\ocW_{s_k}^{-1}(\cE'_k)}   \Big| \fF_{s_{k-1}} \Big] \bigg| \fF_s \bigg] \\
 &&\= E_\oP \bigg[ \, \underset{j=1}{\overset{k-1}{\prod}} \b1_{\ocW_{s_j}^{-1}(\cE'_j)} \psi_{k-1} \big(\ocW_{s_{k-1}} \big)  \Big| \fF_s \bigg]
 \= E_\oP \bigg[ \, \underset{j=1}{\overset{k-2}{\prod}} \b1_{\ocW_{s_j}^{-1}(\cE'_j)} E_\oP \Big[ \b1_{ \ocW_{s_{k-1}}^{-1}(\cE'_{k-1})} \psi_{k-1} \big(\ocW_{s_{k-1}}\big)  \Big| \fF_{s_{k-2}} \Big] \bigg| \fF_s \bigg] \\
 &&\= E_\oP \bigg[ \, \underset{j=1}{\overset{k-2}{\prod}} \b1_{\ocW_{s_j}^{-1}(\cE'_j)} \psi_{k-2} \big(\ocW_{s_{k-2}} \big)  \Big| \fF_s \bigg]
 \= \cds \= E_\oP \Big[ \,   \b1_{ \ocW_{s_1}^{-1}(\cE'_1)} \psi_1 \big(\ocW_{s_1} \big)  \big| \fF_s \Big]
 \= \psi_0 \big(\ocW_s \big) , \q \hb{$\oP-$a.s.}
 \eeas
It follows that
  $ E_\oP \Big[  \b1_{\underset{i=1}{\overset{n}{\cap}}\ocW_{t_i}^{-1}(\cE_i)}  \b1_{\underset{j=1}{\overset{k}{\cap}}\ocW_{s_j}^{-1}(\cE'_j)} \Big| \fF_s \Big]
   \=   \underset{i=1}{\overset{n}{\prod}} \b1_{\ocW_{t_i}^{-1}(\cE_i)} E_\oP \Big[
 \underset{j=1}{\overset{k}{\prod}} \b1_{\ocW_{s_j}^{-1}(\cE'_j)} \Big| \fF_s \Big]
 \=   \underset{i=1}{\overset{n}{\prod}} \b1_{\ocW_{t_i}^{-1}(\cE_i)} \psi_0 \big(\ocW_s \big) $,  $\oP-$a.s.
 By Dynkin's Pi-Lambda Theorem, the Lambda-system
 $
 \big\{ A \ins \cF^{\ocW}_\infty \n : E_\oP \big[ \b1_A \big| \fF_s \big] \= \ol{\beta}$,
    \,  $\oP-$a.s. for some  real$-$valued $\cF^{\ocW}_s-$measurable random  variable  $\ol{\beta}$ on $\oO $\big\}
 contains the Pi-system $ 
 \Big\{ \Big(\underset{i=1}{\overset{n}{\cap}}\ocW_{t_i}^{-1}(\cE_i) \Big) \bigcap \Big(\underset{j=1}{\overset{k}{\cap}}\ocW_{s_j}^{-1}(\cE'_j)\Big) \n : t \ls t_1 \< \cds \< t_n \ls s \= s_0 \< s_1 \< \cds \< s_k \< \infty , \,   \{\cE_i\}^n_{i=1} \cp \{\cE'_j\}^k_{j=1} \sb \sB(\hR^d)  \Big\}$
 and thus includes $\cF^{\ocW}_\infty$.

   Let $\oxi$ be a real-valued, $\fF_s-$measurable random variable that is $\oP-$integrable and let $A \ins \cF^{\ocW}_\infty$.
   There exists a real$-$valued $\cF^{\ocW}_s-$measurable random variable  $\ol{\beta}$  such that
 $ E_\oP \big[ \b1_A \big| \fF_s \big] \= \ol{\beta} $, $\oP-$a.s.
 Since it holds for any   $\cA_s \ins \cF^{\ocW}_s \sb \fF_s$ that
 $ E_\oP [\b1_{\cA_s} \b1_A ] \= E_\oP \Big[\b1_{\cA_s} E_\oP \big[ \b1_A \big| \fF_s \big] \Big]
 \= E_\oP \big[\b1_{\cA_s} \ol{\beta} \, \big]   $,
 we see that $ E_\oP \big[ \b1_A \big| \cF^{\ocW}_s \big] \= \ol{\beta}
 \= E_\oP \big[ \b1_A \big| \fF_s \big] $, $\oP-$a.s. Then the tower property implies that
 \beas
   E_\oP \big[\b1_A \oxi \big]
 & \tn \= & \tn 
  E_\oP \Big[  \oxi   E_\oP \big[\b1_A      \big| \fF_s \big]  \Big]
    \= E_\oP \Big[ \oxi   E_\oP \big[\b1_A      \big| \cF^{\ocW}_s \big] \Big]
 \=    E_\oP \Big[ E_\oP \big[ \, \oxi E_\oP [\b1_A      | \cF^{\ocW}_s ] \big| \cF^{\ocW}_s \big] \Big] \\
  & \tn  \=  & \tn   E_\oP \Big[ E_\oP \big[\b1_A  \big| \cF^{\ocW}_s \big] E_\oP \big[ \, \oxi  \big| \cF^{\ocW}_s \big] \Big]
  \=   E_\oP \Big[ E_\oP \Big[\b1_A E_\oP \big[ \, \oxi  | \cF^{\ocW}_s  ]     \big| \cF^{\ocW}_s \Big]  \Big]
\= E_\oP \Big[  \b1_A E_\oP \big[ \, \oxi \big| \cF^{\ocW}_s \big]        \Big]   .
\eeas
As   $A$ runs through $  \cF^{\ocW}_\infty$, we obtain that
$   E_\oP  \big[ \, \oxi \big| \cF^{\ocW}_\infty \big]
\=  E_\oP  \big[ \, \oxi \big| \cF^{\ocW}_s \big]$, $\oP -$a.s. \qed

\bibliographystyle{siam}
\bibliography{OSMC_bib}

\end{document}